\documentclass[preprint,10pt]{elsarticle} 
\usepackage{amssymb}
\usepackage{moreverb} 
 \usepackage{float}   
\usepackage{listings}
\usepackage[linesnumbered,ruled]{algorithm2e}

\SetCommentSty{mycommfont}

\usepackage[english]{babel}

\usepackage{nomencl}  
\usepackage{leftidx} 
\usepackage{subfigure}
 \usepackage{amsmath, amsthm, amssymb}

\usepackage[left=1.5cm,top=2cm,right=1.5cm,bottom=2.5cm]{geometry} %
\usepackage{color}
\usepackage{listings}
\usepackage{mathrsfs}

\newtheorem{remark}{Remark}[section]

\newtheorem{theorem}{Theorem}[section]

\newtheorem{lemma}[theorem]{Lemma}
\newenvironment{demoa}
{\emph{Proof.}}
{\quad \hfill $\Box$}

\definecolor{mygreen}{rgb}{0,0.6,0}
\definecolor{mygray}{rgb}{0.5,0.5,0.5}
\definecolor{mymauve}{rgb}{0.58,0,0.82}
\usepackage{physics}
\usepackage{amsmath}
\usepackage{tikz}
\usepackage{mathdots}
\usepackage{yhmath}
\usepackage{cancel}
\usepackage{array}
\usepackage{multirow}
\usepackage{amssymb}
\usepackage{gensymb}
\usepackage{tabularx}
\usepackage{extarrows}
\usepackage{booktabs}

\usepackage{pgfplots} 
\usepackage{pgfgantt}
\usepackage{pdflscape}
 
\pgfplotsset{compat=newest} 
\pgfplotsset{plot coordinates/math parser=false}

\usetikzlibrary{fadings}
\usetikzlibrary{patterns}
\usetikzlibrary{shadows.blur}
\usetikzlibrary{shapes}

\usepackage{hyperref}
\usepackage{comment}
\usepackage{lipsum} 
\usepackage{color}

\usepackage{colortbl} 

\usepackage{graphicx,epstopdf}
\DeclareUnicodeCharacter{001A}{}
\begin{document}
\newcommand{\dom}[1]{\Omega^{#1}}   
\newcommand{\domPAR}[1]{\Omega^{#1'}}   

 \newcommand{\inpP}{\boldsymbol{\mu}}
  \newcommand{\iinpP}[1]{\inpP_{#1}}
 \newcommand{\TRAC}{T}
  \newcommand{\Eunit}[1]{\boldsymbol{\widehat{E}_{#1}}}
  \newcommand{\TRACbar}{\bar{T}}

 \newcommand{\F}{\boldsymbol{F}}
  \newcommand{\Agp}{\boldsymbol{A}}
    \newcommand{\AgpW}{\boldsymbol{\bar{A}}}
          \newcommand{\RANDOMm}{\boldsymbol{\Omega}}
\newcommand{\dRiter}[1]{\Delta R^{(#1)}}         
                  \newcommand{\randn}{\textrm{randn}}
 \newcommand{\Ebar}{\bar{\boldsymbol{E}}}
    \newcommand{\Ecir}{\boldsymbol{\breve{E}}} 
  \newcommand{\tolSVDb}[1]{\bar{\epsilon}_{#1}} 
    \newcommand{\tolSVDbP}[1]{\bar{\gamma}_{#1}} 
    \newcommand{\tolSVDbPm}{\boldsymbol{\bar{\gamma}}} 
    \newcommand{\innerG}[3]{\langle #1,#2 \rangle_{#3}}
\newcommand{\Ubar}{\boldsymbol{\bar{U}}}
 \newcommand{\Sbar}{\bar{\boldsymbol{S}}}
 \newcommand{\K}{\boldsymbol{K}}
 
     \newcommand{\dCiter}[1]{\boldsymbol{ C^{(#1)}}}
            \newcommand{\dQiter}[1]{{\boldsymbol{\Delta Q}^{(#1)}}}
\newcommand{\range}[1]{\textrm{range }\!(#1)}

\newcommand{\Qiter}[1]{{\boldsymbol{{Q}^{(#1)}}}}
\newcommand{\Qit}[2]{{\boldsymbol{Q_{#1}^{(#2)}}}}
\newcommand{\Piter}[1]{{\boldsymbol{P^{(#1)}}}}
\newcommand{\Citer}[1]{{\boldsymbol{C^{(#1)}}}}
    \newcommand{\tolSVDbM}{\boldsymbol{\bar{\epsilon}} } 
      \newcommand{\ddB}[2]{\boldsymbol{\Delta B}^{(#2)}_{#1}} 

    \newcommand{\Ared}{\boldsymbol{\bar{A}}}
    \newcommand{\Fmat}{\boldsymbol{F}}
  \newcommand{\hadmP}{\circ} 
\newcommand{\tolNEW}{\varepsilon_{nr}}
\renewcommand{\P}{\boldsymbol{P}}
\newcommand{\Pmon}{\boldsymbol{P}}

\newcommand{\V}{\boldsymbol{V}}
\renewcommand{\v}{\boldsymbol{v}}

\newcommand{\Pbool}[1]{\mathbb{P}_{#1}}
 \newcommand{\f}{\boldsymbol{f}}
  \newcommand{\fE}[1]{\boldsymbol{f}^{#1}}
    \newcommand{\fEt}[1]{\boldsymbol{f}^{#1^T}\!\!}
    \newcommand{\eINT}{\boldsymbol{e}}

    \newcommand{\errorINT}{{r}}
        \newcommand{\eDECM}{ {e}_{decm}}

        \newcommand{\rERRORv}{\boldsymbol{r}}
        \newcommand{\rERROR}{r}

 \newcommand{\SNAPf}{\boldsymbol{H}} 
  \newcommand{\tolF}{\epsilon} 
    \newcommand{\tolB}{\epsilon_b} 

 \newcommand{\BasisU}{\boldsymbol{\phi}}   
 \newcommand{\ndofFE}{N_{dof}}
 \newcommand{\nmodesU}{\nU}
  \newcommand{\nmodesF}{\nMOD}
  \newcommand{\nmodesA}{\nMOD}

 \newcommand{\dFE}{\boldsymbol{d}}
  \newcommand{\qRED}{\boldsymbol{q}}
 \newcommand{\Par}[1]{\left( #1 \right)}

 \newcommand{\nSAMPL}{P}   
  \newcommand{\intG}[3]{\displaystyle \int_{#1}  \,#3  \,\, d#2}  
 \newcommand{\x}{\boldsymbol{x}}  
  \newcommand{\xRED}[1]{\boldsymbol{x}_{#1}} 
    \newcommand{\xREDall}{\X} 

  \newcommand{\Xtil}{\boldsymbol{\tilde{X}}}
    \newcommand{\xG}[2]{\boldsymbol{\bar{x}}^{#1}_{#2}}  
        \newcommand{\xGpar}[1]{\boldsymbol{\bar{x}'}_{\!\!#1}}  

    \newcommand{\setALLGP}{\boldsymbol{ \bar{X} }}  
        \newcommand{\xGall}{\boldsymbol{  {X}}_{\!FE} }  
        
                \newcommand{\xGparALL}{\boldsymbol{  {\bar{X}'}} }  

\newcommand{\xFE}{\xGall}
    \newcommand{\Xfe}{\xGall}

    \newcommand{\setREDP}{\xREDall}  
        \newcommand{\setY}{\boldsymbol{\mathcal{Y} }}  
        \newcommand{\setX}{\boldsymbol{\mathcal{X} }}  

    \newcommand{\setREDPit}[1]{\boldsymbol{\mathcal{S}^{(#1)} }}  

  \newcommand{\nSD}{d}  
    \newcommand{\ndim}{\nSD}  

 \newcommand{\paramSPC}{\mathcal{D}}    
  \newcommand{\paramSAMPLE}{\mathcal{D}^s}    

 \newcommand{\RRn}[2]{\mathbb{R}^{#1 \times #2}}  
 \newcommand{\Rn}[1]{\mathbb{R}^{#1}}  
  \newcommand{\Nn}[1]{\mathbb{N}^{#1}}  

 \newcommand{\nU}{n}  
 \newcommand{\ngausRED}{m}   
 \newcommand{\wRED}[1]{\omega_{#1}}   
  \newcommand{\wREDv}{\boldsymbol{\omega}}   
    \newcommand{\wREDall}{\boldsymbol{\omega}}   

  \newcommand{\wREDvIT}[1]{\boldsymbol{\omega}^{(#1)}}   
    \newcommand{\wREDit}[2]{\omega^{(#1)}_{#2}}   

  \newcommand{\mIT}[1]{m^{(#1)}}   

  \newcommand{\W}{\boldsymbol{W}}   
    \newcommand{\Wfe}{\boldsymbol{W}_{\!FE}}   
    \newcommand{\WfeE}[1]{\boldsymbol{W}_{\!FE}^{(#1)}}   
      \newcommand{\dB}{\boldsymbol{\Delta B}} 

  \renewcommand{\b}{\boldsymbol{b}}   
  \newcommand{\bEXACT}{\boldsymbol{b}_{FE}}   
  \newcommand{\bFE}{\boldsymbol{b}_{FE}}   
  \newcommand{\bDECM}{\boldsymbol{b}_{decm}}   
 \newcommand{\dRmin}{\Delta R_{m}}         
\newcommand{\dRmax}{\Delta R_{M}}         
    \newcommand{\ORTHtol}{{\textrm{ORTHtol}} } 
    \newcommand{\SVDtrun}{{\textrm{SVDtrun}} } 
            \newcommand{\dQ}{\boldsymbol{\Delta Q}}
                        \newcommand{\dH}{\boldsymbol{\Delta H}}

            \newcommand{\dA}{\boldsymbol{\Delta A}}

 \newcommand{\nelem}{N_{el}}  
  \newcommand{\ngausE}{r}  
    \newcommand{\nnodeE}{l}  

  \newcommand{\ngausT}{M}  
  \newcommand{\Mallgauss}{M}  

 \newcommand{\defeq}{\mathrel{\mathop:}=}
    \newcommand{\set}[3]{\{ #1\}_{#2}^{#3}} 
  \newcommand{\card}[1]{\textrm{card}(#1)}  
  \newcommand{\supp}[1]{\textrm{supp}(#1)}  

  \newcommand{\derpar}[2]{\dfrac{\partial #1 }{\partial #2} }
   \newcommand{\dderpar}[3]{\dfrac{\partial^2 #1 }{\partial #2\partial #3} }

    \newcommand{\derparAT}[3]{\left . \dfrac{\partial #1 }{\partial #2}\right|_{#3}  \!\! }

  
    \newcommand{\der}[2]{\dfrac{\textrm{d} #1 }{\textrm{d} #2} }
 \newcommand{\ORDER}[1]{\mathcal{O}(#1)}

 \newcommand{\coldos}[2]{\begin{bmatrix} #1 \\ #2  \end{bmatrix}}
 \newcommand{\matcdos}[4]{\begin{bmatrix} #1 & #2 \\ #3 & #4  \end{bmatrix}}
\newcommand{\matctres}[9]{\begin{bmatrix} #1 & #2 & #3 \\ #4 & #5 & #6 \\ #7 & #8 & #9  \end{bmatrix}}
  \newcommand{\RANK}[1]{\textrm{\texttt{rank}} ( #1 ) }

 \newcommand{\rowdos}[2]{\begin{bmatrix} #1  & #2  \end{bmatrix}}
 \newcommand{\rowtres}[3]{\begin{bmatrix} #1  & #2 & #3 \end{bmatrix}}

\newcommand{\rowcuatro}[4]{\begin{bmatrix} #1  & #2 & #3 & #4  \end{bmatrix}}
\newcommand{\rowocho}[8]{\begin{bmatrix} #1  & #2 & #3 & #4 & #5 & #6 & #7 & #8   \end{bmatrix}}
\newcommand{\rowsiete}[7]{\begin{bmatrix} #1  & #2 & #3 & #4 & #5 & #6 & #7 \end{bmatrix}}

\newcommand{\coltres}[3]{\begin{bmatrix} #1 \\ #2 \\ #3  \end{bmatrix}}

 \newcommand{\colcuatro}[4]{\begin{bmatrix} #1 \\ #2 \\ #3 \\ #4  \end{bmatrix}}
  \newcommand{\colcuatroSP}[5]{\begin{bmatrix} #1 \vspace{#5 cm} \\ #2 \vspace{#5 cm} \\ #3 \vspace{#5 cm}  \\ #4  \end{bmatrix}}

\newcommand{\colcinco}[5]{\begin{bmatrix} #1 \\ #2 \\ #3 \\ #4 \\#5  \end{bmatrix}}
 \newcommand{\rowseis}[6]{\begin{bmatrix} #1 & #2 & #3 & #4 & #5 & #6  \end{bmatrix}}

 \newcommand{\colseis}[6]{\begin{bmatrix} #1 \\ #2 \\ #3 \\ #4 \\#5 \\#6 \end{bmatrix}}
 \newcommand{\colsiete}[7]{\begin{bmatrix} #1 \\ #2 \\ #3 \\ #4 \\#5 \\#6 \\#7 \end{bmatrix}}
 \newcommand{\colocho}[8]{\begin{bmatrix} #1 \\ #2 \\ #3 \\ #4 \\#5 \\#6 \\ #7 \\ #8 \end{bmatrix}}

 \newcommand{\zero}{{\boldsymbol{0}}} %
  \newcommand{\zeros}{{\boldsymbol{0}}} %

  %
  \newcommand{\refeq}[1]{Eq.(\ref{#1})} 
  \newcommand{\refpar}[1]{(\ref{#1})} 
 \newcommand{\citeg}[2]{\cite[p.~#2]{#1}}

  \newcommand{\normd}[1]{\Vert #1 \Vert_2} 
    \newcommand{\normLZ}[1]{\Vert #1 \Vert_{0}} 

   \newcommand{\normF}[1]{\Vert #1 \Vert_{F}}  
     \newcommand{\normGEN}[2]{\Vert #1 \Vert_{#2}}  
     
  \newcommand{\maxARG}[2]{\underset{#1}{\textrm{max}} \, #2 } 
 \newcommand{\minARG}[2]{\underset{#1}{\textrm{min}} \, #2 } 
  \newcommand{\argMIN}[2]{\textrm{arg} \, \underset{#1}{\textrm{min}} \, #2 } 

   \newcommand{\MAX}[1]{ \textrm{max}(#1) } 
   \newcommand{\MIN}[1]{ \textrm{min}(#1) } 

 \newcommand{\coseno}[1]{\textrm{cos}\,#1} 
  \newcommand{\cotan}[1]{\textrm{cot}\,#1} 
 
  \newcommand{\bbee}{\begin{enumerate}}
 \newcommand{\eeee}{\end{enumerate}}
 \newcommand{\bbii}{\begin{itemize}}
 \newcommand{\eeii}{\end{itemize}}

\newcommand{\snapMsDOM}[1]{\textbf{A}^{\!\!\boldsymbol{\sigma}}_{\boldsymbol{#1}}}

    \newcommand{\tolSVr}{{\tol}_{{\lambda}}} 
    \newcommand{\tolSVd}{{\tol}_d} 
        \newcommand{\tolSVs}{{\tol}_{\sigma}} 
                \newcommand{\tolSVf}{{\tol}_{f}} 

    \newcommand{\tolSVD}{{\tol}_{svd}} 
    \newcommand{\tolECM}{{\tol}_{ecm}} 

    \newcommand{\tol}{\epsilon} 
     \newcommand{\tolINT}{\theta_{int}}

\newcommand{\ngausDOM}[1]{m_{gs}}
     \newcommand{\ngausDOMred}[1]{m^{*}_{gs}}

               \newcommand{\coorgDOMloc}[2]{\boldsymbol{x}^{#1}_{\!{#2}}}   
 \newcommand{\nmodesSdom}[1]{r^{#1}}  

\newcommand{\nMOD}{p}
\newcommand{\pNMOD}{p}
\newcommand{\sTRAIL}{q}
\newcommand{\qTRAIL}{q}

      \newcommand{\BdomRED}[1]{\boldsymbol{B}^{*{#1}}}
\newcommand{\WdomRED}[1]{\boldsymbol{\mathcal{W}}^{\lowercase{#1}}}
 \newcommand{\DOFl}{\textbf{L}} 
 \newcommand{\DOFr}{\textbf{R}} 
  \newcommand{\DOFrr}{\textrm{R}} 

  \newcommand{\DOFp}{\textbf{P}} 
    \newcommand{\DOFz}{\textbf{Z}} 
        \newcommand{\DOFf}{\textbf{F}} 
        \newcommand{\DOFy}{\textbf{Y}} 

    \newcommand{\DOFs}{\textbf{S}} 

  \newcommand{\nmodesRdom}[1]{p^{#1}}
    \newcommand{\nmodesRdomBASE}[1]{\bar{p}^{#1}}

          \newcommand{\alphaB}{\boldsymbol{\alpha}}

             \newcommand{\SNAPforce}{\boldsymbol{A^f}}
  \newcommand{\nSNAP}{P}
  \newcommand{\BasisS}[2]{\boldsymbol{\varUpsilon}^{{#1}}_{#2}}
         \newcommand{\stress}{\boldsymbol{\sigma}}
              \newcommand{\q}{\boldsymbol{q}}
               \newcommand{\Dq}{\boldsymbol{\Delta q}}
                \newcommand{\Q}{\boldsymbol{Q}}
                \newcommand{\PsiB}{\boldsymbol{\Psi}}
                 \newcommand{\PhiB}{\boldsymbol{\Phi}}
     \newcommand{\dX}{\boldsymbol{\Delta X}}
          \newcommand{\dw}{\boldsymbol{\Delta w}}
        \newcommand{\s}{\boldsymbol{s}}
                \renewcommand{\S}{\boldsymbol{S}}
                                 \newcommand{\Vbar}{\boldsymbol{\bar{V}}}
      \newcommand{\n}{\boldsymbol{n}}
          \newcommand{\N}{\boldsymbol{N}}
                    \newcommand{\NshapeG}[1]{\boldsymbol{N}^{(#1)}}
  \newcommand{\BshapeG}[2]{\boldsymbol{B}^{(#1)}_{#2}}

          \newcommand{\Nglo}{\boldsymbol{\mathcal{N}}}

     \renewcommand{\u}{\boldsymbol{u}}
          \newcommand{\uT}[1]{\boldsymbol{u}^T\!{(#1)}}

          \renewcommand{\t}{\boldsymbol{u}}
          \newcommand{\w}{\boldsymbol{w}}

           \newcommand{\y}{\boldsymbol{y}}
        \newcommand{\z}{\boldsymbol{z}}
                \newcommand{\zDECM}{\boldsymbol{z}}
                \newcommand{\wDECM}{\boldsymbol{w}^*}
                                \newcommand{\Wdecm}{\boldsymbol{W}^*}

                                \newcommand{\wwDECM}[1]{{w}^*_{#1}}

                \newcommand{\xDECM}[1]{\boldsymbol{x}^{*}_{#1}}
                \newcommand{\Xdecm}{\boldsymbol{X}^*}
                \newcommand{\Xcecm}{\boldsymbol{X}^{cecm}}
                \newcommand{\wCECM}{\boldsymbol{w}^{cecm}}
                \newcommand{\CboolCONV}{\mathcal{C}}

                    \newcommand{\Z}{\boldsymbol{Z}}
          \newcommand{\X}{\boldsymbol{X}}
                    \newcommand{\Y}{{\boldsymbol{Y}\!}}
    \newcommand{\TRUE}{\textrm{\textbf{true}}}
    \newcommand{\FALSE}{\textrm{{false}}}
          \newcommand{\xbar}{\boldsymbol{\bar{x}}}
                    \newcommand{\xBAR}{\boldsymbol{\bar{x}}}

                    \newcommand{\xTILDE}{\boldsymbol{\tilde{x}}}
                    \newcommand{\Xtilde}{\boldsymbol{\tilde{X}}}

          \renewcommand{\L}{\boldsymbol{L}}
           \newcommand{\R}{\boldsymbol{R}}
                   \renewcommand{\r}{\boldsymbol{r}}
                  \newcommand{\rBAR}{\boldsymbol{\bar{r}}}

                  \newcommand{\C}{\boldsymbol{C}}
                   \newcommand{\OmegaB}{\boldsymbol{\Omega}}
                   \newcommand{\xiB}{\boldsymbol{\xi}}
                                      \newcommand{\xPAR}{\boldsymbol{x'}}
                                       \newcommand{\xxPAR}{{x'}}

                                      \newcommand{\xiG}[1]{\boldsymbol{\bar{\xi}}_{#1}}
                                      
    \newcommand{\phiMAP}[1]{\boldsymbol{\varphi}^{#1}}
      \newcommand{\Expectation}[1]{\mathbb{E}\Par{#1}}

  \newcommand{\ddSS}{\displaystyle} 
\newcommand{\cred}[1]{ {\color{red}  #1}}
\newcommand{\cblue}[1]{ {\color{blue}  #1}}
\newcommand{\cgreen}[1]{ {\color{green}  #1}}
\newcommand{\cgray}[1]{ {\color{gray}  #1}}
\newcommand{\corag}[1]{ {\color{orange}  #1}}

   \newcommand{\intGG}[4]{\displaystyle \int_{#1}^{#2}  \,#4  \,\, d#3}
 \newcommand{\RRaam}{\;\;\;{\color{black}\Rightarrow}\;\;\;}
 \newcommand{\RRaa}{{\color{black}\Rightarrow}\;\;\;}
 \newcommand{\rraa}{\rightarrow}
\newcommand{\U}{\boldsymbol{U}}
\newcommand{\Uz}{{{\boldsymbol{U}}^{\zDECM}}}
\newcommand{\Uorth}{\boldsymbol{U}_{\!\bot}}
\newcommand{\UUorth}{{U}_{\!\bot}}

\newcommand{\uORTH}[1]{\boldsymbol{u}_{\!\bot,#1}}

\newcommand{\Vorth}{\boldsymbol{V}_{\!\bot}}
\newcommand{\Sorth}{\boldsymbol{S}_{\!\bot}}
\newcommand{\SSorth}[1]{{{S}_{\!\bot}}_{#1}}

\newcommand{\Ue}[1]{\boldsymbol{U}^{(#1)}}

\newcommand{\Utrail}{\boldsymbol{U_t}}
\newcommand{\Strail}{\boldsymbol{S_t}}
\newcommand{\GbarTRAIL}{\boldsymbol{\bar{G}_t}}
\newcommand{\Gtrail}{\boldsymbol{{G}_t}}

 \newcommand{\Fsnp}{\boldsymbol{\mathcal{F}}}
 \newcommand{\D}{\boldsymbol{D}}
 \renewcommand{\d}{\boldsymbol{d}}

\newcommand{\mRED}{m}
     \newcommand{\strain}{{\boldsymbol{\varepsilon}}} 
 \newcommand{\setPoints}{\boldsymbol{\mathcal{I}}}  
  \newcommand{\ssetPoints}[1]{{\mathcal{I}_{#1}}}  
 \newcommand{\setElems}{\boldsymbol{\mathcal{E}}}  
  \newcommand{\ssetElems}[1]{{\mathcal{E}_{#1}}}  
 \newcommand{\setGauss}{\boldsymbol{\mathcal{G}}}  
  \newcommand{\ssetGauss}[1]{{\mathcal{G}_{#1}}}  

 \newcommand{\setPointsCAND}{\textbf{y}_{\boldsymbol{0}}}  

\newcommand{\setIndices}{\boldsymbol{\bar{z}}}  

\newcommand{\Kll}{\boldsymbol{K_{ll}}}

\newcommand{\Fd}[1]{\boldsymbol{F_{#1}}}

\newcommand{\A}{\boldsymbol{A}}
\newcommand{\Afe}{\boldsymbol{A}_{FE}}
\newcommand{\AfeE}[1]{\boldsymbol{A}_{FE}^{(#1)}}

\newcommand{\AfeW}{\boldsymbol{\bar{A}}}

\newcommand{\Abar}{{\boldsymbol{\bar{A}}}}
\newcommand{\Uw}{{\boldsymbol{\bar{U}}}}
\newcommand{\Ew}{{\boldsymbol{\bar{E}}}}

\newcommand{\Abcs}[1]{\boldsymbol{A}_{#1}}

\newcommand{\Ab}{\textbf{A}}
\newcommand{\Qb}{\textbf{Q}}
\newcommand{\Yb}{\textbf{Y}}
\newcommand{\Eb}{\textbf{E}}
\newcommand{\Gb}{\textbf{G}}
\newcommand{\Gbar}{\boldsymbol{\bar{G}}}
\newcommand{\Eset}{\boldsymbol{E_{\setPoints}}}

        \newcommand{\SVD}[1]{\textrm{\texttt{SVD}} ({#1}) } 
                \newcommand{\MakeOneZero}[1]{\textrm{\texttt{MAKE1ZERO}} ({#1}) } 
                                 \newcommand{\MakeOneZeroLOC}[1]{\textrm{\texttt{MAKE1ZEROloc}} ({#1}) } 
                 \newcommand{\floor}[1]{\textrm{\texttt{floor}} ({#1}) } 
                 \newcommand{\CEIL}[1]{\textrm{\texttt{ceil}} ({#1}) } 

                 \newcommand{\SPARSIF}[1]{\textrm{\texttt{SPARSIF}} ({#1}) } 
                                \newcommand{\SPARSIFloc}[1]{\textrm{\texttt{SPARSIFloc}} ({#1}) } 
                                \newcommand{\EvalBasis}[1]{\textrm{\texttt{EVALBASIS}} ({#1}) } 

                                       \newcommand{\SOLVERES}[1]{\textrm{\texttt{SOLVERES}} ({#1}) } 

\newcommand{\MeshINFO}{\boldsymbol{\mathcal{M}}}
\newcommand{\ElemINFO}{\boldsymbol{\mathcal{E}} }
\newcommand{\Nsteps}{N_{steps}}
\newcommand{\AuxVAR}{\boldsymbol{\mathcal{A}}}

                \newcommand{\SRSVD}[1]{\textrm{\texttt{SRSVD}} ({#1}) } 

        \newcommand{\diag}[1]{\,\textrm{diag} {#1} } 
 \newcommand{\J}{\boldsymbol{J}}
  \newcommand{\Jhat}{\boldsymbol{\widehat{J}}}
  \newcommand{\JhatK}[1]{\boldsymbol{\widehat{J}}^{\,#1}}

 \renewcommand{\a}{\boldsymbol{a}}
\newcommand{\e}{\boldsymbol{e}}
\newcommand{\E}{\boldsymbol{E}}
\newcommand{\B}{\boldsymbol{B}}
  \newcommand{\minusSET}{\,\backslash\,}
 \newcommand{\G}{\boldsymbol{G}}
  \newcommand{\Gw}{\boldsymbol{G_w}}

 \newcommand{\g}{\boldsymbol{g}}
\renewcommand{\H}{\boldsymbol{H}}
\newcommand{\h}{\boldsymbol{h}}
\newcommand{\m}{\boldsymbol{m}}
    \newcommand{\argMAX}[2]{\textrm{arg}\,\underset{#1}{\textrm{max}} \, #2 } 
 \renewcommand{\c}{\boldsymbol{c}}
 \newcommand{\includeCOMMENT}[2]{
  \floatstyle{boxed}
  \restylefloat{figure}
\begin{figure}[H]
\centering
		\includegraphics[width=#2\textwidth]{EPS_textP/#1.eps}
\end{figure}
 \floatstyle{plain}url
  \restylefloat{figure}
 }
 \newcommand{\bbqq}{\begin{quotation}}
\newcommand{\eeqq}{\end{quotation}}
   \newcommand{\figONE}[3]{
 \floatstyle{boxed}
  \restylefloat{figure}
\begin{figure}[H]
\label{fig:#1}
\centering
		\includegraphics[width=#2\textwidth]{EPS_textP/#1.eps}
		 \caption{#3 }
\end{figure}
 \floatstyle{plain}
  \restylefloat{figure}
 }

    \newcommand{\figONEj}[3]{
 \floatstyle{boxed}
  \restylefloat{figure}
\begin{figure}[H]
\label{fig:#1}
\centering
		\includegraphics[width=#2\textwidth]{EPS_textP/#1.jpg}
		 \caption{#3 }
\end{figure}
 \floatstyle{plain}
  \restylefloat{figure}
 }
 
    \newcommand{\figTWO}[3]{
 \floatstyle{boxed}
  \restylefloat{figure}
\begin{figure}[H]
\label{fig:#1}
\centering
		\subfigure[]{\label{fig:mode2a}\includegraphics[width=0.44\textwidth]{EPS_textP/#1.eps}}
		\subfigure[]{\label{fig:mode2a}\includegraphics[width=0.44\textwidth]{EPS_textP/#2.eps}}
		 \caption{#3 }
\end{figure}
 \floatstyle{plain}
  \restylefloat{figure}
 }

 \newcommand{\nmodesFdom}{p} 
 \newcommand{\BasisF}{\boldsymbol{\Lambda}}   

 \newcommand{\ident}{\boldsymbol{I}} 
   \newcommand{\ones}{\boldsymbol{{1}}}

 \newcommand{\nstrain}{s} 
\newcommand{\ngaus}{n_g} 
   \newcommand{\ncol}[1]{\textrm{\texttt{ncol}} ( #1 ) }
   \newcommand{\nrow}[1]{\textrm{\texttt{nrow}} ( #1 ) }

 \newcommand{\SelectOper}[1]{\boldsymbol{\mathcal{P}}_{#1}}
   \renewcommand{\trace}[1]{\textrm{tr}\,(#1)} 
   \newcommand{\absval}[1]{\vert #1 \vert}

  \newcommand{\diagOL}[4]{\textrm{diag}\,( #1,#2,#3,#4)  }

\newcommand{\gradT}[1]{\boldsymbol{\nabla}^{T}   \!\! #1}

\excludecomment{raulCommentNotPdf}

\newcommand{\jahoHIDE}[1]{}

\newcommand{\jahoAI}[1]{ }

 \newcommand{\jaho}[1]{
  $\Box \;$[\textbf{JAHO: }     \emph{\textbf{ \cblue{#1}}} ]$\,\Box\,$    %
  }

\begin{frontmatter}

\title{CECM: A continuous  empirical cubature method  with application to the dimensional hyperreduction of parameterized   finite element  models}

\author[rvt,els]{J.A. Hern\'{a}ndez\corref{cor1}}
\ead{jhortega@cimne.upc.edu}
\author[rvt,elss]{J.R. Bravo}
\author[rvt,elss]{S. Ares de Parga}

\cortext[cor1]{Corresponding author}

\address[rvt]{Centre Internacional de M\`{e}todes Num\`{e}rics en Enginyeria (CIMNE), Barcelona, Spain}

\address[elss]{Universitat Polit\`{e}cnica de Catalunya, Department of Civil and Environmental Engineering (DECA), Barcelona, Spain}

\address[els]{Universitat Polit\`{e}cnica de Catalunya,  E.S. d'Enginyeries Industrial, Aeroespacial i Audiovisual de Terrassa (ESEIAAT), Terrassa, Spain}%

 \begin{abstract}
 
 We present the Continuous Empirical Cubature Method (CECM), a novel algorithm for empirically devising efficient integration rules. The CECM aims to improve existing cubature methods by producing rules that are close to the optimal, featuring far less points than the number of functions to integrate.

The CECM consists on a two-stage strategy. First, a point selection strategy is applied for obtaining an initial approximation to the cubature rule,  featuring as many points as functions to integrate. The second stage consists in a sparsification strategy in which, alongside the indexes and corresponding weights, the spatial coordinates of the points are also considered as design variables. The positions of the initially selected points are changed to render their associated weights to zero, and in this way, the minimum number of points is achieved.

Although originally conceived within the framework of hyper-reduced order models (HROMs), we present the method's formulation in terms of generic vector-valued functions, thereby accentuating its versatility across various problem domains. To demonstrate the extensive applicability of the method, we conduct numerical validations using univariate and multivariate Lagrange polynomials. In these cases, we show the method's capacity to retrieve the optimal Gaussian rule. We also asses the method for an arbitrary exponential-sinusoidal function in a 3D domain, and finally consider an example of the application of the method to the hyperreduction of a multiscale finite element model, showcasing notable computational performance gains.

A secondary contribution of the current paper is the Sequential Randomized SVD (SRSVD) approach for computing the Singular Value Decomposition (SVD) in a column-partitioned format. The SRSVD is particularly advantageous when matrix sizes approach memory limitations.

\end{abstract}


\end{frontmatter}
 
  \section{Introduction}

The present paper is concerned with a classical problem of numerical analysis:   the   approximation of   integrals over 1D, 2D and 3D domains  of parameterized functions as a weighted sum of the values of such functions at a set of  $m$ points $\{\x_1,\x_2 \ldots \x_m\}$: 
\begin{equation}
   \intG{\Omega}{\Omega}{f(\x; \mu)} \approx  \sum_{g=1}^m f(\x_g; \mu) w_g, 
\end{equation}
(with $m$ as small as possible).  This problem is generally known as either   quadrature (for 1D domains) or \emph{cubature} (for higher dimensions), and has  a long pedigree stretching back as far as C.F. Gauss, who devised in 1814 the eponymous quadrature rule for univariate polynomials.

\subsection{Cubature problem in hyperreduced-order models}

The recent development of the so-called \emph{hyperreduced-order models} (HROMs) for   parameterized  finite element (FE) analyses \cite{farhat2014dimensional,hernandez2017dimensional} has   sparked the resurgence of interest in this classical   problem. Indeed, a crucial step in the construction of such HROMs is the solution of the cubature problem associated to the evaluation of the   nonlinear term(s) in the pertinent governing equations. For instance,  in a Galerkin-based structural HROM, the    nonlinear term   is typically the projection of the nodal FE internal forces $\F^h \in \Rn{\ndofFE}$ (here $\ndofFE$ denotes the number of degrees of freedom of the FE model) onto the span of the displacement modes, i.e.: $\F = \BasisU^T \F^h$, $\BasisU  \in \RRn{\ndofFE}{\nmodesU}$ being the matrix of displacement modes.   The basic premise in these HROMs is that  the number of modes is much smaller that the number of FE degrees of freedom  ($\nmodesU << \ndofFE$). This in turn implies that the internal forces per unit volume will also reside in a space of relatively small dimensions  (independent of the size of the underlying FE mesh), and therefore, its integral over the spatial domain will be, in principle, amenable to approximation by an efficient cubature rule, featuring far less points than the original FE-based rule.    The challenge lies in determining the minimum number of cubature points necessary for achieving a prescribed accuracy, as well as their location and  associated \emph{positive} weights.   The requisite of positive   weights arises from the fact that, in a Galerkin  FE framework,  the   Jacobian matrix of the discrete system of  equations is a weighted sum of the contribution  at each FE Gauss point. Thus,  if  the Jacobian matrices at point level are positive definite, the global matrix is only guaranteed to inherit this desirable attribute   if the cubature weights are positive  \cite{hernandez2017dimensional}.

Before delving into the description of  the diverse approaches proposed to date to deal with this cubature problem in the context of HROMs, it proves convenient   to formally  formulate the    problem in terms of a generic  parameterized vector-valued function $\a : \dom{} \times  \paramSPC \rightarrow \Rn{\nU}$.  Let  $\dom{} = \cup_{e=1}^{\nelem} \dom{e}$ be a  finite element partition of the spatial domain $\dom{} \subset \Rn{\nSD}$  ($\nSD =1, 2$ or $ 3$). For   simplicity of exposition,   assume  that  all   elements are isoparametric and of the same order of interpolation, possessing   $\ngausE$ Gauss points each.     Suppose   we are given the values of  the integrand functions   for $\nSAMPL$ instantiations of the input parameters (  $\{\iinpP{i}  \}_{i=1}^{\nSAMPL} =  \paramSAMPLE  \subset \paramSPC$) at all the  {Gauss} points of the discretization.    The integral of the function over $\dom{}$ for each $\iinpP{j}$ ($j=1,2 \ldots \nSNAP$) can be   calculated by the corresponding element  {Gauss} rule as   
  \begin{equation}
  \label{eq:1}
   b_k =   \sum_{e=1}^{\nelem}  \intG{\dom{e}}{\dom{}}{a_i(\x,\iinpP{j})}  =     \sum_{e=1}^{\nelem} \sum_{g=1}^{\ngausE}   {a_i(\xG{e}{g},\iinpP{j})} W_g^e,    \hspace{0.5cm}   k = (j-1) \nU + i; \;\;\;  j = 1,2 \ldots \nSAMPL; \;\;\; i = 1,2 \ldots \nU. 
  \end{equation}
Here,  $\xG{e}{g} \in \dom{e}$ denotes the position of the $g$-th  Gauss point of element $\dom{e}$,  whereas  $W^e_g > 0$ is the product of the   Gauss weight and the Jacobian of the isoparametric transformation at such a point.  Each $b_k$ ($k = 1,2 \ldots \nSNAP \nU $) is therefore  considered as the ``exact'' integral, that is,  the reference value we wish to approximate.     The above expression can be written in a compact  matrix form as
    \begin{equation}
    \label{eq:3sdfsd}
     \bEXACT = \Afe^T \Wfe,
    \end{equation}    
     where $\bEXACT \in \Rn{\nU \nSNAP}$ is the vector of ``exact'' integrals defined in Eq. \refpar{eq:1}, 
     $\Afe$ is the  matrix obtained from evaluating the integrand functions at all the FE Gauss point, $\xGall = \set{\set{\xG{e}{g}}{g=1}{\ngausE}}{e=1}{\nelem}$, while $\Wfe$    designates  the vector of FE weights,  formed by  gathering all the Gauss weights in a single column vector. Each column of $\Afe$ is the discrete representation of a scalar-valued integrand function, and thus  the total number of  columns    is   equal to the  the number of sampling parameters $\nSAMPL$ times the number of integrand functions per parameter, $\nU$.     The number of rows of $\Afe$, on the other hand,  is equal to the total number of integration points ($\ngausT = \nelem \cdot \ngausE$). In terms of   element contributions, matrix $\Afe$ is expressible  as 
     
     \begin{equation}
     \label{eq:8243t--}
      \Afe = \colcuatro{\AfeE{1}}{\AfeE{2}}{\vdots}{\AfeE{\nelem}}_{\nelem \cdot \ngausE \times \nU \nSNAP} \hspace{-0.1cm} \textrm{where}  \hspace{0.5cm}  \AfeE{e} = \begin{bmatrix} 
                  a_1(\xG{e}{1},\iinpP{1})  &  a_2(\xG{e}{1},\iinpP{1})  &  \cdots & a_n(\xG{e}{1},\iinpP{1}) & \cdots &  a_n(\xG{e}{1},\iinpP{\nSAMPL})   \\
                   a_1(\xG{e}{2},\iinpP{1})  &  a_2(\xG{e}{2},\iinpP{1})  &  \cdots & a_n(\xG{e}{2},\iinpP{1}) & \cdots &  a_n(\xG{e}{2},\iinpP{\nSAMPL})   \\
                   \vdots & \vdots & \vdots & \ddots & \vdots & \vdots  \\     
                  a_1(\xG{e}{\ngausE},\iinpP{1})  &  a_2(\xG{e}{\ngausE},\iinpP{1})  &  \cdots & a_n(\xG{e}{\ngausE},\iinpP{1}) & \cdots &  a_n(\xG{e}{\ngausE},\iinpP{\nSAMPL})   \\    
                 \end{bmatrix}_{\ngausE \times \nU \nSNAP}
     \end{equation} 
     (here $\AfeE{e} \in \RRn{\ngausE}{\nU \nSNAP}$ denotes the block matrix corresponding to the $\ngausE$ Gauss points of  element $\dom{e}$). The same notational scheme is used for the vector of FE weights:
     \begin{equation}
     \label{eq:WFE}
      \Wfe = \colcuatro{\WfeE{1}}{\WfeE{2}}{\vdots}{\WfeE{\nelem}}_{\nelem \cdot \ngausE \times 1} \hspace{0.5cm} \textrm{where}  \hspace{1cm} \WfeE{e}  = \colcuatro{W^e_1}{W^e_2}{\vdots}{W^{e}_{\ngausE}}_{\ngausE \times 1}.
     \end{equation}

 The cubature problem consists in finding  a  set of   points $\xREDall \defeq \{\xRED{g}\}_{g=1}^{\ngausRED} $ ($\xRED{g} \in \dom{}$) and associated \emph{positive} weights $\set{\wRED{g}}{g=1}{\ngausRED}$  (with $\ngausRED$ \emph{as small as possible}) such that the vector of ``exact'' integrals $\bFE$  is approximated     to some desired level of accuracy  $0 \le \tolB \le 1$, i.e.:  
    \begin{equation}
  \label{eq:4}
   \normd{ \A^T\!\!(\xREDall) \wREDall -\bEXACT }{} \le   \tolB \normd{\bEXACT   }{}.
  \end{equation}
    Here, $\normd{\bullet}{}$  is the standard Euclidean norm, whereas   $\A(\xREDall)$ and $\wREDall$ denote the matrix of the integrand evaluated at the set of points $\xREDall$ and their associated weights, respectively: 
         \begin{equation}
    \label{eq:9235add}
    \A(\xREDall)  = \colcuatro{\A(\xRED{1})}{\A(\xRED{2})}{\vdots}{\A(\xRED{\ngausRED})}   \defeq  \begin{bmatrix} 
                  a_1(\xRED{1},\iinpP{1})  &  a_2(\xRED{1},\iinpP{1})  &  \cdots & a_n(\xRED{1},\iinpP{1}) & \cdots &  a_n(\xRED{1},\iinpP{\nSAMPL})   \\
                   a_1(\xRED{2},\iinpP{1})  &  a_2(\xRED{2},\iinpP{1})  &  \cdots & a_n(\xRED{2},\iinpP{1}) & \cdots &  a_n(\xRED{2},\iinpP{\nSAMPL})   \\
                   \vdots & \vdots & \vdots & \ddots & \vdots & \vdots  \\     
                    a_1(\xRED{\ngausRED},\iinpP{1})  &  a_2(\xRED{\ngausRED},\iinpP{1})  &  \cdots & a_n(\xRED{\ngausRED},\iinpP{1}) & \cdots &  a_n(\xRED{\ngausRED},\iinpP{\nSAMPL})   \\  
                 \end{bmatrix}_{\ngausRED  \times  \nSAMPL  \nU}     \hspace{0.5cm} \wREDv = \colcuatro{\wRED{1}}{\wRED{2}}{\vdots}{\wRED{\ngausRED}}_{\ngausRED \times 1}.
  \end{equation}    
  
  \begin{remark}
  \label{rem:iwe,sdfsddds}
   A remark concerning notation is in order here. In \refeq{eq:9235add}, $\A(\x)$ ($\x \in \dom{}$) represents a vector-valued function  that returns the $\nU$ entries of the  integrand function $\a(\x)$ for  the $\nSNAP$ samples  of the input parameters, in the form of row matrix (i.e., $\A(\x) \in \RRn{1}{\nU \nSNAP}$). On the other hand,  when the argument of $\A$ is not a single point, but a collection of $\ngausRED$  points $\xREDall \defeq \{\xRED{g}\}_{g=1}^{\ngausRED} $, then $\A(\xREDall)$ represents a matrix with as many rows as points in the set, i.e. $\A(\xREDall) \in \RRn{\ngausRED}{\nU \nSNAP}$. According to this notational convention, the matrix defined in \refeq{eq:8243t--} can be compactly written as $\Afe =  \A(\xGall)$.  
  \end{remark}

\subsection{State-of-the-art on cubature rules for HROMs}
\label{sec:statenfass}

The first attempts  to  solve the above described  cubature problem  in  the context of reduced-order modeling were carried by the  computer graphics community. The cubature scheme proposed by An et al. in Ref. \cite{an2009optimizing} in 2010 for dealing with the evaluation of the internal forces in geometrically nonlinear models   may be regarded as the  germinal paper in this respect. An and co-workers  \cite{an2009optimizing}  addressed the cubature problem \refpar{eq:4} as a \emph{best subset selection problem}    (i.e., the desired set of points is considered a subset of the entire set of Gauss points,  $\setREDP \subset \xGall$).  They proposed  a greedy strategy that incrementally constructs the   set of points  by minimizing the norm of the   residual of the integration   at each iteration,   while enforcing  the positiveness of the weights. Subsequent papers in the computer graphics community (see Ref. \cite{kim2013subspace} and references therein) revolved around the same idea, and focused fundamentally in improving the efficiency of the scheme originally proposed by An et al. \cite{an2009optimizing}---which turned out to be ostensibly inefficient, for it solves a nonnegative-least squares problem, using the standard Lawson-Hanson algorithm \cite{lawson1974solving},   each time a new point enters the set.

Interesting re-interpretations of the cubature problem  came with the works of Von Tycowicz et al. \cite{von2013efficient} and Pan et al. \cite{pan2015subspace} ---still  within computer graphics circles. 
Both  works recognized the analogy between  the discrete cubature problem  and the quest for \emph{sparsest  solution of underdetermined systems  of equations}, a problem which is   common to  many disciplines  such as       signal processing, inverse problems, and genomic
data analysis \cite{donoho2006most}. Indeed, if we regard the vector of reduced weights $\wREDv$  as  a sparse vector of the same length as $\Wfe$, then the \emph{best subset selection} problem can be posed as that of minimizing the nonzero entries  of $\wREDv$:
\begin{equation}
\label{eq:10}
   \minARG{   {\wREDv} \ge \zero}{  \normLZ{ {\wREDv}} }, \hspace{0.2cm} \textrm{subject to }  \normGEN{\Afe^T \wREDv - \bEXACT}{2} \le \tolB   \normd{\bEXACT}
\end{equation}
\noindent where $\normLZ{\cdot}$ stands for  the $\ell_0$ pseudo-norm ---the number of nonzero entries of the vector. It is well-known \cite{boyd2004convex} that this problem is  computationally intractable (NP hard), and therefore, recourse to either suboptimal greedy heuristic or convexification  is to be made. Von Tycowicz et al  \cite{von2013efficient}  adapted the  algorithm proposed originally in Ref. \cite{blumensath2010normalized}  for compressed sensing applications (called \emph{normalized iterative hard thresholding}, abbreviated NIHT) by incorporating the positive constraints, reporting significant improvements in performance with respect to the original NNLS-based algorithm of An et al. \cite{an2009optimizing}. The work by Pan et al.    \cite{pan2015subspace}, on the other hand,  advocated an alternative approach ---also borrowed from the compressed sensing literature, see Ref. \cite{yang2011alternating} ---  based on the \emph{convexification} of problem \refpar{eq:10}. Such a  convexification consists in  replacing the  $\ell_0$ pseudo-norm     by the $\ell_1$ norm ---an idea that, in turn, goes back to the seminal paper by Chen et al. \cite{chen2001atomic}.  In doing so, the problem becomes tractable, and can be solved by standard  Linear Programming (LP) techniques.

Cubature schemes did not enter the  computational engineering   scene   until the appearance in 2014 of the \emph{Energy-Conserving Mesh Sampling and Weighting} (ECSW) \emph{scheme} proposed by C. Farhat  and co-workers \cite{farhat2014dimensional}.  The ECSW is, in essence, a   nonnegative least squares method (NNLS), very much aligned to the original proposal by An et al \cite{an2009optimizing}, although much more algorithmically efficient.  Indeed, Farhat and co-workers realized that the NNLS itself produces sparse approximations, and therefore it  suffices to introduce  a control-error parameter   inside the standard NNLS algorithm ---rather than invoking the NNLS at each greedy iteration, as proposed originally in An's paper \cite{an2009optimizing}. The efficiency of the   ECSW  was tested against other sparsity recovery algorithms by Farhat's team in Ref. \cite{chapman2016accelerated}, arriving at the conclusion that, if equipped with an updatable QR decomposition for calculating the unrestricted least-squares problem of each iteration, the ECSW outperformed existing implementations based on convexification of the original problem. It should be pointed out that, although the ECSW is a mesh sampling procedure, and therefore,    the entities selected by the ECSW are finite elements rather than single Gauss points,  the formulation of the problem   is  rather similar to the one described in the foregoing: the only   differences are that, firstly, each element contribution $\AfeE{e}$ in  \refeq{eq:8243t--} collapses into a single row obtained as the weighted sum of the Gauss points rows; and,  secondly,   the   vector of FE weights $\Wfe$ is replaced by an all-ones vector.

The Empirical Cubature Method, hereafter referred to as Discrete Empirical Cubature Method (DECM),  introduced by the first author  in Ref. \cite{hernandez2017dimensional} for  parametrized  finite element structural problems, also addresses the problem via a  greedy algorithm, in the spirit of An's approach \cite{an2009optimizing}, but   exploits the fact that deriving a cubature rule for  integrating the set of functions contained column-wise in matrix $\Afe$ is equivalent to deriving a cubature rule for a set of \emph{orthogonal bases} for such functions.  Ref. \cite{hernandez2017dimensional} demonstrates that this brings two salient advantages in the points selection process. Firstly, the   algorithm invariably converges to zero integration error when the number of selected points is equal to the number of orthogonal basis functions;  and secondly,  the algorithm need not   enforce the positiveness of the weights at each iteration. Furthermore, Ref. \cite{hernandez2017dimensional}  recognizes that the cubature problem is ill-posed when $\bFE \approx \zero$ ---this occurs, for instance,  in self-equilibrated structural problems, such as computational homogenization \cite{hernandez2014high,oliver2017reduced}---and shows that this can be   overcome by enforcing    the sum of the reduced weights to  be equal to the volume of the domain.    In Ref. \cite{hernandez2020multiscale}, the first author  proposed an  improved version of the original DECM, in which the local least-squares  problem at each iteration are solved  by rank-one updates.

 Another approach also introduced recently in the computational engineering community is the   \emph{Empirical Quadrature Method} (EQM), developed  by A. Patera and co-workers \cite{patera2017lp,yano2019lp,yano2019discontinuous}. It should be noted that  the name similarity with the above described Empirical Cubature Method    is only coincidental, for the EQM is not based on the  nonnegative least squares method, like the ECM, but rather  draws on the previously mentioned $\ell_1$ convexification of   problem \ref{eq:10}. Thus,  in the EQM,   the integration  rule is determined by linear programming techniques, as in   the method advocated in the work by Pan et al.  \cite{pan2015subspace} for computer graphics applications.

\subsection{Efficiency of best subset selection algorithms}
\label{sec:823ds}

\begin{figure}[!ht]
  \centering
  \includegraphics[width=0.99\textwidth]{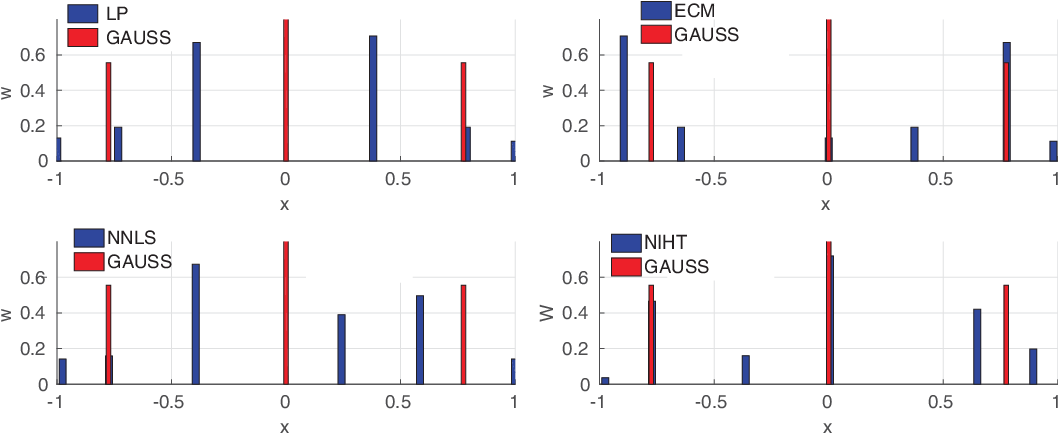}
  \caption{ Location of the points and magnitude of the weights of the integration rules for polynomial of order $\nSAMPL = 5$ in $\dom{}{} = [-1,1]$ provided by the: a) Linear programming-based   strategy (LP);  b) Discrete Empirical Cubature Method  (DECM);   c) Non-negative least-squares (NNLS).  d) Normalized iterative hard thresholding (NIHT).  Note that in the  4 cases, the number of integration points is equal to the number of functions to be integrated (i.e., monomials) $m = \nSAMPL +1 = 6$. The optimal solution is provided by the Gaussian quadrature rule of $m^* = 3$ points, also displayed in each of the four graphs.  \jahoHIDE{Generated by \href{/home/joaquin/Desktop/CURRENT_TASKS/PAPERS_2020_onwards/OptimalECM/DOCS/MATLAB/PaperFigura1.m}{PaperFigura1.m}}   }
     \label{fig:Figure1PB}
 \end{figure}

 The  \emph{best subset selection}    algorithms described in the foregoing vary in the   way the corresponding optimization problem is formulated, and   also in computational performance (depending on the nature and size of the problem under consideration), yet all of them have something in common: none of them is able to provide the optimal solution,  not even when   the optimal integration points are  contained in the set of FE Gauss  points. We have corroborated this claim by examining the number of points provided by all these methods when the integrand is a 1D  {polynomial} in the interval $\dom{}  = [-1,1]$. In   Figure \ref{fig:Figure1PB}, we show, for the case of polynomials  of order $\nSAMPL = 5$,  the  location of the points and the associated weights provided by:\footnote{The NNLS and LP analyses can be carried out by calling standard libraries (here we have used the Matlab functions \emph{lsqnonneg} and \emph{linprog}, respectively),   the ECM algorithm is given in Ref. \cite{hernandez2020multiscale}, whereas for the NIHT we have used the   codes  given in Ref. \cite{kim2016greedy}}  1) the nonnegative least-squares method (NNLS); 2) the linear-programming based method (LP); 3) the Discrete Empirical Cubature Method (DECM); 4) and the normalized iterative hard thresholding (NIHT). We also show in each Figure  the 3-points optimal Gauss rule, which in this case is   $\wRED{1}^*=\wRED{3}^* = 5/9$, $\wRED{2}^*=8/9$, and $x_1^*=-x_3^*=\sqrt{{3}/{5}}$ $x_2^* =  0$. The employed spatial discretization features $\nelem = 1000$ elements, with one Gauss point per element (located at the midpoint), and it was arranged in such a way that 3  of  the corresponding element midpoints  coincide with  the  optimal Gauss points.  It can be seen that, as asserted, none of the the four schemes is able to arrive at the optimal quadrature rule. Rather, the four methods provide  quadrature rules with   $m= \nSAMPL +1 =  6$ points, that is, with as many points as functions to be integrated; in the related literature,  these rules are known as \emph{interpolatory quadrature rules}   \cite{golub2009matrices}.  Different experiments with different  initial discretizations and/or polynomial orders led invariably to the same conclusion (i.e., all of them produce interpolatory rules).

\subsection{Goal and methodology}
\label{sec:goal}

  Having described the   capabilities and limitations of existing  cubature methods in the context of HROMs,   we now focus on the actual goal of the present paper, which is to enhance such methods so that they produce rules close to the optimal ones ---or at least rules featuring far less points than integrand functions. Our proposal in this respect  draws inspiration from   the \emph{elimination} algorithm   advocated,   apparently independently,  by  Bremer et al.  \cite{bremer2010nonlinear}   and Xiao et al. \cite{xiao2010numerical} in the context of the so-called \emph{Generalized Gaussian Rules} (see   Refs.  \cite{ma1996generalized}, \cite{mousavi2010generalized},  \cite{mousavi2011numerical}), which, as its name indicates,  is a  research discipline that seeks to extend the scope of  the quadrature rule originally developed for polynomials  by  C.F. Gauss.   To the best of the authors' knowledge, cross-fertilization between this field  and the field of hyperreduction of parameterized finite element models  has not  yet taken place.   This lack of cross-fertilization may be attributed to the fact that the former   is fundamentally concerned with  parametric families of    functions whose analytical expression is known, while the latter concentrates in huge databases of \emph{empirical} functions  (i.e., functions derived from   \emph{computational experiments}),  whose values are only given    at certain points of the spatial domain ( the Gauss points of the FE mesh).  The present work, thus, appears to be the first  attempt to   combine ideas from these  two  related disciplines. 
  
  The  intuition behind the    elimination algorithm presented in Refs. \cite{bremer2010nonlinear,xiao2010numerical}  goes as follows. Consider, for instance (the same arguments  can be used with either the LP or NIHT approaches  ),   the points and weights provided by the \emph{interpolatory} DECM rule  shown in Figure \ref{fig:Figure1PB}.b.  Observe that the distribution of weights   is rather irregular, being the   difference between the  largest and smallest weights   more pronounced than in the case of the optimal rule ---for instance, the smallest weight  is only 5 \% of the total length of the domain. This suggests  that we may  get rid of some of the points in the initial set, on the grounds that, as their  contribution to the integration error is not significant (relatively small weights), a slight ``readjustment'' of the positions and weights of the remaining points may suffice to return the integration error to zero.     Since we cannot know a priori how many points  in total can be eliminated, this operation must be carried out carefully, removing one point at a time.

\subsubsection{Sparsification problem}
\label{sec:sparsification}
 
Although inspired by this elimination scheme, our approach addresses the problem from a different perspective, more in line with the \emph{sparsification} formulation presented in expression  \refpar{eq:10}, in which the goal is to drive to zero as many weights as possible. 
To understand how our   {sparsification} scheme works, it proves useful to draw a  {physical analogy}  in which the integration points are regarded as \emph{particles} endowed with\emph{ nonnegative masses} (the weights), and which are subject to nonlinear conservation equations (the integration conditions).  At the beginning, the particles have  the positions and masses (all positive) determined by one of the interpolatory cubature rules discussed previously in Section \ref{sec:823ds}. The goal is to, progressively,  drive to zero the mass of as many particles as possible, while keeping the remaining particles within the spatial domain, and with nonnegative masses.   To this end, at each   step, we   reduce  the mass of the particle that  least contributes to the conserved quantities, and then calculate the position and masses of the remaining particles so that the nonlinear conservation equations are satisfied.  

For solving the nonlinear balance equations using standard methods (i.e., Newton's), it is necessary to have a     {continuous} (and differentiable) representation of the integrand functions. In contrast to the cases presented in Refs. \cite{bremer2010nonlinear} and \cite{xiao2010numerical}, in our case, the analytical expression of such functions are in general not available. To overcome this obstacle,  we propose to  construct local polynomial interpolatory functions using the values of the integrand functions at the Gauss points of each finite element traversed by the particles. 
 
  Another crucial difference of our approach with respect to Refs. \cite{bremer2010nonlinear,xiao2010numerical}  is the procedure to solve the  nonlinear equations at each step. Due to the underdetermination of such equations,   there are an infinite number of possible configurations of the system for the majority of the steps.   Both Refs. \cite{xiao2010numerical} and \cite{bremer2010nonlinear} use the pseudo-inverse of the Jacobian matrix, a fact that is equivalent to choosing the (non-sparse) $\ell_2$ minimum-norm solution \cite{boyd2004convex} in each iteration.  By contrast, here we employ  sparse solutions, with as many nonzero entries as functions to be integrated.  The rationale for employing this sparse solution is that, on the one hand,  it minimizes the number of particles that move at each iteration, and consequently, diminish the computational effort of tracking the particles through the mesh;  and, on the other hand,  it reduces the overall error inherent to the recovery of the integrand functions via interpolation.

  It should be stressed that we do not employ a specific strategy for directly enforcing the positiveness of the masses (weights). Rather,   we     force   the constant function to appear in the set of integrand functions; in our physical analogy,  this implies  that one of the balance equations    is the \emph{conservation of mass}.  Since the total mass of the system is to be conserved, reducing the mass of one particle leads to  an increase in the overall  mass of the remaining particles, and this tends to ensure that their masses remain positive.  On the other hand,   when  a particle attempts to leave the domain,   we return it back to its previous position, and proceed with the solution scheme. If convergence is not achieved, or the constraints are massively violated, we simply abandon our attempt of reducing the weight of the current controlled particle, and move to the next particle in the list. The process terminates ---hopefully at the optimum--- when we have tried to make zero the masses of all particles.

We choose  as initial interpolatory  rule  ---over the other methods discussed in Figure \ref{fig:Figure1PB}--- the  Discrete Empirical Cubature Method, DECM. The reason for this choice is twofold.    Firstly,  we have empirically found that   the DECM  gives  points that  tend to be close to the optimal ones ---for instance, in Figure \ref{fig:Figure1PB}.b two of the points calculated by the DECM practically coincide with the    optimal Gauss points $x_1 = \sqrt{3/5}$ and $x_2 =0$.  Secondly,  the DECM  does not  operate directly on the sampling matrix $\Afe$ defined in \refeq{eq:8243t--}, but rather on an  orthogonal basis matrix  for its column space \cite{hernandez2020multiscale}. As a consequence, the cubature problem translates into one of  integrating orthogonal basis functions, and this property  greatly facilitates  the convergence of the nonlinear problem alluded to earlier.   The combination of the DECM followed by the continuous search process will be referred to hereafter as the \emph{Continuous Empirical Cubature Method} (CECM).

\subsection{Sequential  randomized SVD (SRSVD)}
\label{sec:srsvd}

 We use  the ubiquitous Singular Value Decomposition (SVD) to determine the orthogonal basis matrix for the column space of $\Afe$ required by the DECM. Computationally speaking, the SVD is by far the most memory-intensive operation of the entire cubature algorithm. For instance, in a parametric function  of  dimension $\nU = 10$, with $\nSNAP =100$ parametric samples  and a mesh of $\nelem = 10^6$ linear hexahedra elements featuring $\ngaus = 8$ Gauss points each,   matrix  $\Afe$ occupies 64 Gbytes of RAM memory.
To overcome such potential memory bottlenecks,   we have devised a   scheme for computing the SVD in which the matrix is   provided into a column-partitioned format,  with the submatrices being processed one at a time.    In contrast to other partitioned schemes, such as the one proposed by the first author in Ref. \cite{hernandez2017dimensional} or the partitioned Proper Orthogonal Decomposition of Ref. \cite{wang2016approximate}, which compute the SVD of the entire matrix from the individual SVDs of each submatrices, our scheme    addresses the problem in an incremental, sequential fashion: at each increment, the current   basis matrix   (for the column space of $\Afe$) is enriched with the left singular vectors coming from the SVD of  the \emph{orthogonal complement}. The advantage of this sequential approach over the concurrent approaches in Refs. \cite{wang2016approximate,hernandez2017dimensional}  is that  it exploits the existence of linear correlations  accross the blocks. For instance, in a  case in which   all   submatrices are full rank, and besides, a linear combination of the first submatrix (this may happen when analyzing periodic functions), our sequential approach would require performing  a single SVD ---that of the first matrix. By contrast,    the concurrent approaches in Refs. \cite{wang2016approximate,hernandez2017dimensional}, would not only need to calculate the SVD of all the submatrices, but they would not provide any benefit at all in terms of computer memory  (in fact the partitioned scheme would end up being more costly than the standard one-block implementation). Lastly, to accelerate the performance of each SVD on the orthogonal complement of the submatrices, we employ a modified version of the  randomized blocked SVD proposed by Martinsson et al. \cite{martinsson2015randomized}, using as prediction for the rank of a given submatrix that of the previous submatrix in the sequence.

 \subsection{Organization of the paper}

 The paper is organized as follows.   The determination of the orthogonal basis functions and their gradients by using the SVD of the sampling matrix are discussed in Section \ref{sec:initial}.  
  Although   an original contribution of the present work, we have relegated the description of the Sequential Randomized SVD (SRSVD) algorithm to Appendix \ref{sec:SRSVD} (in order not to interrupt the continuity of the presentation of the cubature algorithm,  which constitutes the primary focus of this paper). On the other hand, the   computation of the interpolatory cubature rule by  the Discrete Empirical Cubature Method, DECM, is presented in Section \ref{sec:DECM1}, and the solution of the continuous sparsification problem in Sections \ref{sec:sparseP} and   \ref{sec:zeroingP}. Except for the   DECM, which can be found in the original reference \cite{hernandez2020multiscale}, we provide the pseudo-codes of all the algorithms involved in both the  cubature and the SRSVD. Likewise,  we have
summarized all the implementation steps in Box \ref{box:1} of Section \ref{sec:cecmLAST}.  The logic of the proposed
methodology can be followed without the finer details from the information in this  Box.   

Sections \ref{sec:uni} and \ref{sec:multo}  are devoted to the numerical validation by comparison with the (optimal) quadrature and cubature rules of  univariate and multivariate Lagrange polynomials. The example presented in Section \ref{sec:expsino}, on the other hand, is intended to illustrate the performance of the method in scenarios where the proposed SRSVD becomes essential ---because the integrand matrix exhausts the memory capabilities of the computer at hand.  Finally, the application of the proposed CECM to the   hyperreduction of a multiscale finite element model is explained   in Section \ref{sec:hyper}.

  \section{Orthogonal basis  for the integrand }
  \label{sec:initial} 
  
\subsection{Basis matrix via   SVD}
\label{sec:wSVD}


As pointed out in the foregoing, our cubature method does not operate directly on the integrand sampling matrix $\Afe$, defined in Eq. \refpar{eq:8243t--}, but on a    basis matrix for its column space, denoted henceforth by $\U \in \RRn{\ngausT}{\pNMOD}$. Since   $\U$  will be a linear combination of the columns of $\Afe$, which are in turn the discrete representation of the scalar integrand functions we wish to integrate, it follows  that  the columns of $\U$ themselves   will be the discrete representations of  basis functions for such integrand functions. These basis functions will be denoted hereafter  by   $u_i: \dom{} \rraa \Rn{}$  ($i=1,2 \ldots \pNMOD$).  In analogy to \refeq{eq:8243t--}, we  can write $\U$ in terms of such  basis functions as  
 \begin{equation}
     \label{eq:dalma1}
      \U = \colcuatro{\Ue{1}}{\Ue{2}}{\vdots}{\Ue{\nelem}}_{\ngausT \times \pNMOD} \hspace{0.25cm} \textrm{where}  \hspace{1cm}  \Ue{e} \defeq \colcuatro{\u({\xG{e}{1}})}{\u({\xG{e}{2}})}{\vdots}{\u({\xG{e}{\ngausE}})}_{\ngausE \times \pNMOD}. 
     \end{equation} 
     while 
     \begin{equation}
      \u(\x) \defeq  \rowcuatro{u_1(\x)}{u_2(\x)}{\cdots}{u_{\pNMOD}(\x)}_{1 \times \pNMOD}. 
     \end{equation}
We shall require these basis functions to be $L_2(\dom{})$-orthogonal, i.e.: 
\begin{equation}
 \intG{\dom{}}{\dom{}}{u_i u_j} =   \delta_{ij},  \hspace{0.5cm}  i,j = 1,2 \ldots \pNMOD,
\end{equation}
 $\delta_{ij}$ being the Kronecker delta. By evaluating the above integral using the FE-Gauss rule (as done in Eq \ref{eq:1}), we get: 
\begin{equation}
 \intG{\dom{}}{\dom{}}{u_i u_j} =  \sum_{e=1}^{\nelem} \sum_{g=1}^{\ngausE}  u_i(\xG{e}{g}) W_g^e u_j(\xG{e}{g}) =  \U_i^T \diag{(\Wfe)} \U_j,  \hspace{0.5cm}  i,j = 1,2 \ldots \pNMOD.
\end{equation}
In the preceding equation, $\U_i$ and $\U_j$ represents the $i$-th and $j$-th columns of $\U$, while $\diag{(\Wfe)}$ stands for a diagonal matrix containing the entries of the vector of FE weights $\Wfe$ ( defined in Eq. \ref{eq:WFE}). The above condition can be cast in a compact form as 
\begin{equation}
\label{eq:awebddddd}
 \U^T \diag{(\Wfe)} \U  = \ident,  
\end{equation}
$\ident$ being  the $\pNMOD \times \pNMOD$ identity matrix. The preceding equation reveals that orthogonality in the $L_2{(\Omega)}$ sense for the basis functions $u_i$   translates  into orthogonality for the columns of $\U$ in the sense defined by the following inner product 
\begin{equation}
\label{eq:56m}
 \innerG{\v_1}{\v_2}{W} \defeq \v_1^T \diag{(\Wfe)} \v_2
\end{equation}
($\v_1,\v2 \in \Rn{\ngausT}$). 
 
 In order to determine $\U$ from $\Afe$, we compute first  the (truncated) Singular Value Decomposition of the \emph{weighted}  matrix defined by  
 \begin{equation}
 \label{eq:82kdsdddd}
  \AfeW \defeq \diag{(\sqrt{\Wfe}) } \Afe
 \end{equation}
that is: 
 \begin{equation}
 \label{eq:firsta}
  \AfeW =  \Uw \S \V^T + \Ew,
 \end{equation}
 symbolyzed in what follows as the operation: 
 \begin{equation}
\label{eq:3qewf}
 [\Uw,\S,\V] = \SVD{\AfeW,\tolSVD}.
\end{equation}
Here,  $\Uw \in \RRn{\ngausT}{\pNMOD}$, $\S \in \RRn{\pNMOD}{\pNMOD}$ and $\V \in \RRn{\ndim \nSNAP}{\pNMOD}$ are the matrices of left-singular vectors, singular values and right-singular vectors, respectively.  The matrix of singular values is diagonal with $S_{ii} \ge S_{i-1,i-1}>0$, while the matrices of left-singular and right-singular vectors obey the orthogonality conditions
\begin{equation}
\label{eq:gavdddd}
 \Uw^T \Uw = \ident, \hspace{1cm}   \V^T \V = \ident. 
\end{equation}
Matrix $\Ew$ in \refeq{eq:firsta}, on the other hand, represents the  truncation term, which is controlled by a user-specified tolerance $0 \le \tolSVD \le 1$ such that 
\begin{equation}
\label{eq:adbcsssa}
 \normF{\Ew} \le \tolSVD \normF{\AfeW}
\end{equation}
(here $\normF{\bullet}$ denotes the Frobenius norm).   
The desired basis matrix $\U$ is computed from $\Uw$ as 
\begin{equation}
\label{eq:42dadd}
 \U = \diag(\sqrt{\Wfe})^{-1} \Uw. 
\end{equation}
It can be readily seen  that, in doing so, the $\Wfe$-orthogonality condition defined in \refeq{eq:awebddddd} is   satisfied. Multiplication of both sides of \refeq{eq:firsta} allows us to write 
 \begin{equation}
 \label{eq:first}
  \Afe =  \U \S \V^T + \E,
 \end{equation}
 where 
 \begin{equation}
 \label{eq:asas....}
  \E \defeq  \diag(\sqrt{\Wfe})^{-1} \Ew. 
 \end{equation}
 Notice that, by virtue of the definition of Frobenius norm, and by using the preceding expression,  we have that 
 \begin{equation}
  \normF{\Ew}^2  = \trace{\Ew^T \Ew} = \trace{\E^T \diag{(\Wfe)} \E} = \normGEN{\E}{W}^2,
 \end{equation}
where $\trace{\bullet}$ stands for the trace operator, and $\normGEN{\bullet}{W}$ designates the norm induced by the inner product introduced in Eq. \refpar{eq:56m}. Since the same reasoning can be applied to $\normF{\AfeW}$,  we can alternatively write the truncation condition  \refpar{eq:adbcsssa} as 
\begin{equation}
\label{eq:adbcsss}
 \normGEN{\E}{W} \le \tolSVD \normGEN{\Afe}{W}.
\end{equation}

\begin{remark}
 \label{rem:001}

When $\Afe$ is too large to be processed as a single matrix, we shall  use, rather than the standard SVD \refpar{eq:3qewf},  the sequential randomized SVD  alluded to in the introductory section \refpar{sec:srsvd}:
\begin{equation}
\label{eq:asdksss}
 [\Uw,\S,\V] = \SRSVD{[\AfeW_1,\AfeW_2, \ldots \AfeW_s],\tolSVD} 
\end{equation}
(here $[\Abar_1,\Abar_2, \ldots \Abar_s]$ stands for a partition of the weighted matrix $\AfeW$ ).  The implementation details  are provided in Algorithm \ref{alg:006} of Appendix \ref{sec:SRSVD}.  

\end{remark}

\subsection{Constant function}
\label{sec:constant_funct}

We argued in Section \ref{sec:goal} that the efficiency of the proposed search algorithm relies on one fundamental requirement: the volume of the domain is to be exactly integrated ---i.e., the sum of the cubature weights must be equal to the volume of the domain $V = \intG{\dom{}}{\dom{}}{}$. If the  integrand functions are provided as a collection of  analytical expressions, this can be achieved by  incorporating a constant function in such a collection, with the proviso that   the value for the  constant    should be sufficiently high so that the SVD regards the function as representative within the sample.   

The same reasoning applies when the      only data available is the empirical matrix $\Afe$: in this case, we may make  $\Afe \leftarrow [\Afe,c \ones]$,  where $\ones$ is an all-ones vector and  $c \in \Rn{}$ the aforementioned constant. Alternatively, to make the procedure less contingent upon the employed constant $c$, we may expand, rather than the original matrix $\Afe$, the basis matrix $\U$ itself.   To preserve   column-wise orthogonality, we proceed by first computing the component of  the all-ones vector   orthogonal to   the column space of $\U$  (with respect to the inner product \refpar{eq:56m}  ):   
\begin{equation}
 \v = \ones -\U \U^T \diag{(\Wfe)} \ones  = \ones -\U \U^T  \Wfe. 
\end{equation}
If $\normd{\v} \approx 0$, then no further operation is needed (the column space of $\U$ already contains the all-ones vector); otherwise, we    set $\v \leftarrow \v/\normd{\v}$, and expand $\U$ as $\U \leftarrow [\v,\U]$.  

 \begin{lemma}
 \label{lemma:1}
 If the column space of the basis matrix $\U$ contains the all-ones (constant)  vector, then 
 \begin{equation}
 \label{eq:a3wdddd}
  \E^T \Wfe = \zero, 
 \end{equation}
 that is,   the integrals of the functions whose discrete representation are the truncation matrix $\E$ in  the SVD  \refpar{eq:first}  are all zero. 
  
 \end{lemma}
 
 \begin{demoa}
 \label{proof:1}
  By construction, the truncation term $\E$ admits also a decomposition of the form $\E = \Uorth \Sorth \Vorth^T$, where $\innerG{\Uorth}{\Uorth}{W} = \ident$, $\innerG{\Uorth}{\U}{W} = \zero$ and $\Vorth^T \V = \ident$. Thus, replacing this decomposition into \refeq{eq:a3wdddd}, we arrive at  
   \begin{equation}
 \label{eq:a3wdddd.}
  \E^T \Wfe =\Vorth \Sorth (\Uorth^T \Wfe) = \zero.  
 \end{equation}
   The proof boils down thus to demonstrate that $\Uorth^T \Wfe = \zero$. This follows easily from the condition that  $\innerG{\Uorth}{\U}{W} = \zero$. Indeed, since the all-ones vector pertains to the column space of $\U$, the matrix of trailing modes $\Uorth$ is also orthogonal to the all-ones vector, hence
  \begin{equation}
   \innerG{\Uorth}{\ones}{W} = \zero \RRaam  \Uorth^T \diag{(\Wfe)} \ones = \Uorth^T  \Wfe = \zero. 
  \end{equation}

 \end{demoa}

\subsection{Evaluation of basis functions}
\label{sec:evalfun}

During the   weight-reduction process, it is necessary to repeatedly evaluate the basis functions, as well as their spatial gradient,  at points, in general, different from  the Gauss points of the mesh.

\subsubsection{Integrand given as analytical expression}
\label{sec:integrandANAL}

If the analytical expressions for the integrand functions $\A(\x)$ and their spatial derivatives $\derpar{\A(\x)}{x_i}$ ($i=1,2 \ldots \ndim$) are available, these evaluations can be readily performed by using the singular values and right-singular vectors of  decomposition \refpar{eq:first} as follows: 
 \begin{equation}
 \label{eq:follws1}
 \u(\x) = \A(\x) \V \S^{-1},
\end{equation}
and 
\begin{equation}
\label{eq:4wdasd}
 \derpar{\u(\x)}{x_i} = \derpar{\A(\x)}{x_i} \V \S^{-1}, \hspace{0.5cm}   i = 1\ldots \ndim. 
\end{equation}

\begin{demoa}
Post-multiplication of both sides of \refeq{eq:first} by $\V$ leads to 
\begin{equation}
 \Afe \V =  \U \S \V^T \V  + \E \V  = \U \S \V^T \V  + \Uorth \Sorth \Vorth^T \V,
\end{equation}
where we have used the matrices introduced in the proof of Lemma \refpar{lemma:1}. By virtue of the orthogonality conditions $\V^T \V = \ident$ and $\Vorth^T \V = \zero $, the above equation becomes  $\Afe \V =  \U \S$; postmultiplication by $\S^{-1}$ finally leads to $\U = \Afe \V \S^{-1}$. This equation  holds, not only  for $\Afe = \A(\Xfe)$, but for any $\A(\x)$, as stated in \refeq{eq:follws1}.

\end{demoa}

\begin{remark}
 \refeq{eq:4wdasd} indicates  that the gradient of the  $j$-th basis function   depends inversely on the $j$-th singular value. Negligible singular values, thus, may give rise to inordinately high gradients, causing convergence issues during the nonlinear readjustment problem. To avoid these numerical issues, the SVD truncation threshold $\tolSVD$ (see Expression \ref{eq:adbcsss}) should be set to a sufficiently large value (typically $\tolSVD \ge 10^{-6}$).
\end{remark}

 \subsubsection{Interpolation using Gauss points }
 \label{sec:InterpolationG}
 In general, however, the analytical expression for the integrand functions are not available, and therefore, the preceding equations cannot be employed for retrieving the values of the orthogonal basis functions. This is the case     encountered when dealing with FE-based reduced-order models, where the only information we have at our disposal  is the value of the basis functions at the Gauss points of the  finite elements, represented by submatrices $\Ue{e} \in \RRn{\ngausE}{\pNMOD}$ ($e=1,2 \ldots \nelem$) in Eq. \refpar{eq:dalma1}. 
 
 In a FE-based reduced-order model,  at element level, the integrand functions are, in general, a nonlinear function of the employed nodal shape functions.  It appears reasonable, thus,  to use also polynomial interpolatory functions to estimate the values of the basis functions at other points of the element     using as interpolatory points, rather than the nodes of the element,  their Gauss points. If we denote by $\NshapeG{e}: \dom{e} \rraa \RRn{1}{r}$   the $\ngausE$ interpolatory functions (arranged as a row matrix), then we can write  
   \begin{equation}
   \label{eq:42ewe}
   \u(\x) = \NshapeG{e}(\x) \Ue{e}  , \hspace{1cm} \x  \in \dom{e} .
  \end{equation}
 Likewise,   the spatial derivatives can be determined as
    \begin{equation}
    \label{eq:72,dwwww}
   \derpar{\u(\x)}{x_i} =  \BshapeG{e}{i}\!(\x) \Ue{e}, \hspace{1cm} \x \in \dom{e},  \;\;\; i = 1 \ldots \ndim
  \end{equation}
  where  
  \begin{equation}
  \label{eq:w3gwe}
    \BshapeG{e}{i} \defeq \derpar{\NshapeG{e}}{x_i},  \;\;\; i = 1 \ldots \ndim. 
  \end{equation}
  
 The level of accuracy of this estimation will depend  on the number of Gauss point per element with respect to  the order  of  the nodal shape functions, as well as  the distorsion of the physical domain $\dom{e}$ with respect to the parent domain ---which is the cause of the aforementioned nonlinearity.  It may be argued that if the element has no distorsion, the evaluation of the integrand via \refeq{eq:42ewe} will be exact if the proper number of Gauss points is used. For instance,  in a  small-strains structural problem, if the element is a 4-noded bilinear quadrilateral, with no distorsion (i.e., a rectangle), and the term to be integrated is the virtual internal work, then the integrand is represented exactly by a quadratic\footnote{Virtual work is the product of virtual strains (which are linear in a 4-noded rectangular ) and stresses (which are therefore also linear)} polynomial. Such a polynomial possesses $(2+1)^2 = 9$ monomials, and therefore a $3 \times 3$ Gauss rule would   be needed.   Notice that this is an  element integration rule with one point  more per  spatial direction than the integration standard rule   for  bilinear quadrilateral elements ($2 \times 2$).

  The expression for $\NshapeG{e}$ and $\BshapeG{e}{i}$ in terms of the coordinates of the Gauss points $\{ \xG{e}{1},\xG{e}{2} \ldots \xG{e}{\ngausE} \}$ can be  obtained by the standard procedure used in  deriving FE shape functions, see e.g. Ref. \cite{liu2003finite}. Firstly,   we introduce the mapping $\phiMAP{e}: \dom{e} \rraa \domPAR{e}$   defined by
     \begin{equation}
     \label{eq:3weasd}
    [\xPAR]_i  = \dfrac{1}{L_i} [\x - \xG{e}{0}]_i, \hspace{1cm} i = 1 \ldots \ndim 
   \end{equation}
   where $[\bullet]_i$ symbolyzes   the $i$-th   component of the argument,  $\xG{e}{0} = (\sum_{g=1}^{\ngausE} \xG{e}{g})/\ngausE$ is the centroid of the Gauss points, and $L_i$ is a scaling length defined by
   \begin{equation}
    L_i = \maxARG{g=1\ldots \ngausE}{( \absval{[\xG{e}{g} - \xG{e}{0}]_i})}, \hspace{1cm} i=1 \ldots \ndim.
   \end{equation}
   The expression of the shape functions in terms of the  scaled   positions of the Gauss points  $\xGparALL = \{\xGpar{1},\xGpar{2} \ldots \xGpar{\ngausE} \}$, where $\xGpar{g} = \phiMAP{e}(\xG{e}{g})$, is given by
   \begin{equation}
   \label{eq:8932e}
    \NshapeG{e}(\xPAR) = \Pmon(\xPAR) \Pmon^{-1}\!(\xGparALL). 
  \end{equation}
  Here,  $\Pmon(\xPAR) \in \RRn{1}{\ngausE}$ is the row matrix containing the monomials up to the   order corresponding to the number and distribution of Gauss points  at point $\xPAR = \phiMAP{e}(\x)$; for instance,   for the case of a 2D $q \times q$ rule, where $\ngausE = q^2$, this row matrix adopts the form
   \begin{equation}
    \P(\xPAR) \defeq [\xxPAR_1^0 \xxPAR_2^0, \xxPAR_1^1 \xxPAR_2^0 \cdots \xxPAR_1^i \xxPAR_2^j \cdots \xxPAR_1^{q-1} \xxPAR_2^{q-1}]_{1 \times \ngausE}. 
   \end{equation}
   The other matrix appearing in \refeq{eq:8932e}, $\Pmon(\xGparALL) \in \RRn{\ngausE}{\ngausE}$, known as the \emph{moment matrix} \cite{liu2003finite}, is formed by stacking the result of applying the preceding mapping to the set of scaled Gauss points $\xGparALL$.   Provided that the element is not overly distorted (no negative Jacobians in the original isoparameteric transformation), the invertibility of $\Pmon(\xGparALL)$ is guaranteed thanks to the   coordinate transformation \refeq{eq:3weasd} ---which ensures that the   coordinates of all  points range between -1 and 1, therefore avoiding scaling issues in the inversion.

    As for the gradient of the shape functions in \refeq{eq:w3gwe}, by applying the chain rule, we get that
   \begin{equation}
   \label{eq:kommon1}
    \BshapeG{e}{i}    = \derpar{\P(\xPAR)}{\xxPAR_j} \derpar{\xxPAR_j}{x_i} \Pmon^{-1}\!(\xGparALL)  = \dfrac{1}{L_i}  \derpar{\P(\xPAR)}{\xxPAR_i} \Pmon^{-1}\!(\xGparALL). 
   \end{equation}

  \section{Discrete Empirical Cubature Method (DECM)}
  \label{sec:DECM1}

Once  the orthogonal basis matrix $\U$ has been computed  by the weighted SVD outlined in Section \ref{sec:wSVD},   the next step consists in determining an \emph{interpolatory cubature rule } (featuring as many points as   functions to be integrated) for the basis functions $\u: \dom{} \rraa \RRn{1}{\pNMOD}$. As pointed out in Section \ref{sec:goal}, we   employ for this purpose the  Empirical Cubature Method, proposed by the first author in Ref. \cite{hernandez2017dimensional}, and further refined in Ref.  \cite{hernandez2020multiscale}. We call it here \emph{Discrete} Empirical Cubature Method, DECM, to emphasize that the cubature points are selected among the Gauss points of the mesh.  The DECM,   symbolized in whats follows as the operation
 \begin{equation}
 \label{eq:2dasdfsad}
  [\zDECM,\wDECM] \leftarrow \textrm{DECM}(\U,\Wfe),
 \end{equation} 
takes as inputs  the basis matrix $\U \in \RRn{\ngausT}{\pNMOD}$ and the vector of positive FE weights $\Wfe \in \Rn{\ngausT}$; and returns  a set of $\pNMOD$ indexes $\zDECM \subset \{1,2 \ldots \ngausT\}$ and a vector of \emph{positive} weights $\wDECM$ such that 
\begin{equation}
\label{eq:4sds}
 \U(\z,:)  \wDECM   = \b. 
\end{equation}
Here,  $\U(\z,:)$  denotes, in the so-called ``colon'' notation \cite{golub2012matrix} (the one used by Matlab), the submatrix of $\U$ formed by the  rows corresponding to indexes $\z$,  while $\b \in \Rn{\pNMOD}$ is the vector of ``exact'' integrals  of the basis functions, that is: 
\begin{equation}
\label{eq:sdas---}
  \b = \intG{\dom{}}{\dom{}}{\u^T} = \U^T \Wfe. 
 \end{equation}

 The points associated to the selected rows will be denoted hereafter by  $\Xdecm = \{\xDECM{1},\xDECM{2} \ldots \xDECM{\pNMOD} \}$ (  $\Xdecm \subset \Xfe$). Hence, according to the notational  convention introduced in Remark \ref{rem:iwe,sdfsddds}, $\Pbool{\zDECM} \U$ may be alternatively  expressed as 
 \begin{equation}
 \label{eq:emploued1}
  \Pbool{\zDECM} \U =  \u(\Xdecm) =  \colcuatro{\u(\xDECM{1})}{\u(\xDECM{2})}{\vdots}{\u(\xDECM{\pNMOD})}  = \begin{bmatrix}
                                                                           u_1(\xDECM{1}) & u_2(\xDECM{1}) & \cdots & u_\pNMOD(\xDECM{1})  \\
                                                                              u_1(\xDECM{2}) & u_2(\xDECM{2}) & \cdots & u_\pNMOD(\xDECM{2})  \\
                                                                             \vdots & \vdots & \ddots & \vdots \\ 
                                                                                u_1(\xDECM{\pNMOD}) & u_2(\xDECM{\pNMOD}) & \cdots & u_\pNMOD(\xDECM{\pNMOD})  \\
                                                                           \end{bmatrix}_{\pNMOD \times \pNMOD}
\end{equation}

\begin{figure}[!ht]
  \centering
   \input{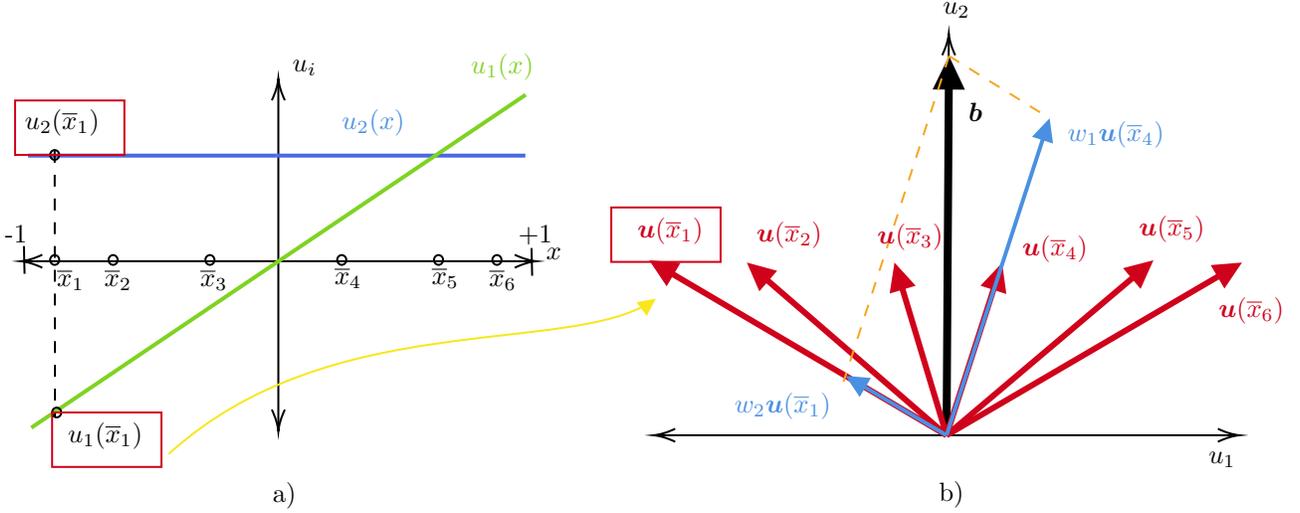}
   \caption{Performance of the Discrete Empirical Cubature Method (DECM) \cite{hernandez2020multiscale}, for the case of the  orthonormal functions $u_1 = \sqrt{3/2} x$ and  $u_2 = \sqrt{1/2}$   in the interval $\dom{} = [-1,1]$. The number of finite elements is $\nelem = 1$, and the number of Gauss points $\ngausT = \ngausE = 6$. a) Graph of the functions along with the location of the Gauss points $\xG{}{6} = -\xG{}{1} = 0.9325$ ,$\xG{}{5} = -  \xG{}{2} = 0.6612$, $\xG{}{4}  = -\xG{}{3}  = 0.2386$.   b) Representation in the plane $u_1u_2$ of each vector $\u(\xG{}{g}) = [u_1(\xG{}{g}),u_2(\xG{}{g})]^T$ ($g=1,2 \ldots 6$), along with the vector of exact integrals $\b = [0,{\sqrt{2}}]^T$.    The DECM chooses points $\xG{}{4}$ and $\xG{}{1}$, giving weights equal to $\wwDECM{1} = 1.5925$ and $\wwDECM{2} = 0.4075$.    }
     \label{fig:Figure2PB}
 \end{figure}

\begin{remark}
\label{remark:1}

It should be stressed that the solution to problem \refeq{eq:4sds} is  not unique. Rather, the number of possible solutions grows combinatorially with the ratio between the  total number of Gauss points and the number of functions  ($\ngausT/\pNMOD$).  The situation is illustrated in Figure \ref{fig:Figure2PB}, where we graphically explain how the DECM works for the case of $\ngausT = 6 $ Gauss points and  polynomial functions up to order 1 (it can  be readily shown that the orthogonal functions in this case are  $u_1 = \sqrt{3/2} x$ and $u_2 = \sqrt{1/2}$, displayed in Figure \ref{fig:Figure2PB}.a). The problem, thus,  boils down   to select $\pNMOD = 2$ points out of $\ngausT = 6$, such that the resulting weights are positive. In Figure \ref{fig:Figure2PB}.b, we plot each   $\u(\xG{}{g})$ $(g = 1,2 \ldots 6)$ along with  the vector of exact integrals, which in this case is equal to $\b = [0,\sqrt{2}]^T $.    It follows from this representation that, out of the ${\ngausT \choose \pNMOD} = {6 \choose 2} = 15$ possible combinations, only  9 pairs are valid solutions. The DECM\footnote{The first vector $\u(\xG{}{4})$ is chosen because is the one which is most positively parallel to $\b$ (notice that, because of symmetry,  it might have chosen $\u(\xG{}{3})$ as well). Then it orthogonally projects $\b$ onto $\u(\xG{}{4})$, giving $\s$, and then search for the vector which is more positively parallel to the residual $\b - \s$, which in this case is $\u(\xG{}{1})$.   } chooses   $\xDECM{1}= \xG{}{4}$ and $\xDECM{2} = \xG{}{1}$, which is the solution that yields the largest ratio between highest and lowest weight. Other possible solutions are, for instance,  pairs $\{\xG{}{1}, \xG{}{6}\}$  and $\{\xG{}{2},\xG{}{6}\}$ ---observe that in both cases, vector $\b$ lies in the cone\footnote{The cone positively spanned by a set of vectors is the set of all possible positive linear combinations of such vectors \cite{davis1954theory}. } ``positively spanned'' by $\{\u(\xG{}{1}), \u(\xG{}{6})\}$ and $\{\u(\xG{}{2}),\u(\xG{}{6})\}$, respectively.    
\end{remark}

 The reader interested in the points selection algorithm behind the  DECM is referred\footnote{It should be pointed that the notation employed in Ref. \cite{hernandez2020multiscale} is different from the one used here. 
 The input of the DECM in Ref \cite{hernandez2020multiscale} (which is called therein simply the \emph{ECM}) is not $\U$, but the transpose of the weighted matrix $\Uw$ (defined in Eq. \ref{eq:firsta}). Likewise the error threshold appearing in Ref. \cite{hernandez2020multiscale}  is to be set to zero in order to produce an interpolatory rule.  } to Ref. \cite{hernandez2020multiscale}, Appendix A, Algorithm 7.

\subsection{Relation between SVD  truncation error and the DECM integration error}
\label{eq:ressss--}

Let us examine now the error incurred in approximating the ``exact'' integrals $\bFE = \Afe^T \Wfe$ by the DECM cubature rule.  This error may be expressed as 
\begin{equation}
 \eDECM \defeq \normd{ \Afe^T \Wdecm  - \Afe^T \Wfe}
\end{equation}
 where $\Wdecm \defeq \Pbool{\zDECM}^T \wDECM$ (a vector of the same length as  $\Wfe$, but with nonzero entries only at the indexes specified by $\zDECM$). Inserting decomposition \refpar{eq:first} in the preceding equation, we get 
\begin{equation}
\label{eq:thus}
 \eDECM = \normd{  \Par{\V \S (\U^T \Wdecm  - \U^T \Wfe)} + (\E^T \Wdecm -\E^T \Wfe  ) }
\end{equation}
From condition \refpar{eq:4sds}, it follows that the term involving the basis matrix $\U$   vanishes; besides,  since by  construction the column space of $\U$ contains the all-ones vector, we have that, by virtue of Lemma \refpar{lemma:1}, $\E^T \Wfe = \zero$. Thus,   \refeq{eq:thus} boils down to 
 \begin{equation}
 \label{eq:mu3ewd}
 \eDECM = \normd{     \E^T \Wdecm  }. 
\end{equation}

The truncation term $\E$ in the above condition is controlled by the SVD tolerance $\tolSVD$ appearing   in  inequality \ref{eq:adbcsss}. Thus, the integration error $\eDECM$ may be lowered to any  desired level by decreasing the SVD tolerance $\tolSVD$. Numerical experience shows that for most problems $\eDECM/ \normd{\bFE}$  is slightly above, but of the same order of magnitude, as  $\tolSVD$.

\jahoHIDE{To gain further insight, go to Matlab's file: }

\jahoHIDE{\input{FutileAttemptsBoundAbove}}

  \section{Global sparsification problem}
  \label{sec:sparseP}
  \subsection{Formulation}

We now concentrate our attention on the  \emph{sparsification} problem outlined   in Section \ref{sec:sparsification}. The design variables in this optimization problem will be a vector of $\pNMOD$   weights   
\begin{equation}
  \w = [w_1,w_2 \ldots w_\pNMOD]^T, \hspace{1cm}  w_g \ge 0,
\end{equation}
(  recall that $\pNMOD$ is the number of orthogonal basis functions we wish to integrate), and  the position of the associated points within the domain: 
\begin{equation}
 \X = \{\x_1,\x_2 \ldots \x_{\pNMOD} \},  \hspace{1cm}  \x_g \in \dom{}. 
\end{equation}
With a minor abuse of notation, we shall also use  $\X$  to denote the    variable formed by stacking   the position of the points into a column matrix, i.e.: $\X = [\x_1^T,\x_2^T \ldots \x_{\pNMOD}]^T$.  On the other hand, we define the \emph{integration residual} as 
\begin{equation}
 \r \defeq \u^T\!(\X) \w - \b
\end{equation}
that is, as the difference between the approximate and the exact integrals of the basis functions.  In the preceding equation,     $\u(\X)$ designates the matrix formed by stacking the rows of all $\u(\x_g) \in \RRn{1}{\pNMOD}$ ($g=1,2 \ldots \pNMOD $) into a single matrix, i.e.: 
\begin{equation}
 \u(\X) \defeq  \colcuatro{\u(\x_1)}{\u(\x_2)}{\vdots}{\u(\x_{\pNMOD})}  = \begin{bmatrix}
                                                                           u_1(\x_1) & u_2(\x_1) & \cdots & u_\pNMOD(\x_1)  \\
                                                                              u_1(\x_2) & u_2(\x_2) & \cdots & u_\pNMOD(\x_2)  \\
                                                                             \vdots & \vdots & \ddots & \vdots \\ 
                                                                                u_1(\x_\pNMOD) & u_2(\x_\pNMOD) & \cdots & u_\pNMOD(\x_\pNMOD)  \\
                                                                           \end{bmatrix}_{\pNMOD \times \pNMOD}
\end{equation}
With the preceding definitions at hand, the sparsification problem can be formulated as  follows:  
\begin{equation}
\label{eq:spars1}
\begin{aligned}
\min_{ \w, \X } \quad & {  \normLZ{ {\w}} }\\
\textrm{s.t.} \quad &  \r = \u^T\!(\X) \w - \b = \zero \\
  &   \w \ge \zero    \\
  &   \X \subset \dom{}
\end{aligned}
\end{equation}

Recall that $ \normLZ{ {\w}}$ stands for the number of nonzero entries of $\w$. Thus, the goal in the preceding optimization problem  is to find the \emph{sparsest} vector of \emph{positive} weights, along with their associated positions within the domain, that render the \emph{integration residual}  $\r$  equal to zero.

\begin{remark}
 The differences between  the    sparsification problem presented above  and the one described in the introductory section (see Problem \ref{eq:10}) are three. Firstly, in the preceding problem, $\w \in \Rn{\pNMOD}$, i.e., the number of weights   is equal to the number of basis functions to be integrated $\pNMOD$, whereas in Problem \ref{eq:10}, this number is equal to the total number of FE Gauss points $\ngausT$ (it is assumed that $\pNMOD << \ngausT$).  Secondly, the integration residual in Problem \ref{eq:spars1}     appears in the form of an   equality constraint, while in Problem \ref{eq:10}   appears as an inequality constraint.   And thirdly, and most importantly, in Problem \ref{eq:spars1}, the positions of the cubature points are considered design variables ---in constrat to the situation encountered in Problem \ref{eq:10}, in which the points are forced to coincide with the FE Gauss points, and thus,  the only design variables are the weights.     
\end{remark}



%
%

\subsection{Proposed sparsification algorithm}

   \begin{algorithm}[!ht]
   \DontPrintSemicolon
   
   \SetKwFunction{FMain}{SPARSIFglo}
    \SetKwProg{Fn}{Function}{}{}
                  \SetKwFunction{ceil}{ceil}

       
      \SetKwFunction{LENGTH}{length}
            \SetKwFunction{SSPARSIF}{SPARSIF}

            \SetKwFunction{SORT}{sort}
            \SetKwFunction{ANYOUT}{ANYOUT}
   \SetKwFunction{RMVPNT}{RMV1PNT}

       \SetKwFunction{FUNGRAD}{FUNGRAD}
%

\newcommand{\Xold}{\boldsymbol{\bar{X}}}
\newcommand{\xOLD}{\boldsymbol{\bar{x}}}

\newcommand{\wOLD}{\boldsymbol{\bar{w}}}
\newcommand{\wOLDe}{\bar{w}}
\newcommand{\iORDER}{\boldsymbol{i}}

\newcommand{\Eall}{\boldsymbol{\mathcal{E}}}
\newcommand{\MSH}{\boldsymbol{\mathcal{M}}}
        \newcommand{\CONV}{{\mathcal{C}}}
        \newcommand{\BREAK}{\textrm{{break}}}

    \Fn{ $[\X,\w,\ElemINFO]$  $\leftarrow$ \FMain{${\Xdecm}$,$\wDECM$,$\Nsteps$,$\AuxVAR$,$\MeshINFO$} }{
    
    \KwData{ $\Xdecm = \{\xDECM{1},\xDECM{2} \ldots \xDECM{\pNMOD} \} (\xDECM{g} \in \dom{}$) and $\wDECM = [{\wwDECM{1}},{\wwDECM{2}},{\ldots},{\wwDECM{\pNMOD}}]$  ($\wwDECM{g} >0$): DECM cubature rule (obtained by function \refpar{eq:2dasdfsad}).  $\Nsteps >1$:   number of steps  used to solve the nonlinear   problem (in the second stage) associated to the residual constraint   $\r = \zero$.  $\AuxVAR$: Remaining variables controlling the solution of this nonlinear problem. $\MeshINFO$: data structures containing   variables needed to evaluate the residual $\r$ at any point $\x \in \dom{}$  (such as  the basis matrix itself $\U$ and the vector of exact integrals   $\b$, among others).  }    
     
    \KwResult{   \emph{Sparsest} weight vector $\w = [{w_{1}},{w_{2}},{\ldots},{w_{\pNMOD}}]$  ($w_{g} \ge 0$) and associated positions $\X = \{\x_{1},\x_{2} \ldots \x_{\pNMOD} \} (\x_{g} \in \dom{}$)  such that $\r = \u^T(\X) \w -\b = \zero$.  }  
    
      $\ElemINFO^{old} \leftarrow \emptyset$   \tcp{Initializations}    
      $N \leftarrow 1$  \tcp{First stage. Number of steps employing in solving the nonlinear problem $\r = 0$   is set to 1.}

        $[\X^{new},\w^{new},\ElemINFO^{new}]$  $\leftarrow$ \SSPARSIF{${\Xdecm}$,$\wDECM$,$N$,$\AuxVAR$,$\MeshINFO$,$\ElemINFO^{old}$} \tcp{\SSPARSIF{} is described in Algorithm \ref{alg:002}  \label{alg:001_a} }

           $[\X,\w,\ElemINFO]$  $\leftarrow$ \SSPARSIF{${\X^{new}}$,$\w^{new}$,$\Nsteps$,$\AuxVAR$,$\MeshINFO$,$\ElemINFO^{new}$}  \tcp{Second stage.  Number of steps  equal to the value specified in the input arguments (typically $\Nsteps \sim 20$)  \label{alg:001_b}}

       }

   \caption{ Proposed two-stage  procedure for solving  the  sparsification problem \refpar{eq:spars1}.   \label{alg:001}} 
   \end{algorithm}


The proposed approach for arriving at the solution of the preceding problem is to construct a sequence of $\pNMOD$-points cubature  rules  
\begin{equation}
 \{\X^{0},\w^{0} \}, \{\X^{1},\w^{1} \} \ldots \{\X^{i},\w^{i} \} \ldots \{\X^{\pNMOD - \mRED},\w^{\pNMOD-\mRED} \} 
\end{equation}
such that  
 \begin{equation}
  \normLZ{\w^{k+1}} =   \normLZ{\w^{k}} -1  
 \end{equation}
 that is, such that  each weight vector in the sequence has one non-zero less than the previous one. The first element in the sequence will be taken as the     cubature rule provided by the DECM (see Section \ref{sec:DECM1}):
 \begin{equation}
  \X^0 = \Xdecm, \hspace{1cm} \w^0 = \wDECM. 
 \end{equation}
 The algorithm proceeds from this initial point by the recursive  application of  an operation consisting in driving the weight of one single point to zero (while forcing the remaining points and weights to obey the constraints appearing in Problem \ref{eq:spars1}); this step will be symbolized hereafter as the function 
 \begin{equation}
 \label{eq:nccasdsd}
   [\CboolCONV,\X^{new},\w^{new},\ElemINFO^{new}] \leftarrow  \MakeOneZero{\X^{old},\w^{old},\Nsteps,\AuxVAR,\MeshINFO,\ElemINFO^{old}}.
 \end{equation}
This function takes as inputs a given cubature rule $\{\X^{old},\w^{old}\}$, and tries to  return a cubature rule $\{\X^{new},\w^{new}\}$ with at least one less nonzero weight. The success of this operation is indicated by the output Boolean variable $\CboolCONV$ ($\CboolCONV=false$ if it fails in producing a sparser cubature rule).

The other inputs in \refpar{eq:nccasdsd} are the following:  1) $\Nsteps$:   number of steps  used to solve the nonlinear   problem associated to the residual constraint $\r = \zero$. 2) $\AuxVAR$: Remaining variables controlling the solution of this nonlinear problem (such as the convergence tolerance for the residual).
3) $\MeshINFO$ and $\ElemINFO^{old}$ are data structures containing the variables needed to evaluate the residual $\r$ at any point $\x \in \dom{}$. $\MeshINFO$ encompasses those variables that do not change   during the execution (such as  the basis matrix itself $\U$, the vector of exact integrals   $\b$, see Eq. \refpar{eq:sdas---}, the    nodal coordinates of the FE mesh, the  connectivity table, the Gauss coordinates,   and in the case of analytical evaluation, the  product $\V_{\! s} \defeq \V \S^{-1}$ appearing in Eqs. \ref{eq:follws1} and Eqs. \ref{eq:4wdasd}).  $\ElemINFO^{old}$, on the other hand, comprises  element variables that are computed   \emph{on demand}\footnote{In the proposed  algorithm, we compute the necessary interpolation variables for each element of the mesh dynamically, as they are needed, rather than precomputing them all at once.  Indeed, each time the position of the points is updated, we check which elements of the mesh contain the updated points.  If all the elements have been previously visited,  we use the  information stored in $\ElemINFO^{old}$ to perform the interpolation; otherwise, we compute the required interpolation variables for the new elements and update  $\ElemINFO^{old}$ into  $\ElemINFO^{new}$  with the new data. }, such as the the inverse of the moment matrix in Eq. \ref{eq:8932e} and the scaling factors in Eq. \ref{eq:3weasd} (required for the interpolation described in Section \ref{sec:InterpolationG}).

   \begin{algorithm}[!ht]
   \DontPrintSemicolon
   
   \SetKwFunction{FMain}{SPARSIF}
    \SetKwProg{Fn}{Function}{}{}
                  \SetKwFunction{ceil}{ceil}

       
      \SetKwFunction{LENGTH}{length}
            \SetKwFunction{SORT}{sort}
            \SetKwFunction{ANYOUT}{ANYOUT}
   \SetKwFunction{RMVPNT}{RMV1PNT}

       \SetKwFunction{FUNGRAD}{FUNGRAD}
%

\newcommand{\Xold}{\boldsymbol{\bar{X}}}
\newcommand{\xOLD}{\boldsymbol{\bar{x}}}

\newcommand{\wOLD}{\boldsymbol{\bar{w}}}
\newcommand{\wOLDe}{\bar{w}}
\newcommand{\iORDER}{\boldsymbol{i}}

\newcommand{\Eall}{\boldsymbol{\mathcal{E}}}
\newcommand{\MSH}{\boldsymbol{\mathcal{M}}}
        \newcommand{\CONV}{{\mathcal{C}}}
        \newcommand{\BREAK}{\textrm{{break}}}
                \newcommand{\RETURN}{\textrm{\textbf{return}}}

    \Fn{ $[\X,\w,\ElemINFO^{new}]$  $\leftarrow$ \FMain{${\X^{old}}$,$\w^{old}$,$N$,$\AuxVAR$,$\MeshINFO$,$\ElemINFO^{old}$} }{
    
    \KwData{ $\X^{old}  \subset \dom{}$,   $\w^{old} \in  \Rn{p}$. $N$ and $\AuxVAR$:   Variables controlling the solution of the nonlinear problem $\r = \zero$. $\MeshINFO$ and  $\ElemINFO^{old}$: data structures containing   variables needed to evaluate the residual $\r$ at any point $\x \in \dom{}$.    }    
     
    \KwResult{   \emph{Sparsest} weight vector $\w \in \Rn{p}$  and associated positions $\X$   such that $\r = \u^T(\X) \w -\b = \zero$.  $\ElemINFO^{new}$: updated structure data with information necessary for performing element interpolation.  }  
    
    $ \CboolCONV \leftarrow true $

      \While{$ \CboolCONV = true $}{

      $[\CboolCONV,\X^{new},\w^{new},\ElemINFO^{new}] \leftarrow   \MakeOneZero{\X^{old},\w^{old},N,\AuxVAR,\MeshINFO,\ElemINFO^{old}}$ \tcp{Described   in Algorithm \ref{alg:003}.      \label{alg:002_a}  }

        \lIf{ $\CboolCONV$ = $\FALSE$  }{$\RETURN$ \tcp{ $\MakeOneZero{}$ has failed to produce a cubature rule with one nonzero weight less.   \label{alg:002_zero}  }} 
        \vspace{-0.2cm}
                 \lIf{ $w^{new}_i \ge \zero$, $\forall i$  }{ $\X \leftarrow \X^{new}$ ;  $\w \leftarrow \w^{new}$ \tcp{The  weights solution  $\w$ can only have positive weights \label{alg:002_b}  }  }      
          \vspace{-0.2cm}
        $ \X^{old} \leftarrow \X^{new}$; $\w^{old} \leftarrow \w^{new}$; $\ElemINFO^{old} \leftarrow \ElemINFO^{new}$   \tcp{Update positions, weights and element interpolation data.}        
       }

       }
  
   \caption{ Sparsification process, given an initial cubature rule $\{\X^{old},\w^{old}\}$ and a number of steps $N$ (  invoked in Line \ref{alg:001_b} of Algorithm \ref{alg:001})  \label{alg:002}. } 
   \end{algorithm}


Due to its greedy or ``myopic'' character,    the DECM  tends to produce weights distribution in  which most of the weights are relatively small in comparison with the total volume of the domain. We have empirically observed that the      readjustment problem associated to the elimination of these   small weights is moderately nonlinear, and in general, one step suffices to ensure convergence. However, as the algorithm advances in the sparsification process, the weights to be zeroed become larger, and, as a consequence,  the  readjustment problem becomes more nonlinear.  In this case, to ensure convergence, it is necessary to   reduce the weights progressively.  To account for this fact, we have devised the two-stage procedure described  in Algorithm \ref{alg:001}. In the first stage (see Line \ref{alg:001_a}), the sparsification process (sketched  in turn in Algorithm \ref{alg:002}) is carried out by decreasing the weight of each chosen weight in one single step. In the second stage, see Line \ref{alg:001_b}, we take the cubature rule produced in the first stage, and try to further decrease the number of nonzero weights by using a higher number of steps $\Nsteps > 1$.

 \section{Local sparsification problem}
   \label{sec:zeroingP}
    
After presenting the global sparsification procedure, we now focus on fleshing out  the details of the fundamental building block of such a procedure, which is the above mentioned function $\MakeOneZero{}$, appearing in Line \ref{alg:002_a} of  Algorithm  \ref{alg:002}. 

The procedural steps are  described in the pseudocode of Algorithm \ref{alg:003}. Given a cubature rule $\{\X^{old},\w^{old}\}$, with $\normLZ{\w^{old}} = m$   ($2 \le m \le \pNMOD$), we seek   a new cubature rule $\{\X,\w\}$ with $\normLZ{\w} = m-1$.   Notice that  there are $m$ different routes for eliminating a nonzero weight --as many as nonzero weights. It may be argued that the higher the contribution of a given    point to the residual $\r$, the higher the  difficulty of  converging to  feasible solutions using as initial point the   cubature rule $\{\X^{old},\w^{old}\}$. To account for this fact, we  sort  the indexes of the points with  nonzero weights   in ascending order according to its contribution to the residual    (which is $s_i = w_i^{old} \normd{\u(\x_i^{old})}$, see Line \ref{alg:003_line1}).   The actual subroutine that performs the zeroing operation is $\SOLVERES{}$ in Line \ref{alg:003_line3}.   If this subroutine fails to determine a feasible solution in which the chosen weight  is set to zero, then the operation is repeated with the next point in the sorted list, and so on until   arriving at the desired sparser solution (if such a solution exists at all). 

   \begin{algorithm}[!ht]
   \DontPrintSemicolon
   
   \SetKwFunction{FMain}{MAKE1ZERO}
    \SetKwProg{Fn}{Function}{}{}
                  \SetKwFunction{ceil}{ceil}

       
      \SetKwFunction{LENGTH}{length}
            \SetKwFunction{SORT}{sort}
            \SetKwFunction{ANYOUT}{ANYOUT}
   \SetKwFunction{RMVPNT}{RMV1PNT}

       \SetKwFunction{SOLVERESfun}{SOLVERES}
       
              \SetKwFunction{EVALBASIS}{EVALBASIS}

%

\newcommand{\Xold}{\boldsymbol{\bar{X}}}
\newcommand{\xOLD}{\boldsymbol{\bar{x}}}

\newcommand{\wOLD}{\boldsymbol{\bar{w}}}
\newcommand{\wOLDe}{\bar{w}}
\newcommand{\iORDER}{\boldsymbol{i}}

\newcommand{\fNZ}{\textbf{F}}

\newcommand{\Eall}{\boldsymbol{\mathcal{E}}}
\newcommand{\MSH}{\boldsymbol{\mathcal{M}}}
        \newcommand{\CONV}{{\mathcal{C}}}
        \newcommand{\BREAK}{\textrm{{break}}}
                \newcommand{\AND}{\textrm{\textbf{AND}}}

    \Fn{ $[\CboolCONV,\X,\w,\ElemINFO]$  $\leftarrow$ \FMain{${\X^{old}}$,$\w^{old}$,$N$,$\AuxVAR$,$\MeshINFO$,$\ElemINFO$} }{
    
    \KwData{  Similar to inputs    in Algorithm \ref{alg:002}.   }

    \KwResult{ $\X \in \dom{}$ and $\w \in \Rn{p}$  such that $\r = \u^T(\X) \w -\b = \zero$.  If $\CboolCONV = true$, then  $\normLZ{\w} = \normLZ{\w^{old}}-1$ (i.e., the output vector weight $\w$ has one more zero that the input weight $\w^{old}$); otherwise, $\w \leftarrow \w^{old}$, $\X \leftarrow \X^{old}$  }  
    
     $\fNZ \leftarrow  \textrm{Indexes nonzero entries of } \w^{old}$  \tcp{$\fNZ \subset \{1,2 \ldots p \}$}

     $[\u^{old}_{\fNZ}, \bullet,\bullet]$ $\leftarrow$ \EVALBASIS($\X^{old}_{\fNZ}$,$\MeshINFO$,$\ElemINFO$)  \tcp{Evaluate basis functions $\u$  at points $\X^{old}_{\fNZ} = \{\X^{old}_{\fNZ(1)},\X^{old}_{\fNZ(2)} \ldots \X^{old}_{\fNZ(m)}  \}$  (using the procedure described in Section \ref{sec:evalfun})  . $\u^{old}_{\fNZ}$ is a $m \times p$ matrix, $m$ being the number of nonzero entries of $\w^{old}$   \label{alg:003_line2}}

      $s_i \leftarrow w^{old}_i \normd{\u^{old}_{\fNZ}(i,:)}$   ($i=1,2 \ldots m$)  \tcp{ $\u^{old}_{\fNZ}(i,:)$ denotes the $i$-th row of $\u^{old}_{\fNZ}$ \label{alg:003_line1} }
      
      $\iORDER \leftarrow$  \textrm{Sort } $\{s_1,s_2 \ldots s_m \} \textrm{ in ascending order; return    indexes describing the arragement}$       
      
       $\CboolCONV \leftarrow false$;   $\;  j \leftarrow 1 $   \tcp{Initializations}
      \While{  $j \le m$ $\AND$ $\CboolCONV=false$ \label{alg:003g} }{
      
      $\DOFrr \leftarrow \f(j)$ \tcp{Index  of  candidate weight  to be zeroed  ($\DOFrr \in \{1,2 \ldots p \}$)}
      $[\CboolCONV,\X,\w,\ElemINFO] \leftarrow$   \SOLVERESfun($\DOFrr,{\X^{old},\w^{old},N,\AuxVAR,\MeshINFO,\ElemINFO})$ \tcp{See Algorithm \ref{alg:004}  \label{alg:003_line3}   }

        {$j \leftarrow j +1$ \tcp{ If $\CboolCONV = false$, trying next candidate.    }}

       }

       }
  
   \caption{ Given a cubature rule $\{\X^{old},\w^{old}\}$, this algorithm tries to determine  a cubature rule $\{\X,\w\}$ with one additional zero weight (invoked in Line \ref{alg:002_zero} of Algorithm \ref{alg:002}).  \label{alg:003}   } 
   \end{algorithm}



   \begin{algorithm}[!ht]
   \DontPrintSemicolon
   
   \SetKwFunction{FMain}{SOLVERES}
    \SetKwProg{Fn}{Function}{}{}
                  \SetKwFunction{ceil}{ceil}

       
      \SetKwFunction{LENGTH}{length}
            \SetKwFunction{SORT}{sort}
            \SetKwFunction{ANYOUT}{ANYOUT}
   \SetKwFunction{RMVPNT}{RMV1PNT}
   \SetKwFunction{EVALBASIS}{EVALBASIS}
   
          \SetKwFunction{NEWTONRmodFUN}{NEWTONRmod}


%

\newcommand{\Xold}{\boldsymbol{\bar{X}}}
\newcommand{\xOLD}{\boldsymbol{\bar{x}}}

\newcommand{\wOLD}{\boldsymbol{\bar{w}}}
\newcommand{\wOLDe}{\bar{w}}
\newcommand{\iORDER}{\boldsymbol{i}}

\newcommand{\fNZ}{\textbf{f}}

\newcommand{\Eall}{\boldsymbol{\mathcal{E}}}
\newcommand{\MSH}{\boldsymbol{\mathcal{M}}}
        \newcommand{\CONV}{{\mathcal{C}}}
        \newcommand{\BREAK}{\textrm{{break}}}
                \newcommand{\AND}{\texttt{\textbf{AND}}}

    \Fn{ $[\CboolCONV,\X,\w,\ElemINFO ]$  $\leftarrow$ \FMain{$\DOFrr$,${\X^{old}}$,$\w^{old}$,$N$,$\AuxVAR$,$\MeshINFO$,$\ElemINFO$} }{
    
   \KwData{ $\DOFrr \in \{1,2 \ldots p \}$: Candidate index. Remaining inputs:  similar to inputs of function in Algorithm \ref{alg:003}.    }      
     
    \KwResult{ $\X \in \dom{}$ and $\w \in \Rn{p}$  such that $\r = \u^T(\X) \w -\b = \zero$.  If $\CboolCONV = true$, then  $\normLZ{\w} = \normLZ{\w^{old}}-1$ and $w_{\DOFrr}  = 0$  }  
    
   
   $n \leftarrow 1$;  $\;\;$ $\CboolCONV \leftarrow true$; $\;\;$ $w^{ref} \leftarrow w^{old}_{\DOFrr}$; $\;\;$ $\X \leftarrow \X^{old}$; $\;\;$ $\w \leftarrow \w^{old}$  \tcp{Initializations}
   
   \While{$n \le N$  $\AND$  $\CboolCONV = true$  }
   {     
   $w_{\DOFrr} = w^{ref}(1-n/N)$  \label{alg:004_line_w}

    $[\CboolCONV,\X,\w,\ElemINFO,$\DOFp$] \leftarrow$   \NEWTONRmodFUN($\X,\w$,$\DOFrr$,$\AuxVAR$,$\MeshINFO$,$\ElemINFO$) \tcp{Modified Newton-Raphson, see Algorithm \ref{alg:005}  \label{alg:004_line3}   }
   
   $n \leftarrow n +1$ 
   }
   
   \lIf{$\CboolCONV = false$ }{$\X \leftarrow \X^{old}$; $\;\;$ $\w \leftarrow \w^{old}$  \tcp{No feasible solution  found. }}     
       }
  
   \caption{ Given a cubature rule $\{\X^{old},\w^{old}\}$ and an index  $\DOFrr \in \{1,2 \ldots p\}$, where $w^{old}_{\DOFrr} \neq 0$,  this algorithm determines  a cubature rule $\{\X,\w\}$ with one additional zero at point $\DOFrr$ (  $w_{\DOFrr} = 0$)   \label{alg:004}.   } 
   \end{algorithm}


\subsection{Modified Newton-Raphson algorithm}

We now move to the above mentioned subroutine $\SOLVERES{}$, appearing in Line \ref{alg:003_line3} of Algorithm \ref{alg:003} ---and with pseudo-code explained in Algorithm  \ref{alg:004}. This subroutine is devoted to the calculation of the position  and weights of the remaining  points when the weight of the  chosen ``control'' point $\DOFrr$   ($\DOFrr \in \{1,2 \ldots \pNMOD \}$, $w_{\DOFrr}^{old} \neq 0$)  is set to zero ---by solving the nonlinear equation corresponding to the integration conditions $\r(\X,\w) = \u^T(\X)\w- \b = \zero$.

   \begin{algorithm}[!ht]
   \DontPrintSemicolon
   
   \SetKwFunction{FMain}{NEWTONRmod}
    \SetKwProg{Fn}{Function}{}{}
                  \SetKwFunction{ceil}{ceil}

       
      \SetKwFunction{LENGTH}{length}
            \SetKwFunction{SORT}{sort}
            \SetKwFunction{ANYOUT}{ANYOUT}
   \SetKwFunction{RMVPNT}{RMV1PNT}

       \SetKwFunction{FUNGRAD}{FUNGRAD}
              \SetKwFunction{EVALBASIS}{EVALBASIS}
              \SetKwFunction{SVDloc}{SVD}

       \SetKwFunction{length}{length}
%

\newcommand{\Xold}{\boldsymbol{\bar{X}}}
\newcommand{\xOLD}{\boldsymbol{\bar{x}}}

\newcommand{\wOLD}{\boldsymbol{\bar{w}}}
\newcommand{\wOLDe}{\bar{w}}
\newcommand{\iORDER}{\boldsymbol{i}}

\newcommand{\Eall}{\boldsymbol{\mathcal{E}}}
\newcommand{\MSH}{\boldsymbol{\mathcal{M}}}
        \newcommand{\CONV}{{\mathcal{C}}}
        \newcommand{\BREAK}{\textrm{\textbf{break}}}
                \newcommand{\AND}{\textrm{\textbf{AND}}}

    \Fn{  $[\CboolCONV,\X,\w,\ElemINFO] \leftarrow$   \FMain{$\X^{old},\w^{old}$,$\DOFrr$,$\AuxVAR$,$\MeshINFO$,$\ElemINFO$} }{
    
    \KwData{ $\{\X^{old},\w^{old}\}$: Cubature rule previous  step. $\DOFrr \in \{1,2 \ldots p\}$: Index controlled point.     $\AuxVAR = \{K_{max},\epsilon_{NR}, N_{neg} \}$, $K_{max}$: Maximum number of iterations; $\epsilon_{NR}$: Tolerance convergence residual; $N_{neg}$: Maximum number of negative weights allowed during iterations;  $\MeshINFO$ and  $\ElemINFO$: data structures containing   variables needed to evaluate the residual $\r$ at any point $\x \in \dom{}$ }    
     
    \KwResult{ $\CboolCONV = true $: The Newton-based iterative algorithm has converged to a feasible solution $\{\X,\w\}$ of the equation $\r = \u(\X)^T \w -\b = \zero$, where $w_{\DOFrr} = w_{\DOFrr}^{old}$ is given.      }  
    
        $\DOFp \leftarrow \emptyset$ \tcp{Indexes of points with fixed position but unknown weights \label{alg:005_0} }  
    
    $\DOFl \leftarrow \textrm{Indexes nonzero entries of } \w^{old}, \textrm{excluding } \DOFrr $

     $k \leftarrow 1$;  $\CboolCONV \leftarrow false$; $\DOFs \leftarrow \DOFl$ ; $\epsilon_{svd} \leftarrow 10^{-10}$ \tcp{Initializations}

      \While{ $k \le K_{max}$  $\AND$  $\CboolCONV = false$ }
      {
      
        $[\u(\X^{old}),\grad{\u}(\X^{old}),\ElemINFO]$ $\leftarrow$ \EVALBASIS($\X^{old}$,$\MeshINFO$,$\ElemINFO$)  \tcp{Determine basis functions  and gradients at $\X^{old}$   \label{alg:005_eval}   }

       $\r \leftarrow \u^T(\X^{old}) \w^{old} - \b$  \tcp{Integration residual \label{alg:005_res} }
        
        \eIf{$\normd{\r} \le \epsilon_{NR}  $ }
        {$\CboolCONV \leftarrow true$;  $\X \leftarrow \X^{old}$; $\w \leftarrow \w^{old}$ \tcp{Converged solution}}
        {
        
     $\J_{\S} \leftarrow  [\J_{\X_{\S}},\J_{\w_{\S}}]  $ \tcp{Jacobian matrix (indexes $\S$). $\J_{\X}$ and $\J_{\w}$ are defined in Eqs. \ref{eq:45345} and \ref{eq:45345a}}
     $E_{feas} \leftarrow false$   \tcp{$E_{feas} \leftarrow false$ while   tentative solution   not feasible} 
 
       \While{ $k \le K_{max}$  $\AND$ $E_{feas} = false$ \label{alg:005_while} }{
        
        $\Jhat{} \leftarrow  [\J_{\X_{\DOFl}},\J_{\w_{\DOFs}}]  $ \tcp{$\DOFl \subseteq \DOFs$ ($\DOFl$: indexes unknown weights and positions)  }
        $[\U_J,\S_J,\G^T]  \leftarrow  $    \SVDloc{$\Jhat{}$,$\epsilon_{svd}$}  \tcp{  $\Jhat{} \approx \U_J \S_J \G $ with relative error $\epsilon_{svd} = 10^{-10}$  \label{alg:005_svd}}
        $n_{dofs} \leftarrow (\nSD +1) \length{\DOFl} + \length{\DOFp} $   \tcp{Number of unknowns ($\nSD$: nº of spatial dimensions)} 
        
        \eIf{ $n_{dofs} $ $<$ \length{$\S_J$} }
        { $\BREAK$ \tcp{Overdetermined system (no solution). Exiting internal   loop without convergence}
   }
        {
        $\c \leftarrow - \S_J^{-1}\U^T_J \r $  \tcp{ $\U_J \S_J \G  \Dq = - \r \RRaam \G \Dq = -\S_J^{-1} \U_J^T \r $}
        
        $\Dq \leftarrow \textrm{Sparse solution of } \G \Dq = \c \textrm{ obtained via QR pivoting} $  \tcp{In Matlab:  $\Dq = \G \backslash \c$  \label{alg:005_qr} }
        
        $\X_{\DOFl} \leftarrow \X^{old}_{\DOFl} + \Dq_{x}$; $\; \w_{\DOFs} \leftarrow \w^{old}_{\DOFs} + \Dq_{w}$ \tcp{$\Dq_{x} = \Dq(1:(\nSD +1) \length{\DOFl})$,$\Dq_{w}$:remaining entries \label{alg:005_upd}}
                
        $m_{neg} \leftarrow   \textrm{Number of negative entries of } \w$   
        
        \lIf{$m_{neg} \ge N_{neg}$}{$\BREAK$    \label{alg:005_neg}  }   
         $\DOFy  \leftarrow \textrm{Indexes points  outside the  domain}  (\, \X_{\DOFy(i)} \notin \dom{})  $   \label{alg:005_out}  
        
        \eIf{ $\DOFy \neq \emptyset$  }
        {   $\DOFp \leftarrow \DOFp \cup \DOFy$;  $\; \DOFl \leftarrow \DOFl \setminus \DOFy   $  $\;\; \DOFs \leftarrow \DOFl \cup \DOFp$ \tcp{ Repeat iteration with   $\X_{\DOFp}$ fixed in previous position    \label{alg:005_rep}}
        
        }
        {$\X^{old} \leftarrow \X; \;\;  \w^{old} \leftarrow \w$;  $ \;\; E_{feas} = true$ \tcp{All points inside, update and exit internal loop}

        }
           
           $k \leftarrow k +1$  \label{alg:005_end}

        }
        
        }
                \lIf{$E_{feas} = false$}{$\BREAK$ (no feasible solution found, exiting without convergence)}{}

        }

       }

       }
  
   \caption{ Modified Newton-Rapshon algorithm for solving the constrained nonlinear equation $\r= \u^T(\X)\w - \b = \zero$, using as initial guess $\{\X^{old},\w^{old}\}$. The unknowns are the position   and weights  of the nonzero entries of $\w^{old}$, except for $\w^{old}_{\DOFrr}$, which is given   (this function is invoked in Line \ref{alg:004_line3} of  Algorithm \ref{alg:004}).    } 
     \label{alg:005}
   \end{algorithm}


To facilitate convergence, the weight $w_{\DOFrr}$ is gradually reduced at a rate dictated by the number of steps $N$ (so that  $w_{\DOFrr} = w_{\DOFrr}^{old}(1-n/N)$ at step $n$, see Line \ref{alg:004_line_w}). 

Suppose we have converged to the solution $\{\X_{(n-1)},\w_{(n-1)}\}$   and we want to determine the solution for the next step $n$ using  a   Newton-Raphson iterative scheme, modified so as to account for the constraints that the points must remain within the domain, and that the weights should be positive (although this latter constraint will be relaxed, as explained in what follows).  The pseudo-code of this modified Newton-Raphson scheme is described in turn in Algorithm \ref{alg:005}.  The integration residual at iteration $k \le K_{max}$ is computed in Line \ref{alg:005_res}. This residual admits the following decomposition in terms of unknown and known variables:    
\begin{equation}
\label{eq:4,s...dd}
 \r = \u^T(\X)\w -\b =  \u^T(\X_{\DOFl}) \w_{\DOFl}  + \u^T(\X_{\DOFp}) \w_{\DOFp} + \u^T(\X_{R}) w_{R} -\b. 
\end{equation}
Here, $\DOFl \subset\{1,2 \ldots \pNMOD\}$ denotes the set of points whose positions and weights are unknown, while $\DOFp  \subset\{1,2 \ldots \pNMOD\}$   is the set in which the positions are fixed, but the weights are unknown.  At the first iteration , $\DOFp = \emptyset$ (see Line \ref{alg:005_0}).  The unknown weights will be collectively denoted hereafter by $\w_{\DOFs} = [\w_{\DOFl}^T,\w_{\DOFp}^T]^T$, and the vector of unknowns (including positions and weights ) by $\q \defeq [\X_{\DOFl}^T,\w_{\DOFs}^T]^T$. 

If the Euclidean norm of the  residual is not below the prescribed error tolerance ( Line \ref{alg:005_res}),  we compute, as customary in Newton-Rapshon procedures,  a  correction  $\Dq  = [ \Delta \X_{\DOFl}^T, \Delta \w_{\DOFs}^T]^T$   by obtaining \emph{one} solution of the  \emph{underdetermined} linear equation  
\begin{equation}
\label{eq:underdeter}
   \JhatK{}  \Dq^{} = - \r.
\end{equation}
Here,   $\JhatK{}$ stands for the  block matrix of the Jacobian matrix $\J \in \RRn{\pNMOD}{(\nSD + 1)\pNMOD}$ formed by the rows  corresponding to the indexes of the unknown positions $\X_{\DOFl}$ and the unknown weights $\w_{\DOFs}$   ,   i.e: 
\begin{equation}
 \Jhat \defeq  \rowdos{\J_{\X_{\DOFl}}}{\J_{\w_{\DOFs}}},
\end{equation}
where  
  \begin{equation}
  \label{eq:45345}
  \J_{\X} \defeq \derpar{\r}{\X}  =   \begin{bmatrix}
            w_{1} \gradT{u_1(\x_1)} &    w_{2} \gradT{u_1(\x_2)}  &  \ldots &   w_p \gradT{u_1(\x_p)}      \\ 
             w_{1} \gradT{u_2(\x_1)} &    w_{2} \gradT{u_2(\x_2)}  &  \ldots &   w_p \gradT{u_2(\x_p)}   \\
             \vdots & \vdots & \ddots & \vdots    \\
               w_{1} \gradT{u_p(\x_1)} &    w_{2} \gradT{u_p(\x_2)}  &  \ldots &   w_p \gradT{u_p(\x_p)}  
          \end{bmatrix}_{p \times \nSD \, p} \hspace{-0.5cm},  
\end{equation}
 and  
   \begin{equation}
  \label{eq:45345a}
  \J_{\w} \defeq \derpar{\r}{\w} = \u^T(\X)  =   \begin{bmatrix}
            {u_1(\x_1)} &     {u_1(\x_2)}  &  \ldots &    {u_1(\x_p)}      \\ 
             {u_2(\x_1)} &     {u_2(\x_2)}  &  \ldots &    {u_2(\x_p)}   \\
             \vdots & \vdots & \ddots & \vdots    \\
               {u_p(\x_1)} &     {u_p(\x_2)}  &  \ldots &    {u_p(\x_p)}  
          \end{bmatrix}_{p \times  p} \hspace{-0.5cm}.
\end{equation}
 Recall that the gradients of the basis functions   can be determined by \refeq{eq:4wdasd}, if   the analytical expressions of the integrand functions are available, or by interpolation via Eq. \refpar{eq:72,dwwww} ---using the values of the basis funtions at the Gauss points of the element containing the corresponding point. These operations are encapsulated in the  function $\EvalBasis{}$,  invoked in Line \ref{alg:005_eval}.
%
%
%

Once we have computed   $\Dq$ from \refeq{eq:underdeter}, we update the positions of the points and the weights  in Line \ref{alg:005_upd}. Since the basis functions are only defined inside the domain $\dom{}$ (this is one of the constraints appearing in the sparsification problem \ref{eq:spars1}),   it is necessary to first  identify (  Line \ref{alg:005_out}) and    
  then correct the positions of those points that happen to fall outside the domain. The identification is made by       determining which finite elements contain the points  in their new positions; for the sake of computational efficiency, the search is limited to a patch of elements   centered at the element containing the point at the previous iteration, and located within a radius   $\normd{\Delta \X_{I}}$  ($I \in \DOFl$)---the    mesh connectivities, stored in the data structure $\MeshINFO$, greatly expedites this search task.  If it happens that a given point is not inside any element ($\X_{I} \notin \dom{}$ for some $I \in  \DOFl$), then  we set  $\X_{I} \leftarrow \X_{I}^{old}$,    
  $\DOFp \leftarrow \DOFp \cup I$, and $\DOFl \leftarrow \DOFl \setminus I$ (see Line \ref{alg:005_rep}). Notice that this amounts to ``freezing'' the position of this critical point at the value of the previous iteration during the remaining iterations of the current step\footnote{This can be done because system \refpar{eq:underdeter} is underdetermined, and therefore, one can constrain some points not to move and still find a solution. It should be noticed that these constrained points are freed at the beginning of each step, see Line \ref{alg:005_0}.    }.  This operation is to be repeated until all the points lie within the domain ---ensuring this is the job of the internal \emph{while} loop starting in Line \ref{alg:005_while}.


    The other constraint defining a feasible solution  in the sparsification problem \ref{eq:spars1} is the positiveness of the weights. However,   we argued  in   Section \ref{sec:sparsification} that, since the volume is exactly integrated, the tendency when one of the weights is reduced is that the remaining weights  increase  to compensate for the loss of  volume. Furthermore, according to the sorting criterion employed in Line \ref{alg:003_line1} of Algorithm \ref{alg:003},  negative weights are the first to be zeroed in each local sparsification step, and, thus,   tend to dissapear as the algorithm progresses. For these reasons, the solution procedure does not incorporate any specific strategy   for enforcing positiveness of the weights ---rather, we limit ourselves to keep the number of negative weights below a user-prescribed threshold $M_{neg}$ during the iterative procedure (Line \ref{alg:005_neg} in Algorithm \ref{alg:005}).   {Nevertheless, as a precautionary measure, Line \ref{alg:002_b} in Algorithm \ref{alg:002} prevents negative weights from appearing in the final solution.}

\subsection{Properties of Jacobian matrix and maximum sparsity}

It only remains to addresss the issue of how to solve the system of linear equations \ref{eq:underdeter}. Solving this system of equations is worthy of special consideration because of two reasons: firstly,   the system is    \emph{underdetermined} (more unknowns than equations), and, secondly, the Jacobian matrix $\Jhat{}$   may become \emph{rank-deficient} during the final  steps of the sparsification process, specially in 2D and 3D problems.

\subsubsection{Rank deficiency}
That the Jacobian matrix $\Jhat{}$ may become rank-deficient can be readily demonstrated by analyzing the case of the integration of polynomials in Cartesian domains.  A polynomial of order $t$ gives rise to $\pNMOD = (t + 1)^{\nSD}$  ($\nSD = 1,2$ or $3$) integration conditions (as many as monomials). If one assumes that  the Jacobian matrix $\Jhat{} $ remains full rank during the entire sparsification process,  then it follows that the number of optimal points  one can get under such an assumption, denoted henceforth by $m_{eff}$, is when $\Jhat{}$ becomes square\footnote{Because if there are less unknowns than equations, there are no solution to the equation $\r = \zero$.} (or underdetermined with less than $\nSD$ surplus unknowns); this condition yields  
\begin{equation}
\label{eq:effective}
 m_{eff} = \CEIL{\dfrac{\pNMOD}{\nSD +1}} = \CEIL{\dfrac{(t+1)^{\nSD}}{\nSD +1}} 
\end{equation}
where  $\CEIL{}$  rounds its argument to the nearest integer greater or equal than itself.  For 1D polynomials ($\nSD  = 1$), it is readily seen that $m_{eff}$ coincides with the number of points of the well-known (optimal) Gauss quadrature rule; for instance, for $t= 3$ (cubic polynomials), the above equation gives $m_{eff} = {4}/{2} = 2$ integration points.  This implies that in this 1D case, the Jacobian matrix \emph{does} remain full rank during the process, as presumed.  However, this   does not hold in the 2D and 3D cases. For instance, for 3D cubic polynomials ($t = 3$, $\nSD = 3$), the above equation yields $m_{eff} = (1+3)^{3}/(1+3) = 16$ points, yet it is well  known that  8-points tensor product rule ($2 \times 2 \times 2$)  can integrate exactly cubic polynomials in any cartesian domain. This means that, in this 3D case, from the rule with 16 nonzero weights, to the cubature rule with 8 nonzero weights, the Jacobian matrix $\Jhat{}$ must remain \emph{necessarily} rank-deficient.  

To account for this potential rank-deficiency, we determine the truncated SVD of $\Jhat{}$   (with error threshold $\epsilon_{SVD} = 10^{-10}$ to avoid near-singular cases) in Line \ref{alg:005_svd} of Algorithm \ref{alg:005}: $\JhatK{} \approx \U_J \S_J \G$. Replacing $\JhatK{}$ by this decomposition in \refeq{eq:underdeter}, and  pre-multiplying both sides of the  resulting equation by $\S_J^{-1} \U_J^T$, we obtain, by virtue of the  property $\U_J^T \U_J = \ident$: 
\begin{equation}
 \G \Dq = \c,
\end{equation}
where $\G$ denotes the transpose of the orthogonal matrix  of right-singular vectors of $\JhatK{}$, while $\c = - \S_J^{-1} \U_J^T \r$, $\S_J$ being the diagonal  matrix of singular values, and $\U_J$ the matrix of left singular vectors. 

\subsubsection{Underdetermination and sparse solutions}

 Let us discuss now the issue of underdeterminacy. It is easy to show that  the preceding system of equations remains underdetermined  during the entire sparsification process, with a degree of underdeterminacy (surplus of  unknowns over number of equations) decaying at each sparsification step until the optimum is reached, when $\G$ becomes as square as possible.  For instance, at the  very first step of the   process,  in a problem with $\pNMOD$ basis functions,  
 $\Jhat{}$ is by  construction\footnote{On the grounds that $\u(\Xdecm)$ is also full rank   because otherwise \refeq{eq:4sds} would not hold. } full rank   (i.e., there are $\pNMOD$ linearly independent equations), while the  number of unknowns is   $(\pNMOD-1) (1 + \nSD)$ ($m = \pNMOD-1$ points with $\nSD$ unknowns coordinates associated to the position of each point and one unknown associated to its weight).   Thus, the solution space  in this case is of dimension  $( \nSD(\pNMOD-1) -1) $.
 
 To update the  weights and the positions, we need to pick up \emph{one} solution from this vast space.  The standard approach 
    in Newton's method for underdetermined systems    (and also the method favored in the literature on   generalized Gaussian quadratures \cite{mousavi2010generalized,xiao2010numerical,bremer2010nonlinear} ) is to use the  \emph{least $\ell_2 $-norm solution}, which is simply $\Dq \defeq \G^{+} \c$, where $\G^{+} = \G^T (\G \G^T)^{-1}$ is the pseudo-inverse of $\G$ (notice that in our case\footnote{In this regard, it should be pointed out that References \cite{mousavi2010generalized,xiao2010numerical,bremer2010nonlinear} calculate  the pseudo-inverse of the Jacobian matrix as $\Jhat^{+} = \Jhat^{T}(\Jhat{} \Jhat^{T})^{-1}$, thus ignoring the fact that, as we have argued in the foregoing,   $\Jhat{}$ might become rank-deficient, and therefore, $\Jhat{} \Jhat{}^T$ cannot be inverted.} $\G^+ = \G^T$).   However, we do not use this approach here because the resulting solution tends to be \emph{dense},  and this implies that the positions of all the cubature points have to updated at all iterations. This is a significant disadvantage in our interpolatory framework, since     updating the position of one point entails an interpolation error of greater or lesser extent depending on the  functions being interpolated and the distance from the FE Gauss points.   Thus, it  would be beneficial for the overall accuracy of the method to determine solutions that \emph{minimize the number of positions being updated at each iteration} ---incidentally, this would also help to reduce the computational effort associated to the spatial search carried out in Line \ref{alg:005_out} of Algorithm \ref{alg:005}. This requisite natural calls for \emph{solution methods that promote sparsity}. For this reason, we use here (see Line \ref{alg:005_qr} in Algorithm \ref{alg:005})  the QR factorization with column pivoting  (QRP) proposed in Golub et al. \cite{golub2012matrix} (page 300, Algorithm 5.6.1), which furnishes a solution with as many nonzero entries as equations\footnote{In Matlab,  this QRP solution is the one obtained in using the ``backslash'' operator (or $\textrm{\texttt{mldivide}}(\G,\c)$). }.   An alternative strategy would be to determine the \emph{least $\ell_1$-norm solution}, which, as discussed in Section \ref{sec:statenfass}, also promotes sparsity \cite{boyd2004convex,chen2001atomic}. However, computing this solution  would involve  addressing  a convex, nonquadratic optimization problem at each iteration, and this would require  considerably more effort and sophistication than the simple QRP method employed in Line \ref{alg:005_qr}.

 \subsection{ Summary}
    \label{sec:cecmLAST}
   By way of conclusion, we summarize in Box \ref{box:1}, all the operations required to determine an optimal cubature rule using as initial data the location of the FE Gauss points, their corresponding weights, and  the values at such Gauss points of the functions we wish to efficiently integrate.

   \floatstyle{boxed}
\newfloat{BOX}{!ht}{box}[section]
\floatname{BOX}{Box}
\begin{BOX}  
 \bbee
 \item   Given the coordinates of the nodes of the finite element mesh, the array of element connectivities, and the position of the  Gauss points  for each element in the parent domain, determine the location of such  points in the physical domain: $\xGall = \{ \xG{1}{1}, \xG{1}{2}, \ldots \xG{e}{i} \ldots \}$. Likewise, compute the vector of (positive) finite element weights ($\Wfe \in \Rn{M}$) for each of these points as the product of the corresponding Gauss weights and the Jacobian of the transformation from the parent domain to the physical domain.   
 \item Determine the values at all Gauss points $\xGall$ of the  parameterized  function $\a : \dom{} \times  \paramSPC \rightarrow \Rn{\nU}$ we wish to efficiently integrate for the chosen   parameters $\{\inpP_{1},\inpP_2 \ldots \inpP_{P}\} \subset \paramSPC$, and store the result in matrix $\Afe \in \RRn{M}{P n}$  (see Eq. \ref{eq:8243t--}). In the case of hyperreduced-order models, the analytical expression of the integrand functions  is normally not available as an explicit function of the input parameters, and  constructing matrix $\Afe$  entails solving the corresponding governing equations for the chosen input parameters. If the matrix proves to be too large to fit into main memory, it should be partitioned into column blocks ${\Afe}_1$, ${\Afe}_2$, $\ldots$ ${\Afe}_p$. Such blocks need not be loaded into main memory all at once.
 \item Compute the weighted matrix $\Abar \defeq \diag{\sqrt{\Wfe}} \A $  (see Eq \ref{eq:82kdsdddd}) (alternatively, one may directly store in Step 2, rather than $\Afe$, $\Abar$ itself; this is especially convenient when the matrix is treated in a partitioned fashion, because it avoids loading the submatrices twice).
 \item  Determine the SVD of $\Abar$ ($\Abar \approx \Ubar \S \V$), with relative truncation tolerance $\epsilon_{svd}$ equal to the desired   error threshold for the integration (see Eq. \ref{eq:3qewf}).  
 If the matrix   is relatively small, one can use directly standard  SVD implementations (see function \texttt{SVDT} in Algorithm \ref{alg:007} of Appendix \ref{sec:SRSVD}). If the matrix is large but still fits  into main memory without compromising the machine performance, the incremental randomized SVD   proposed in Appendix \ref{sec:SRSVD} (function \texttt{RSVDinc}, described in Algorithm \ref{alg:009})  may be used instead. Lastly, if the matrix does not fit into main memory, and is therefore provided in a partitioned format,   the Sequential Randomized SVD described also in Appendix \ref{sec:SRSVD} (function $\texttt{SRSVD()}$ in Algorithm \ref{alg:006}) is to be used. 
 \item Determine a $\Wfe$-orthogonal basis matrix for the range of $\Afe$ by making $\U = \diag{(\sqrt{\Wfe})}^{-1} \Ubar$. Following the guidelines outlined in Section \ref{sec:constant_funct},  augment $\U$ with one additional column if necessary so that the column space of $\U\in \RRn{M}{p}$ contains the constant function. 
 \item  Apply the Discrete Empirical Cubature Method (DECM, see Section \ref{sec:DECM1})    to compute a set of indexes $\zDECM$ and positive  weights $\wDECM$ such that $\U(\zDECM,:)^T \wDECM = \U^T \Wfe$ (see Eq. \ref{eq:2dasdfsad}). 
 
 \item Using as initial solution the weights $\wDECM$ obtained by the DECM, as well as the corresponding positions $\Xdecm$,   
 solve the \emph{sparsification   } problem \ref{eq:spars1} by means of function \texttt{SPARSIFglo}  in Algorithm \ref{alg:001}: 
 
 \begin{equation*}
  [\X,\w,\ElemINFO] \leftarrow \texttt{SPARSIFglo} ({{\Xdecm},\wDECM,\Nsteps,\AuxVAR,\MeshINFO})
 \end{equation*}
 The desired cubature rule is given by $\w^{cecm} = [w_{g_1},w_{g_2}, \ldots w_{g_m} ]$ ($w_{g_i} >0$) and $\X^{cecm}  = \{\x_{g_1},\x_{g_2}, \ldots \x_{g_m} \} $, where $g_1,g_2 \ldots g_m$ denote the indexes of the nonzero entries of the (sparse) output weight vector $\w$.

 \eeee
 
 \caption{ Algorithmic steps involved in the  proposed \emph{Continuous Empirical Cubature Method} (CECM). }    \label{box:1}
\end{BOX}

   \section{Numerical assessment}
A repository containing both the Continuous Empirical Cubature Method (CECM), as well as the Sequential Randomized SVD (SRSVD), allowing to reproduce the following examples is publicly available at \url{https://github.com/Rbravo555/CECM-continuous-empirical-cubature-method}

    \label{sec:numass1}

\subsection{  Univariate polynomials}
\label{sec:uni}

  \begin{figure}[!ht]
  \centering
  \subfigure[]{\label{fig:FIG_3_a}\includegraphics[width=0.49\textwidth]{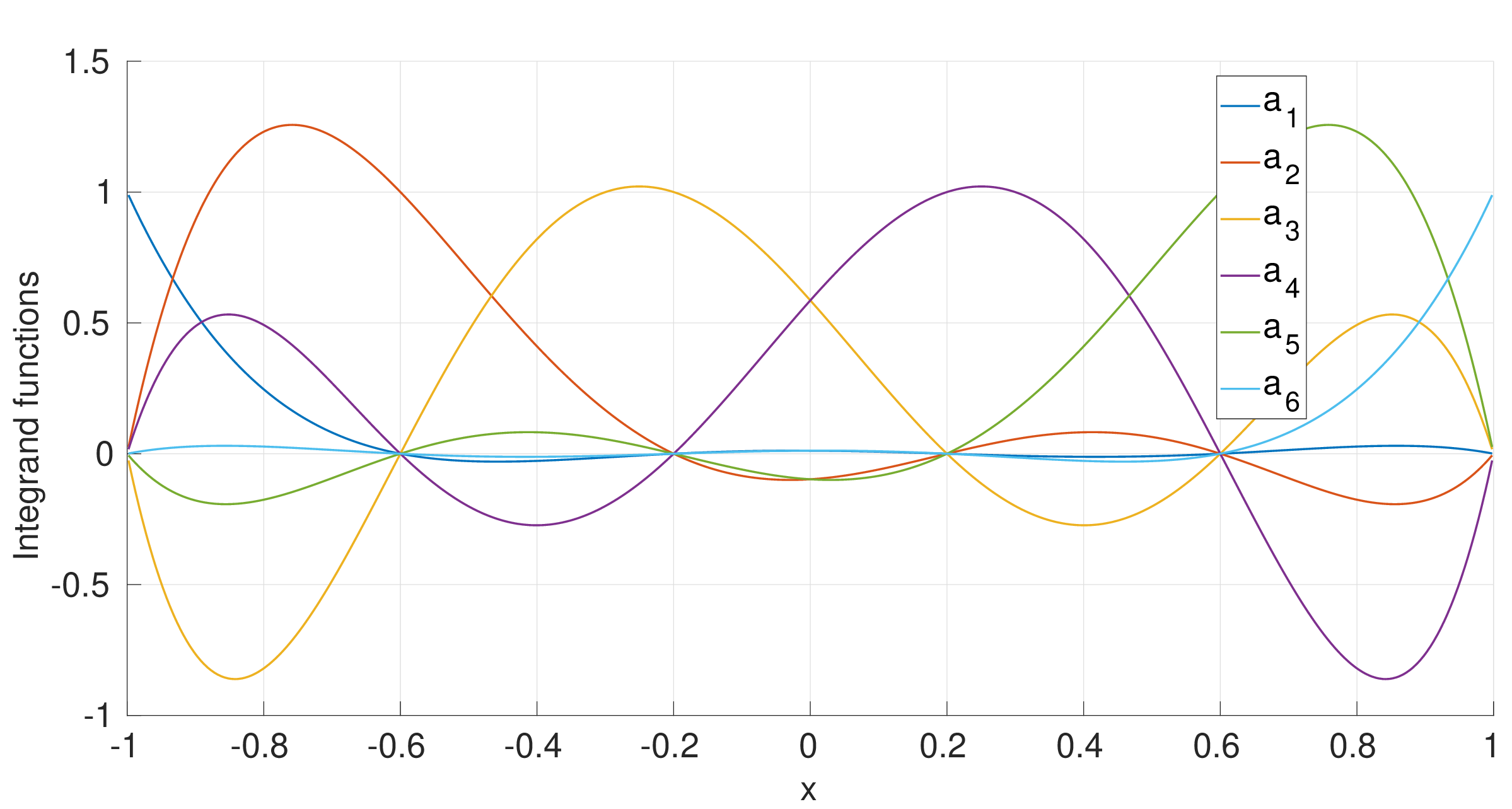}}   
     \subfigure[]{\label{fig:FIG_3_b}\includegraphics[width=0.49\textwidth]{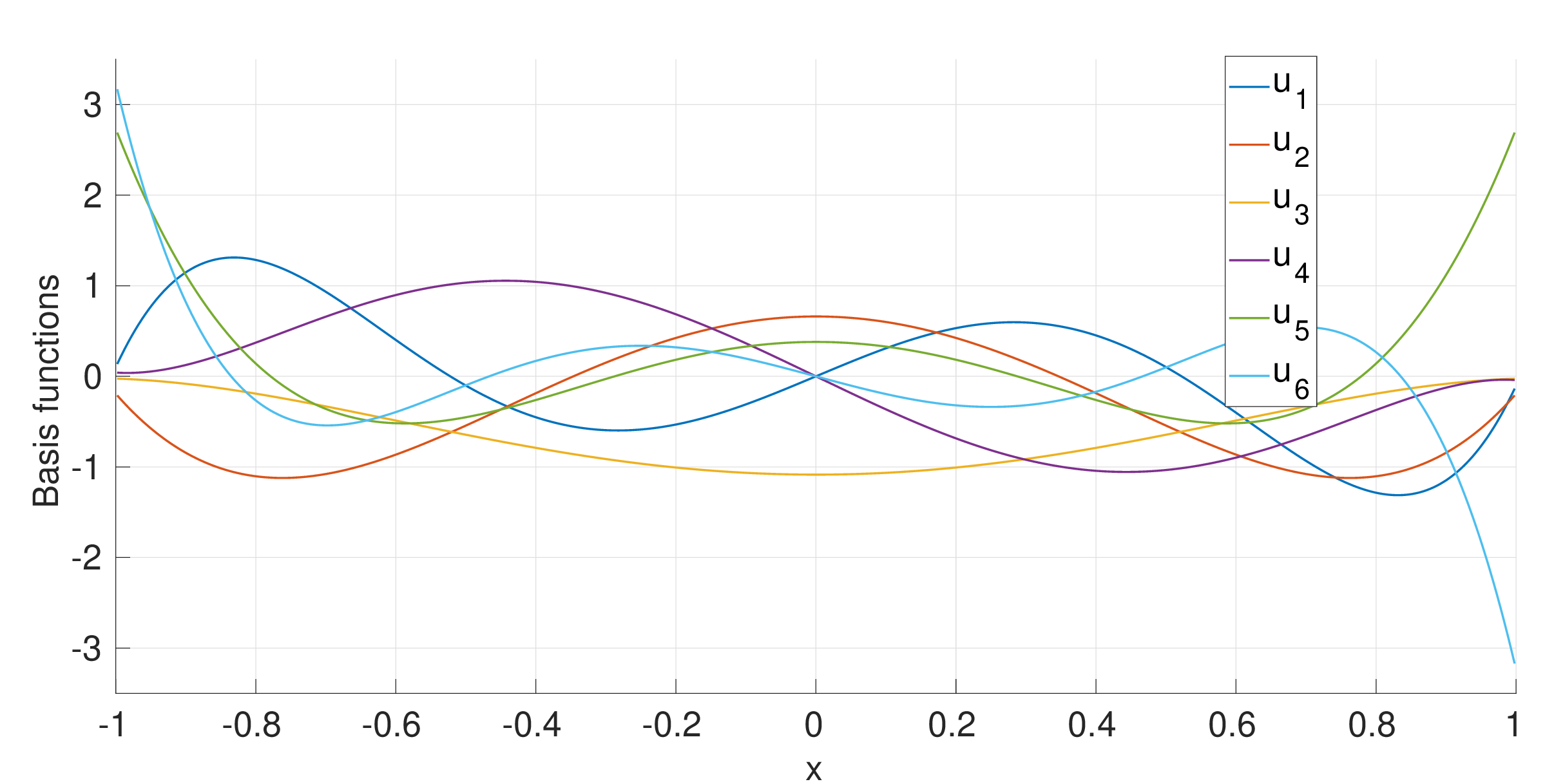}}   
   \caption{ Integration of  univariate Lagrange polynomials of degree    $p = 5$.  a) Original integrand functions $a_1,a_2 \ldots a_6$.   b) Orthogonal   basis functions $u_1$, $u_2$ \ldots $u_6$ derived from the weighted SVD described in step 5 of Box \ref{box:1}). The CECM operates on these basis functions rather than on the original integrand functions in Figure \ref{fig:FIG_3_a}.     }
  \label{fig:FIG_3}
\end{figure}

We begin the assessement of the proposed methodology by examining the example used for motivating the proposal: the integration of univariate polynomials  in the domain $\dom{} = [-1,1]$ .   The employed finite element mesh features $\nelem = 200$ equally-sized elements, and $\ngausE = 4$ Gauss points per element, resulting in a total of $\ngausT = 800$ Gauss points.  Given a   degree $p>0$, and a set of $P= p+1$ equally space nodes $x_1,x_2 \ldots x_{P}$,  we seek the optimal integration rule for the Lagrange polynomials  
 \begin{equation} 
   \label{eq:res1}
    a_i(x)  = \prod_{j=1, i\neq j}^{P}  \dfrac{x-x_j}{x_i-x_j}, \hspace{1cm}  i=1,2 \ldots P
   \end{equation}
 (graphically represented in Figure \ref{fig:FIG_3_a} for degree up to $p=5$).

    \begin{figure}[!ht]
  \centering
  \subfigure[Iteration 1]{\label{fig:DECM1}\includegraphics[width=0.329\textwidth]{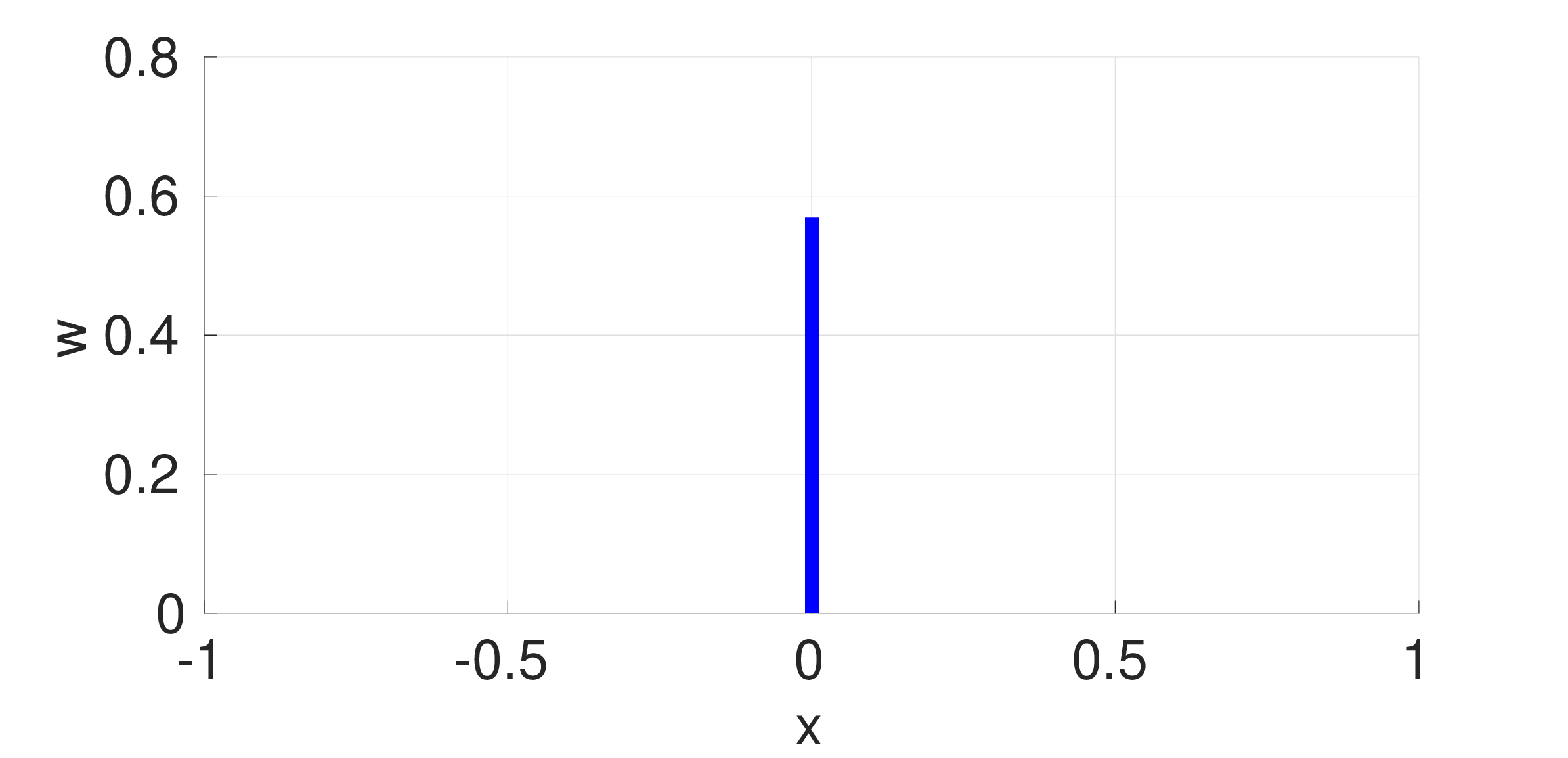}}   
     \subfigure[Iteration 2]{\label{fig:DECM2}\includegraphics[width=0.329\textwidth]{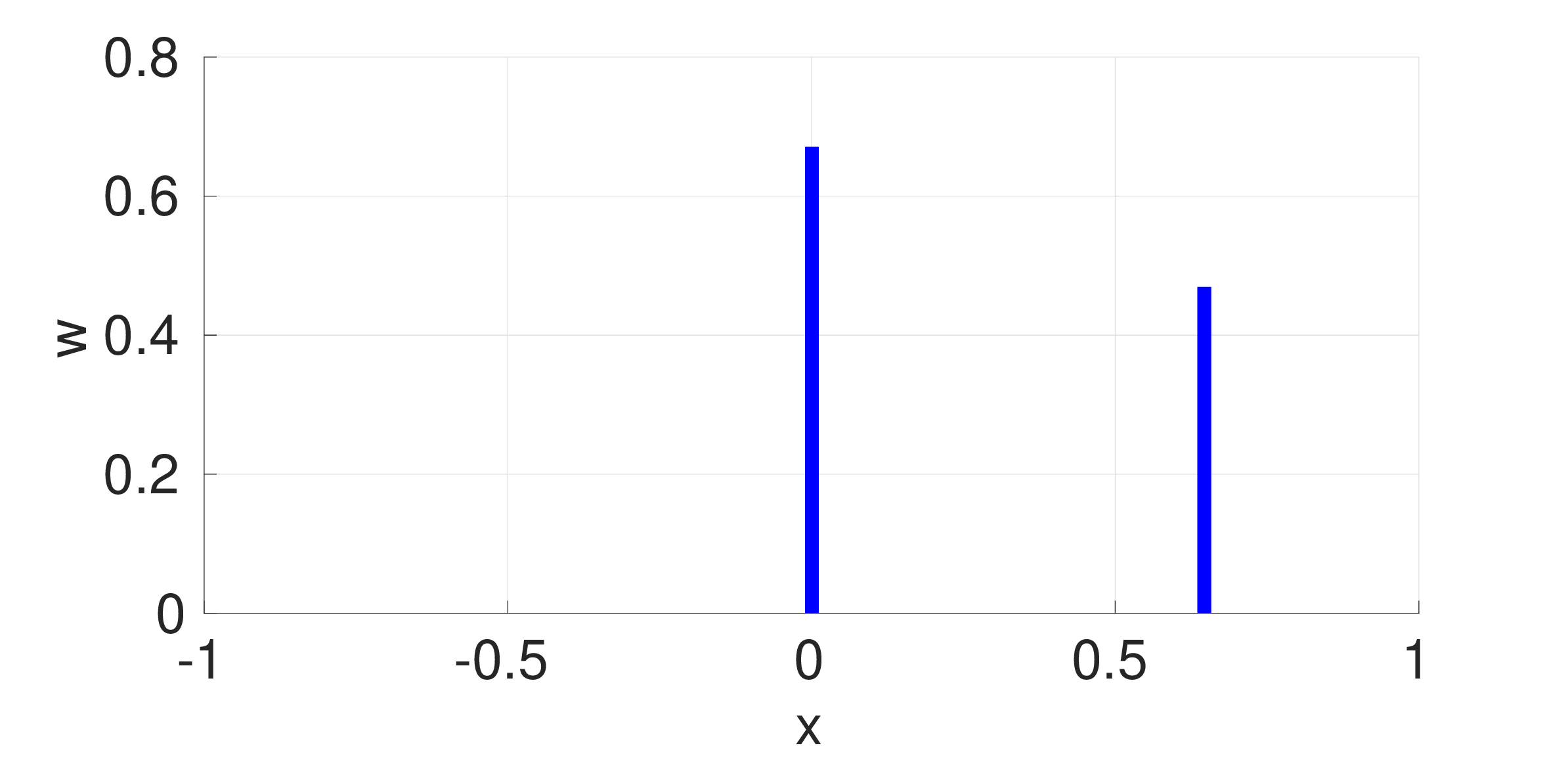}}   
     \subfigure[Iteration 3]{\label{fig:DECM3}\includegraphics[width=0.329\textwidth]{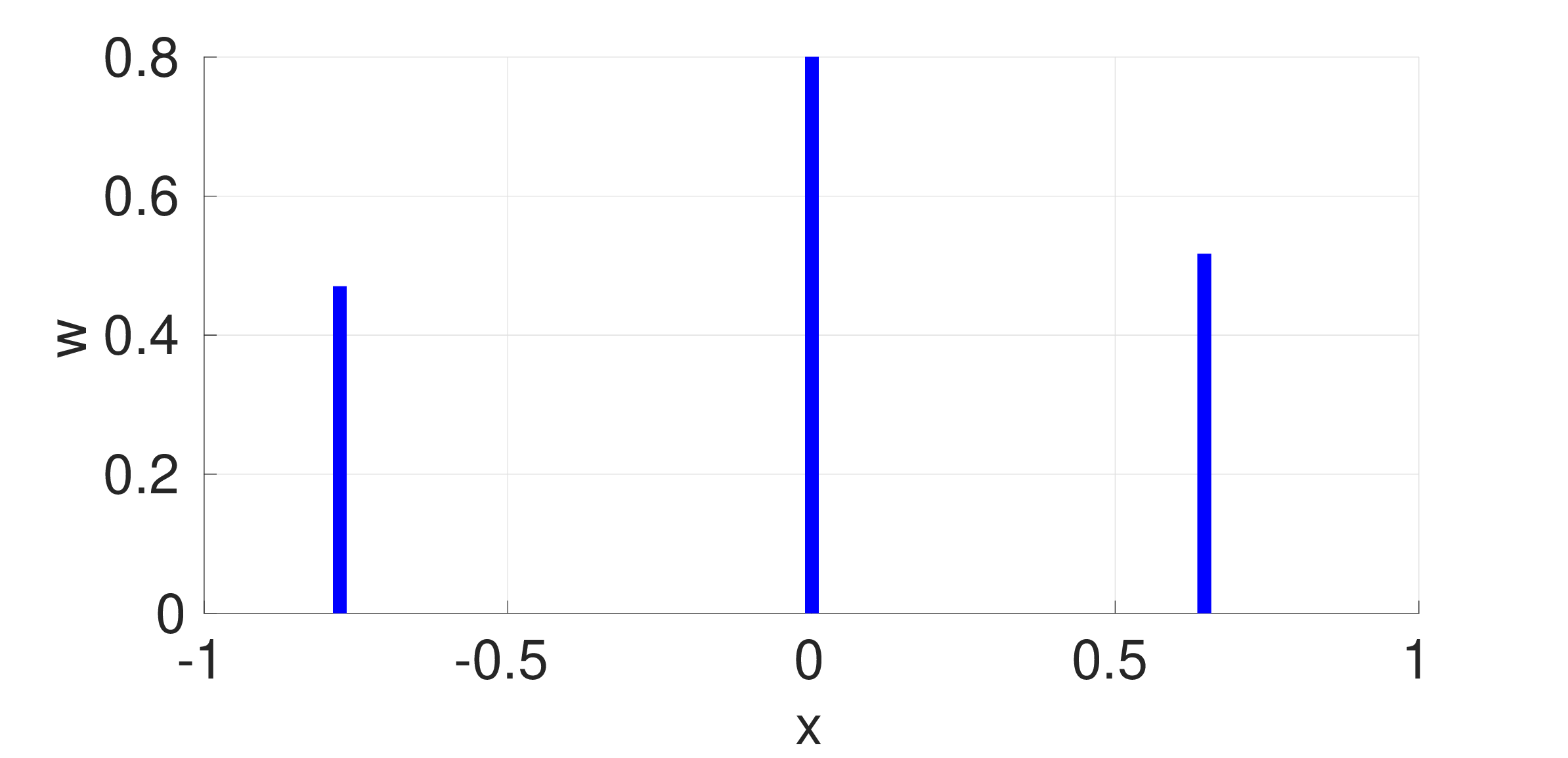}}   
     \subfigure[Iteration 4]{\label{fig:DECM4}\includegraphics[width=0.329\textwidth]{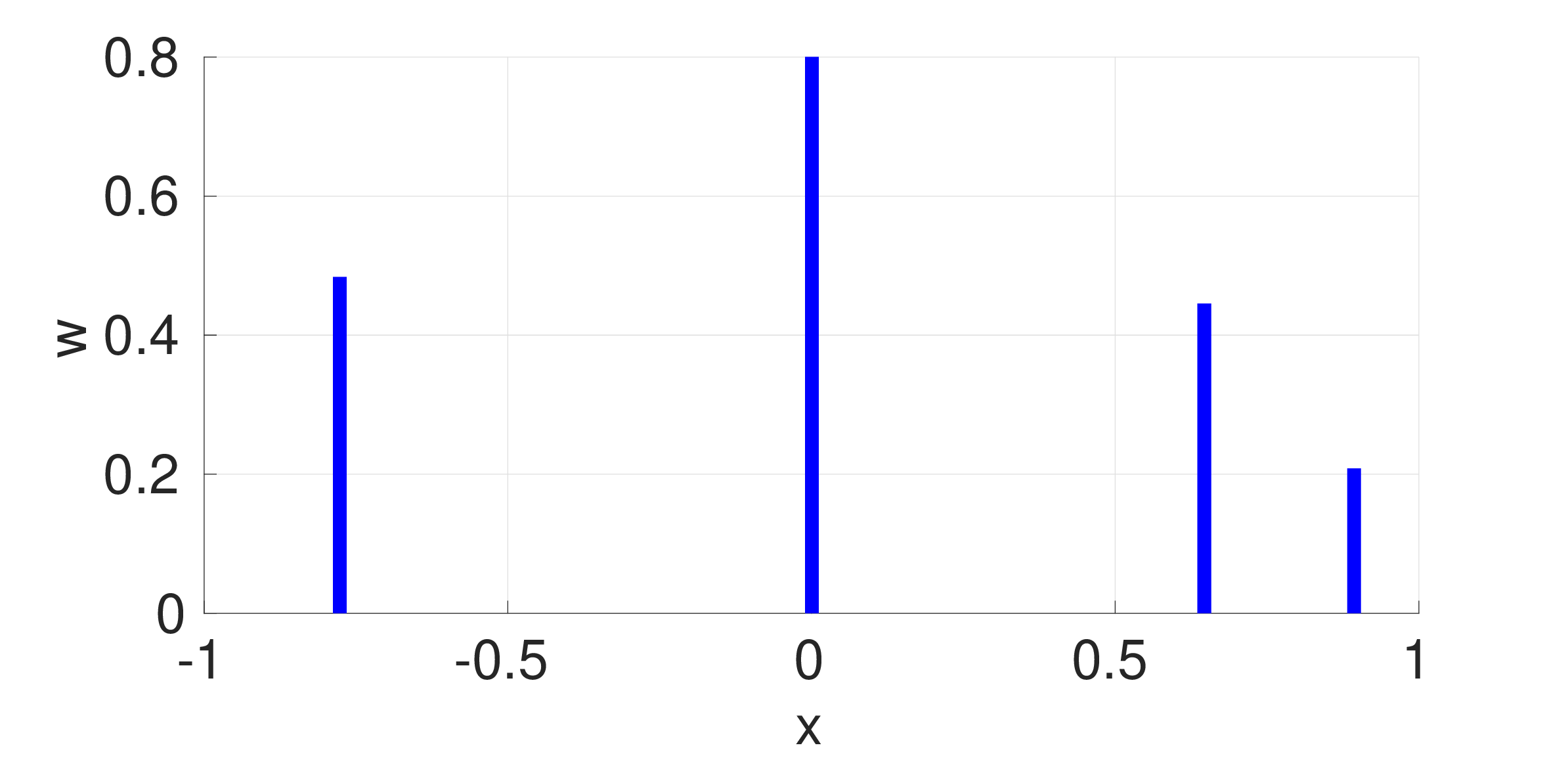}}   
     \subfigure[Iteration 5]{\label{fig:DECM5}\includegraphics[width=0.329\textwidth]{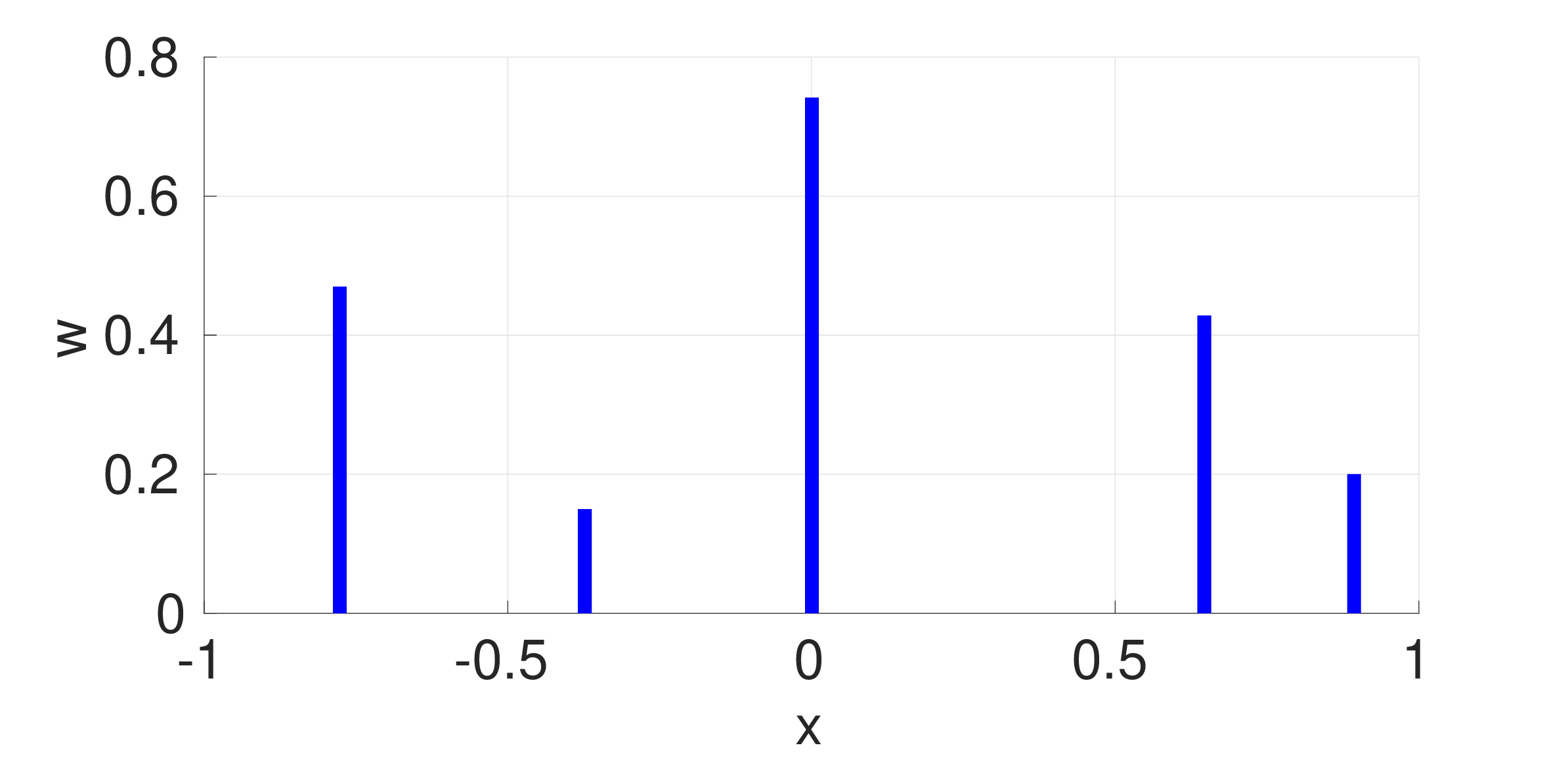}}   
     \subfigure[Iteration 6. DECM rule (6 points)]{\label{fig:DECM6}\includegraphics[width=0.329\textwidth]{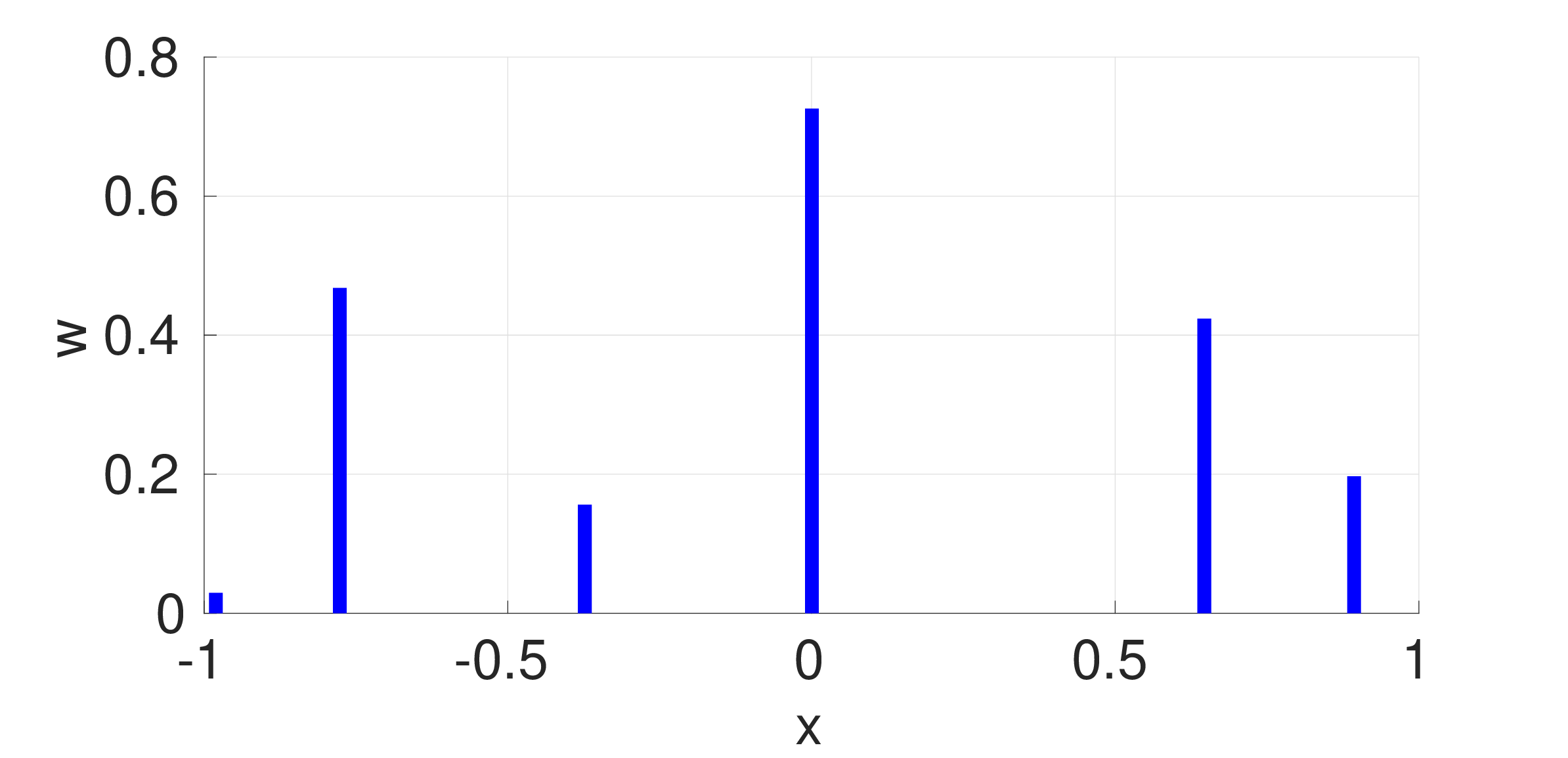}}   
   \caption{ Location ($x$) and corresponding weights ($w$) selected by the DECM (see Box \ref{box:1}, step 6), at each iteration, for the case of polynomials of  degree $p = 5$. The final  integration rule (iteration 6, graph  \ref{fig:DECM6}) is the starting point for the subsequent sparsification problem (illustrated, in turn, in Figure \ref{fig:FIG_5}).       }
  \label{fig:FIG_4}
\end{figure}

\jahoHIDE{Generated with   \href{MATLAB_EXAMPLES_CODE/EXAMPLES_1D/CECM_1D.m}{EXAMPLES\_1D/CECM\_1D.m}, input InputDataFile ='DATA\_Lag\_poly\_5'  }
 The value of these Lagrange polynomials at the $\ngausT$ Gauss points are stored in the matrix $\Afe \in \RRn{\ngausT}{P}$, which is then subjected to the weighted SVD (step 5 in Box \ref{box:1}) for determining   $L_2$-orthogonal basis functions $u_1$, $u_2$ \ldots  $u_P$ ---plotted   in Figure \ref{fig:FIG_3_b} for the case $p=1,2 \ldots 5$.  The truncation tolerance in the SVD of expression \refpar{eq:3qewf} is set in this case  to $\tolSVD = 0$, for we seek quadrature rules that \emph{exactly} integrate any polynomial up to the specified degree.  
 The resulting basis matrix  $\U$ (obtained by Eq. \ref{eq:42dadd} from the left singular vectors of the above mentioned SVD), along with the full-order weight vector $\Wfe \in \Rn{\ngausT}$,  are then used to determine an interpolatory quadrature rule by means of the DECM (step 6 in Box \ref{box:1}). As commented in Section \ref{sec:DECM1}, this algorithm selects one point at each iteration, until arriving at as many points as basis functions. By way of illustration, we show in Figure \ref{fig:FIG_4} the iterative sequence leading to the interpolatory quadrature rule with $6$ points (i.e., for polynomials up to degree 5).    
 
%
%
%
%
%
%
%
%
%
%
%

\jahoHIDE{See \url{ MATLAB_EXAMPLES_CODE/EXAMPLES_1D/Generate_table.m} }

\begin{table}[!ht]
\footnotesize
\centering
\begin{tabular}{|c|c|c|c|c|c|c|c|c|}
\cline{1-4} \cline{6-9}
\rowcolor[gray]{0.8}
degree &
  positions &
  weights &
  error (Gs.) &
 \cellcolor{white}  \multirow{9}{*}{} &
  degree &
  positions &
  weights &
  error (Gs.) \\ \cline{1-4} \cline{6-9} 
1 &
  -7.31662253127291E-17 &
  2 &
  2.2504e-16 &
   &
  \multirow{5}{*}{9} &
  -0.906179845938664 &
  0.236926885056189 &
  \multirow{5}{*}{4.8426e-16} \\ \cline{1-4} \cline{7-8}
\multirow{2}{*}{2} &
  -0.718298153787432 &
  0.784973499347808 &
  \multirow{2}{*}{} &
   &
   &
  -0.538469310105683 &
  0.478628670499367 &
   \\ \cline{2-3} \cline{7-8}
 &
  0.464059849765363 &
  1.21502650065219 &
   &
   &
   &
  3.18949244021535E-16 &
  0.568888888888889 &
   \\ \cline{1-4} \cline{7-8}
\multirow{2}{*}{3} &
  -0.577350269189626 &
  1 &
  \multirow{2}{*}{1.6653e-16} &
   &
   &
  0.538469310105683 &
  0.478628670499366 &
   \\ \cline{2-3} \cline{7-8}
 &
  0.577350269189626 &
  1 &
   &
   &
   &
  0.906179845938664 &
  0.236926885056189 &
   \\ \cline{1-4} \cline{6-9} 
\multirow{3}{*}{4} &
  -0.966877924131567 &
  0.25396505608136 &
  \multirow{3}{*}{} &
   &
  \multirow{6}{*}{10} &
  -0.940907514603222 &
  0.151215583523548 &
  \multirow{6}{*}{} \\ \cline{2-3} \cline{7-8}
 &
  -0.266817254357718 &
  1.00671682129073 &
   &
   &
   &
  -0.694289509608625 &
  0.336602757683989 &
   \\ \cline{2-3} \cline{7-8}
 &
  0.695455188587597 &
  0.739318122627913 &
   &
   &
   &
  -0.287589748592848 &
  0.462395955495934 &
   \\ \cline{1-4} \cline{7-8}
\multirow{3}{*}{5} &
  -0.774596669241483 &
  0.555555555555556 &
  \multirow{3}{*}{8.4549e-16} &
  \multirow{13}{*}{} &
   &
  0.194969778013914 &
  0.482737964706169 &
   \\ \cline{2-3} \cline{7-8}
 &
  5.19179955929866E-17 &
  0.88888888888889 &
   &
   &
   &
  0.637318996155301 &
  0.382659779209725 &
   \\ \cline{2-3} \cline{7-8}
 &
  0.774596669241483 &
  0.555555555555556 &
   &
   &
   &
  0.927198931106149 &
  0.184387959380636 &
   \\ \cline{1-4} \cline{6-9} 
\multirow{4}{*}{6} &
  -0.837102793435635 &
  0.405516777593141 &
  \multirow{4}{*}{} &
   &
  \multirow{6}{*}{11} &
  -0.932469514203152 &
  0.17132449237917 &
  \multirow{6}{*}{1.0484e-15} \\ \cline{2-3} \cline{7-8}
 &
  -0.245834821655188 &
  0.717131675886906 &
   &
   &
   &
  -0.661209386466264 &
  0.360761573048139 &
   \\ \cline{2-3} \cline{7-8}
 &
  0.45930568155797 &
  0.62534095948649 &
   &
   &
   &
  -0.238619186083197 &
  0.467913934572691 &
   \\ \cline{2-3} \cline{7-8}
 &
  0.906836939036126 &
  0.252010587033462 &
   &
   &
   &
  0.238619186083197 &
  0.46791393457269 &
   \\ \cline{1-4} \cline{7-8}
\multirow{4}{*}{7} &
  -0.861136311594053 &
  0.347854845137454 &
  \multirow{4}{*}{5.8993e-16} &
   &
   &
  0.661209386466265 &
  0.360761573048139 &
   \\ \cline{2-3} \cline{7-8}
 &
  -0.339981043584856 &
  0.652145154862546 &
   &
   &
   &
  0.932469514203152 &
  0.17132449237917 &
   \\ \cline{2-3} \cline{6-9} 
 &
  0.339981043584856 &
  0.652145154862545 &
   &
   &
  \multirow{7}{*}{12} &
  -0.967242104087041 &
  0.088484427217067 &
  \multirow{7}{*}{} \\ \cline{2-3} \cline{7-8}
 &
  0.861136311594053 &
  0.347854845137454 &
   &
   &
   &
  -0.803941159117123 &
  0.240851330885801 &
   \\ \cline{1-4} \cline{7-8}
\multirow{5}{*}{8} &
  -0.998731268512389 &
  0.081227685280535 &
  \multirow{5}{*}{} &
   &
   &
  -0.493300094782076 &
  0.37168382888409 &
   \\ \cline{2-3} \cline{7-8}
 &
  -0.719380135473419 &
  0.445922523698483 &
   &
   &
   &
  -0.083660679344391 &
  0.434164516926175 &
   \\ \cline{2-3} \cline{7-8}
 &
  -0.166434516034323 &
  0.623254258303059 &
   &
  \multirow{3}{*}{} &
   &
  0.346567257597574 &
  0.411907795624614 &
   \\ \cline{2-3} \cline{7-8}
 &
  0.446668026062019 &
  0.562352205951177 &
   &
   &
   &
  0.712911972240284 &
  0.308489844409304 &
   \\ \cline{2-3} \cline{7-8}
 &
  0.885864638161693 &
  0.287243326766745 &
   &
   &
   &
  0.943166559093504 &
  0.144418256052949 &
   \\ \cline{1-4} \cline{6-9} 
\end{tabular}
\caption{ Quadrature rules computed by the CECM for univariate polynomials of degree up to $p = 12$. The rightmost column represents the relative deviations with respect to the optimal Gaussian rules (for polynomials of odd degree),  calculated as  $e^2 = (\normd{\X_{cecm}-\X_{gauss}}^2 + \normd{\w_{cecm}-\w_{gauss}}^2$)/($\normd{\X_{gauss}}^2 + \normd{\w_{gauss}}^2$).     }
\label{tab:1}
\end{table}

\jahoHIDE{See \url{MATLAB_EXAMPLES_CODE/EXAMPLES_1D/Table_2columns.tgn}, generated by https://www.tablesgenerator.com/latex_tables}

The final  step in the process is the   sparsification of the vector of DECM weights to produce the final CECM rule (step 7 in Box \ref{box:1}). Table \ref{tab:1} shows the location and weights obtained in this sparsification for polynomials up to degree $p=12$. The parameters used in this process are: $K_{max}$= 40,$\epsilon_{NR} = 10^{-8}$, $N_{neg}=5$ and $\Nsteps = 20$ (the definition of these parameters is given  in Algorithm \ref{alg:005}), and  we use  analytical evaluation of the integrand and their derivatives through formulas \refpar{eq:follws1} and \refpar{eq:4wdasd}, respectively.  

\begin{figure}[!ht]
  \centering
  \subfigure[Initial DECM rule ($m=6$ points) ]{\label{fig:CECM1}\includegraphics[width=0.30\textwidth]{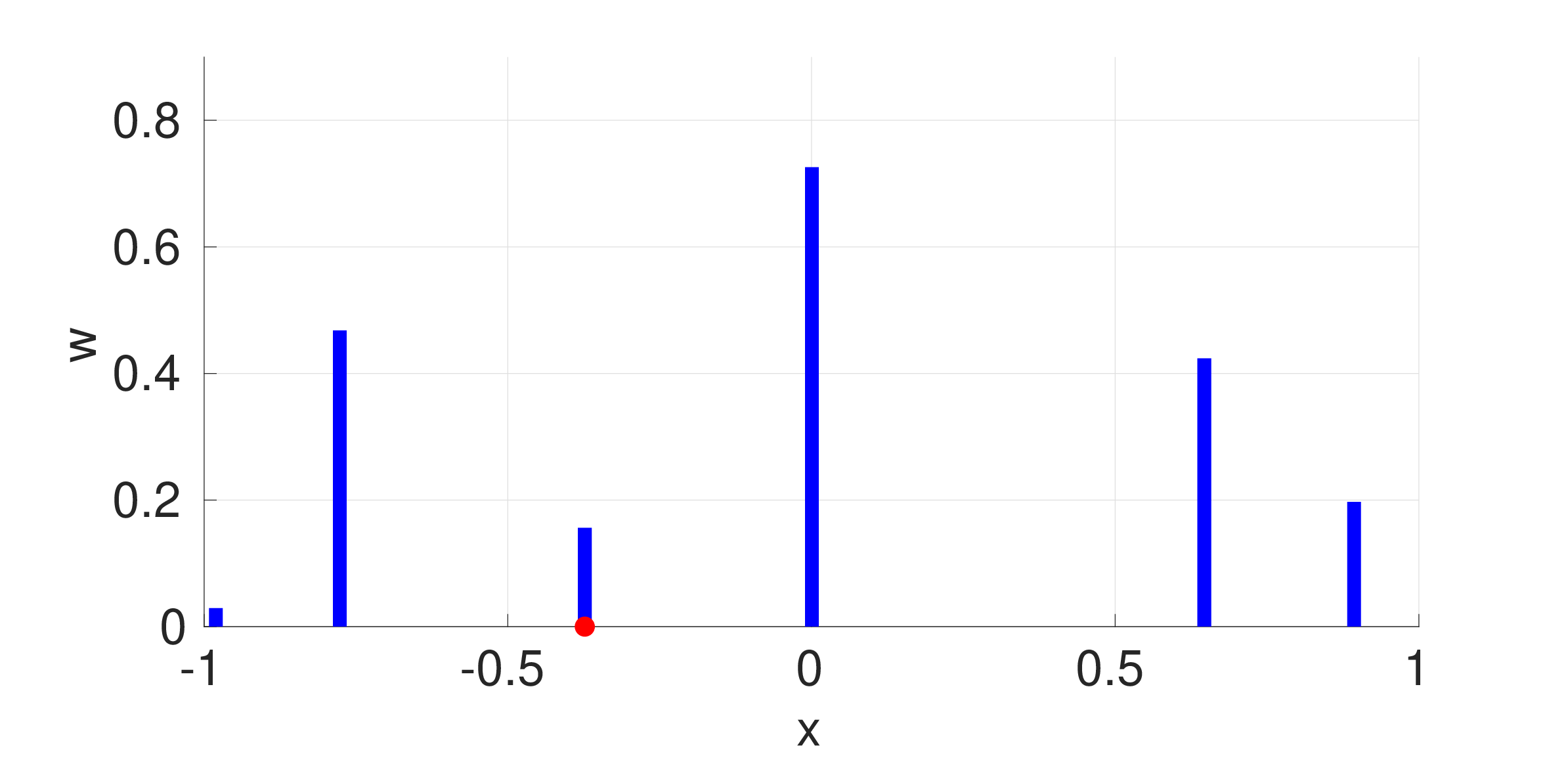}}   
     \subfigure[$m= 5$ points ($t=1,k=4$)]{\label{fig:CECM2}\includegraphics[width=0.30\textwidth]{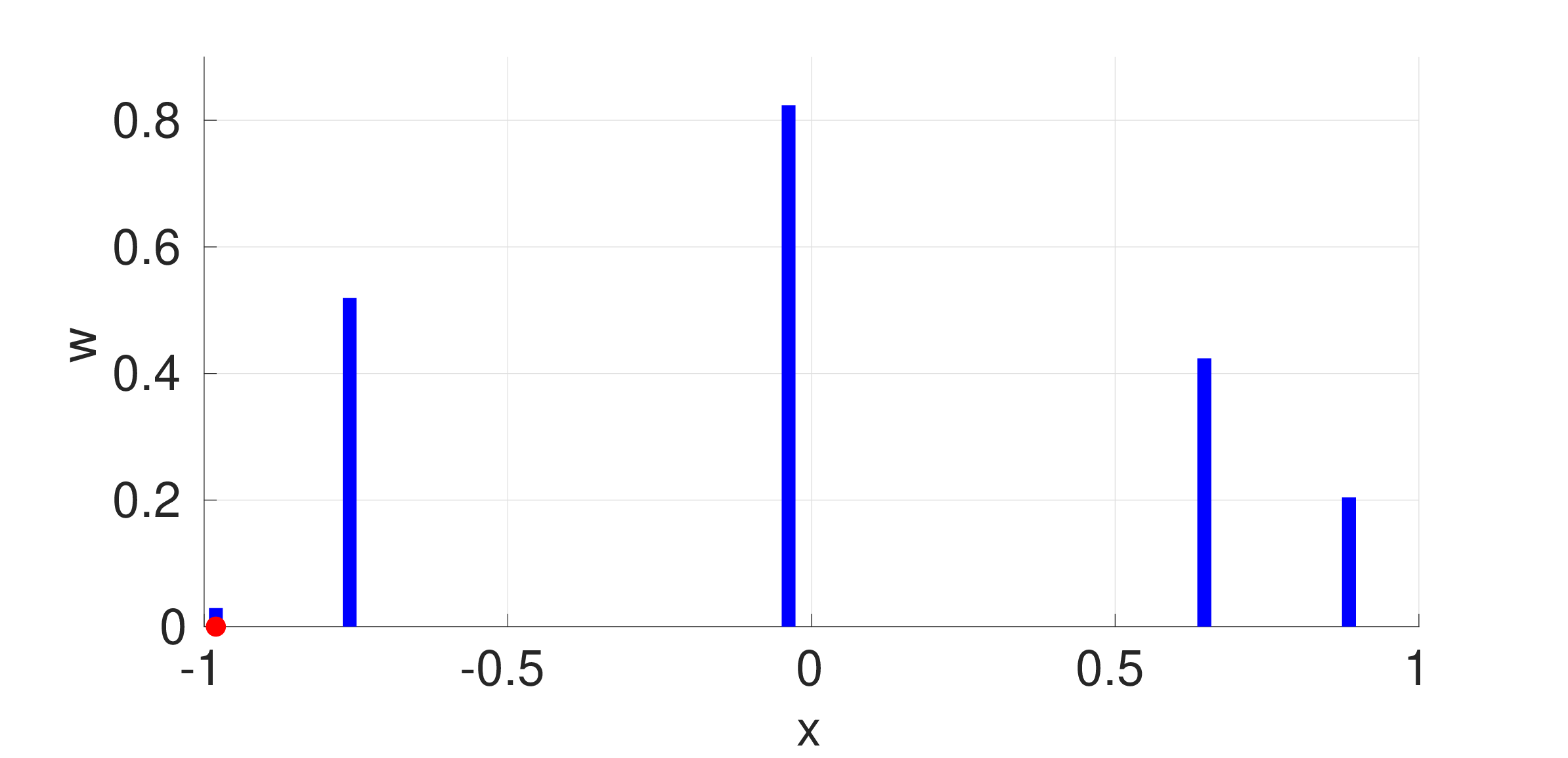}}   
     \subfigure[$m=4$ points ($t=1,k=4$) ]{\label{fig:CECM3}\includegraphics[width=0.30\textwidth]{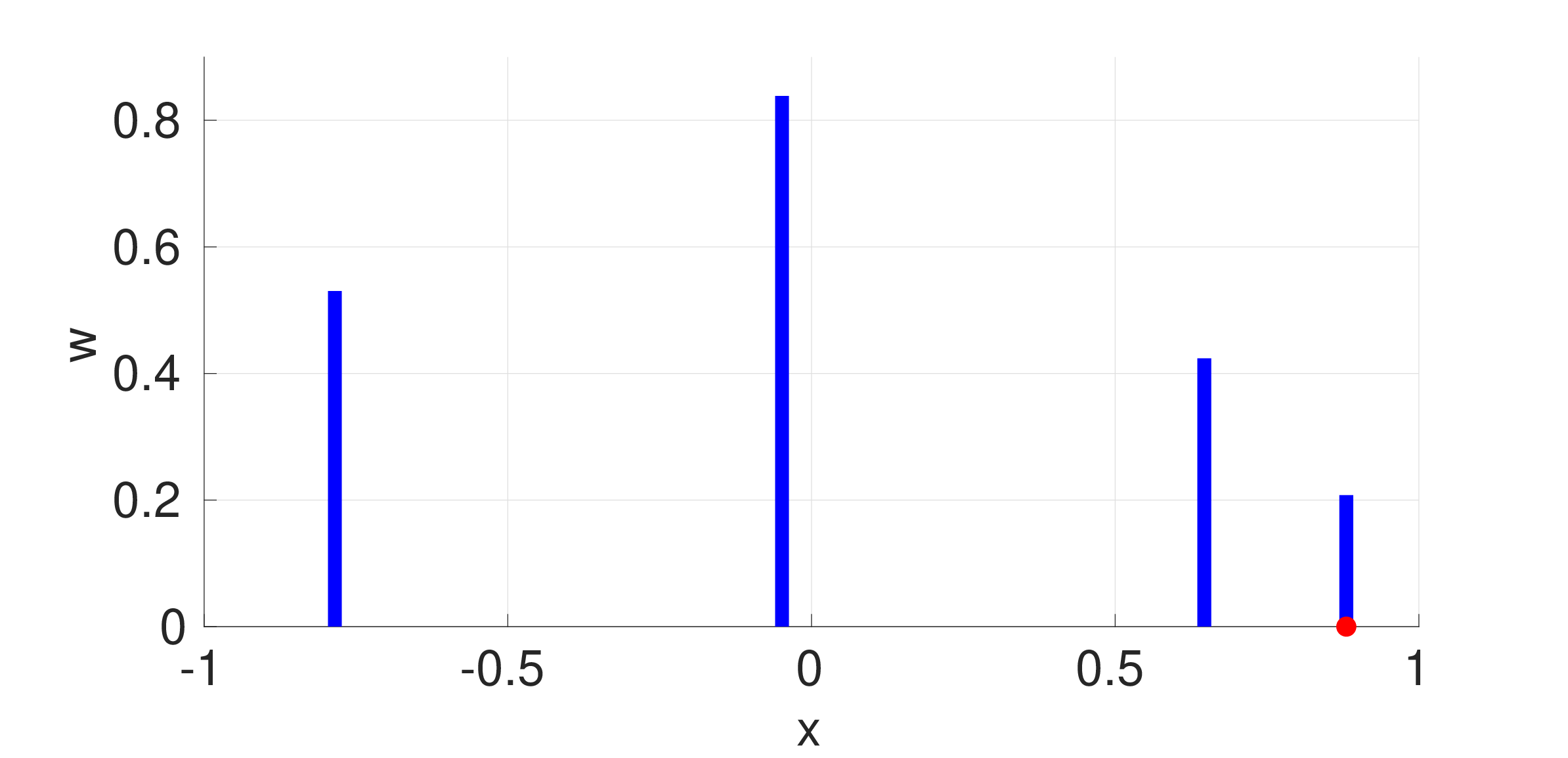}}   
     \subfigure[Final CECM rule, $m=3$ points ($t=1,k=6$) ]{\label{fig:CECM4}\includegraphics[width=0.4\textwidth]{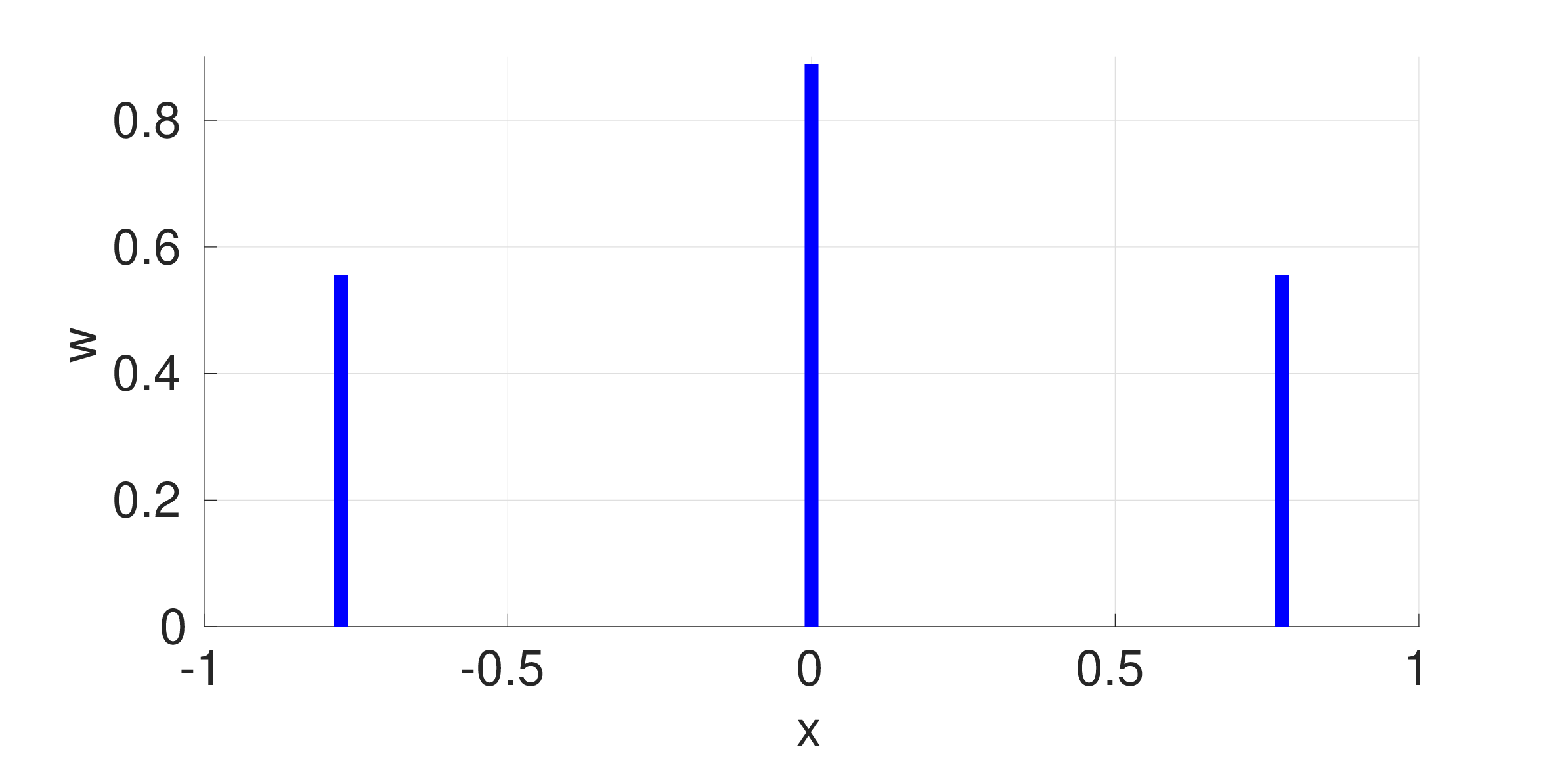}}   
   \caption{ Location   and   weights of the quadrature rules generated during  the sparsification process for the case of univariate polynomials of  degree $p = 5$ in $\Omega = [-1,1]$. Variable $t$ represents the number of passess over the loop that selects the weights to be zeroed in Algorithm \ref{alg:003} (see line \ref{alg:003g}), whereas $k$ indicates the total number of iterations of the modified Newton-Raphson scheme in  Algorithm \ref{alg:005}.     The initial  integration rule, displayed in Figure \ref{fig:CECM1} is the one corresponding to the last iteration of the DECM, displayed previously in  Figure \ref{fig:DECM6}.  The red point in each graph indicates the point whose weight is to be zeroed in the following step.  The final quadrature rule is the one shown in graph \ref{fig:CECM4}.   The   values of the coordinates and weights are given in Table \ref{tab:1}, in which one can see that this quadrature rule is  indeed the optimal 3-points Gauss rule.    }
  \label{fig:FIG_5}
\end{figure}

It can be inferred from the information displayed  in Table  \ref{tab:1}   that,   for polynomials of even degree, the CECM provides rules whose number of points is equal to $(p+2)/2$, whereas for polynomials of odd order, the number of points is equal to $(p+1)/2$ in all cases. Thus, for instance, both CECM rules derived from polynomials of degree 4 and 5   possess 3 points; notice that the rule for the polynomials of degree 4 is asymmetric, whereas the one for polynomials of degree 5   is symmetrical.     Furthermore,  comparison of this symmetrical rule with the  corresponding Gauss rule with the same number of points\footnote{A procedure for determining Gauss rules of arbitrary number of points is given in  Ref. \cite{golub2009matrices}, page 86. }  reveals that they are identical (relative error below $10^{-15}$). The same trend is observed for the remaining CECM rules for polynomials of odd degree.  
To further corroborate this finding, we extended the   study to cover the cases of polynomials from $p=13$   to $p=25$, and the result was invariably the same.  Thus, we can assert that, at least for univariate polynomials, \emph{the proposed CECM is able to arrive at the optimal cubature rule}, that is, the rule with the \emph{minimal number of   points }. To gain further insight into the performance of the method, we present in Figure \ref{fig:FIG_5} the sequence of rules produced during   the sparsification process (from the 6-points DECM   rule (Figure \ref{fig:CECM1}) to the optimal 3-points (Gauss) rule of Figure \ref{fig:CECM4}.

 \subsection{Multivariate polynomials}
 \label{sec:multo}

 \jahoHIDE{ See SummaryExamples2023.mlx} 
 
Let us now extend the preceding assessment to the  integration of    \emph{multivariate Lagrange polynomials} in 2D and 3D cartesian domains --- for which it is known that the optimal rules are tensor product of  univariate Gauss rules \cite{golub2009matrices}. More specifically, we shall focus here  on bivariate and trivariate Lagrange polynomials on biunit squares ($\Omega = [-1,1] \times [-1,1]$) and cubes  ($\Omega = [-1,1] \times [-1,1] \times [-1,1]$), respectively.

\begin{table}[]
\centering
\footnotesize
\begin{tabular}{cccccc|c|cc|c|c|}
\cline{1-5} \cline{7-11} 
\multicolumn{1}{|c|}{\multirow{2}{*}{d}} &
  \multicolumn{2}{c|}{positions} &
  \multicolumn{1}{c|}{\multirow{2}{*}{weights}} &
  \multicolumn{1}{c|}{\multirow{2}{*}{error (Gs.)}} &
  \cellcolor{white}\multirow{29}{*}{} &
  d &
  \multicolumn{2}{c|}{positions} &
  \multirow{2}{*}{weights} &
  \multirow{2}{*}{error (Gs.)} \\ \cline{2-3} \cline{7-9}
\multicolumn{1}{|c|}{} &
  \multicolumn{1}{c|}{x} &
  \multicolumn{1}{c|}{y} &
  \multicolumn{1}{c|}{} &
  \multicolumn{1}{c|}{} &
   &
   &
  \multicolumn{1}{c|}{x} &
  y &
   &
   \\ \cline{1-5} \cline{7-11} 
\multicolumn{1}{|l|}{1} &
  \multicolumn{1}{l|}{0.000000000} &
  \multicolumn{1}{l|}{-0.000000000} &
  \multicolumn{1}{l|}{4.000000000} &
  \multicolumn{1}{l|}{1.1104e-15} &
   &
  \multirow{16}{*}{6} &
  \multicolumn{1}{l|}{0.835410647} &
  0.241939775 &
  0.295169258 &
  \multirow{16}{*}{} \\ \cline{1-5} \cline{8-10}
\multicolumn{1}{|l|}{\multirow{4}{*}{2}} &
  \multicolumn{1}{l|}{0.706474840} &
  \multicolumn{1}{l|}{0.648819610} &
  \multicolumn{1}{l|}{0.707815861} &
  \multicolumn{1}{l|}{\multirow{4}{*}{}} &
   &
   &
  \multicolumn{1}{l|}{0.238387911} &
  0.241939775 &
  0.521300430 &
   \\ \cline{2-4} \cline{8-10}
\multicolumn{1}{|l|}{} &
  \multicolumn{1}{l|}{0.505956954} &
  \multicolumn{1}{l|}{-0.513753481} &
  \multicolumn{1}{l|}{1.262661689} &
  \multicolumn{1}{l|}{} &
   &
   &
  \multicolumn{1}{l|}{0.839282581} &
  -0.465525643 &
  0.250410863 &
   \\ \cline{2-4} \cline{8-10}
\multicolumn{1}{|l|}{} &
  \multicolumn{1}{l|}{-0.471826192} &
  \multicolumn{1}{l|}{0.648819610} &
  \multicolumn{1}{l|}{1.059826914} &
  \multicolumn{1}{l|}{} &
   &
   &
  \multicolumn{1}{l|}{0.255295202} &
  -0.465525643 &
  0.443797308 &
   \\ \cline{2-4} \cline{8-10}
\multicolumn{1}{|l|}{} &
  \multicolumn{1}{l|}{-0.658817575} &
  \multicolumn{1}{l|}{-0.513753481} &
  \multicolumn{1}{l|}{0.969695536} &
  \multicolumn{1}{l|}{} &
   &
   &
  \multicolumn{1}{l|}{0.463842461} &
  0.836215301 &
  0.255076972 &
   \\ \cline{1-5} \cline{8-10}
\multicolumn{1}{|l|}{\multirow{4}{*}{3}} &
  \multicolumn{1}{l|}{0.577350269} &
  \multicolumn{1}{l|}{0.577350269} &
  \multicolumn{1}{l|}{1.000000000} &
  \multicolumn{1}{l|}{\multirow{4}{*}{2.0914e-15}} &
   &
   &
  \multicolumn{1}{l|}{0.909197293} &
  0.836215301 &
  0.100975680 &
   \\ \cline{2-4} \cline{8-10}
\multicolumn{1}{|l|}{} &
  \multicolumn{1}{l|}{0.577350269} &
  \multicolumn{1}{l|}{-0.577350269} &
  \multicolumn{1}{l|}{1.000000000} &
  \multicolumn{1}{l|}{} &
   &
   &
  \multicolumn{1}{l|}{0.528704543} &
  -0.910089745 &
  0.157360912 &
   \\ \cline{2-4} \cline{8-10}
\multicolumn{1}{|l|}{} &
  \multicolumn{1}{l|}{-0.577350269} &
  \multicolumn{1}{l|}{0.577350269} &
  \multicolumn{1}{l|}{1.000000000} &
  \multicolumn{1}{l|}{} &
   &
   &
  \multicolumn{1}{l|}{0.952103717} &
  -0.910089745 &
  0.044168747 &
   \\ \cline{2-4} \cline{8-10}
\multicolumn{1}{|l|}{} &
  \multicolumn{1}{l|}{-0.577350269} &
  \multicolumn{1}{l|}{-0.577350269} &
  \multicolumn{1}{l|}{1.000000000} &
  \multicolumn{1}{l|}{} &
   &
   &
  \multicolumn{1}{l|}{-0.471295457} &
  0.241939775 &
  0.451263578 &
   \\ \cline{1-5} \cline{8-10}
\multicolumn{1}{|l|}{\multirow{9}{*}{4}} &
  \multicolumn{1}{l|}{0.235019870} &
  \multicolumn{1}{l|}{0.887041325} &
  \multicolumn{1}{l|}{0.346272588} &
  \multicolumn{1}{l|}{\multirow{9}{*}{}} &
   &
   &
  \multicolumn{1}{l|}{-0.242986893} &
  0.836215301 &
  0.293371414 &
   \\ \cline{2-4} \cline{8-10}
\multicolumn{1}{|l|}{} &
  \multicolumn{1}{l|}{0.928559333} &
  \multicolumn{1}{l|}{0.872658207} &
  \multicolumn{1}{l|}{0.111090783} &
  \multicolumn{1}{l|}{} &
   &
   &
  \multicolumn{1}{l|}{-0.444666035} &
  -0.465525643 &
  0.391121725 &
   \\ \cline{2-4} \cline{8-10}
\multicolumn{1}{|l|}{} &
  \multicolumn{1}{l|}{0.235019870} &
  \multicolumn{1}{l|}{0.192355601} &
  \multicolumn{1}{l|}{0.936139040} &
  \multicolumn{1}{l|}{} &
   &
   &
  \multicolumn{1}{l|}{-0.205611449} &
  -0.910089745 &
  0.185331125 &
   \\ \cline{2-4} \cline{8-10}
\multicolumn{1}{|l|}{} &
  \multicolumn{1}{l|}{0.928559333} &
  \multicolumn{1}{l|}{0.175016914} &
  \multicolumn{1}{l|}{0.279848771} &
  \multicolumn{1}{l|}{} &
   &
   &
  \multicolumn{1}{l|}{-0.913220650} &
  0.241939775 &
  0.173212013 &
   \\ \cline{2-4} \cline{8-10}
\multicolumn{1}{|l|}{} &
  \multicolumn{1}{l|}{-0.703200858} &
  \multicolumn{1}{l|}{0.889364290} &
  \multicolumn{1}{l|}{0.251544803} &
  \multicolumn{1}{l|}{} &
   &
   &
  \multicolumn{1}{l|}{-0.899639238} &
  -0.465525643 &
  0.166229135 &
   \\ \cline{2-4} \cline{8-10}
\multicolumn{1}{|l|}{} &
  \multicolumn{1}{l|}{-0.703200858} &
  \multicolumn{1}{l|}{0.195018638} &
  \multicolumn{1}{l|}{0.687833785} &
  \multicolumn{1}{l|}{} &
   &
   &
  \multicolumn{1}{l|}{-0.836453352} &
  0.836215301 &
  0.165982881 &
   \\ \cline{2-4} \cline{8-10}
\multicolumn{1}{|l|}{} &
  \multicolumn{1}{l|}{0.235019870} &
  \multicolumn{1}{l|}{-0.713942525} &
  \multicolumn{1}{l|}{0.682449449} &
  \multicolumn{1}{l|}{} &
   &
   &
  \multicolumn{1}{l|}{-0.828149253} &
  -0.910089745 &
  0.105227960 &
   \\ \cline{2-4} \cline{7-11} 
\multicolumn{1}{|l|}{} &
  \multicolumn{1}{l|}{0.928559333} &
  \multicolumn{1}{l|}{-0.718476493} &
  \multicolumn{1}{l|}{0.203099967} &
  \multicolumn{1}{l|}{} &
   &
  \multirow{16}{*}{7} &
  \multicolumn{1}{l|}{0.861136312} &
  0.861136312 &
  0.121002993 &
  \multirow{16}{*}{5.7779e-16} \\ \cline{2-4} \cline{8-10}
\multicolumn{1}{|l|}{} &
  \multicolumn{1}{l|}{-0.703200858} &
  \multicolumn{1}{l|}{-0.713255984} &
  \multicolumn{1}{l|}{0.501720814} &
  \multicolumn{1}{l|}{} &
   &
   &
  \multicolumn{1}{l|}{0.339981044} &
  0.861136312 &
  0.226851852 &
   \\ \cline{1-5} \cline{8-10}
\multicolumn{1}{|l|}{\multirow{9}{*}{5}} &
  \multicolumn{1}{l|}{0.774596669} &
  \multicolumn{1}{l|}{0.774596669} &
  \multicolumn{1}{l|}{0.308641975} &
  \multicolumn{1}{l|}{\multirow{9}{*}{5.9957e-16}} &
   &
   &
  \multicolumn{1}{l|}{0.861136312} &
  0.339981044 &
  0.226851852 &
   \\ \cline{2-4} \cline{8-10}
\multicolumn{1}{|l|}{} &
  \multicolumn{1}{l|}{-0.000000000} &
  \multicolumn{1}{l|}{0.774596669} &
  \multicolumn{1}{l|}{0.493827160} &
  \multicolumn{1}{l|}{} &
   &
   &
  \multicolumn{1}{l|}{0.339981044} &
  0.339981044 &
  0.425293303 &
   \\ \cline{2-4} \cline{8-10}
\multicolumn{1}{|l|}{} &
  \multicolumn{1}{l|}{0.774596669} &
  \multicolumn{1}{l|}{-0.000000000} &
  \multicolumn{1}{l|}{0.493827160} &
  \multicolumn{1}{l|}{} &
   &
   &
  \multicolumn{1}{l|}{-0.339981044} &
  0.861136312 &
  0.226851852 &
   \\ \cline{2-4} \cline{8-10}
\multicolumn{1}{|l|}{} &
  \multicolumn{1}{l|}{0.000000000} &
  \multicolumn{1}{l|}{-0.000000000} &
  \multicolumn{1}{l|}{0.790123457} &
  \multicolumn{1}{l|}{} &
   &
   &
  \multicolumn{1}{l|}{0.861136312} &
  -0.339981044 &
  0.226851852 &
   \\ \cline{2-4} \cline{8-10}
\multicolumn{1}{|l|}{} &
  \multicolumn{1}{l|}{-0.774596669} &
  \multicolumn{1}{l|}{0.774596669} &
  \multicolumn{1}{l|}{0.308641975} &
  \multicolumn{1}{l|}{} &
   &
   &
  \multicolumn{1}{l|}{-0.339981044} &
  0.339981044 &
  0.425293303 &
   \\ \cline{2-4} \cline{8-10}
\multicolumn{1}{|l|}{} &
  \multicolumn{1}{l|}{0.774596669} &
  \multicolumn{1}{l|}{-0.774596669} &
  \multicolumn{1}{l|}{0.308641975} &
  \multicolumn{1}{l|}{} &
   &
   &
  \multicolumn{1}{l|}{0.339981044} &
  -0.339981044 &
  0.425293303 &
   \\ \cline{2-4} \cline{8-10}
\multicolumn{1}{|l|}{} &
  \multicolumn{1}{l|}{-0.774596669} &
  \multicolumn{1}{l|}{0.000000000} &
  \multicolumn{1}{l|}{0.493827160} &
  \multicolumn{1}{l|}{} &
   &
   &
  \multicolumn{1}{l|}{-0.861136312} &
  0.861136312 &
  0.121002993 &
   \\ \cline{2-4} \cline{8-10}
\multicolumn{1}{|l|}{} &
  \multicolumn{1}{l|}{-0.000000000} &
  \multicolumn{1}{l|}{-0.774596669} &
  \multicolumn{1}{l|}{0.493827160} &
  \multicolumn{1}{l|}{} &
   &
   &
  \multicolumn{1}{l|}{-0.339981044} &
  -0.339981044 &
  0.425293303 &
   \\ \cline{2-4} \cline{8-10}
\multicolumn{1}{|l|}{} &
  \multicolumn{1}{l|}{-0.774596669} &
  \multicolumn{1}{l|}{-0.774596669} &
  \multicolumn{1}{l|}{0.308641975} &
  \multicolumn{1}{l|}{} &
   &
   &
  \multicolumn{1}{l|}{0.861136312} &
  -0.861136312 &
  0.121002993 &
   \\ \cline{1-5} \cline{8-10}
\multicolumn{5}{l}{\multirow{5}{*}{}} &
  \multirow{5}{*}{} &
   &
  \multicolumn{1}{l|}{-0.861136312} &
  0.339981044 &
  0.226851852 &
   \\ \cline{8-10}
\multicolumn{5}{l}{} &
   &
   &
  \multicolumn{1}{l|}{0.339981044} &
  -0.861136312 &
  0.226851852 &
   \\ \cline{8-10}
\multicolumn{5}{l}{} &
   &
   &
  \multicolumn{1}{l|}{-0.861136312} &
  -0.339981044 &
  0.226851852 &
   \\ \cline{8-10}
\multicolumn{5}{l}{} &
   &
   &
  \multicolumn{1}{l|}{-0.339981044} &
  -0.861136312 &
  0.226851852 &
   \\ \cline{8-10}
\multicolumn{5}{l}{} &
   &
   &
  \multicolumn{1}{l|}{-0.861136312} &
  -0.861136312 &
  0.121002993 &
   \\ \cline{7-11} 
\end{tabular}
\caption{Cubature rules computed by the CECM for bivariate polynomials   of degree up to $p = 7$ (for both variables) in in $\Omega = [-1,1] \times [-1,1]$. The rightmost column represents the relative deviations with respect to the optimal  Gauss product rules (for polynomials of odd degree).  }
\label{tab:2}
\end{table}

\jahoHIDE{Generated by \url{/home/joaquin/Desktop/CURRENT_TASKS/PAPERS_2020_onwards/OptimalECM/MATLAB_EXAMPLES_CODE/EXAMPLES_2D/Generate_table2D.m}. See table Table\_2.tgn. Created using \url{https://www.tablesgenerator.com/}}

  Given a degree $p$, and a set of $p+1$ equally spaced nodes for each direction, let us define the monomials: 
    \begin{equation} 
    \label{eq:res100}
     \Gamma_i(x)  = \prod_{h=1, i\neq h}^{P}  \dfrac{x-x_h}{x_i-x_h}, \hspace{1cm}  i=1,2 \ldots (p+1)
    \end{equation}
       \begin{equation} 
    \label{eq:res101}
     \Gamma_j(y)  = \prod_{h=1, j\neq h}^{P}  \dfrac{y-y_h}{y_j-y_h}, \hspace{1cm}  j=1,2 \ldots (p+1)
    \end{equation}
       \begin{equation} 
    \label{eq:res102}
     \Gamma_k(z)  = \prod_{h=1, k\neq h}^{P}  \dfrac{z-z_h}{z_k-z_h}, \hspace{1cm}  k=1,2 \ldots (p+1).
    \end{equation}
 The expression for the $P = (p+1)^2$ integrand functions for the case of bivariate polynomials is given by     
    \begin{equation}
     a_l(x,y) =  \Gamma_i(x)   \Gamma_j(y),  \hspace{1cm}    l = (j-1) (p+1) +  i, \;\;\; i,j= 1,2 \ldots ( p+1), 
    \end{equation}
  whereas  for    trivariate polymomials, the $P = (p+1)^3$ integrand functions adopt the expression
     \begin{equation}
     a_l(x,y,z) =  \Gamma_i(x)   \Gamma_j(y)   \Gamma_k(z),  \hspace{1cm}    l = (k-1)(p+1)^2 + (j-1)(p+1) + i , \;\;\; i,j,k= 1,2 \ldots  (p+1).  
    \end{equation}

           \begin{figure}[!ht]
  \centering
  \subfigure[DECM rule, $m = 16$ points]{\label{fig:2CECM1}\includegraphics[width=0.30\textwidth]{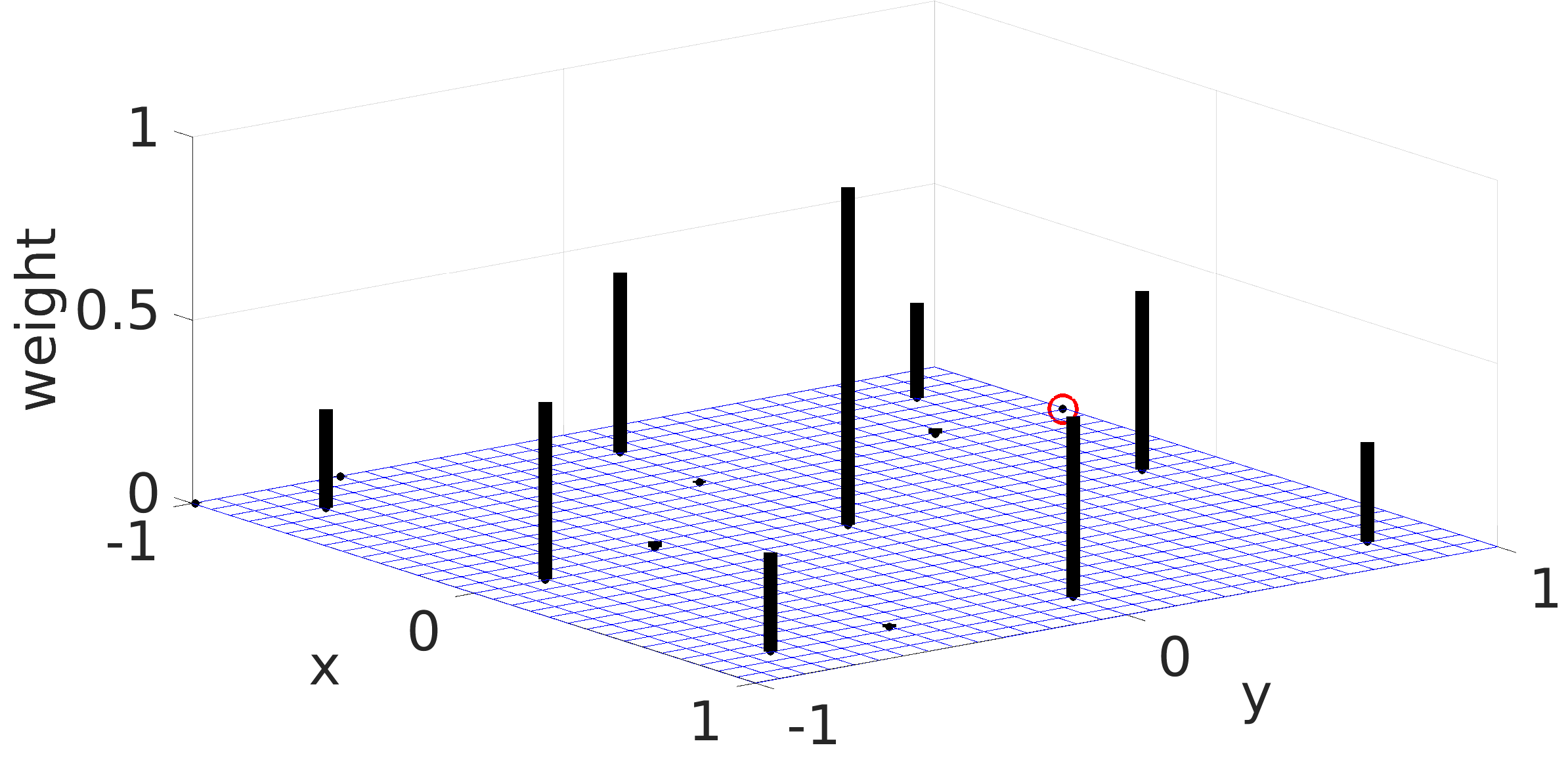}}   
     \subfigure[$m=15$ points ($t=1$,$k=3$)]{\label{fig:2CECM2}\includegraphics[width=0.30\textwidth]{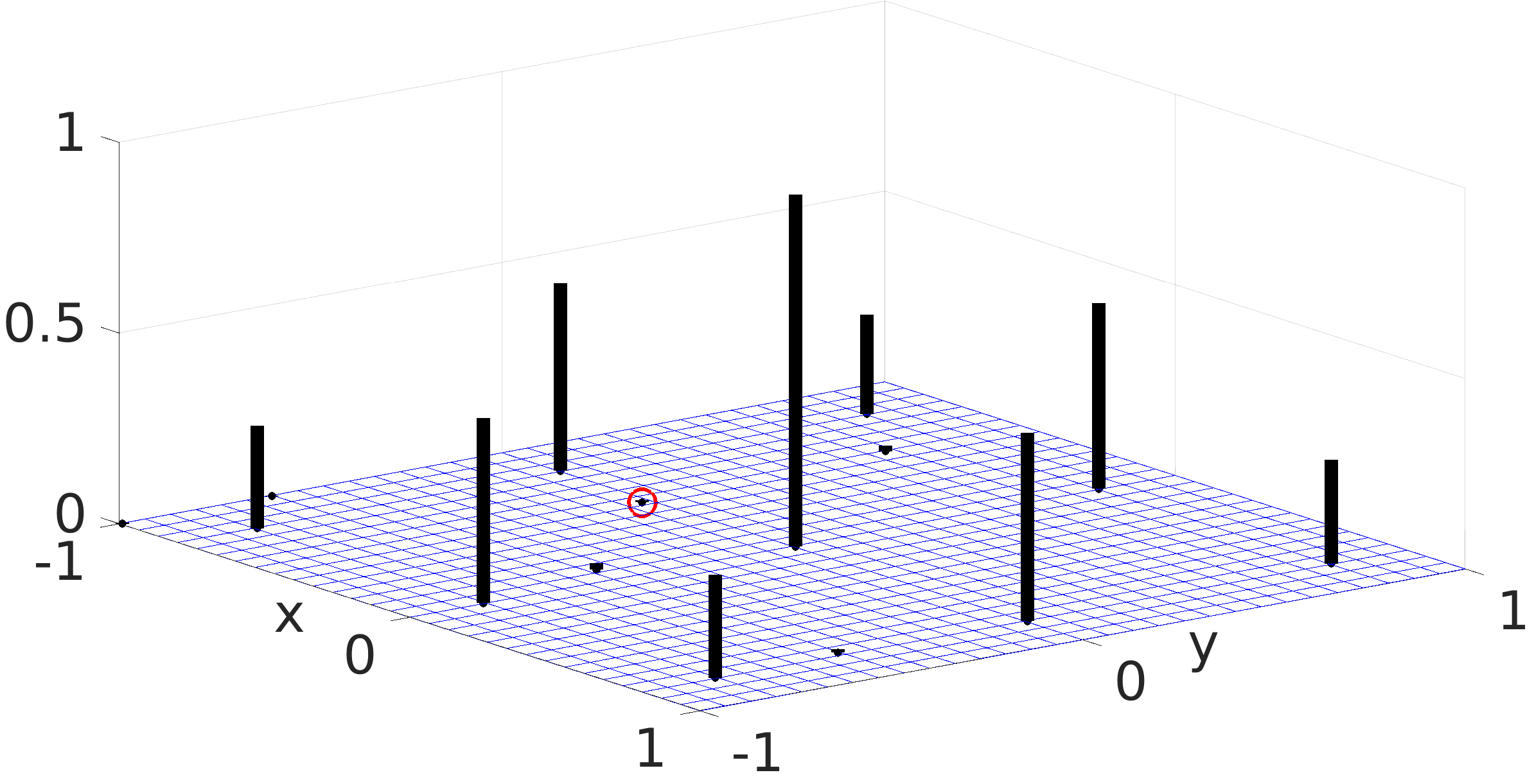}}   
     \subfigure[$m=14$ points ($t=1$,$k=3$)]{\label{fig:2CECM3}\includegraphics[width=0.30\textwidth]{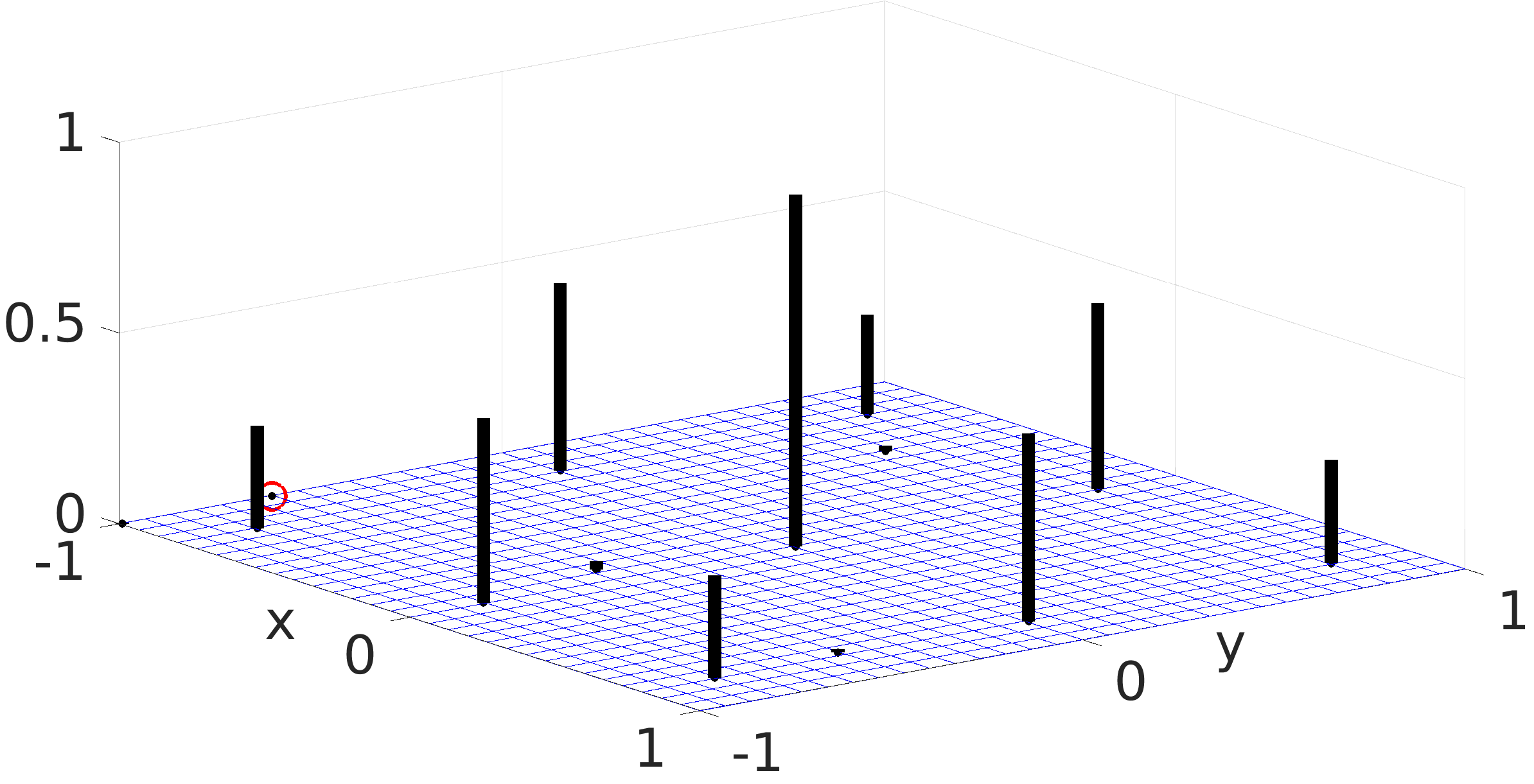}}   
     \subfigure[$m=13$ points ($t=1$,$k=3$)]{\label{fig:2CECM4}\includegraphics[width=0.30\textwidth]{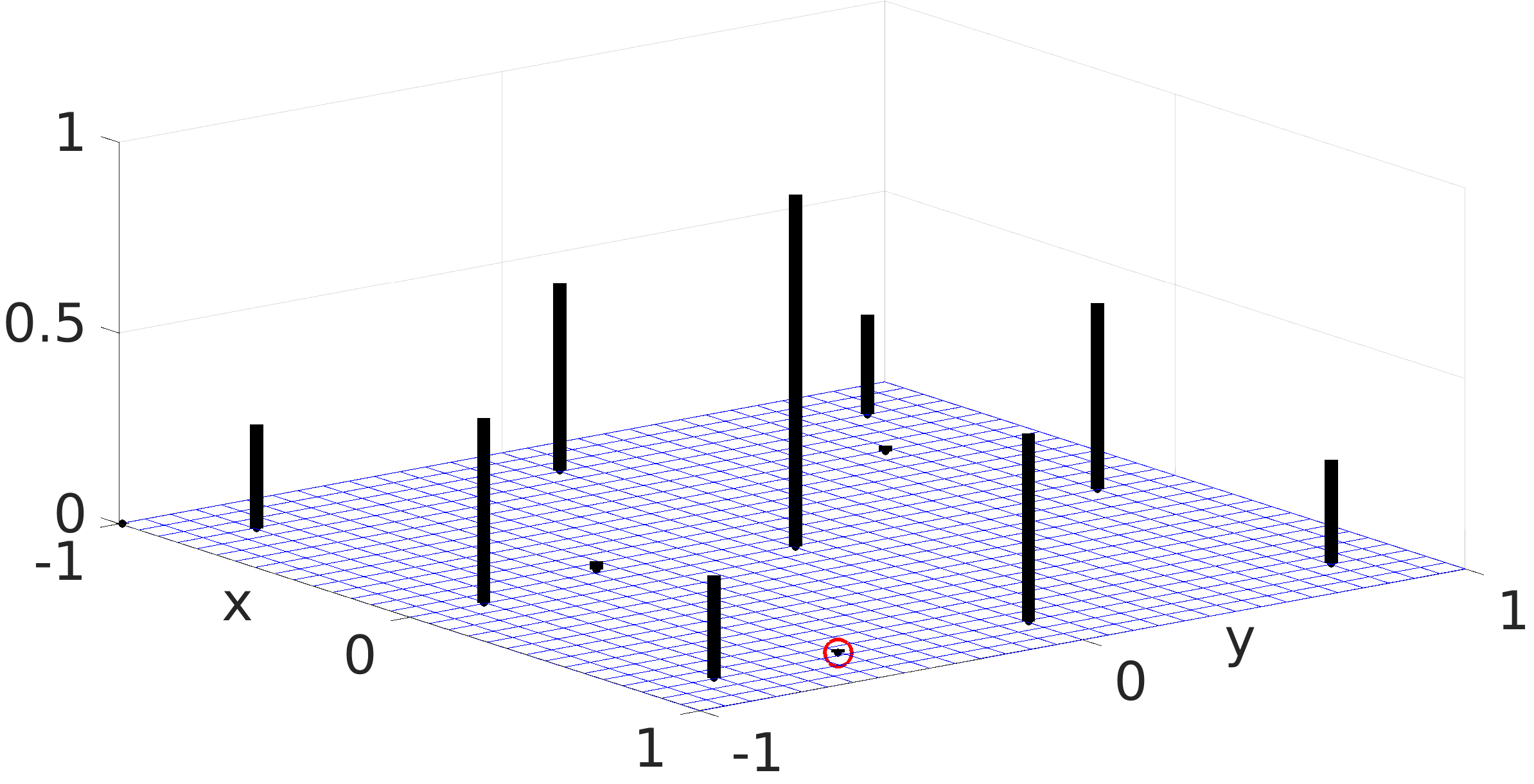}}   
     \subfigure[$m=12$ points ($t=1$,$k=4$)]{\label{fig:2CECM5}\includegraphics[width=0.30\textwidth]{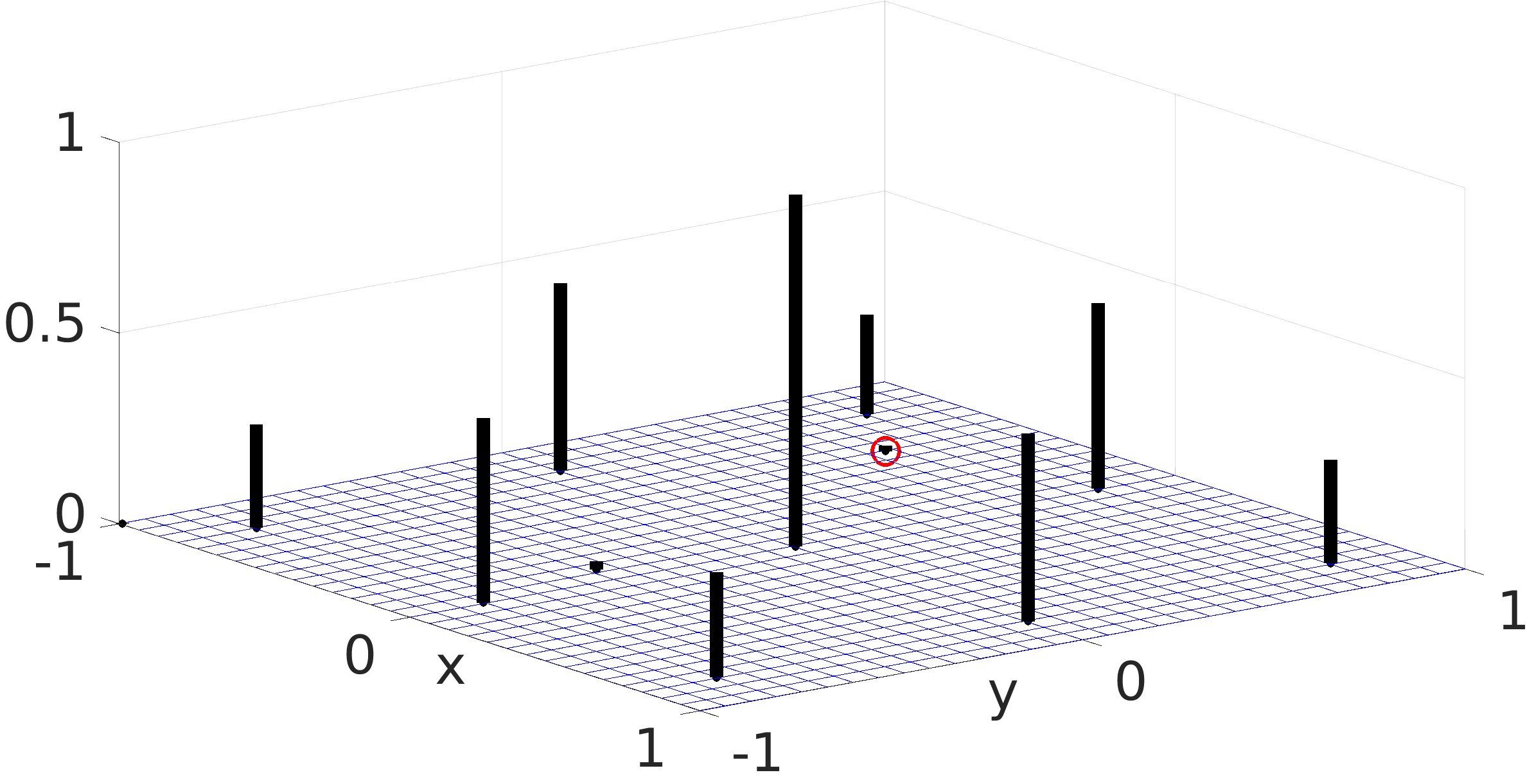}}   
     \subfigure[$m=11$ points ($t=1$,$k=4$)]{\label{fig:2CECM6}\includegraphics[width=0.30\textwidth]{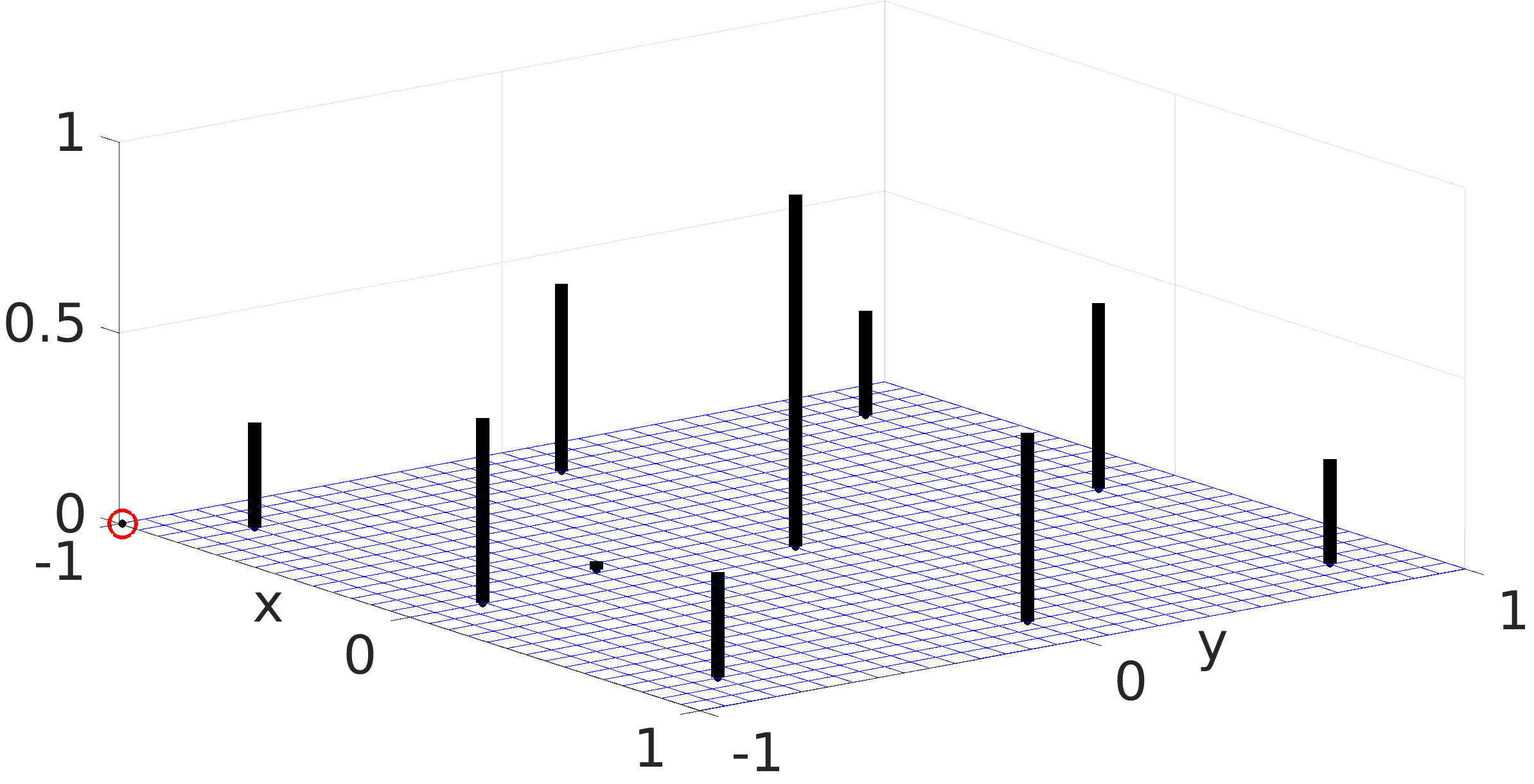}}  
      \subfigure[$m=10$ points  ($t=1$,$k=4$) ]{\label{fig:2CECM7}\includegraphics[width=0.30\textwidth]{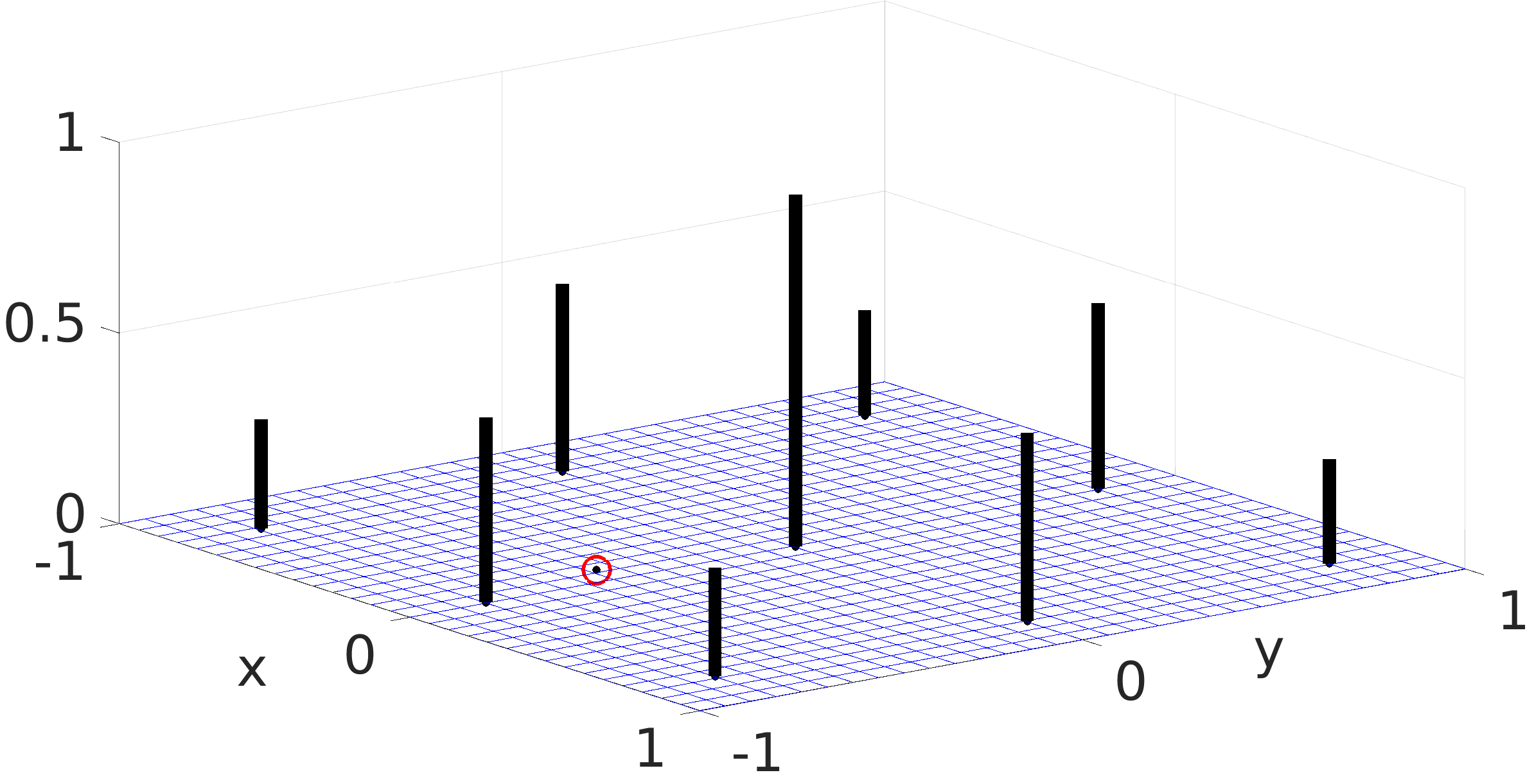}}   
     \subfigure[$m=9$ points  ($t=1$,$k=3$)]{\label{fig:2CECM8}\includegraphics[width=0.30\textwidth]{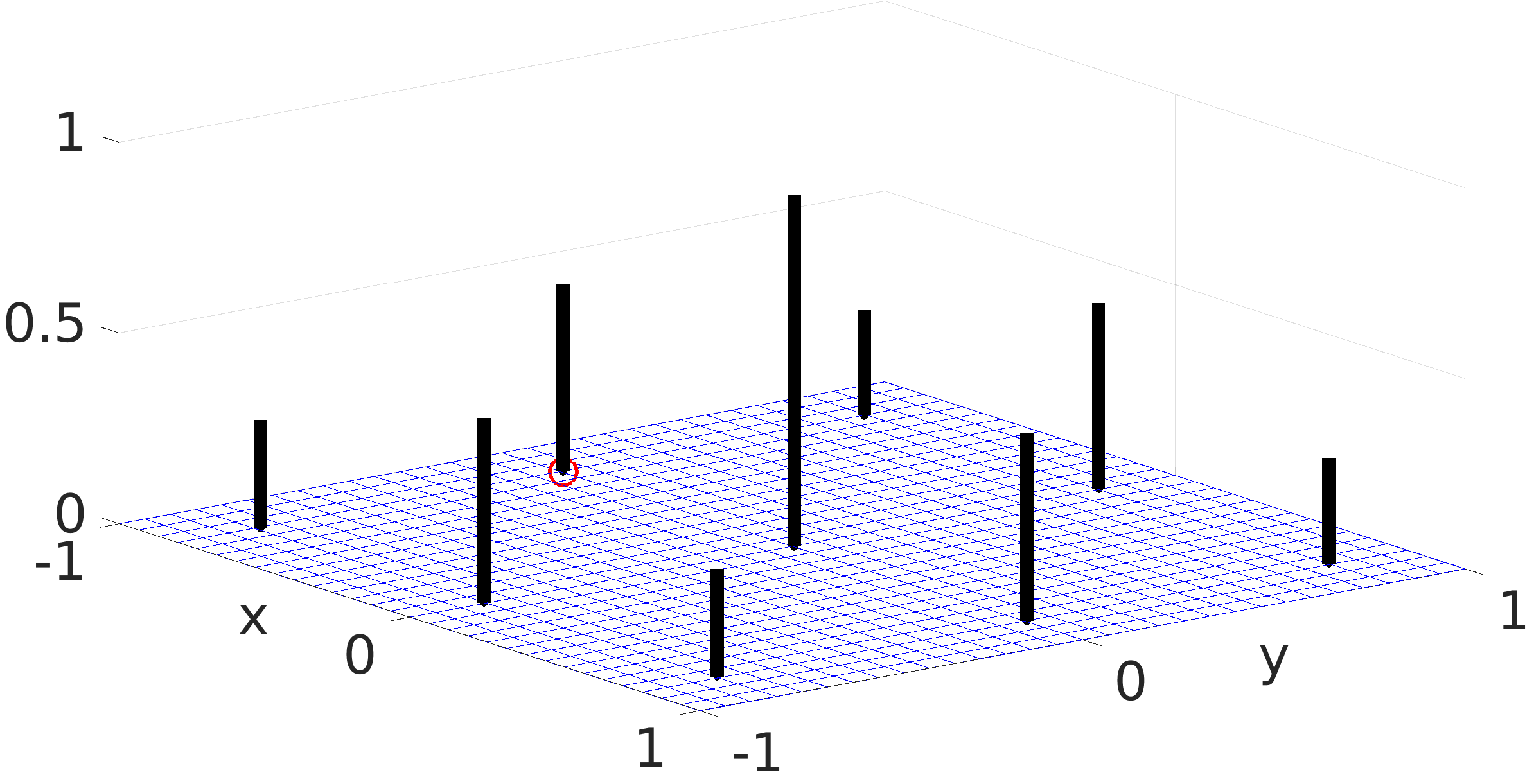}}   
     \subfigure[$m=8$ points  ($t=5$,$k=6$)]{\label{fig:2CECM9}\includegraphics[width=0.30\textwidth]{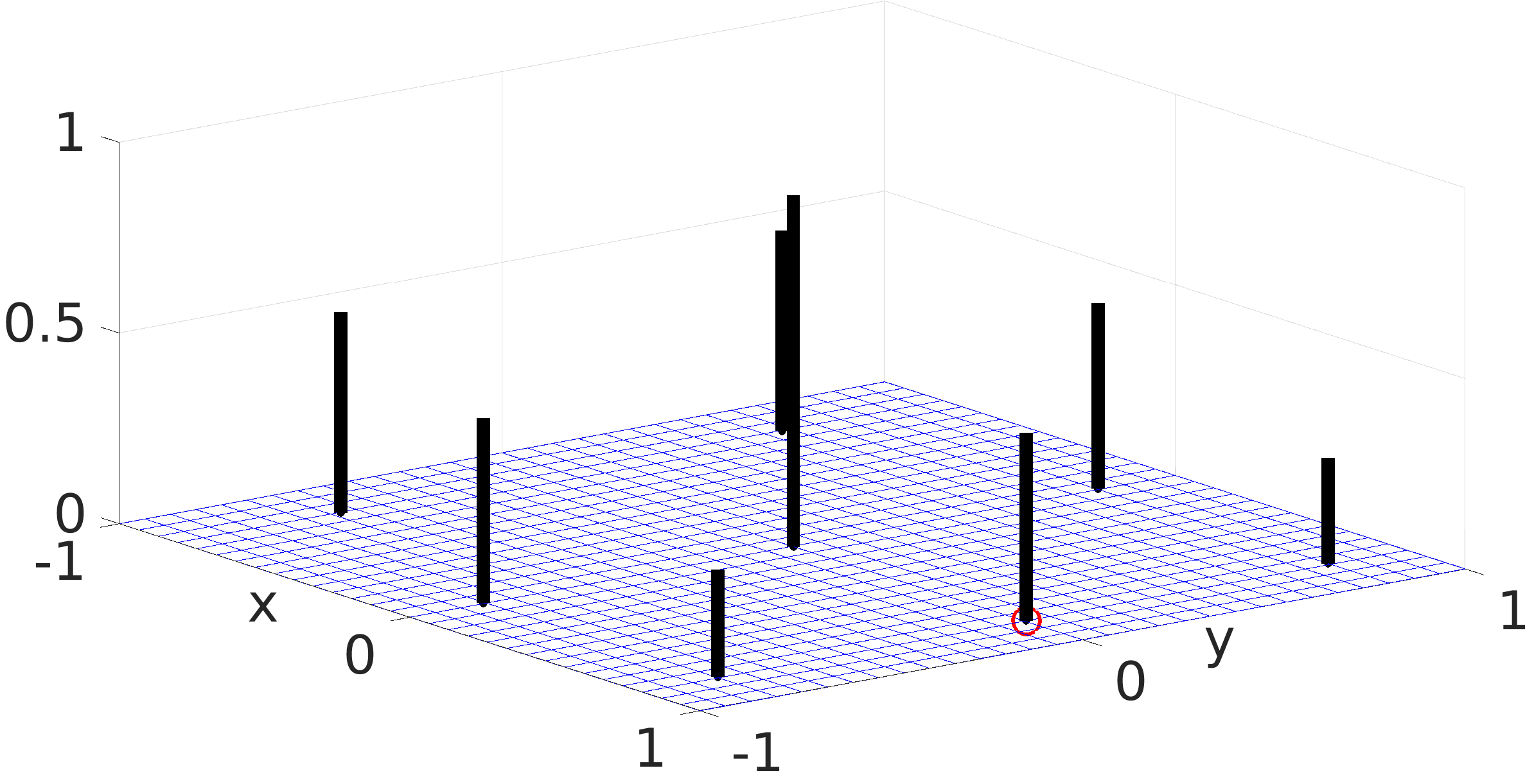}} 
       \subfigure[$m=7$ points  ($t=3$,$k=6$)]{\label{fig:2CECM10}\includegraphics[width=0.30\textwidth]{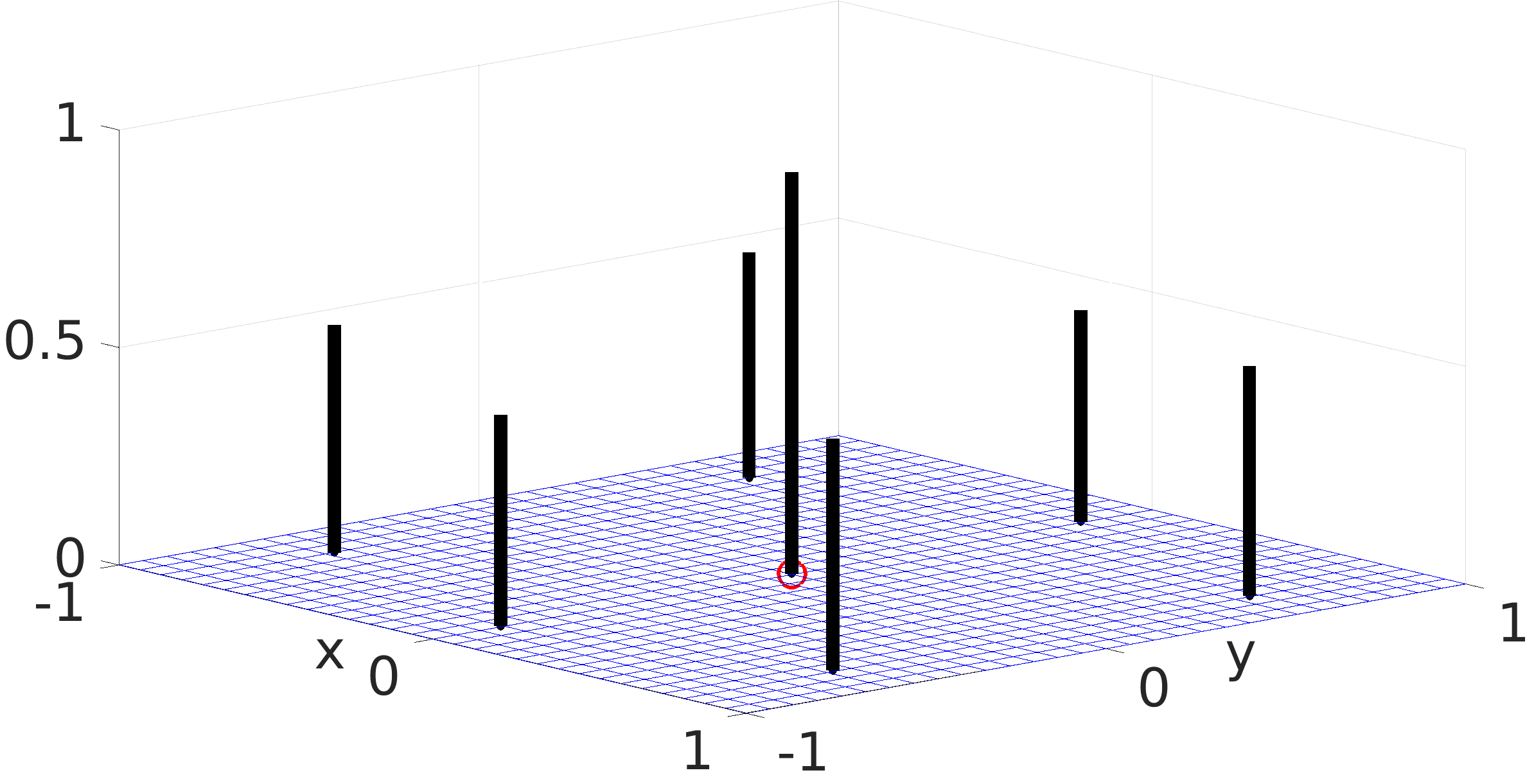}}   
     \subfigure[$m=6$ points ($t=3$,$k=5$)]{\label{fig:2CECM11}\includegraphics[width=0.30\textwidth]{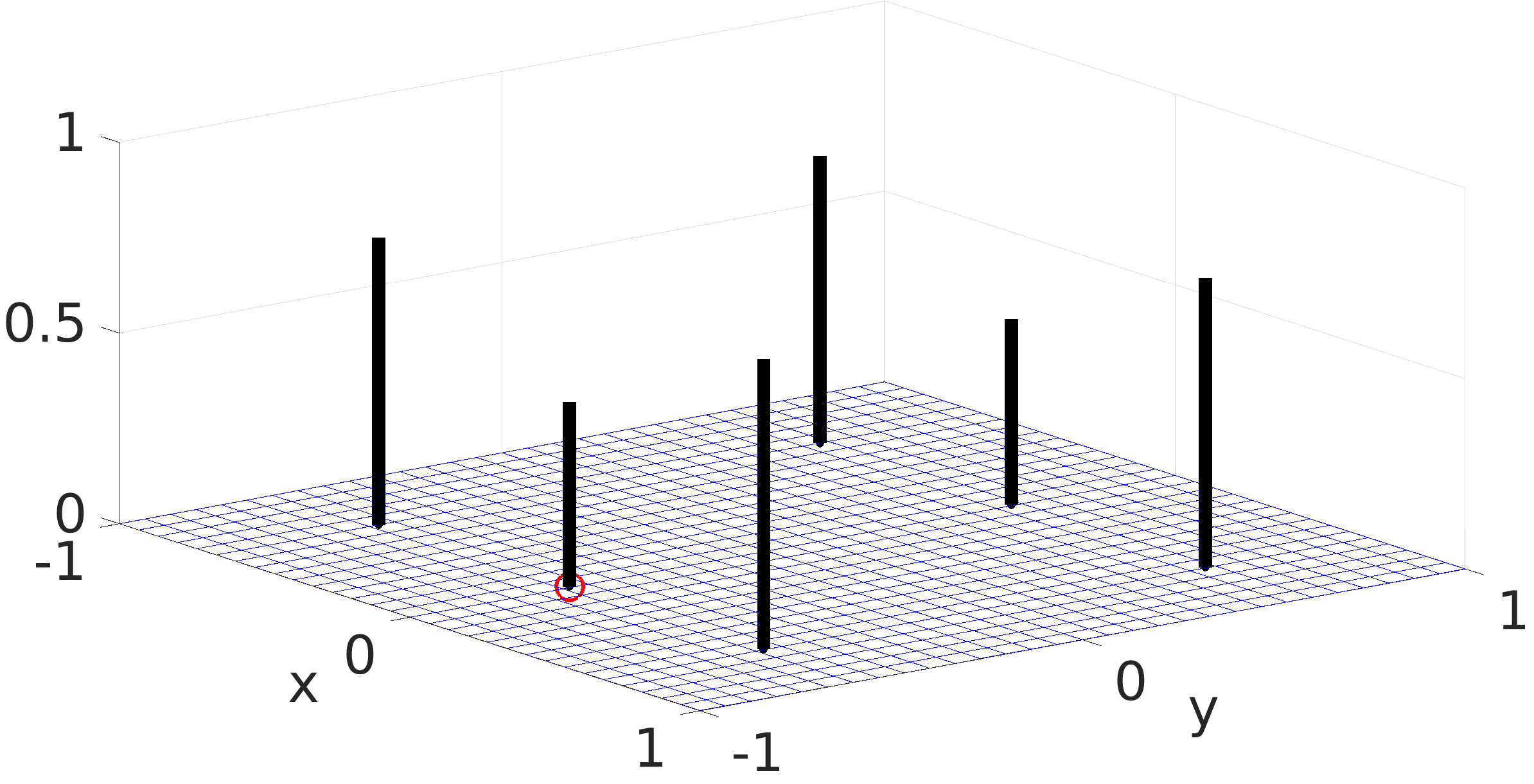}}   
     \subfigure[$m=5$ points ($t=1$, $k=5$)]{\label{fig:2CECM12}\includegraphics[width=0.30\textwidth]{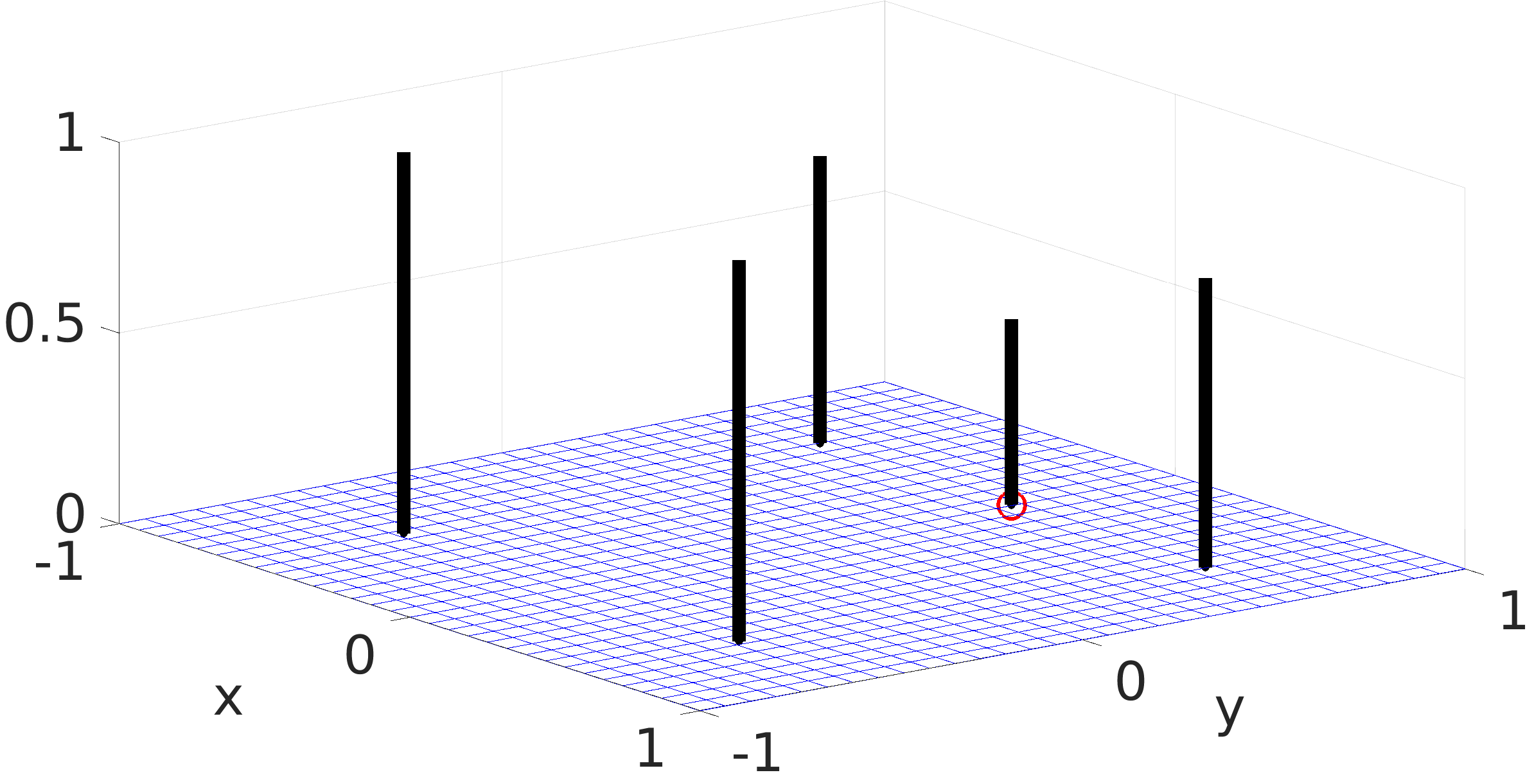}}
       \subfigure[Final CECM rule, $m=4$ points  ($t=1$,$k=5$)]{\label{fig:2CECM13}\includegraphics[width=0.4\textwidth]{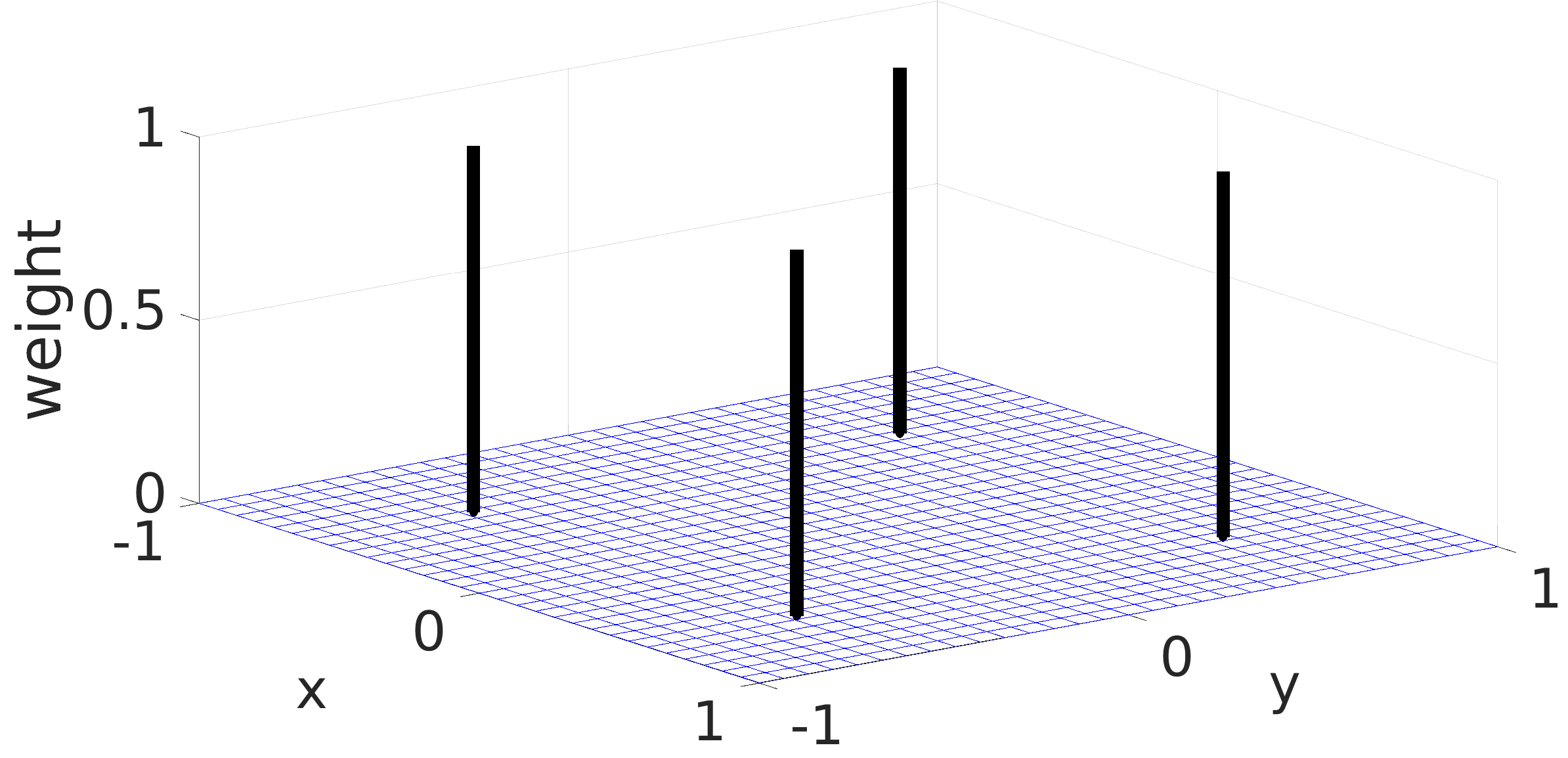}}
   \caption{ Locations  and   weights   of the cubature rules generated during  the sparsification process for the case of bivariate polynomials of  degree $p = 3$ in $\Omega = [-1,1] \times [-1,1]$.   The red circle in each graph indicates the point whose weight is to be zeroed in the following step. Variables $t$ and $k$, on the other hand, have the same interpretation as in Figure \ref{fig:FIG_5}.  The initial DECM rule has $(p+1)^2 = 16$ points (see Figure \ref{fig:2CECM1}), while  the final   rule features $(p+1)^2/2^2 = 4$ points, see graph \ref{fig:2CECM12}. The  exact values of the coordinates and weights are given in Table \ref{tab:2}  ( it can be seen that the CECM rule   coincides with  the standard $2 \times 2$  product Gauss   rule).        }
  \label{fig:FIG_6}
\end{figure}
        
    We use structured meshes  of  $20 \times 20$   quadrilateral elements for the square, and $20 \times 20 \times 20 $ hexahedra elements for the cube, with  Gauss integration rules for each element of $2 \times 2$ and $2 \times 2 \times 2$ points, respectively. The   parameters governing the performance of the CECM are the same employed in the univariate case.  We examined the rules computed by the CECM for degrees up to $p=12$ for both 2D and 3D cases (the same degree for all variables). For reasons of space limitation, we only display the coordinates and weights up to $p = 7$   for the 2D case (see Table \ref{tab:2}) , and in Table \ref{tab:3} up to $p=4$ for the 3D case. The comparison with the product Gauss rules contained in both   tables  reveals the very same  pattern  observed in the case of univariate polynomials: for even degrees, the CECM produces asymmetrical rules, whereas for odd degrees, the CECM produces symmetrical rules identical to the corresponding product Gauss rules (featuring in both cases $(p+1)^d/2^d$ points, where $d=2,3$).   Although not shown here, the same trend was observed for the remaining polynomial degrees.

\begin{table}[!ht]
\centering
\footnotesize
\begin{tabular}{|c|ccc|c|c|}
\hline
\multirow{2}{*}{d} &
  \multicolumn{3}{c|}{positions} &
  \multirow{2}{*}{weights} &
  \multirow{2}{*}{error (Gs.)} \\ \cline{2-4}
 &
  \multicolumn{1}{c|}{x} &
  \multicolumn{1}{c|}{y} &
  z &
   &
   \\ \hline
1 &
  \multicolumn{1}{c|}{8.38519099255208E-17} &
  \multicolumn{1}{c|}{-1.37013188830199E-16} &
  2.54794197642291E-16 &
  8.00000000000022 &
  2.7534e-14 \\ \hline
\multirow{8}{*}{2} &
  \multicolumn{1}{l|}{0.743706246140823} &
  \multicolumn{1}{l|}{0.588904961511756} &
  0.484154086744496 &
  0.86563029280071 &
  \multirow{8}{*}{} \\ \cline{2-5}
 &
  \multicolumn{1}{l|}{0.743706246140821} &
  \multicolumn{1}{l|}{-0.566022287327392} &
  0.831399263412485 &
  0.499063612121599 &
   \\ \cline{2-5}
 &
  \multicolumn{1}{l|}{-0.44820563907193} &
  \multicolumn{1}{l|}{0.741835057873505} &
  0.830045419614278 &
  0.613961492138332 &
   \\ \cline{2-5}
 &
  \multicolumn{1}{l|}{0.743706246140824} &
  \multicolumn{1}{l|}{0.588904961511757} &
  -0.688486047024208 &
  0.608724673042735 &
   \\ \cline{2-5}
 &
  \multicolumn{1}{l|}{0.743706246140825} &
  \multicolumn{1}{l|}{-0.56602228732739} &
  -0.400930513175059 &
  1.03489533941434 &
   \\ \cline{2-5}
 &
  \multicolumn{1}{l|}{-0.44820563907193} &
  \multicolumn{1}{l|}{-0.449336183017346} &
  0.830045419614395 &
  1.01362448933892 &
   \\ \cline{2-5}
 &
  \multicolumn{1}{l|}{-0.448205639071929} &
  \multicolumn{1}{l|}{0.739388192080586} &
  -0.401584450027124 &
  1.274239461537 &
   \\ \cline{2-5}
 &
  \multicolumn{1}{l|}{-0.448205639071928} &
  \multicolumn{1}{l|}{-0.450823176382297} &
  -0.401584450027067 &
  2.08986063960634 &
   \\ \hline
\multirow{8}{*}{3} &
  \multicolumn{1}{l|}{0.577350269189625} &
  \multicolumn{1}{l|}{0.577350269189626} &
  0.577350269189626 &
  0.999999999999999 &
  \multirow{8}{*}{4.3425e-16} \\ \cline{2-5}
 &
  \multicolumn{1}{l|}{0.577350269189626} &
  \multicolumn{1}{l|}{0.577350269189626} &
  -0.577350269189626 &
  0.999999999999999 &
   \\ \cline{2-5}
 &
  \multicolumn{1}{l|}{-0.577350269189626} &
  \multicolumn{1}{l|}{0.577350269189626} &
  0.577350269189626 &
  1 &
   \\ \cline{2-5}
 &
  \multicolumn{1}{l|}{-0.577350269189626} &
  \multicolumn{1}{l|}{0.577350269189626} &
  -0.577350269189626 &
  1 &
   \\ \cline{2-5}
 &
  \multicolumn{1}{l|}{0.577350269189626} &
  \multicolumn{1}{l|}{-0.577350269189626} &
  0.577350269189626 &
  1 &
   \\ \cline{2-5}
 &
  \multicolumn{1}{l|}{0.577350269189626} &
  \multicolumn{1}{l|}{-0.577350269189626} &
  -0.577350269189626 &
  1 &
   \\ \cline{2-5}
 &
  \multicolumn{1}{l|}{-0.577350269189626} &
  \multicolumn{1}{l|}{-0.577350269189626} &
  0.577350269189626 &
  1 &
   \\ \cline{2-5}
 &
  \multicolumn{1}{l|}{-0.577350269189626} &
  \multicolumn{1}{l|}{-0.577350269189626} &
  -0.577350269189626 &
  1 &
   \\ \hline
\multirow{27}{*}{4} &
  \multicolumn{1}{l|}{0.691575606960029} &
  \multicolumn{1}{l|}{0.160521574170478} &
  -0.261954268935127 &
  0.701290700192769 &
  \multirow{27}{*}{} \\ \cline{2-5}
 &
  \multicolumn{1}{l|}{-0.282918324688049} &
  \multicolumn{1}{l|}{-0.217978724118607} &
  -0.27678081615893 &
  1.00337289978755 &
   \\ \cline{2-5}
 &
  \multicolumn{1}{l|}{0.691575606960029} &
  \multicolumn{1}{l|}{-0.282053603748311} &
  0.696630178723833 &
  0.561750852289478 &
   \\ \cline{2-5}
 &
  \multicolumn{1}{l|}{-0.28291832468805} &
  \multicolumn{1}{l|}{-0.211998475048984} &
  0.69305243630094 &
  0.733763793522678 &
   \\ \cline{2-5}
 &
  \multicolumn{1}{l|}{-0.282918324688049} &
  \multicolumn{1}{l|}{0.707430491792361} &
  -0.27678081615893 &
  0.734672939267592 &
   \\ \cline{2-5}
 &
  \multicolumn{1}{l|}{0.691575606960029} &
  \multicolumn{1}{l|}{-0.722360355870564} &
  -0.261954268935126 &
  0.50651426374409 &
   \\ \cline{2-5}
 &
  \multicolumn{1}{l|}{0.691575606960028} &
  \multicolumn{1}{l|}{0.691783591238766} &
  0.696630178723832 &
  0.412586162346216 &
   \\ \cline{2-5}
 &
  \multicolumn{1}{l|}{0.691575606960028} &
  \multicolumn{1}{l|}{0.861535244017794} &
  -0.261954268935126 &
  0.294025738847668 &
   \\ \cline{2-5}
 &
  \multicolumn{1}{l|}{-0.28291832468805} &
  \multicolumn{1}{l|}{0.708932471043466} &
  0.693052436300939 &
  0.536831048244394 &
   \\ \cline{2-5}
 &
  \multicolumn{1}{l|}{0.691575606960028} &
  \multicolumn{1}{l|}{0.067249317677481} &
  -0.960515811309942 &
  0.175016341839065 &
   \\ \cline{2-5}
 &
  \multicolumn{1}{l|}{-0.282918324688048} &
  \multicolumn{1}{l|}{-0.910801772049135} &
  -0.27678081615893 &
  0.330495727456156 &
   \\ \cline{2-5}
 &
  \multicolumn{1}{l|}{-0.282918324688049} &
  \multicolumn{1}{l|}{0.26535993911279} &
  -0.980548984419782 &
  0.246322371421189 &
   \\ \cline{2-5}
 &
  \multicolumn{1}{l|}{-0.28291832468805} &
  \multicolumn{1}{l|}{-0.904960687303036} &
  0.693052436300938 &
  0.248651857974508 &
   \\ \cline{2-5}
 &
  \multicolumn{1}{l|}{0.691575606960028} &
  \multicolumn{1}{l|}{-0.988158764971828} &
  0.696630178723831 &
  0.128498060521655 &
   \\ \cline{2-5}
 &
  \multicolumn{1}{l|}{-0.989432805330447} &
  \multicolumn{1}{l|}{-0.053487256448397} &
  -0.23226227865251 &
  0.203195059703459 &
   \\ \cline{2-5}
 &
  \multicolumn{1}{l|}{0.691575606960029} &
  \multicolumn{1}{l|}{0.804707428230073} &
  -0.960515811309941 &
  0.096280486097567 &
   \\ \cline{2-5}
 &
  \multicolumn{1}{l|}{0.691575606960028} &
  \multicolumn{1}{l|}{-0.750168356265182} &
  -0.960515811309941 &
  0.118969763496464 &
   \\ \cline{2-5}
 &
  \multicolumn{1}{l|}{-0.989432805330449} &
  \multicolumn{1}{l|}{-0.053487256448396} &
  0.703880237838895 &
  0.148998882881412 &
   \\ \cline{2-5}
 &
  \multicolumn{1}{l|}{-0.28291832468805} &
  \multicolumn{1}{l|}{-0.695806658208985} &
  -0.980548984419782 &
  0.180892299291507 &
   \\ \cline{2-5}
 &
  \multicolumn{1}{l|}{-0.989432805330447} &
  \multicolumn{1}{l|}{0.754804289030939} &
  0.072783264075214 &
  0.124784457736787 &
   \\ \cline{2-5}
 &
  \multicolumn{1}{l|}{-0.282918324688049} &
  \multicolumn{1}{l|}{0.96494897614471} &
  -0.980548984419783 &
  0.062699664209236 &
   \\ \cline{2-5}
 &
  \multicolumn{1}{l|}{-0.989432805330446} &
  \multicolumn{1}{l|}{-0.797964254050424} &
  -0.226012314072442 &
  0.11467159716885 &
   \\ \cline{2-5}
 &
  \multicolumn{1}{l|}{-0.989432805330448} &
  \multicolumn{1}{l|}{0.754804289030939} &
  0.807516321196852 &
  0.067775700765236 &
   \\ \cline{2-5}
 &
  \multicolumn{1}{l|}{-0.989432805330447} &
  \multicolumn{1}{l|}{-0.053487256448394} &
  -0.925565909634688 &
  0.06232167903138 &
   \\ \cline{2-5}
 &
  \multicolumn{1}{l|}{-0.989432805330448} &
  \multicolumn{1}{l|}{-0.797964254050423} &
  0.705427653813099 &
  0.084038633103628 &
   \\ \cline{2-5}
 &
  \multicolumn{1}{l|}{-0.989432805330448} &
  \multicolumn{1}{l|}{0.754804289030939} &
  -0.748349786829047 &
  0.085270559050749 &
   \\ \cline{2-5}
 &
  \multicolumn{1}{l|}{-0.989432805330448} &
  \multicolumn{1}{l|}{-0.797964254050424} &
  -0.918958907578276 &
  0.036308460008693 &
   \\ \hline
\end{tabular}
\caption{Cubature rules computed by the CECM for trivariate polynomials   of degree up to $p = 4$  in   $\Omega = [-1,1] \times [-1,1] \times [-1,1]$. The rightmost column represents the relative deviations with respect to the optimal  Gauss product rules (for polynomials of odd degree).  }
\label{tab:3}
\end{table}

\begin{figure}[!ht]
  \centering
  \subfigure[DECM rule, $m=64$ points]{\label{fig:3CECM_64}\includegraphics[width=0.30\textwidth]{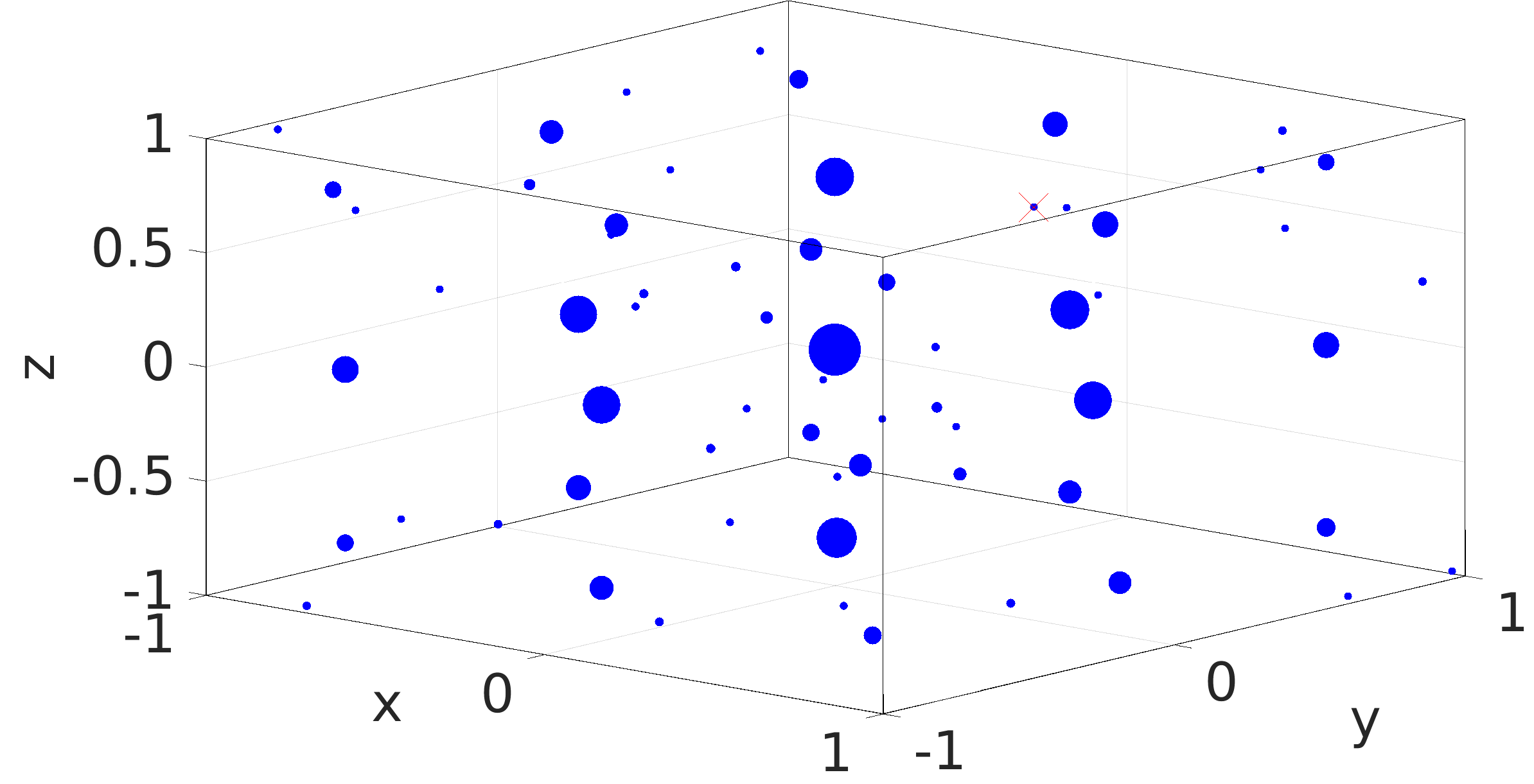}}   
  \subfigure[  $m=60$ points ($t=1,k=3$)]{\label{fig:3CECM_64a}\includegraphics[width=0.30\textwidth]{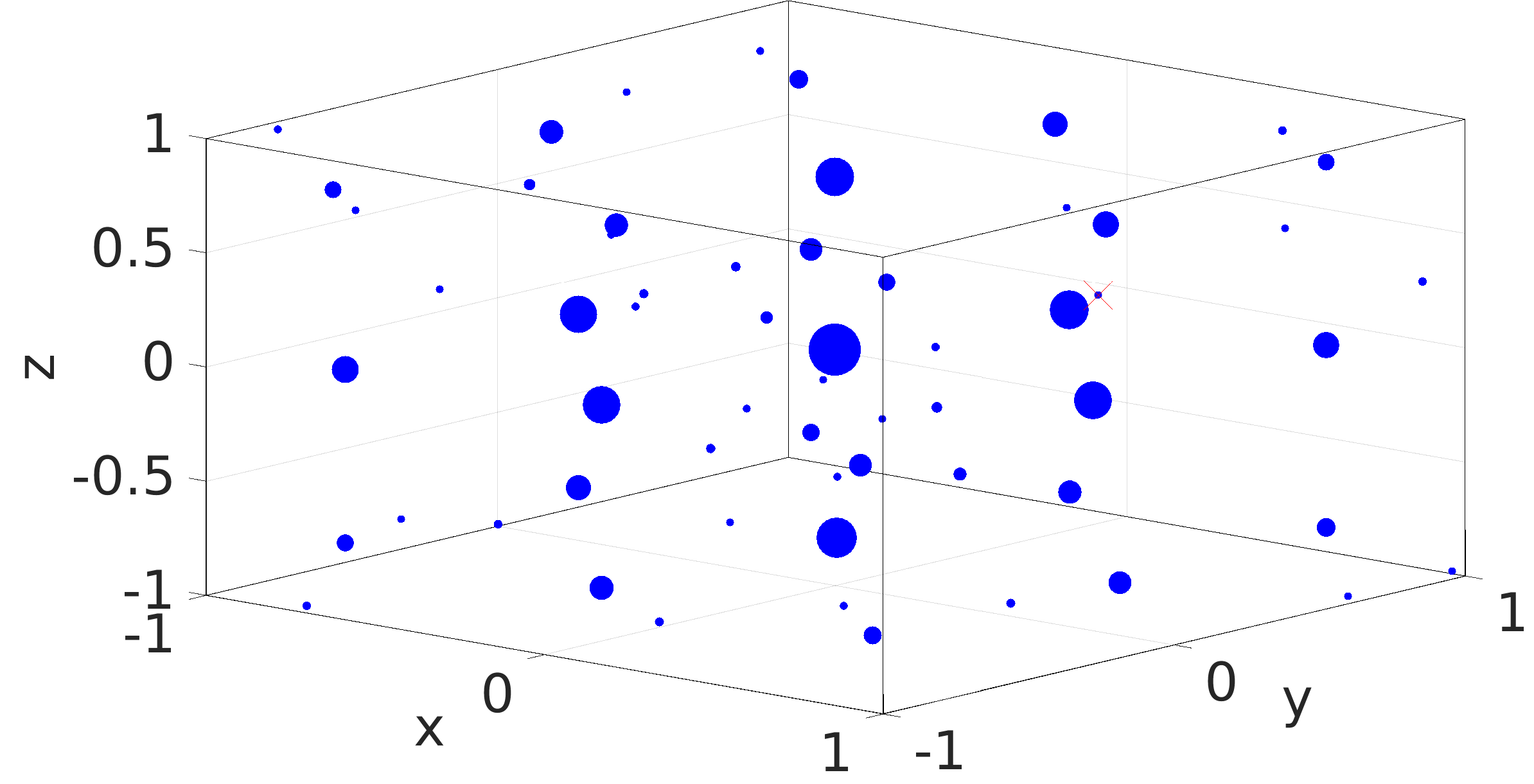}} 
    \subfigure[  $m=55$ points ($t=1,k=3$)]{\label{fig:3CECM_55}\includegraphics[width=0.30\textwidth]{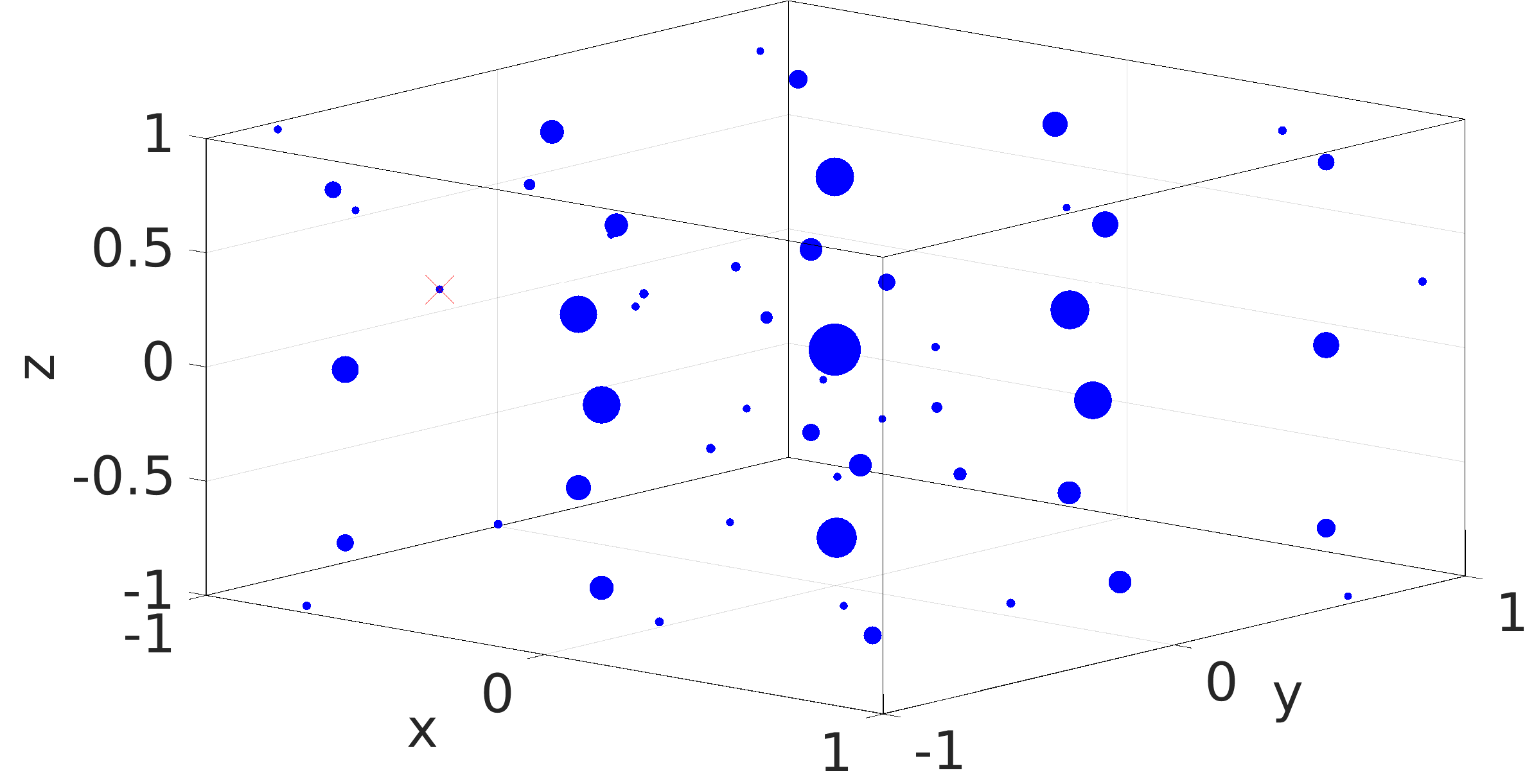}} 
    \subfigure[  $m=50$ points ($t=1,k=4$)]{\label{fig:3CECM_50}\includegraphics[width=0.30\textwidth]{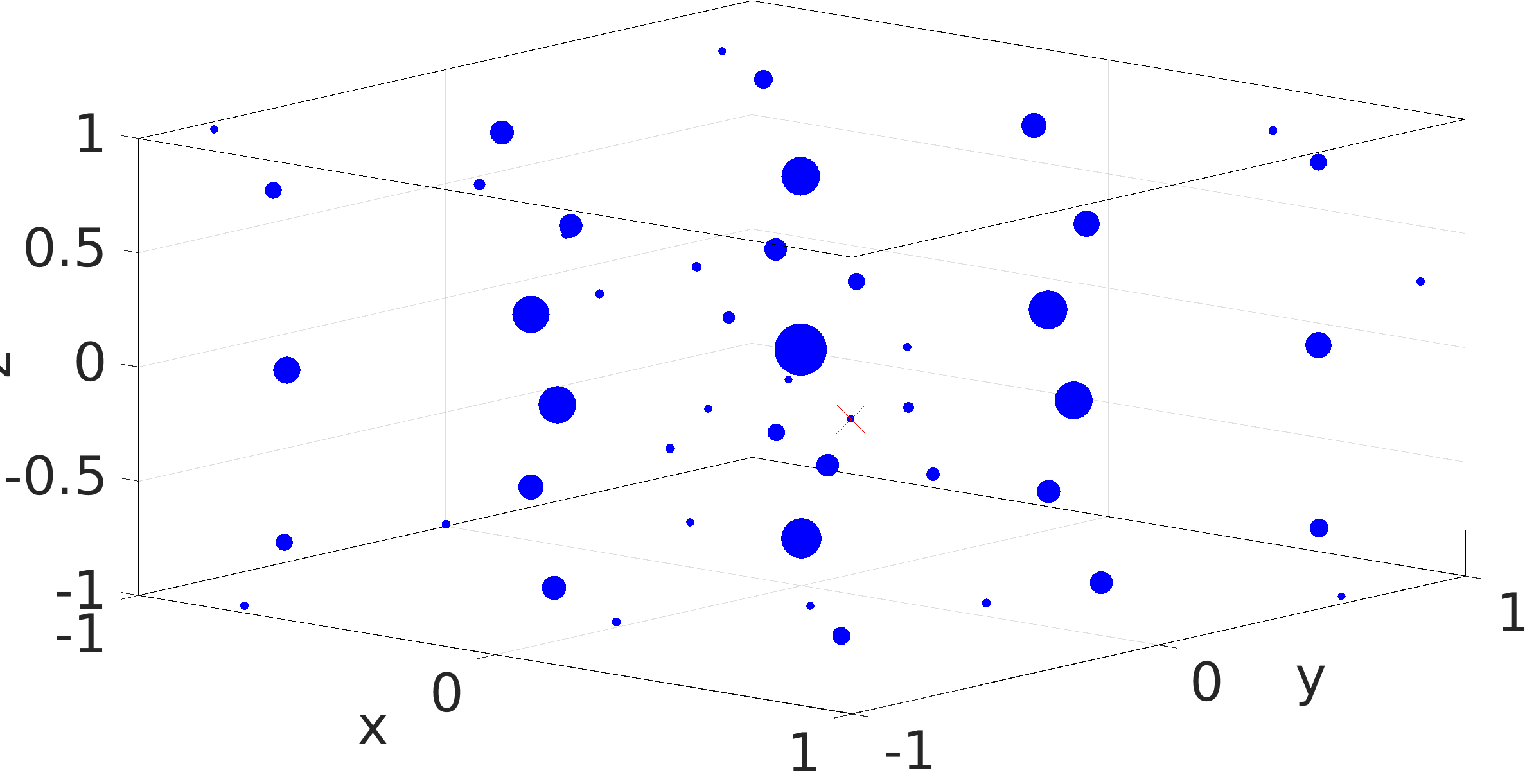}} 
    \subfigure[  $m=45$ points ($t=1,k=4$) ]{\label{fig:3CECM_45}\includegraphics[width=0.30\textwidth]{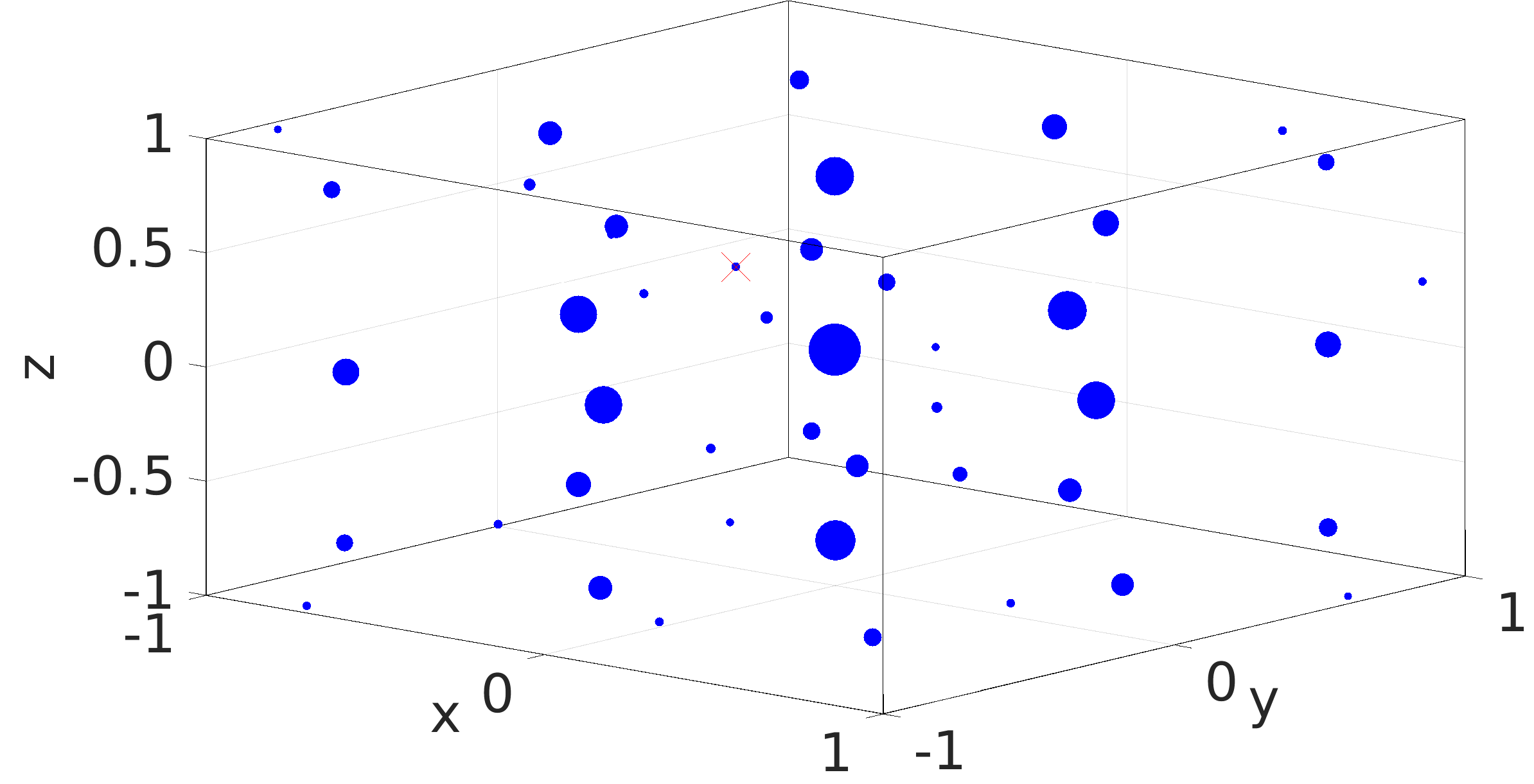}} 
    \subfigure[  $m=40$ points ($t=1,k=4$)]{\label{fig:3CECM_40}\includegraphics[width=0.30\textwidth]{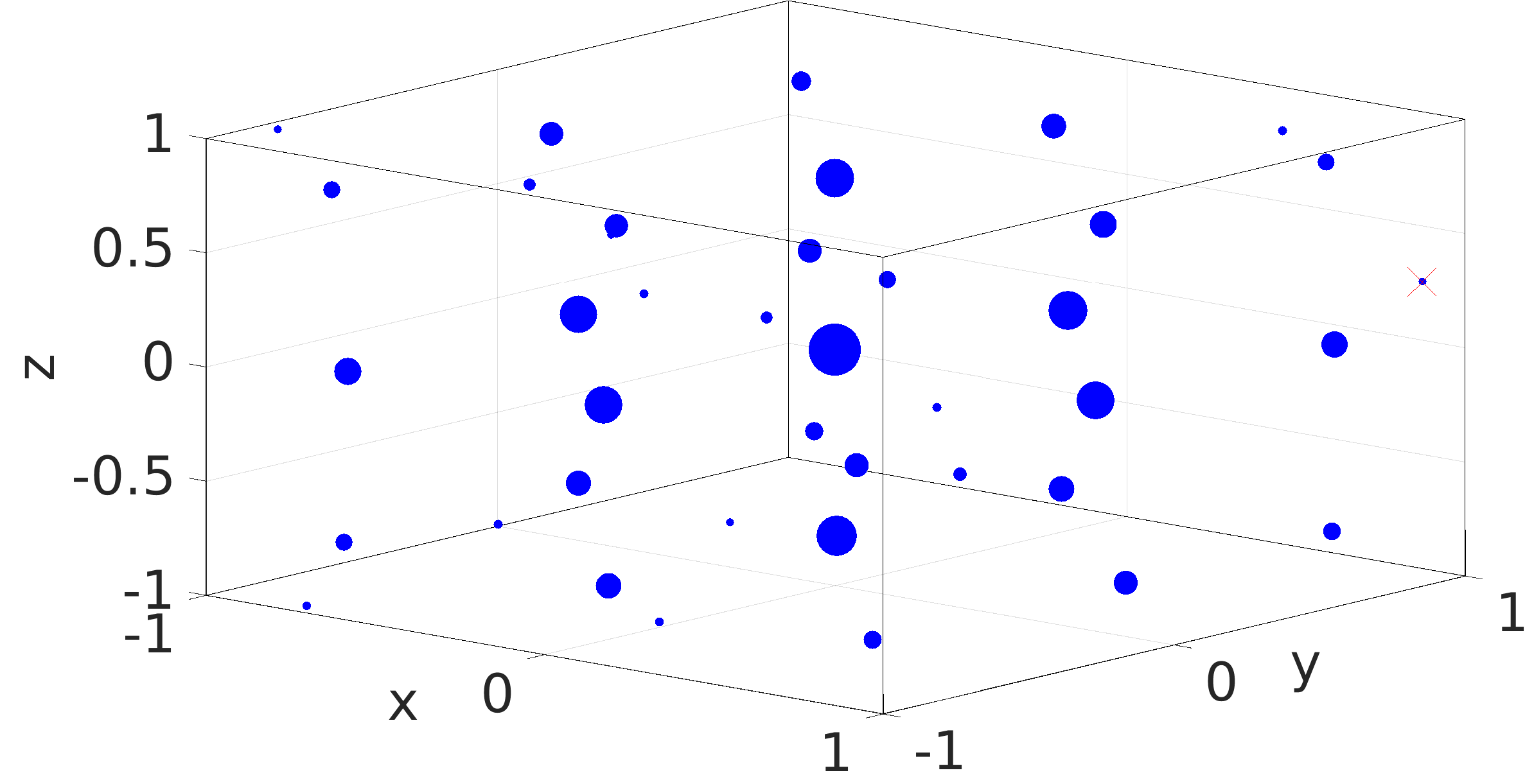}} 
    \subfigure[  $m=35$ points ($t=1,k=5$) ]{\label{fig:3CECM_35}\includegraphics[width=0.30\textwidth]{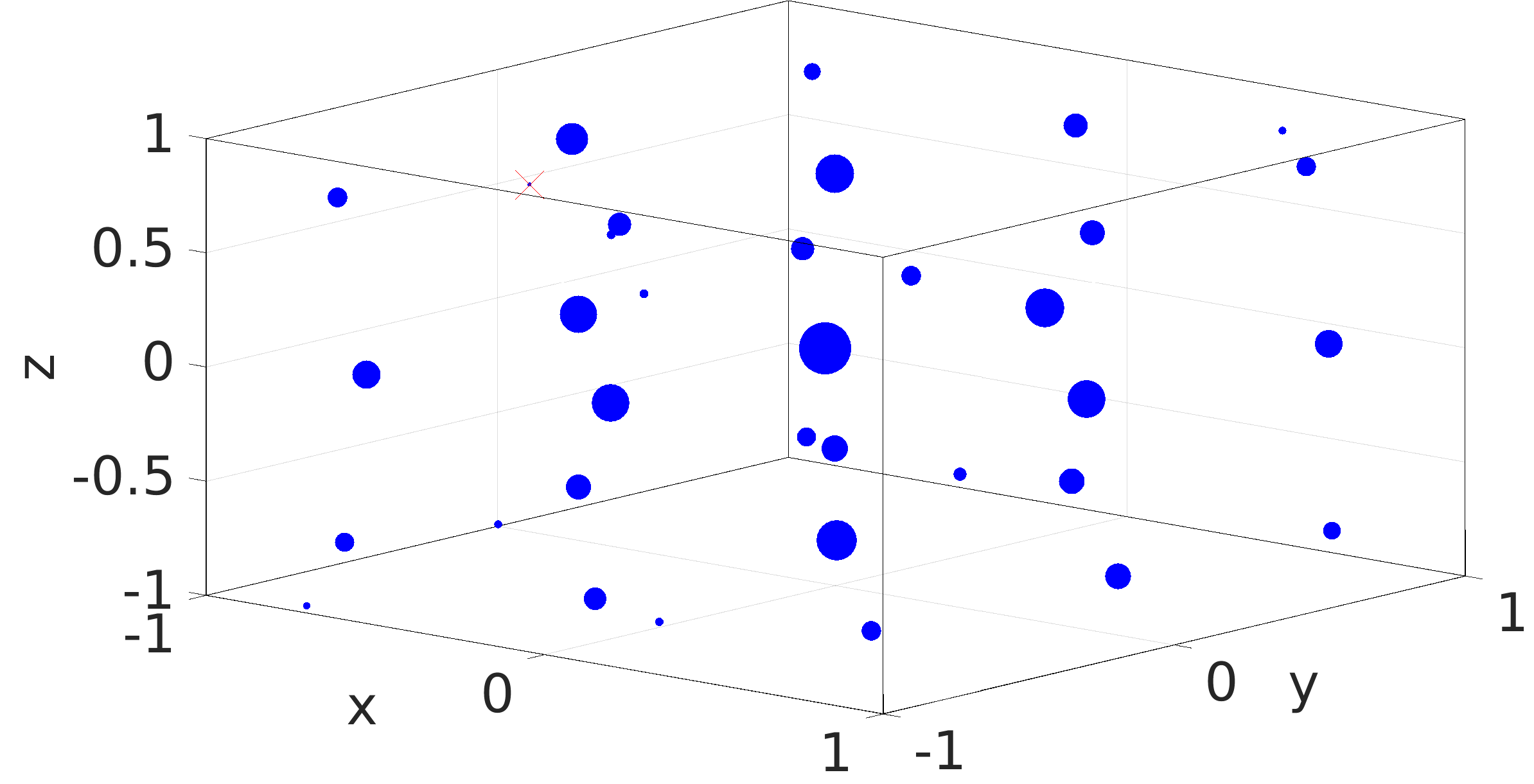}} 
    \subfigure[  $m=30$ points ($t=1,k=5$) ]{\label{fig:3CECM_30}\includegraphics[width=0.30\textwidth]{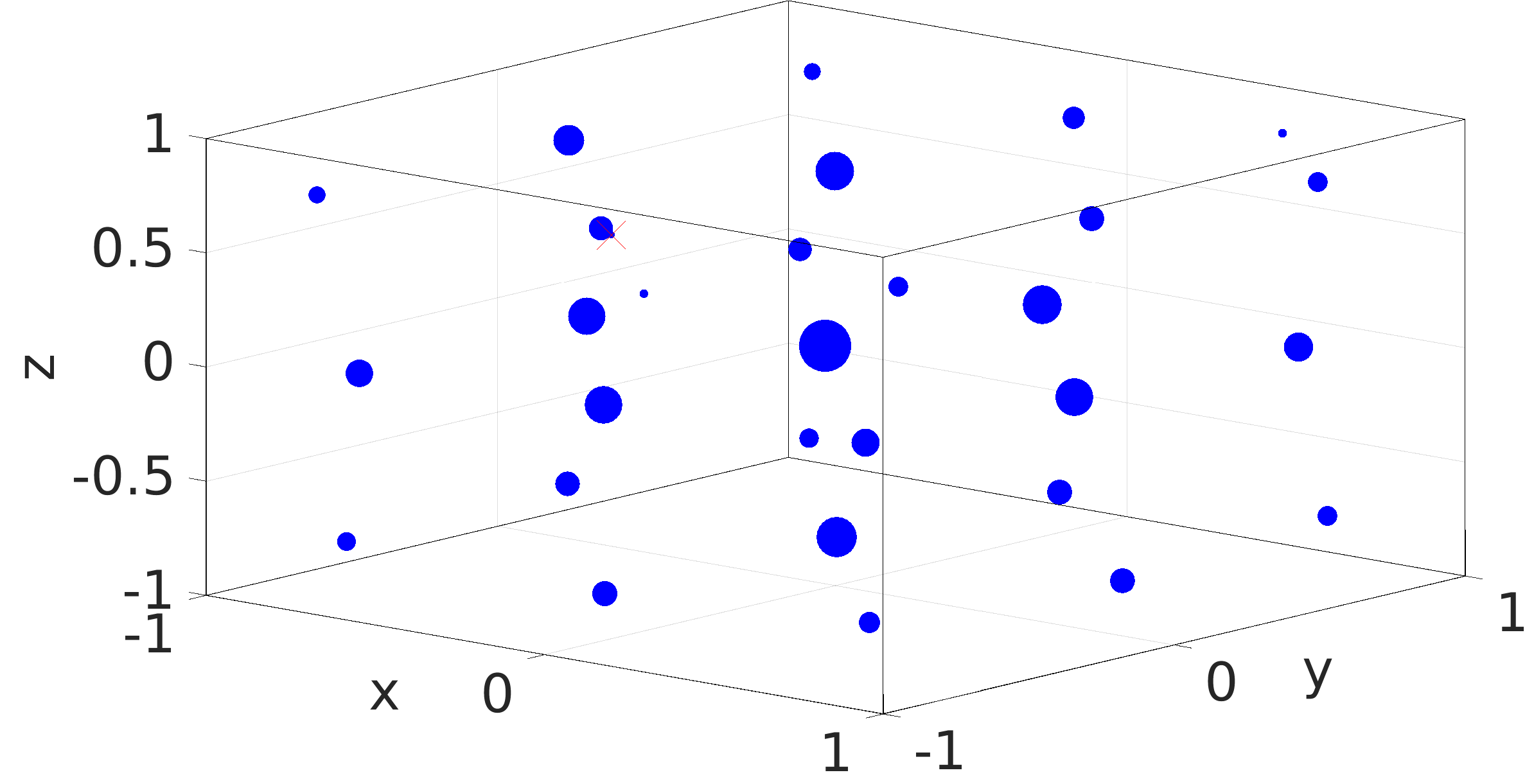}} 
    \subfigure[  $m=25$ points ($t=2,k=6$)]{\label{fig:3CECM_25}\includegraphics[width=0.30\textwidth]{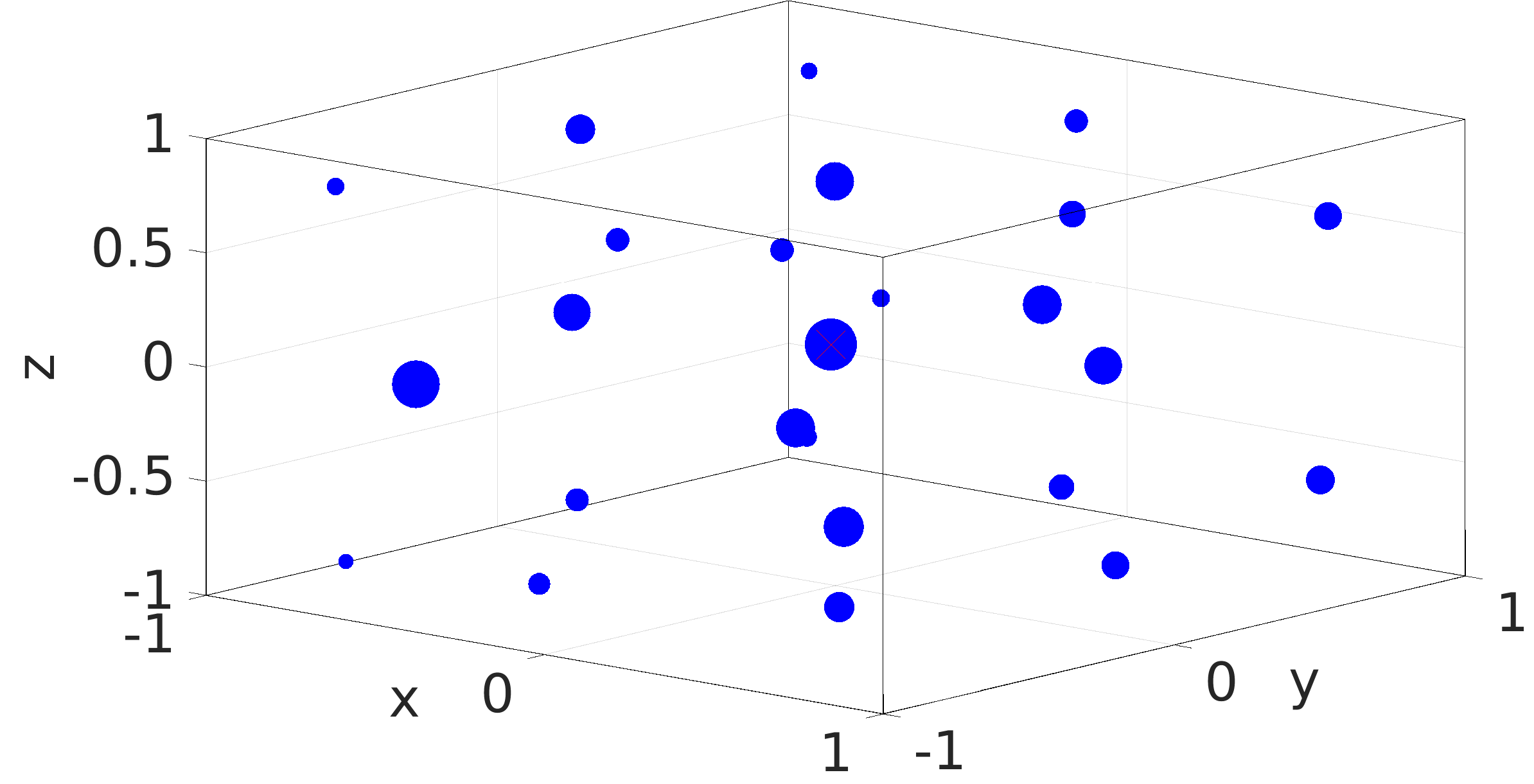}} 
    \subfigure[  $m=20$ points  ($t=1,k=16$)]{\label{fig:3CECM_20}\includegraphics[width=0.30\textwidth]{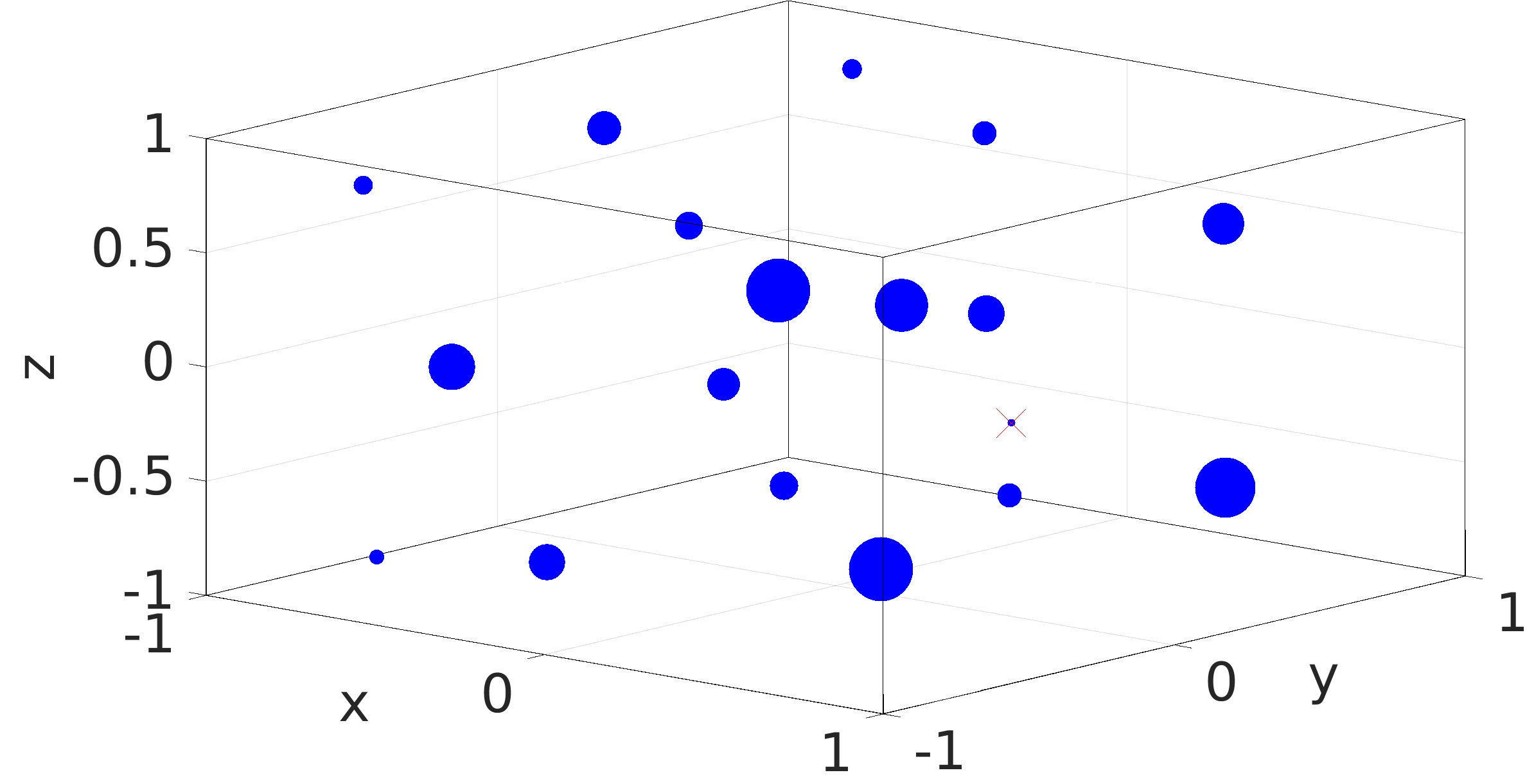}} 
    \subfigure[  $m=15$ points   ($t=1,k=1$) ]{\label{fig:3CECM_15}\includegraphics[width=0.30\textwidth]{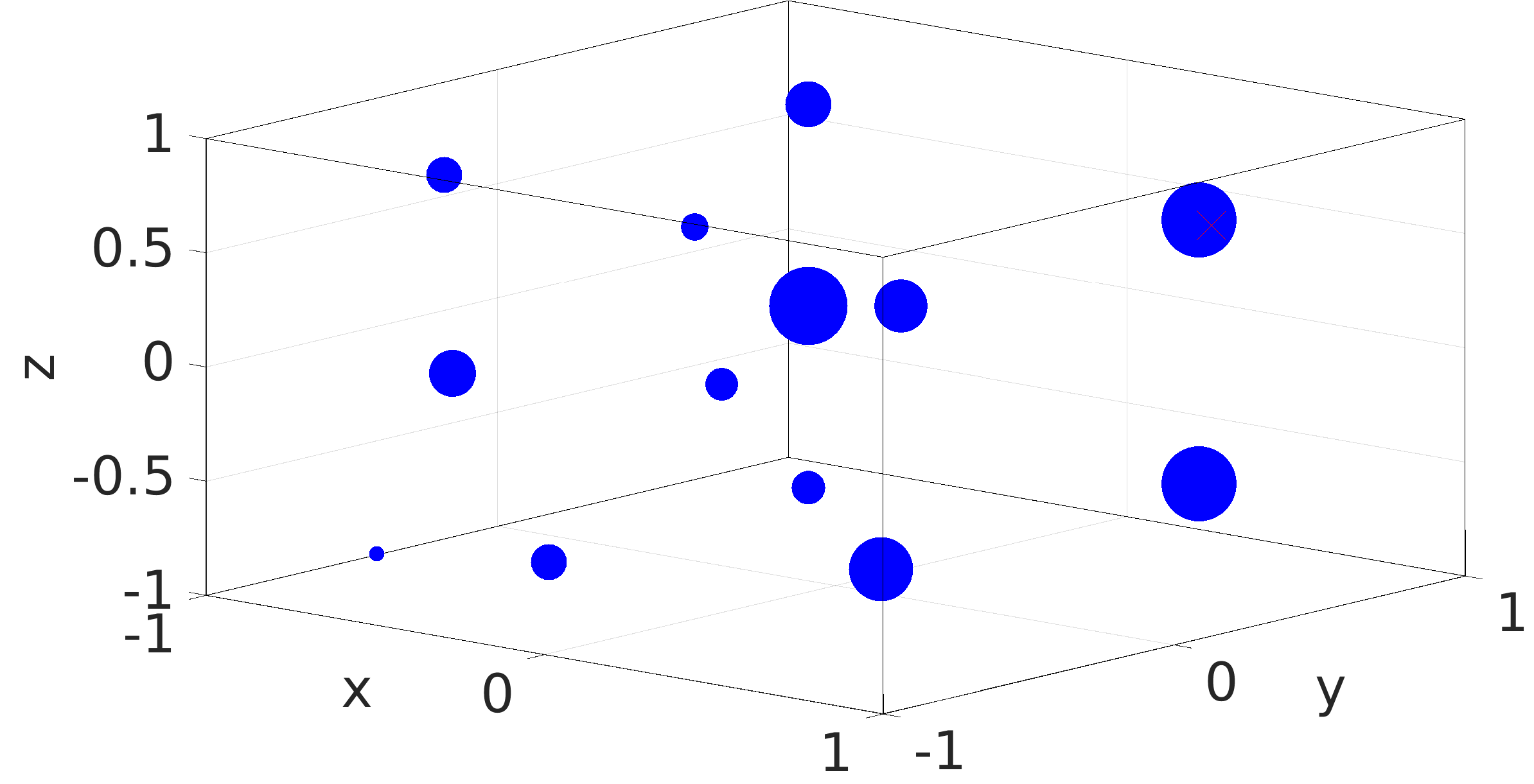}} 
    \subfigure[  $m=10$ points   ($t=1,k=1$) ]{\label{fig:3CECM_10}\includegraphics[width=0.30\textwidth]{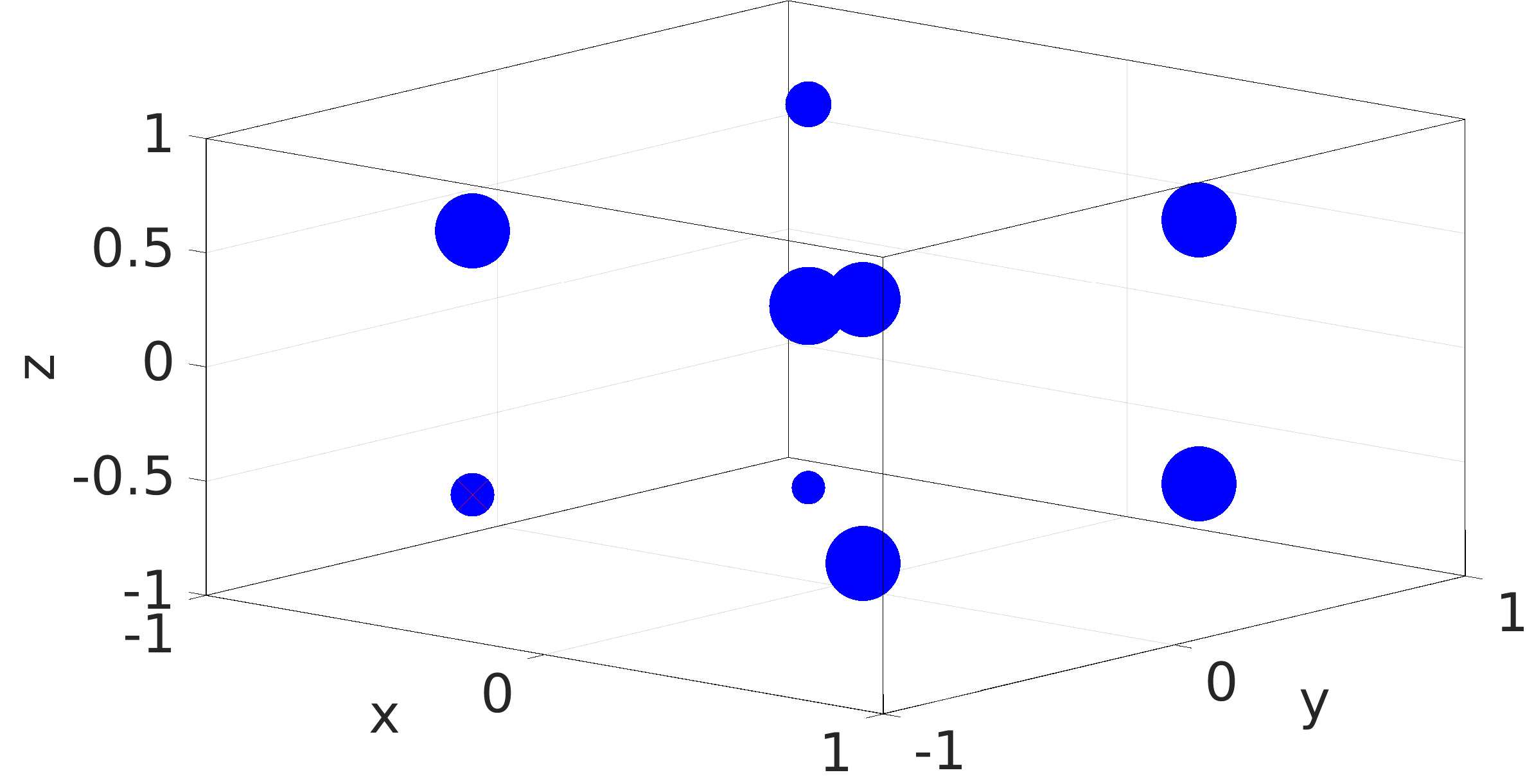}} 
    \subfigure[ Final CECM rule:  $m=8$ points   ($t=2,k=13$)]{\label{fig:3CECM_8}\includegraphics[width=0.40\textwidth]{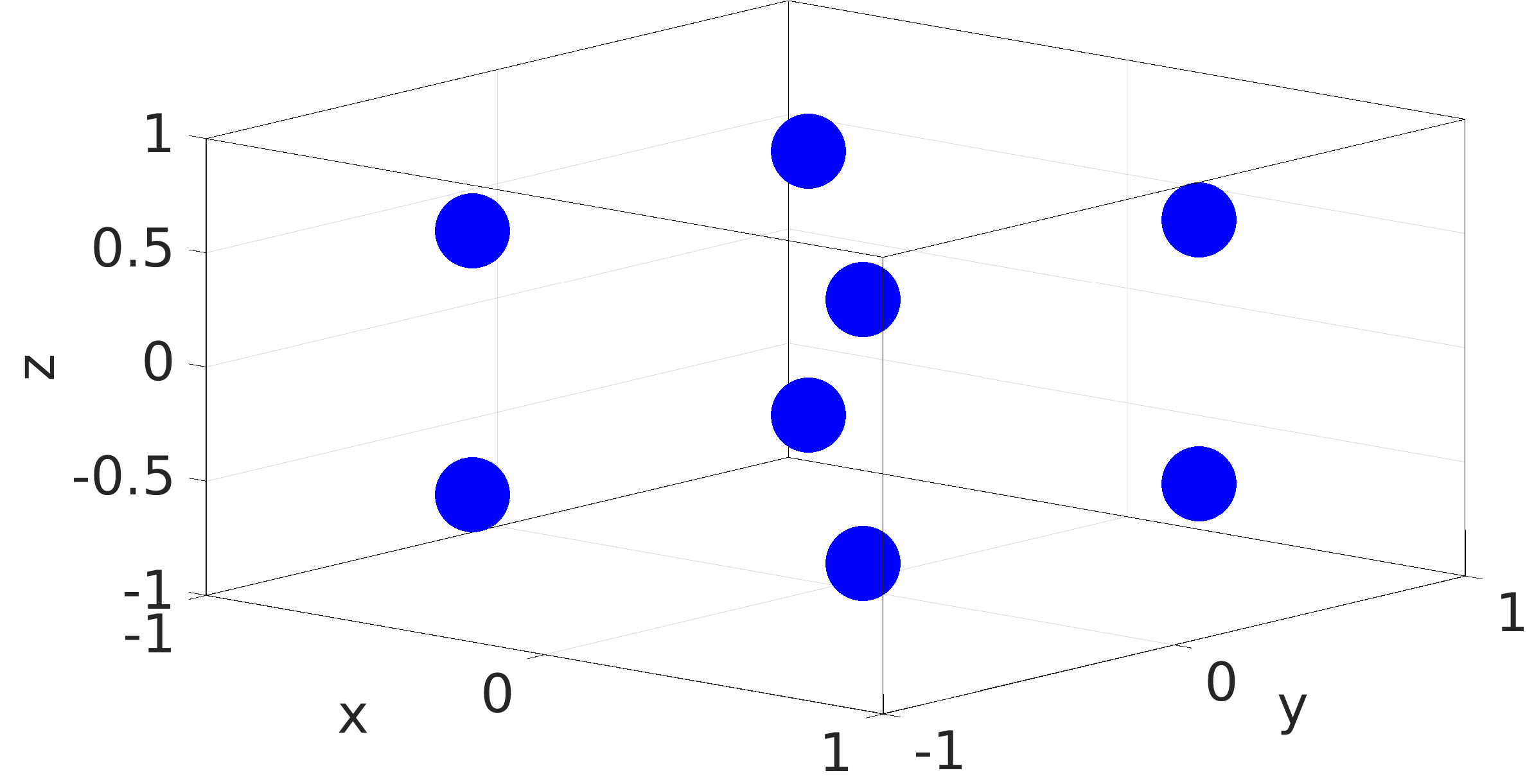}} 
   \caption{  Locations  and   weights   of the cubature rules generated during  the sparsification process for the case of trivariate polynomials of  degree $p = 3$ in $\Omega = [-1,1] \times [-1,1] \times [-1,1]$.   Variables $t$ and $k$  have the same interpretation as in Figures \ref{fig:FIG_5} and \ref{fig:FIG_6}.  The initial DECM rule has $(p+1)^3 = 64$ points (see Figure \ref{fig:3CECM_64}), and  the final   rule   $(p+1)^3/2^3 = 8$ points, see graph \ref{fig:3CECM_8}. The  exact values of the coordinates and weights are given in Table \ref{tab:3}, wherein it can be seen that the computed rule  does  coincide with  the standard $2 \times 2 \times 2$ product Gauss   rule.      }
  \label{fig:FIG_7}
\end{figure}

      These results provide  further confirmation of the ability of the proposed sparsification algorithm to arrive at the \emph{integration rules with minimal number of points}.  Figures \ref{fig:FIG_6} and \ref{fig:FIG_7} depict the sparsification process for the case $p = 3$ in 2D and 3D, respectively. As done previously in Figure \ref{fig:FIG_5} for the univariate case,  we show in each graph    the number of trials required  to find the point whose weight is to be zeroed (i.e., the number of times the method passess over the loop in line \ref{alg:003g} of Algorithm \ref{alg:003}), as well as the     number of iterations required  for zeroing the chosen weight (in the modified Newton-Raphson scheme of  Algorithm \ref{alg:005}).  Whereas in the case of univariate polynomials,   displayed previously in Figure \ref{fig:FIG_5}, the method  succesfully determines the weights to be zeroed on the first trial, in the multivariate case several trials are necessary in some cases, especially when the algorithm approaches the optimum. For instance, to produce the rule with 8 points in the bivariate case, see Figure \ref{fig:2CECM9}, the algorithm tries  $t=5$ different points until finding the appropriate combination.   
      Closer examination of the causes for this increase of  iterative effort indicates that the most common cause is the violation of the constraint that the points must remain within the domain.

 \jahoHIDE{See \href{EXAMPLES_2D/CECM_2D.m}{CECM\_2D.m},  InputDataFile = 'DATAlagrange2D'; stop  at   CODE\_CECM/ContinuousECM/HistoryPointsALL2D.m  }

   \jahoHIDE{See \href{EXAMPLES_3D/CECM_3D.m}{CECM\_3D.m} }
 
 \jahoHIDE{See also \href{/home/joaquin/Desktop/CURRENT_TASKS/MATLAB_CODES/TESTING_PROBLEMS_FEHROM/ContinuousEmpiricalCubatureM/Paper_hernandez2021ecm/09_ASSESSMENTpaper/FINAL2023_CECM/SummaryExamples2023.mlx}{SummaryExamples2023.mlx}}

    \jahoHIDE{ \url{MATLAB_EXAMPLES_CODE/EXAMPLES_3D/CECM_3D.m}, InputDataFile ='DATAlagrange3D\_p3'; Go to:  \url{ MATLAB_EXAMPLES_CODE/CODE_CECM/ContinuousECM/HistoryPointsPlot3D.m}, and set 
     AUXVAR.DATALOC.PlotEvolutionWeights\_MarkerSizeMin = 20 ;
    AUXVAR.DATALOC.PlotEvolutionWeights\_MarkerSizeMax = 200 ;
    AUXVAR.DATALOC.ColorMarkerPOINTS = [0,0,1] ;
    
    }

   \subsection{Exponential-sinusoidal   function  }
   \label{sec:expsino}
   
      \jahoHIDE{See \href{MATLAB_EXAMPLES_CODE/SummaryExamples2023.pdf}{SummaryExamples2023.pdf}}
   
   \jahoHIDE{  \href{/home/joaquin/Desktop/CURRENT_TASKS/MATLAB_CODES/TESTING_PROBLEMS_FEHROM/ContinuousEmpiricalCubatureM/Paper_hernandez2021ecm/09_ASSESSMENTpaper/FINAL_IMPLEM_CECM/EXAMPLES_3D/CECM_3D.m}{CECM\_3D.m}   with  input DATAexpCOS_nopart.m, }
   
   \jahoHIDE{SEe also \href{MATLAB_EXAMPLES_CODE/EXAMPLE_partSVD/TESTING_SVDpart.pdf}{TESTING\_SVDpart.pdf}}

We next study the derivation of a cubature rule for the following parameterized, vector-valued function:
   \begin{equation}
   \label{eq:param333}
   \begin{split}
    \a:  & \, \Omega \times \paramSPC \rightarrow \Rn{6},    
       \\& \,\,  (\x,\boldsymbol{\mu})   \mapsto      (a_1,a_2, \ldots a_6)
    \end{split}
   \end{equation}
   where 
     \begin{equation}
    \Omega =  [-1,1] \times [-1,1] \times [-1,1],  \hspace{1cm}   \textrm{(spatial domain)}
   \end{equation}
    \begin{equation}
   \hspace{0.5cm}   \paramSPC =  [1,\pi] \times [1,\pi],   \hspace{2.5cm}   \textrm{(parametric domain)}    
   \end{equation}
   \begin{equation}
   \begin{split}
  &  a_1  =  B(x_1) C(x_1,\mu_1) E(x_1,\mu_1) + 1     
 \\&   a_2  = B(x_2) C(x_2,\mu_1) E(x_2,\mu_1)  + 1
  \\&   a_3  = B(x_1) C(x_1,\mu_1) E(x_2,\mu_1)  +  1
    \\&   a_4  = B(x_2) C(x_2,\mu_1) E(x_1,\mu_1)  + 1 
        \\&   a_5  = B(x_1) C(x_1,\mu_1) E(x_3,\mu_2)  + 1 
            \\&   a_6  = B(x_3) C(x_3,\mu_2) E(x_2,\mu_1)   + 1
    \end{split}
   \end{equation}
   and 
   \begin{equation}
    B(r) = 1-r, \hspace{1cm}  C(r,s) = \coseno{3 \pi s (r+1)}, \hspace{1cm} E(r,s) = e^{(-1+r)s} . 
   \end{equation}

We use a structured spatial mesh of $30 \times 30 \times 30$ hexahedra elements,   each element being equipped with a product Gauss rule of    $3 \times 3 \times 3$ points. Unlike the case of polynomials discussed in the foregoing, where we knew beforehand which was the space of functions to be integrated ---the SVD only played a secondary, orthogonalizing role therein---, in this problem we have   to delineate first  the space in which the integrand lives.
This task naturaly confronts us with the   question of how dense should be the sampling of the parametric space so that the column space of the corresponding integrand matrix $\Afe{}$ becomes representative of this  linear space. 
We address here this question   by gradually increasing the number of sampled points in parametric space, applying the SVD with a fixed user-prescribed   truncation tolerance to the corresponding integrand matrix (here we use $\epsilon_{SVD}=10^{-4}$) , and  then examining when   the rank of the approximation (number of retained singular values) appears to converge to a maximum value.  Since there are only two parameters here, it is computationally affordable\footnote{  Higher parameter dimensions may require more sophisticated sampling strategies, such as the greedy adaptive procedure advocated in Ref. \cite{bui2007goal} for reduced-order modeling purposes. } to conduct this exploration by uniformly sampling the parametric space. 

 \begin{table}[]
\centering
\begin{tabular}{cccccccccc}
 &
   &
   &
   &
   &
   &
   &
   &
   &
   \\ \cline{4-10} 
\multicolumn{3}{c|}{} &
  \multicolumn{2}{c|}{SVD} &
  \multicolumn{4}{c|}{SRSVD} &
  \multicolumn{1}{c|}{ERROR SING. VAL.} \\ \hline
\multicolumn{1}{|c|}{$n_{samp}$} &
  \multicolumn{1}{c|}{$n_{col}$} &
  \multicolumn{1}{c|}{Size (GB)} &
  \multicolumn{1}{c|}{Time (s)} &
  \multicolumn{1}{c|}{Rank} &
  \multicolumn{1}{c|}{$N_{part}$} &
  \multicolumn{1}{c|}{$N_{iter}^{avg}$} &
  \multicolumn{1}{c|}{Time (s)} &
  \multicolumn{1}{c|}{Rank} &
  \multicolumn{1}{c|}{} \\ \hline
\multicolumn{1}{|c|}{4} &
  \multicolumn{1}{c|}{96} &
  \multicolumn{1}{c|}{0.56} &
  \multicolumn{1}{c|}{2.1} &
  \multicolumn{1}{c|}{36} &
  \multicolumn{1}{c|}{1} &
  \multicolumn{1}{c|}{3} &
  \multicolumn{1}{c|}{2.5} &
  \multicolumn{1}{c|}{36} &
  \multicolumn{1}{c|}{5.05E-15} \\ \hline
\multicolumn{1}{|c|}{6} &
  \multicolumn{1}{c|}{216} &
  \multicolumn{1}{c|}{1.26} &
  \multicolumn{1}{c|}{5.4} &
  \multicolumn{1}{c|}{53} &
  \multicolumn{1}{c|}{1} &
  \multicolumn{1}{c|}{3} &
  \multicolumn{1}{c|}{5.9} &
  \multicolumn{1}{c|}{53} &
  \multicolumn{1}{c|}{5.26E-15} \\ \hline
\multicolumn{1}{|c|}{8} &
  \multicolumn{1}{c|}{384} &
  \multicolumn{1}{c|}{2.24} &
  \multicolumn{1}{c|}{11.4} &
  \multicolumn{1}{c|}{70} &
  \multicolumn{1}{c|}{1} &
  \multicolumn{1}{c|}{2} &
  \multicolumn{1}{c|}{6.6} &
  \multicolumn{1}{c|}{70} &
  \multicolumn{1}{c|}{4.86E-15} \\ \hline
\multicolumn{1}{|c|}{11} &
  \multicolumn{1}{c|}{726} &
  \multicolumn{1}{c|}{4.23} &
  \multicolumn{1}{c|}{28.2} &
  \multicolumn{1}{c|}{90} &
  \multicolumn{1}{c|}{3} &
  \multicolumn{1}{c|}{1.67} &
  \multicolumn{1}{c|}{11.9} &
  \multicolumn{1}{c|}{90} &
  \multicolumn{1}{c|}{5.93E-14} \\ \hline
\multicolumn{1}{|c|}{16} &
  \multicolumn{1}{c|}{1536} &
  \multicolumn{1}{c|}{8.96} &
  \multicolumn{1}{c|}{84.7} &
  \multicolumn{1}{c|}{122} &
  \multicolumn{1}{c|}{4} &
  \multicolumn{1}{c|}{1.67} &
  \multicolumn{1}{c|}{21.3} &
  \multicolumn{1}{c|}{122} &
  \multicolumn{1}{c|}{3.66E-14} \\ \hline
\multicolumn{1}{|c|}{22} &
  \multicolumn{1}{c|}{2904} &
  \multicolumn{1}{c|}{16.94} &
  \multicolumn{1}{c|}{234.0} &
  \multicolumn{1}{c|}{131} &
  \multicolumn{1}{c|}{9} &
  \multicolumn{1}{c|}{1.22} &
  \multicolumn{1}{c|}{35.3} &
  \multicolumn{1}{c|}{131} &
  \multicolumn{1}{c|}{2.62E-13} \\ \hline
\multicolumn{1}{|c|}{31} &
  \multicolumn{1}{c|}{5766} &
  \multicolumn{1}{c|}{33.63} &
  \multicolumn{1}{c|}{*} &
  \multicolumn{1}{c|}{*} &
  \multicolumn{1}{c|}{16} &
  \multicolumn{1}{c|}{1.06} &
  \multicolumn{1}{c|}{65.2} &
  \multicolumn{1}{c|}{133} &
  \multicolumn{1}{c|}{*} \\ \hline
\end{tabular}
\caption{   Comparison of the performance of the proposed   Sequential Randomized SVD (see Algorithm \ref{alg:006} in Appendix \ref{sec:SRSVD}) with respect to the standard SVD in determining an approximate orthogonal basis matrix (truncation tolerance $\epsilon_{svd} = 10^{-4}$) for the column space of the integrand matrix of the vector-valued function ( $6$ components)  defined in Eq. \refpar{eq:param333}. The number of spatial integration points is  $\ngausT = (30 \cdot 3)^3 = 729000$, and the function is sampled at the points of   uniform grids in parameter space of varying size ($n_{samp} \times n_{samp}$, see first column).  The number of columns of the integrand matrix is therefore  $n_{col} = 6 n_{samp}^2 $, and its size (in gigabytes) equal to $8 \cdot 10^{-9} n_{col}  \ngausT $. 
For the matrix of 33.63 GB (last row) there is no information on either the computing time  nor the rank (number of basis vectors) for the standard truncated SVD (using the builtin Matlab function $\texttt{svd}$, see Algorithm \ref{alg:007} in Appendix \ref{sec:SRSVD}), because the computation exhausted the memory capabilities of the employed 64 GB RAM computer. $N_{part}$ denotes the number of partitions of the integrand matrix in the case of the SRSVD, and $N_{iter}^{avg}$ the average number of iterations employed by the incremental randomized orthogonalization (see Algorithm \ref{alg:010} in Appendix \ref{sec:SRSVD}) for all the partitions. The rightmost column represents the relative difference between the singular values computed by both methods ($ \normd{\S_{svd}-\S_{srsvd}}/\normd{\S_{svd}}$). 
}
\label{tab:4}
\end{table}



We show in  Table \ref{tab:4} the result of this convergence study using both the standard SVD and the proposed Sequential Randomized SVD (described in Appendix \ref{sec:SRSVD}). The study has been devised so that the  size of the integrand matrix doubles at each refinement step. Likewise, the block partition of the integrand matrix in the case of the SRSVD has been taken so that the size of each block matrix is approximately 2 GB.  This convergence study reveals that the dimension of the  linear space in which the integrand lies (for the prescribed tolerance) is around 130.  The study  also serves to highlight the advantages of the proposed SRSVD in terms of both computing time and memory requirements:  for the matrix of size $16.96$ GB, the SRSVD is almost 7 times faster than the SVD, and  for the largest matrix of $33.63$ GB, the standard SVD cannot handle the operation because it exhausts the memory capabilities of the employed computer (which has 64 GB RAM\footnote{The code is implemented in Matlab, and exectuted in an Intel(R) Core(TM) i7-8700 CPU, 3.20GHz with 64 Gb RAM (Linux platform) }); the SRSVD, by constrast,   returns  the result in approximately 1 minute.         The reasons of this clear outperformance of the SRSVD over the SVD are further discussed in Section \ref{sec:numst} of Appendix \ref{sec:SRSVD}.

          \begin{figure}[!ht]
  \centering
  \subfigure[Initial DECM rule, $m=133$ points]{\label{fig:iniSIN}\includegraphics[width=0.49\textwidth]{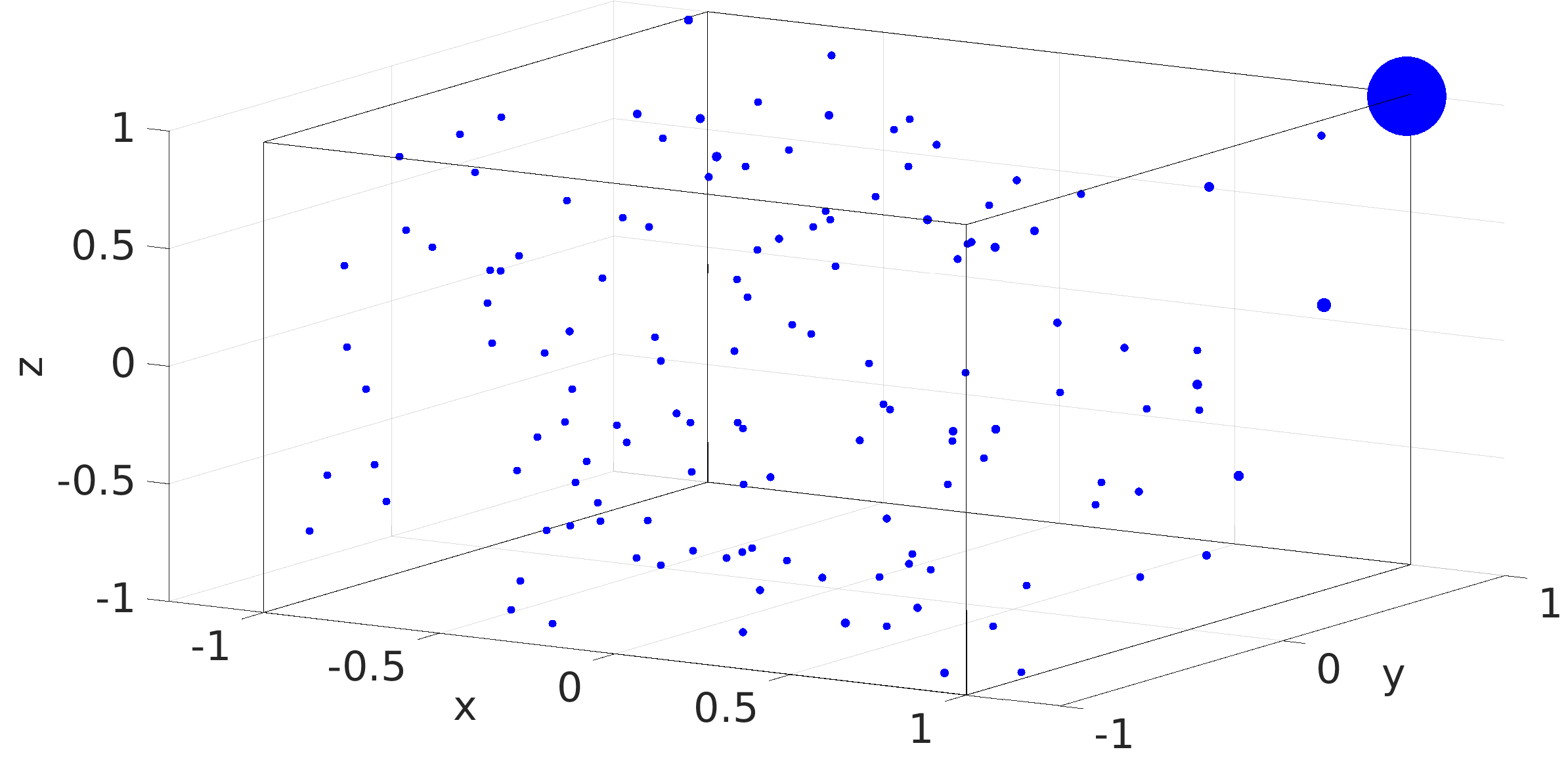}}   
   \subfigure[ ]{\label{fig:trials}\includegraphics[width=0.49\textwidth]{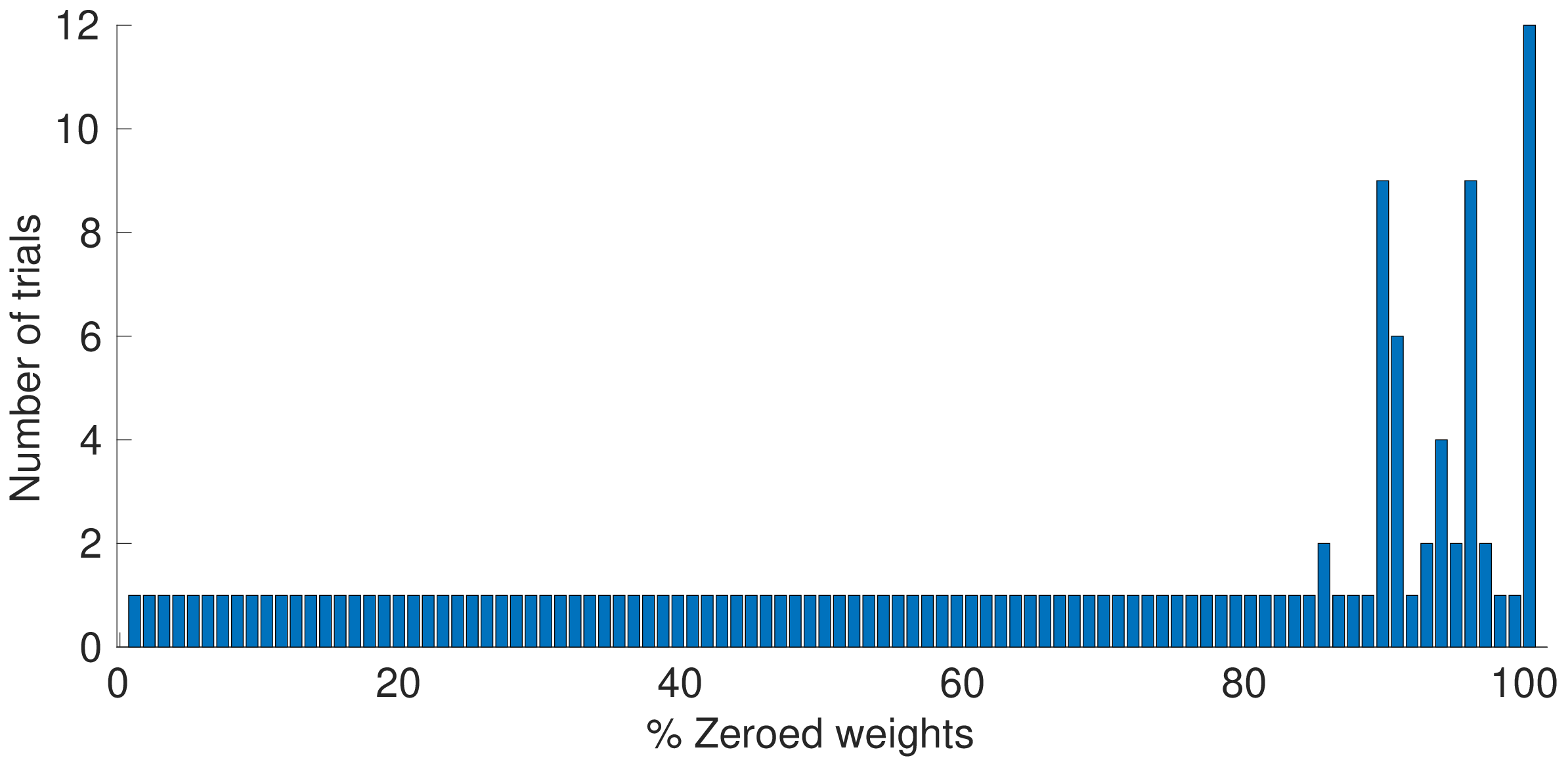}}  
    \subfigure[ Final CECM rule:  $m=38$ points  ]{\label{fig:finSIN}\includegraphics[width=0.49\textwidth]{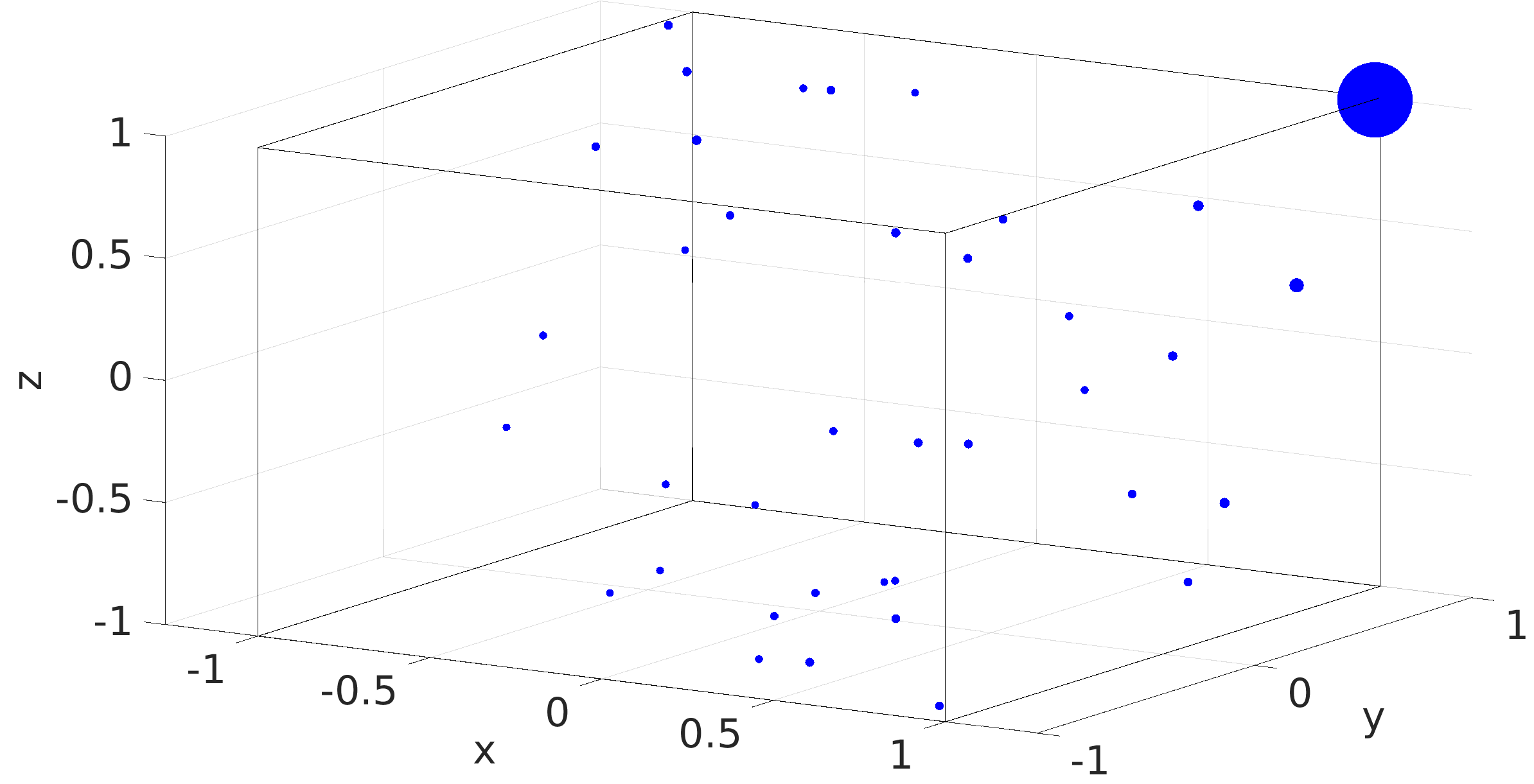}} 
     \subfigure[  ]{\label{fig:iter}\includegraphics[width=0.49\textwidth]{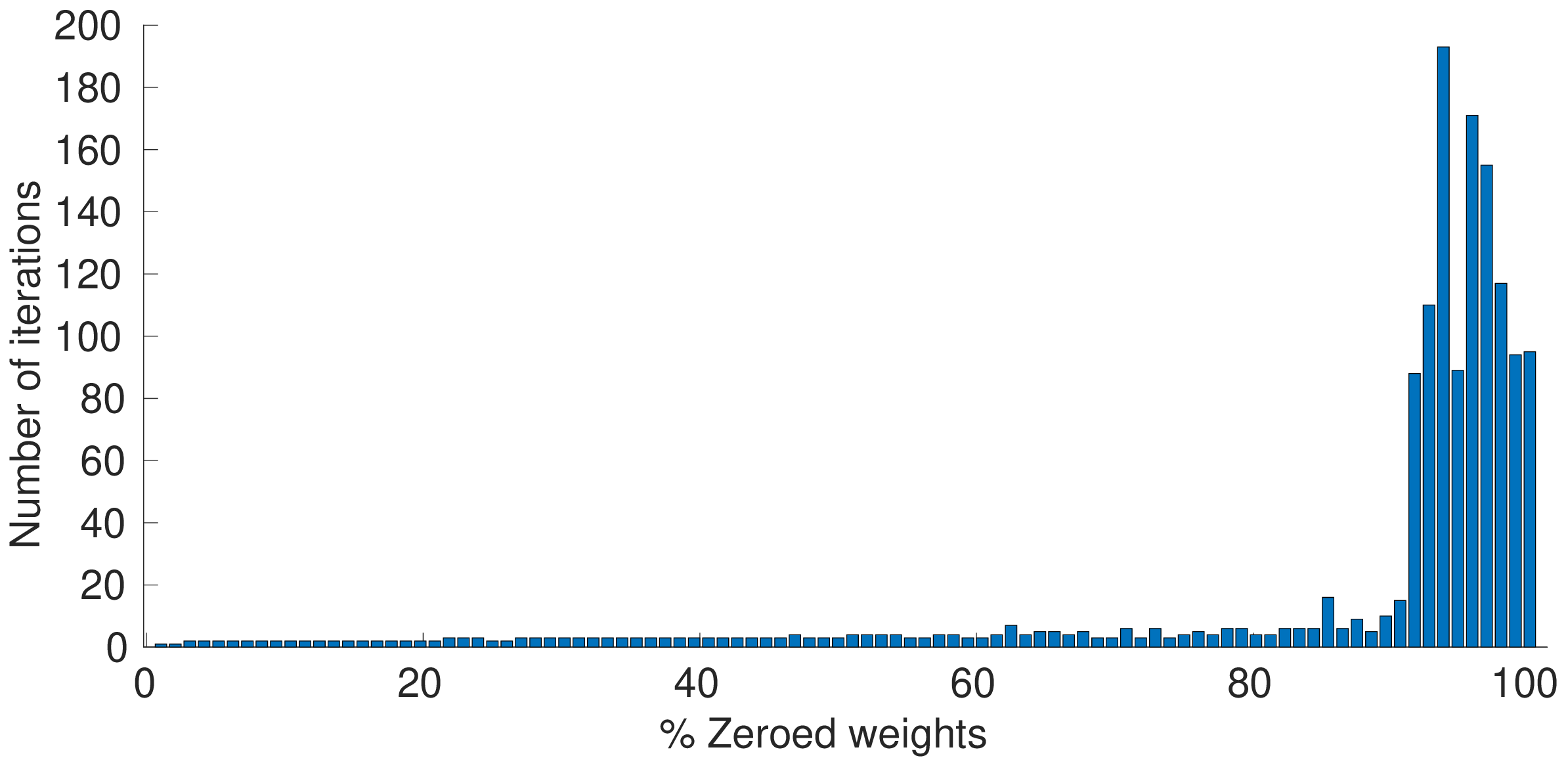}} 
   \caption{ Cubature of  the parameterized function defined in Eq. \ref{eq:param333} in $\Omega = [-1,1] \times [-1,1] \times [-1,1]$. The orthogonal basis vectors are determined by the SRSVD using a parametric grid of $31 \times 31$ points (see last row of Table \ref{tab:3}). a) Initial interpolatory DECM rule.   b) Number of passes over the loop that selects the weights to be zeroed in Algorithm \ref{alg:003} (see line \ref{alg:003g}) as a function of the percentage of zeroed weights.        c) Final CECM  rule.      d)   Total number of iterations of the modified Newton-Raphson scheme in  Algorithm \ref{alg:005} as a function of the percentage of zeroed weights.       }
  \label{fig:FIG_8}
\end{figure}

  Figure \ref{fig:iniSIN} shows the DECM rule determined using the 133 left singular vectors provided by the SRSVD for the parametric grid of $31 \times 31$ points (see last row of Table \ref{tab:4}), while Figure \ref{fig:finSIN} displays the final 38-points CECM rule obtained after the sparsification process ---the reduction factor is approximately 3.4. The variables controlling the sparsification are the same employed in the polynomial case.   
Further  information about the sparsification process are displayed in  Figures \ref{fig:trials} and \ref{fig:iter}. The number of  trials taken by the algorithm to find the weight to be zeroed (versus the percentage of zeroed weights) is shown in Figure \ref{fig:trials}, whereas   Figure  \ref{fig:iter} represents the total number of accumulated nonlinear  iterations (also versus the  percentage of zeroed weights).      It can be seen in  Figure \ref{fig:trials} that   approximately 90 \% of the weights are zeroed on the first trial; it is only in the last 10 \% that the number of trials increases considerably. The same behavior is observed in terms of   accumulated  nonlinear iterations in Figure \ref{fig:iter}: while the first 90 \% of the points are zeroed  in   5 iterations on average, for the  last    10 \% of points, the  number of iterations required for this very task raises  sharply   (close to 200 iterations in some cases). This is   not only due to the increase of the number of attempts to zeroed the weights reflected in Figure \ref{fig:trials}, but also because, at this juncture of the sparsification process, the weights of the points are relatively large and, to ensure converge,  the  
  problem of driving the integration residual to zero  is   solved in more than one step\footnote{The sparsification process enters the  second  stage in line \ref{alg:001_b} of Algorithm \ref{alg:001}}     (we take here $N_{steps} =20$).  
  
  Obviously, this uneven distribution of iterations during the sparsification process translated  into an equally uneven computing time distribution:      zeroing the  first 90 \%  of weights   took less than 1 minute, whilst the remaining 10 \%   required about 8 minutes.

   \subsection{Hyperreduction of multiscale finite element models}
   \label{sec:hyper}
   \jahoHIDE{\href{/home/joaquin/Desktop/CURRENT_TASKS/MATLAB_CODES/TESTING_PROBLEMS_FEHROM/ContinuousEmpiricalCubatureM/Paper_hernandez2021ecm/09_ASSESSMENTpaper/MultiscaleHROM/READMEmulti.pdf}{MultiscaleHROM/READMEmulti.pdf}}

     \begin{figure}[!ht]
  \centering
  \includegraphics[width=0.99\textwidth]{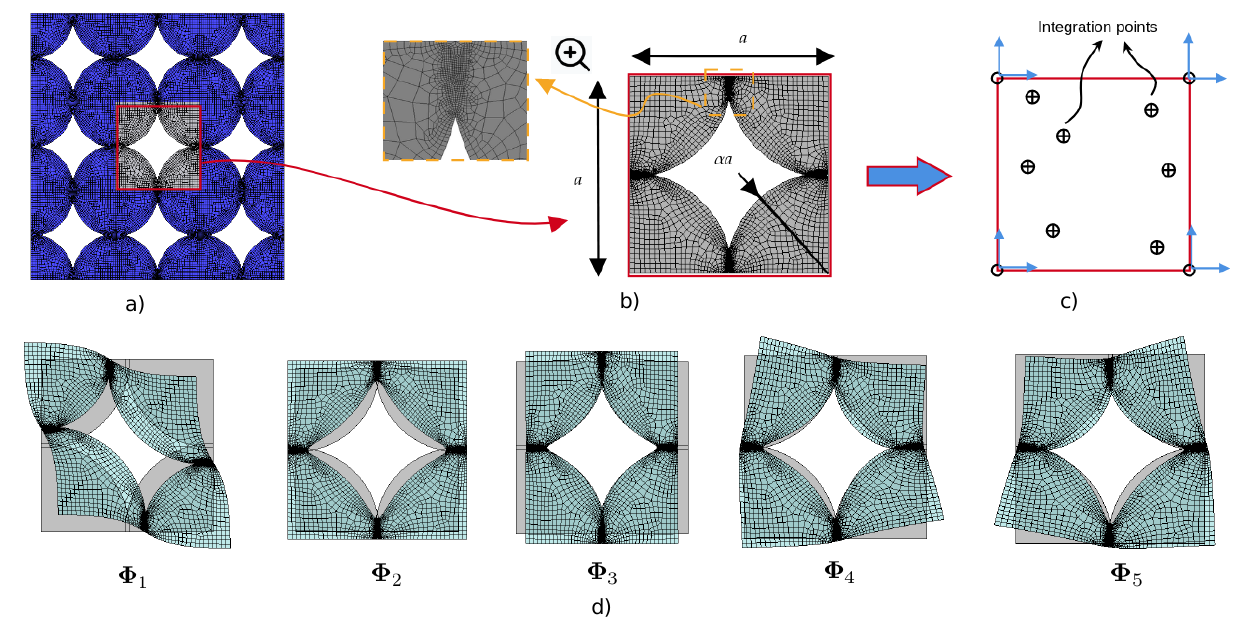}
  \caption{ Assessment of the performance of the CECM in the hyperreduction of multiscale finite elements.  a) Periodic structure under study.  b) Fine-scale mesh of the unit cell for which the coarse-scale representation is required ($a= 0.195$ m, $\alpha = 0.5135$). The total number of Gauss points is  $M = 42471$.  c) Coarse-scale representation of the unit cell, possessing 8 degrees of freedom and a number of integration points to be determined by the CECM.  d) Deformational modes  of the unit cell. The integral to be tackled by the CECM is the projection of the fine-scale nodal internal forces onto the span of these modes (see Eq. \ref{eq:internal}).     }
     \label{fig:9}
 \end{figure}
 
   \begin{figure}[!ht]
  \centering
  \includegraphics[width=0.99\textwidth]{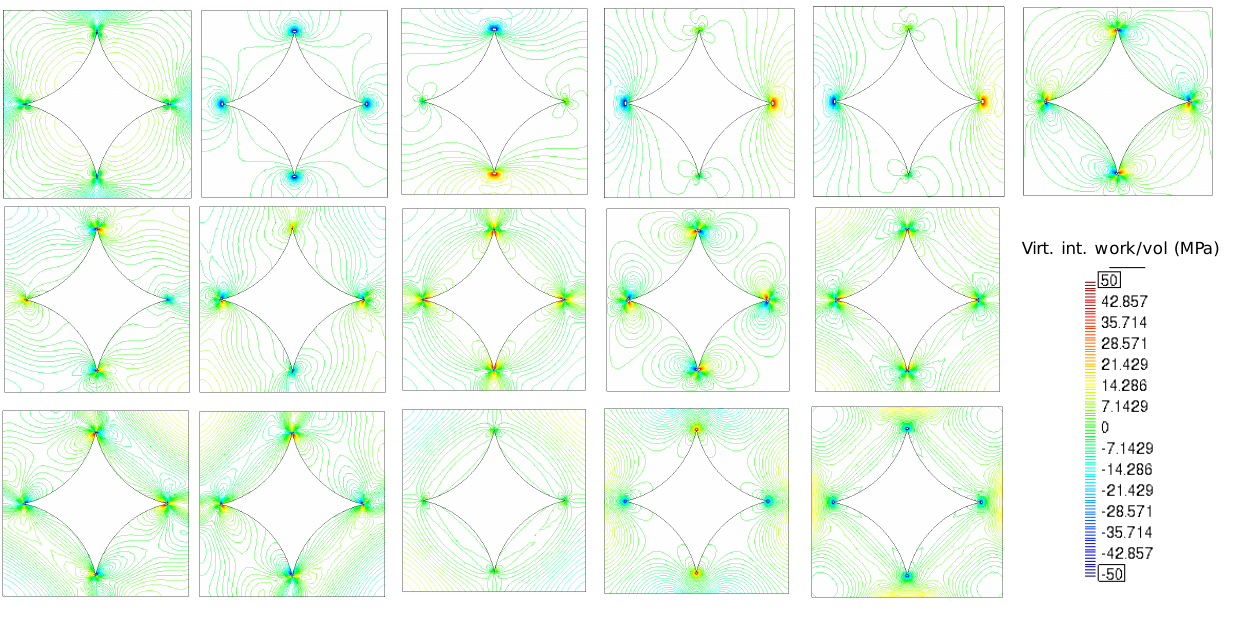}
  \caption{ Countor lines corresponding to the 16 basis functions for the integrand   in Eq. \ref{eq:internal}, which represents the virtual work per unit volume  (for the case in which the reduced basis $\PhiB$ is formed by the $n=5$ modes shown previously in Figure \ref{fig:9}.d)    }
     \label{fig:10}
 \end{figure}

        \begin{figure}[!ht]
  \centering
  \subfigure[Initial DECM rule, $m = 16$ points]{\label{fig:unitD16}\includegraphics[width=0.325\textwidth]{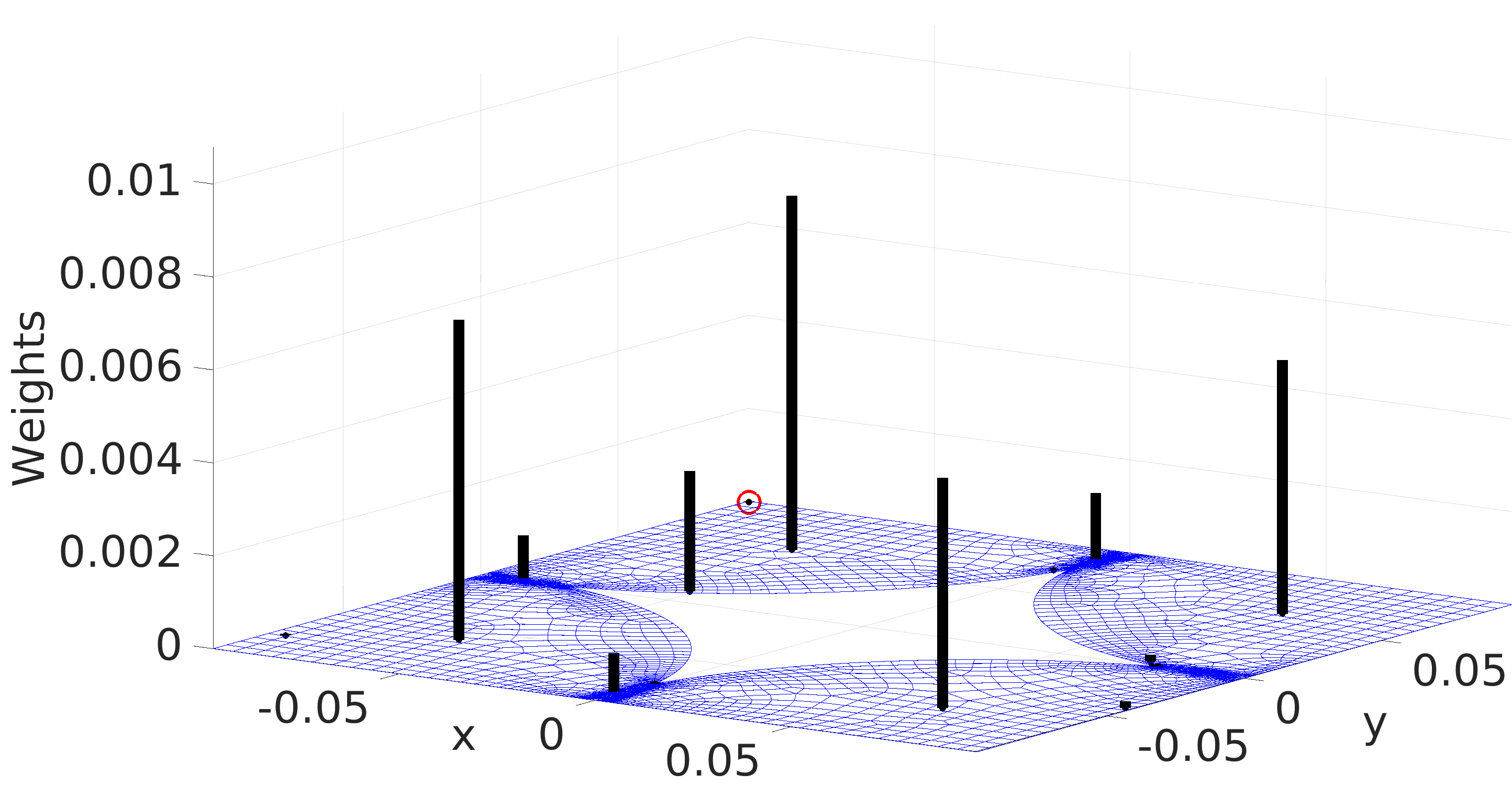}}   
   \subfigure[ $m = 15$ points ($t=1,k= 3$)]{\label{fig:unitD15}\includegraphics[width=0.325\textwidth]{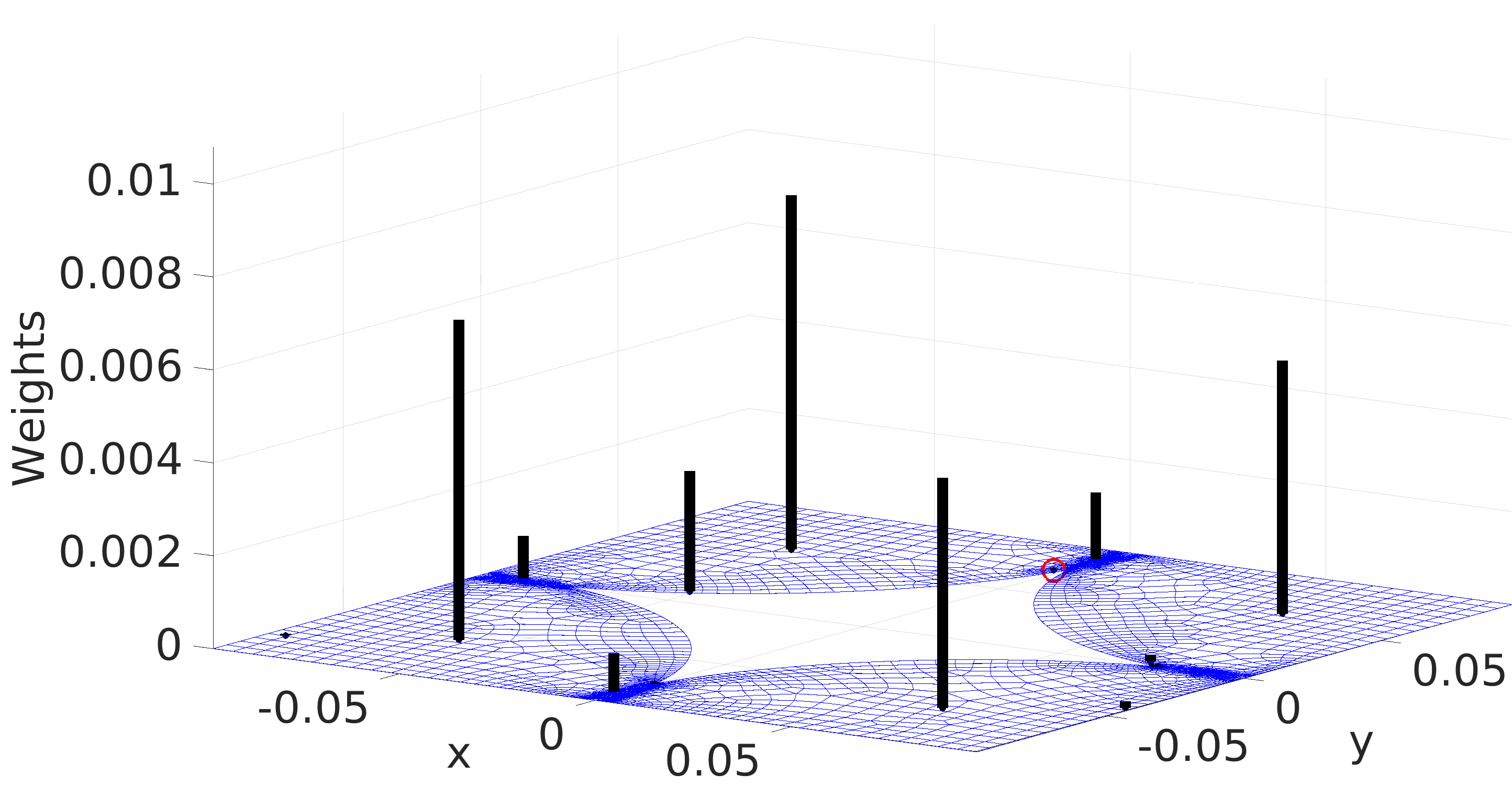}}  
    \subfigure[$m = 14$ points ($t=1,k= 3$) ]{\label{fig:unitD14}\includegraphics[width=0.325\textwidth]{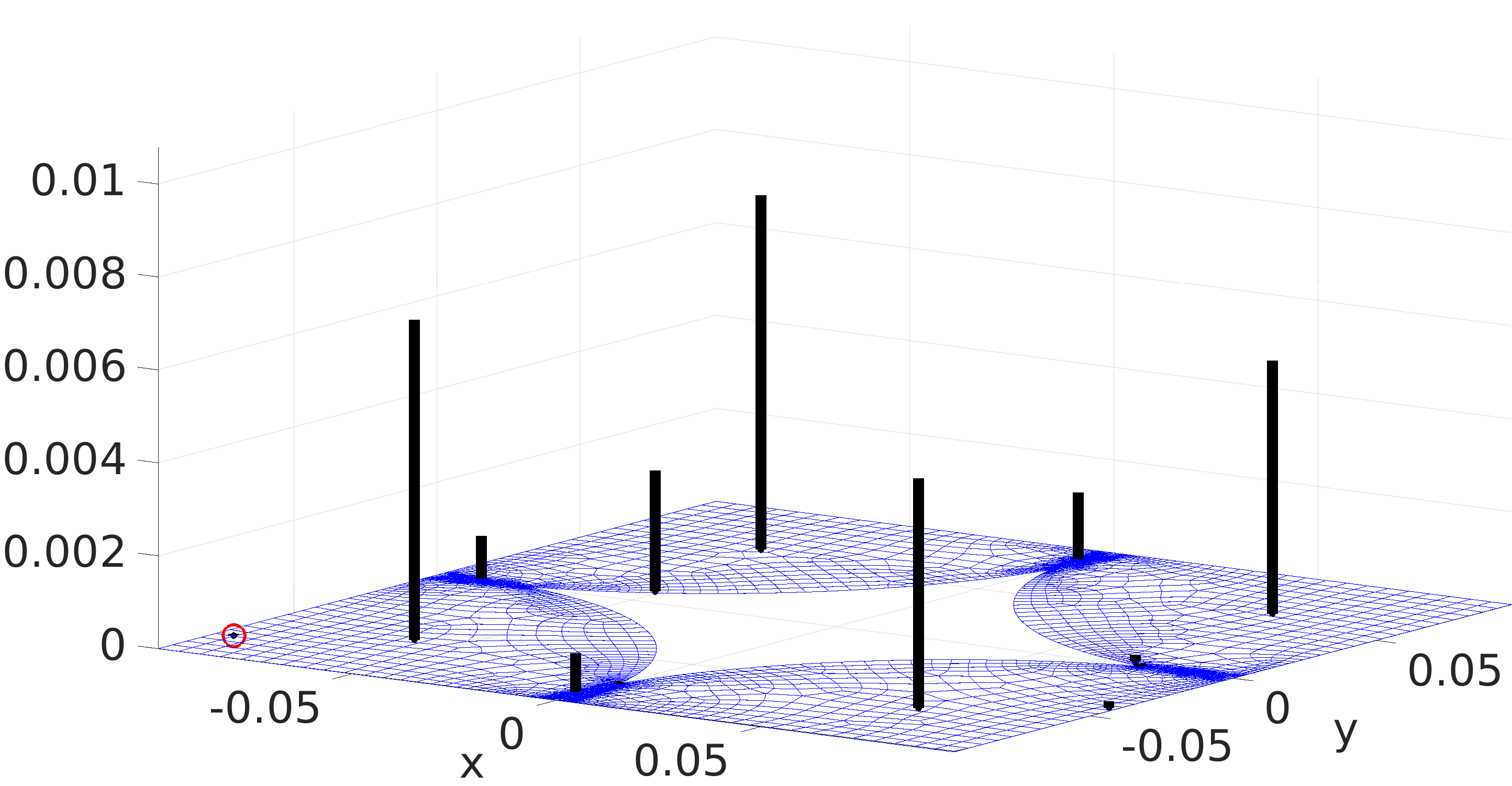}}  
     \subfigure[$m = 13$ points ($t=1,k= 5$) ]{\label{fig:unitD13}\includegraphics[width=0.325\textwidth]{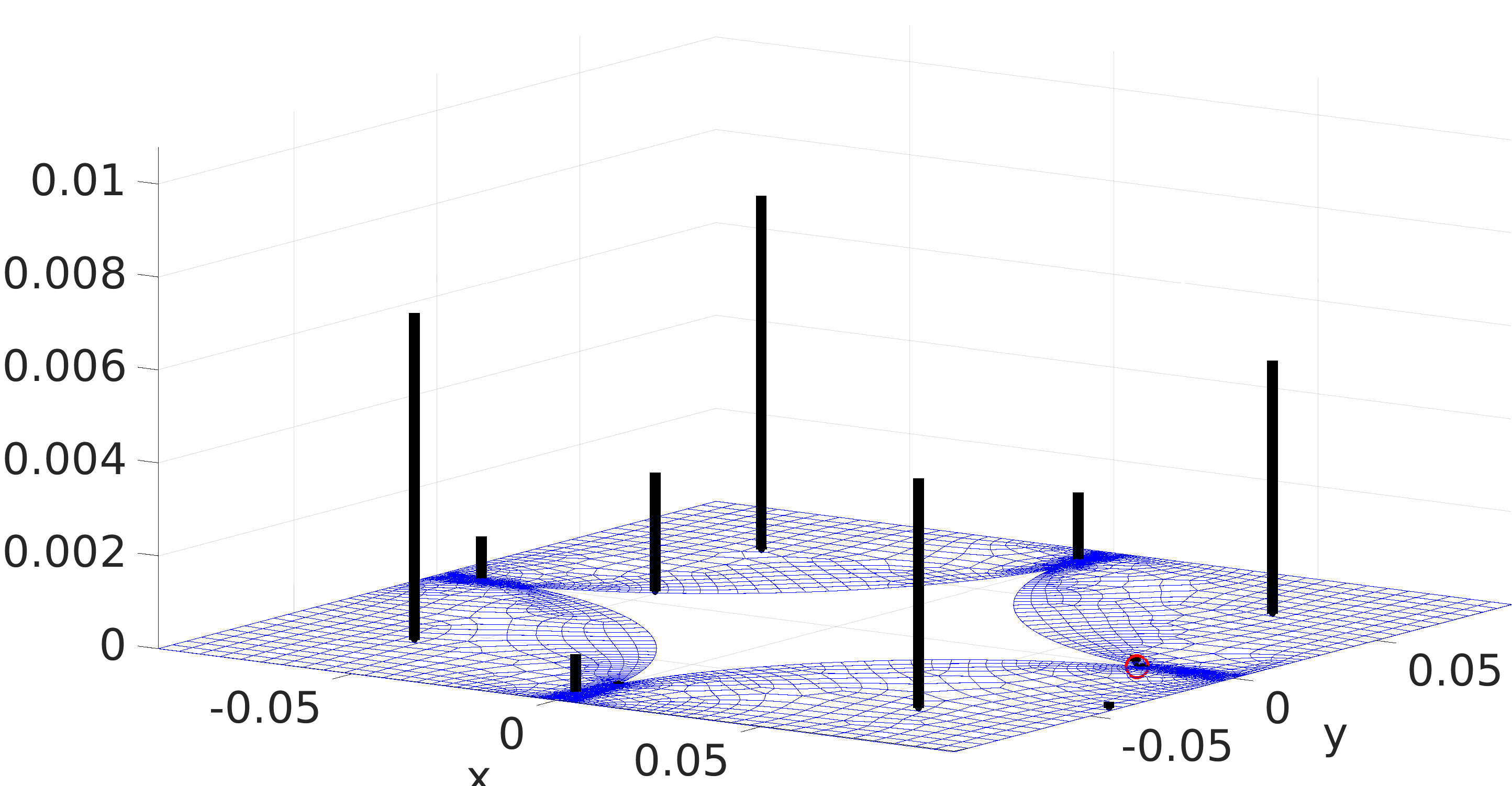}}  
      \subfigure[$m = 12$ points ($t=1,k= 3$)]{\label{fig:unitD12}\includegraphics[width=0.325\textwidth]{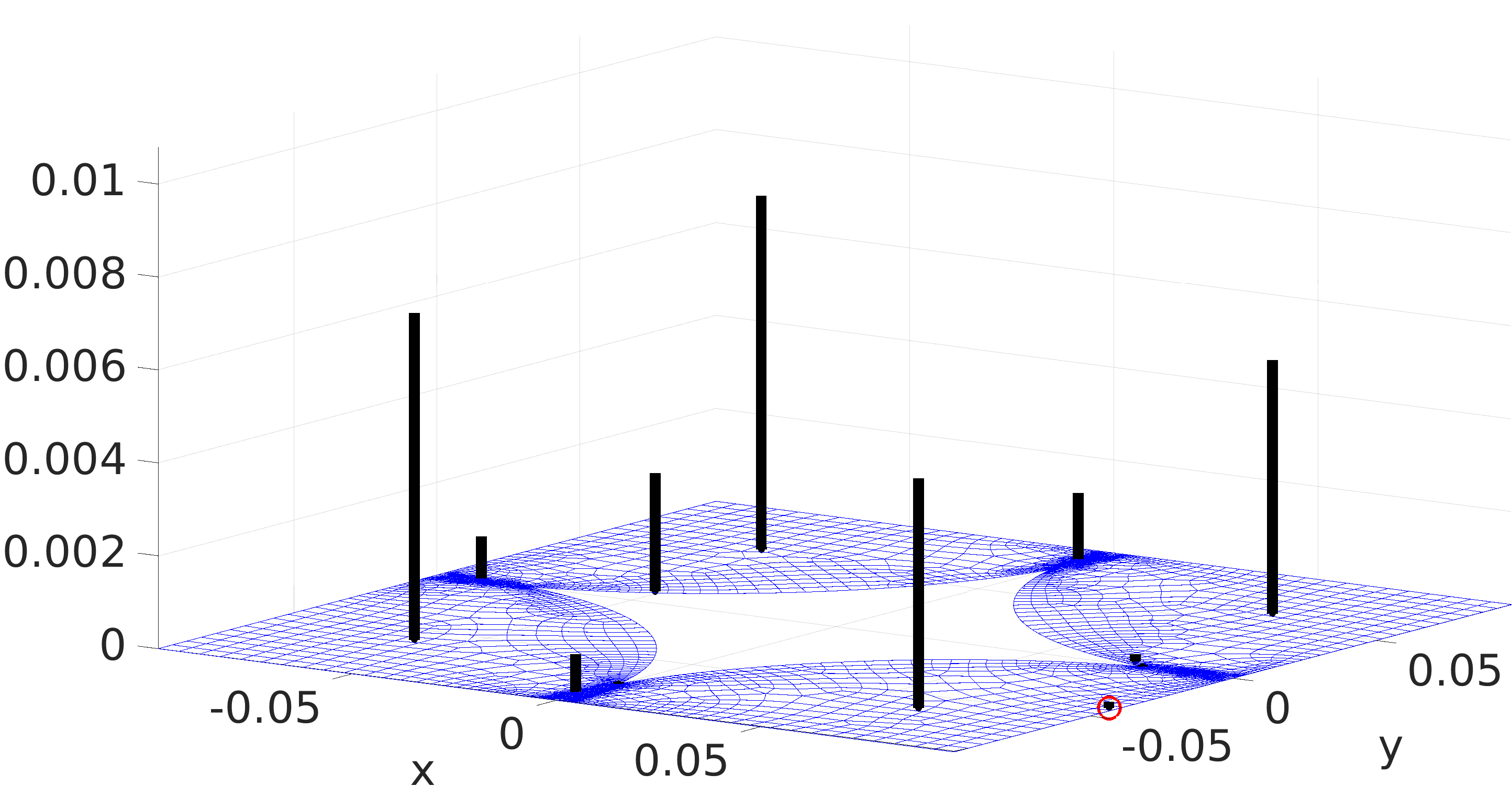}}  
       \subfigure[$m = 11$ points ($t=1,k= 5$)]{\label{fig:unitD11}\includegraphics[width=0.325\textwidth]{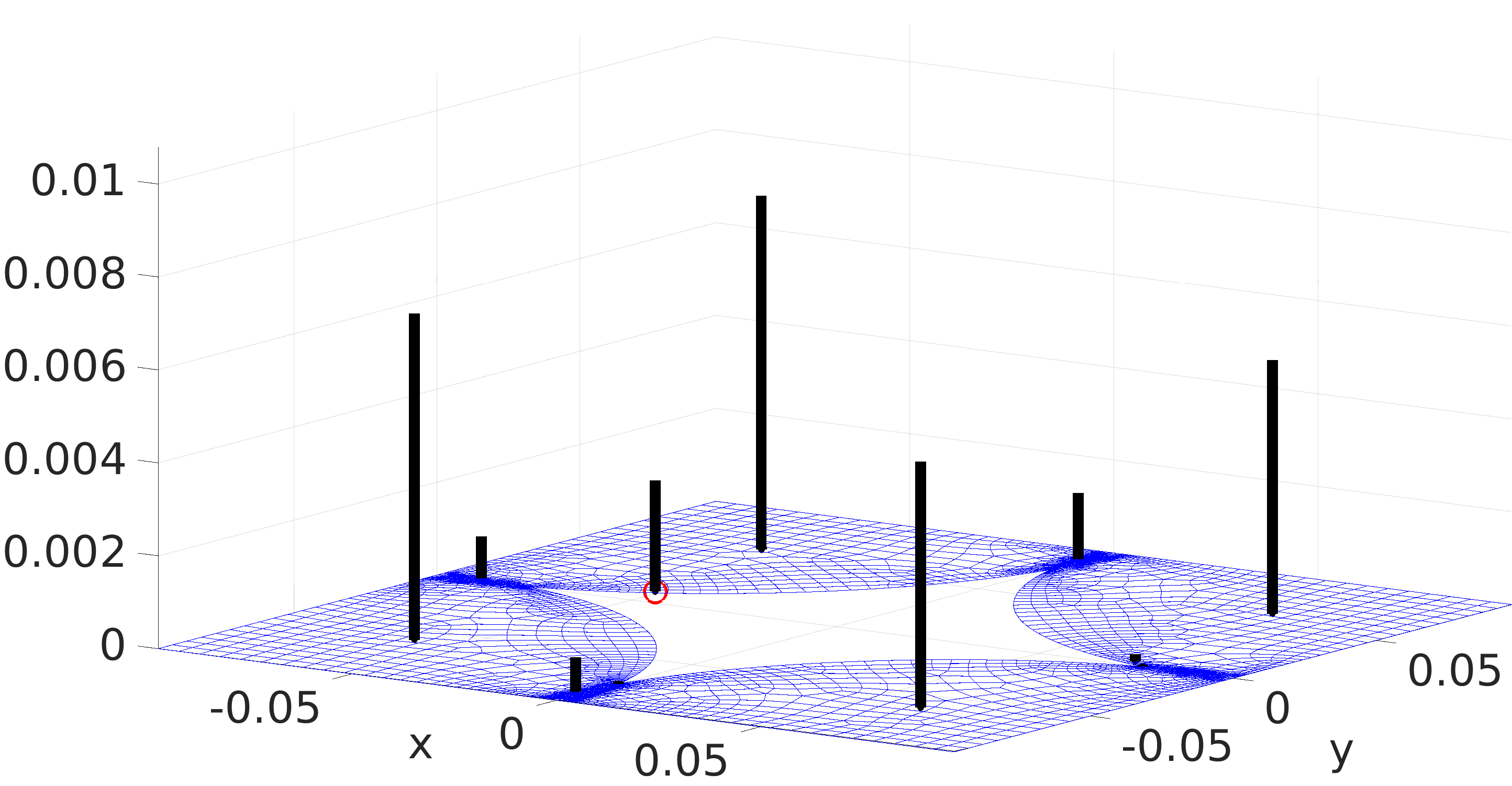}}  
   \subfigure[$m = 10$ points ($t=1,k= 12$)]{\label{fig:unitD10}\includegraphics[width=0.325\textwidth]{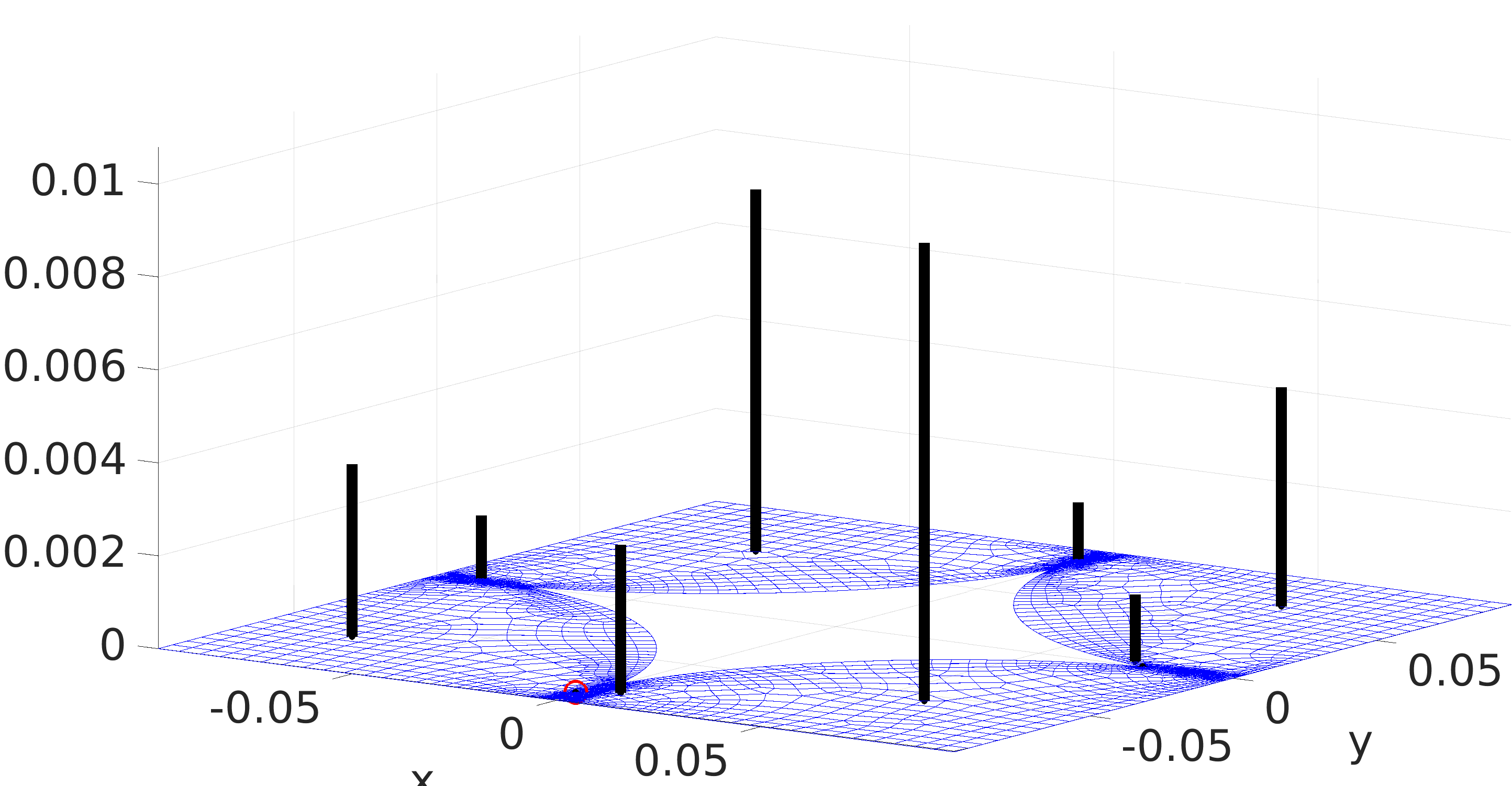}} 
    \subfigure[$m = 9$ points ($t=1,k= 7$) ]{\label{fig:unitD9}\includegraphics[width=0.325\textwidth]{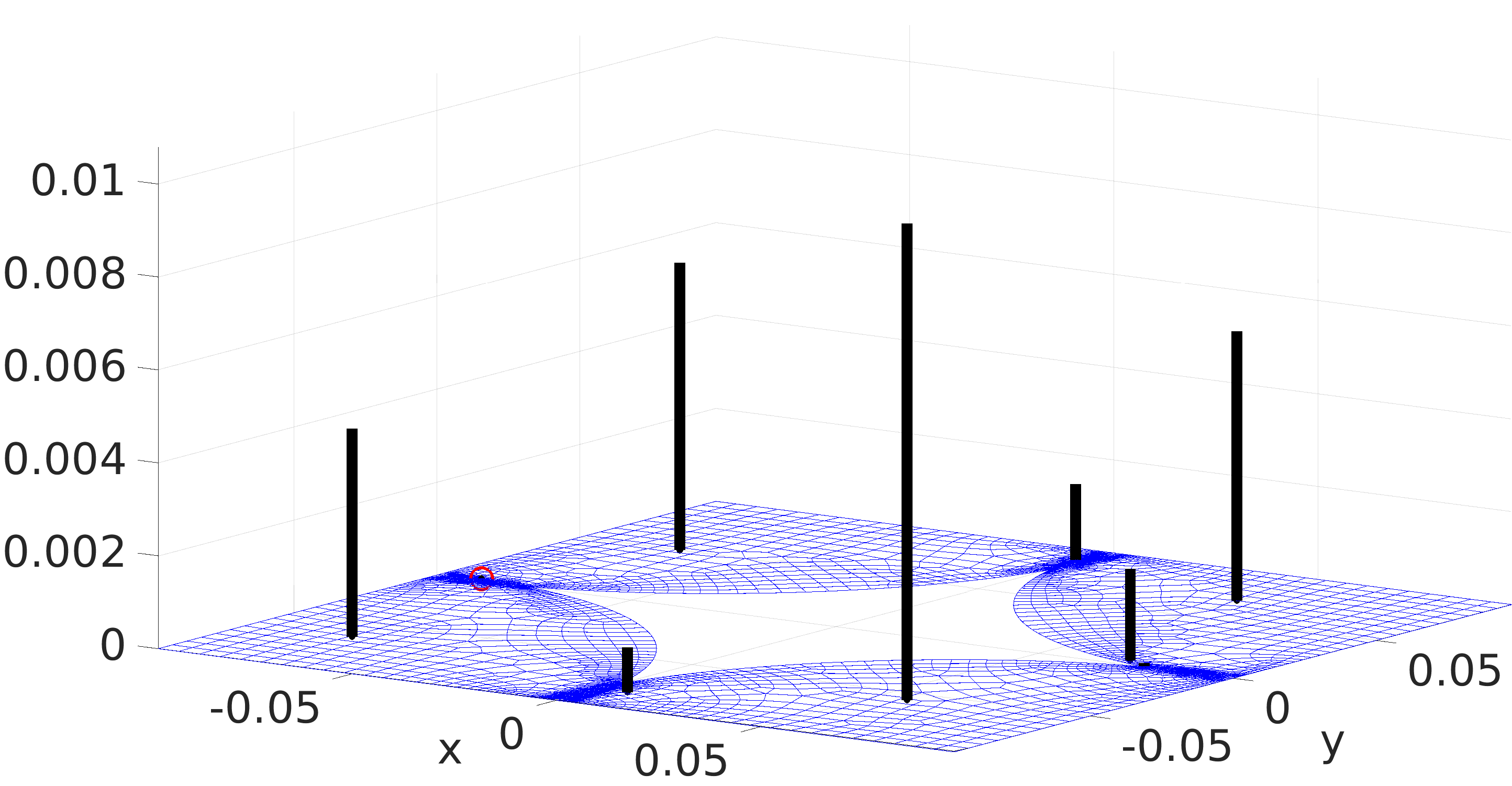}} 
        \subfigure[$m = 8$ points  ($t=1,k= 13$)]{\label{fig:unitD8}\includegraphics[width=0.325\textwidth]{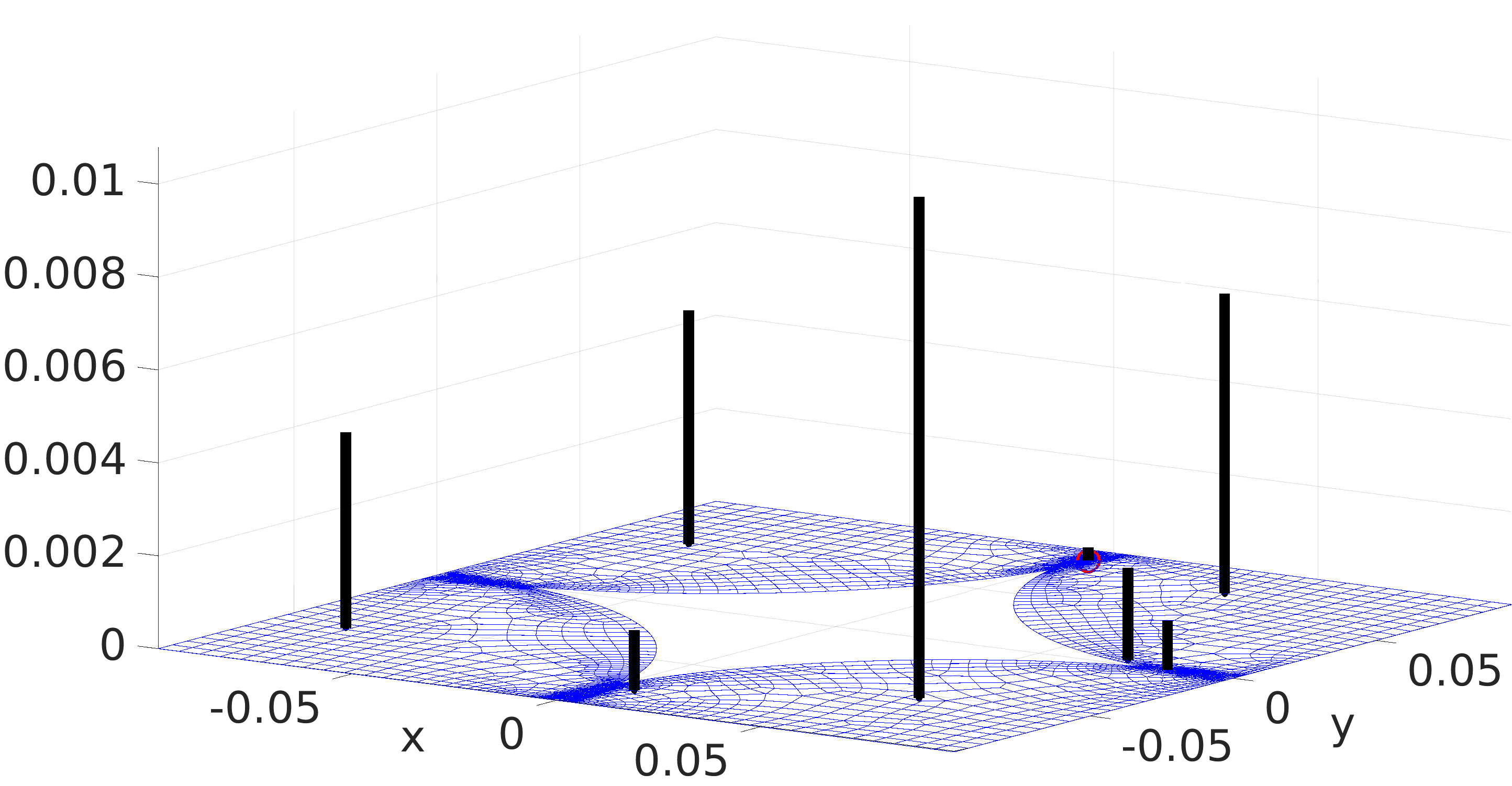}}
                \subfigure[$m = 7$ points ($t=1,k= 5$)]{\label{fig:unitD7}\includegraphics[width=0.325\textwidth]{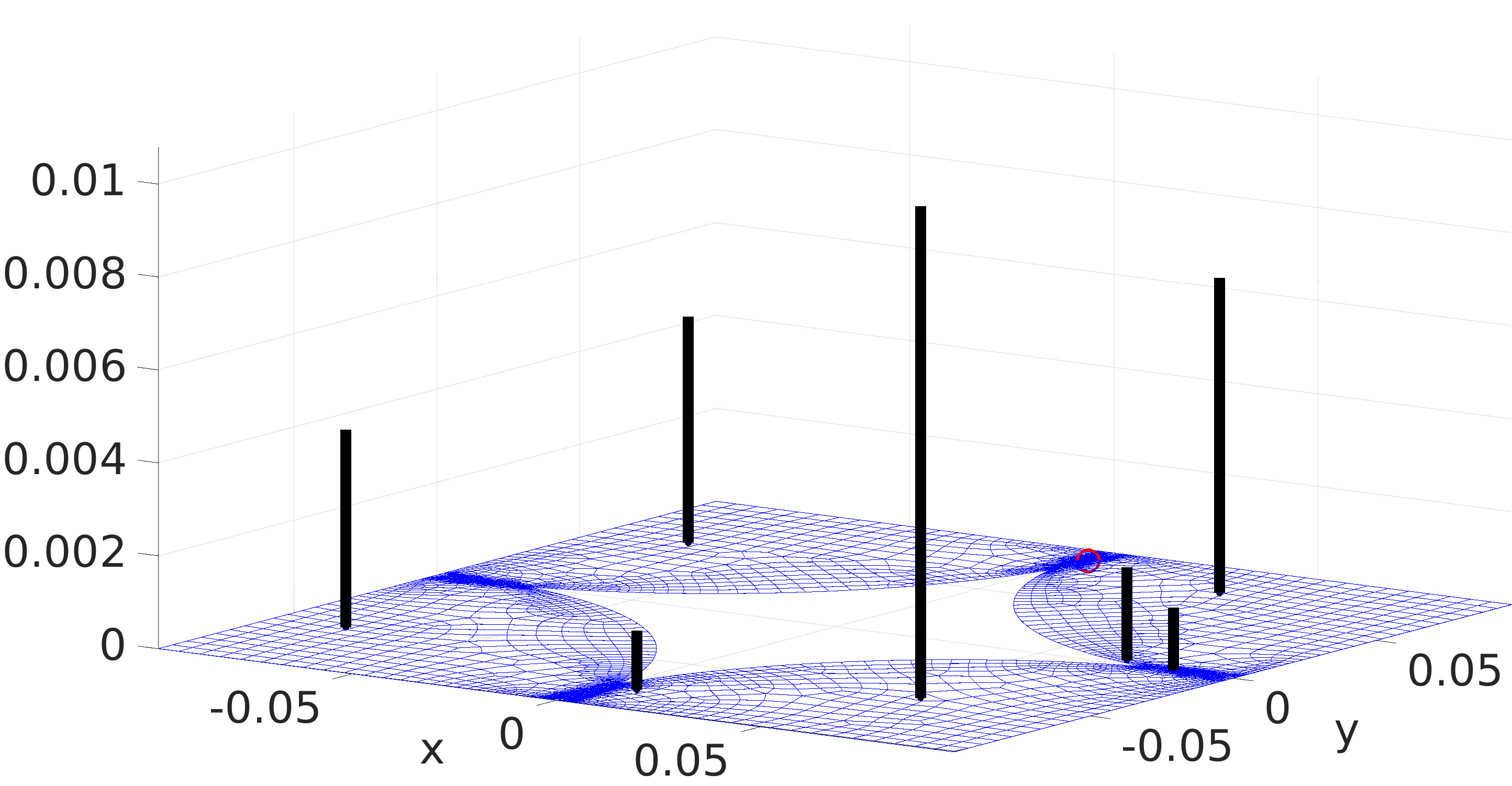}}
                                \subfigure[  CECM rule, $m = 6$ points ($t=1,k= 97$)]{\label{fig:unitD6}\includegraphics[width=0.325\textwidth]{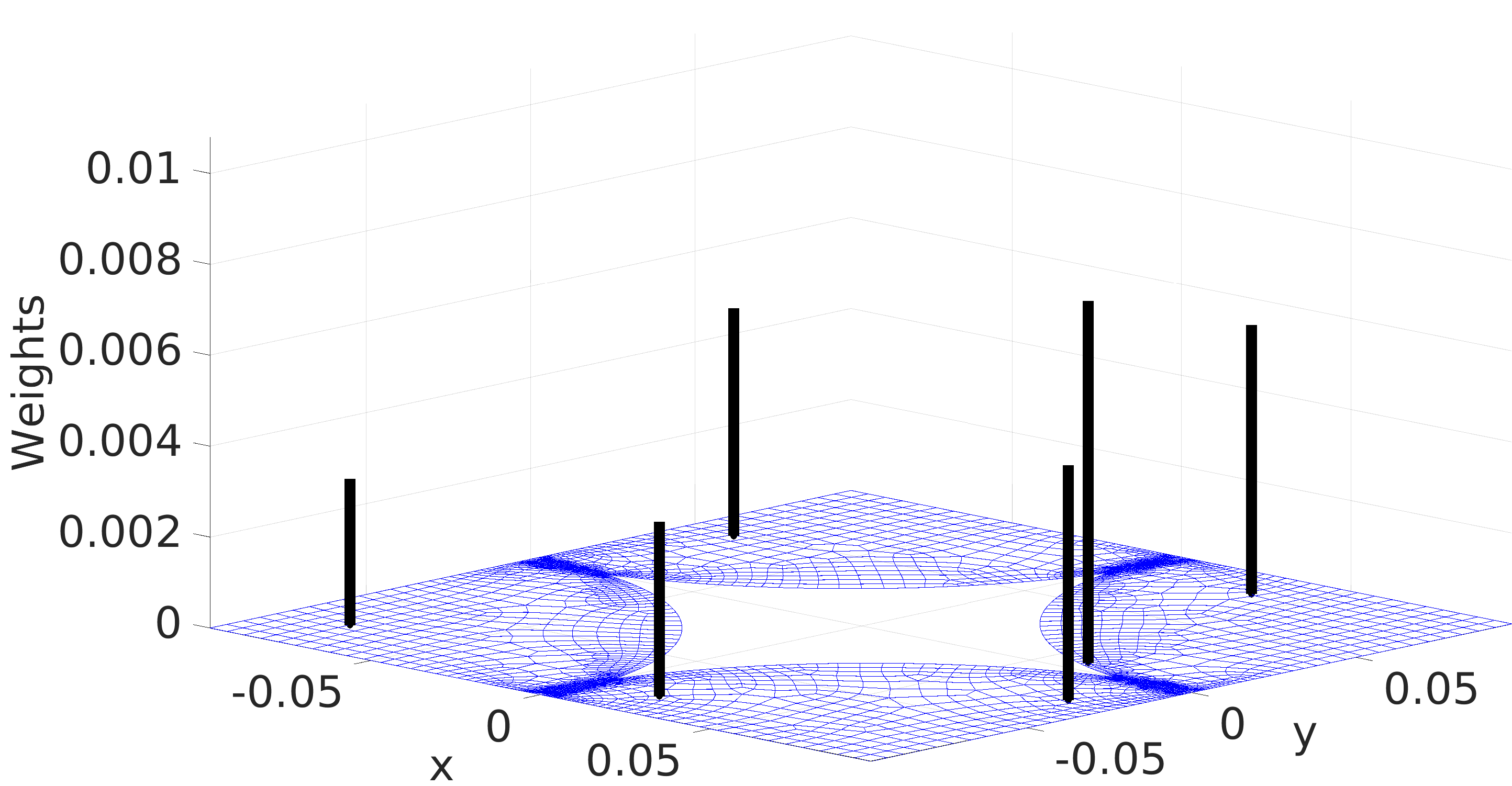}}
  \caption{   Locations  and   weights   of the cubature rules generated during  the sparsification process for the case of the virtual internal work per unit volume in Eq. \ref{eq:internal} (integrand functions shown previously in Figure \ref{fig:10}).   The red circle in each graph indicates the point whose weight is to be zeroed in the following step. Variables $t$ and $k$, on the other hand, have the same interpretation as in Figure \ref{fig:FIG_5}.  The initial DECM rule has $16$ points (see Figure \ref{fig:unitD16}), while  the final   rule features  $m=6$ points, see Figure \ref{fig:unitD6}.  The total reduction in number of integration points is   $M/m =42471/6 = 7078.5$.     }
  \label{fig:11}
\end{figure}

We conclude the validation of the proposed CECM by illustrating its use  in the hyperreduction of finite element  models.  More specifically, attention is concentrated in the derivation of low-dimensional surrogate models in the context of  {multiscale } finite element methods.  The employed multiscale methodology is the \emph{Empirical Interscale FE }  (EIFE) method  developed by the first author and co-worker in Refs.  \cite{giuliodori2023multiscale} (for beam elements) and  Ref. \cite{hernandez2023multiscale}  (for  solid elements). The chosen example is the  modeling  of the periodic structure displayed in Figure \ref{fig:9}.a, made of tiling copies of the unit cell displayed in turn in Figure \ref{fig:9}.b.  The FE mesh of this unit cell, also depicted in Figure \ref{fig:9}.b, is formed by     19356 nodes, and $\nelem =   4719$ quadratic quadrilateral elements,  with $\ngausE = 9$ Gauss points per element (hence $M = \nelem \ngausE =   42471$). The areas in which stress concentration are likely to appear are densely meshed to   avoid both highly distorted elements and    pronounced interlemental jumps ---both issues are detrimental, not only to the accuracy of the FE analysis per se, but also to the  accuracy of the final cubature rule.

The goal in the EIFE method is to replace the  fine-scale representation of the unit cell in Figure \ref{fig:9}.b by a surrogate coarse-scale element such as the one shown in Figure \ref{fig:9}.c.  The coarsening process involves two sequential stages;  the first stage involves    the reduction of number of  degrees of freedom (DOFs), and the second  stage involves an additional    reduction in complexity in  terms  of    number of integration points for each term of the governing equations (this is why the second stage is called    ``hyperreduction'').  The particular details on how the first reduction of number of DOFs is carried out are of no concern here ---the reader interested in such details is referred to Ref. \cite{hernandez2023multiscale}. It suffices to say that the process involves  running firstly  FE analyses in   domains comprising several cells (such as the one in Figure \ref{fig:9}.a) under appropriate prescribed displacements, and then extracting characteristic deformational patterns of the unit cell in the center via the SVD. In the case studied here, in which the coarse model has 8 DOFs and the FE analyses are conducted, for  simplicity of exposition\footnote{The procedure can be applied to any constitutive stress-strain law. },  in the  linear elastic regime (Young's Modulus $E= 70000 $ MPa and Poisson's ration $\nu=0.3$, in plane strain), the resulting deformational patterns  are the $n=5$ modes whose deformed shapes are   displayed in Figure \ref{fig:9}.d.  


Our interest lies in the hyperreduction stage, in which one has to devise efficient cubature rules for each one of the integrals appearing in the (reduced) governing equations. We focus here in the   internal force term, but a similar procedure can be followed for the integrals associated to    body forces and surface traction.  The reduced internal forces are given by the projection onto the space spanned by the deformational modes $\PhiB$ of the   FE nodal internal forces $\boldsymbol{F_{int}}$: 
\begin{equation}
\label{eq:internal}
 \PhiB^T  \boldsymbol{F_{int}}  = \intG{\Omega}{\Omega}{\PhiB^T  \B_{FE}^T \stress}. 
\end{equation}
Here, $\B_{FE}$ denotes the standard strain-displacement FE matrix (in its globally supported format), whereas $\stress$ is  the   stress vector. Our parametric integrand is therefore equal to  $\a = (\B_{FE} \PhiB)^T \stress$, the work per unit volume done over virtual strains of the form $\B_{FE} \PhiB_i$ ($i=1,2 \ldots n$) by  the  stresses $\stress$ caused in turn by strains also of  the form $\strain = (\B_{FE} \PhiB)  \q$ (Galerkin projection).  It follows then that   $\iinpP{} = \q$, that is,  in this problem,  the amplitude of the deformational modes in the expression for the stresses plays the role of input parameters. Since we are assuming linear elastic behavior,   the number of possible stress states is equal to $n=5$ as well. This implies that the integrand matrix $\Afe{}{}$   is formed by  $n^2 = 25$  columns, which are all the possible combinations of stress modes times virtual strain modes.  The SVD (step 4 in Box \ref{box:1} ) of this matrix ( using a  fairly low truncation tolerance, $\epsilon_{svd} = 10^{-10}$, to eliminate numerical errors) reveals that there are redundant work modes: out of the $n^2 = 25$ work modes, only 15 are linearly independent\footnote{ It can be readily seen that  the existence of these 10 redundant work patterns is nothing but the consequence of the symmetry of the elastic problem ---Betti's reciprocity theorem \cite{reddy2017energy}. In general, if there are $n$ deformational modes, the number of independent work modes will be equal to $(n+1)n/2$.  }. By augmenting the integrand basis with an additional vector for accounting for the integration of the volume (step 5 in Box \ref{box:1}  ), we end up with 16 basis functions for the integrand; the contourn lines of these 16 functions are depicted in Figure \ref{fig:10}.

The 16-points interpolatory DECM  rule corresponding to these 16 functions is displayed in turn in Figure \ref{fig:11}.a. The fact that the weights of  8 of the 16 points are comparatively small gives an indication that the number of points can be further reduced by the proposed sparsification process. The result of this     sparsification process is displayed in Figures \ref{fig:11}.b to  \ref{fig:11}.k.  The number of trials and accumulated nonlinear iterations are also shown in the captions of these Figures.  As observed with the 3D analytical function in Section \ref{sec:expsino}, the number of iterations required for the last 10 \% of the weights   (which in this case is just the last step)   requires significantly more iterations  than the previous 90 \%.     Nevertheless, the total computing time is not signficant ---less than 6 seconds for the entire procedure, including the computation of the DECM rule. To assess the final integration error ,  we compare the reduced stiffness matrix
\begin{equation}
 \PhiB^T  \K  \PhiB  = \intG{\Omega}{\Omega}{\PhiB^T  (\B_{FE}^T \C \B_{FE})  \PhiB^T}. 
\end{equation}
(here $\K$ is the nodal stiffness matrix of the unit cell, and $\C$ the corresponding elasticity matrix) computed by both the original 42471-points Gauss element rule and the 6-points CECM rule. The difference turns out to be below 0.005 \%. This result demonstrates that,  one the one hand, that the computed  6-points CECM effectively encodes the physics of the coarsened unit cell, and on the other hand, that the proposed cubature algorithm is able to deal with complex domains in which the integrand function is only defined at some points of such domain ---the FE Gauss points.  In fact, it should be remarked that, since   the tolerance employed in the SVD is negligible,    this   small error is to be exclusively attributed to the  fitting procedure outlined in Section \ref{sec:InterpolationG} for constructing approximations of the integrand at element level .

\label{sec:conclusions}
\section{Conclusions}

In this paper we presented the Continuous Empirical Cubature Method (CECM), a novel algorithm designed to enhance the efficiency of numerical integration rules. The CECM is a two-stage algorithm whose first stage consists on the application of a point selection strategy for obtaining an interpolatory rule ---featuring as many points as functions to integrate. We have used for this purpose the Discrete Empirical Cubature Method (DECM) \cite{ hernandez2017dimensional,hernandez2020multiscale}. Then, for the second stage, we have applied a sparsification strategy whose aim is rendering to zero the associated weight of as many points as possible. To this end, the locations of the initially selected points were changed following a modified Newton-Raphson algorithm whose details have been outlined within the text (The code has also been made available at \url{https://github.com/Rbravo555/CECM-continuous-empirical-cubature-method}).  

The versatility of the method was highlighted in the numerical assessment section, showcasing its effectiveness across a diverse range of problems. For the case of univariate and multivariate Lagrange polynomials, the CECM was able to recover the optimal Gaussian rule whenever the number of function to integrate was odd. For even number of functions to integrate, the number of points still coincides with the Gauss rule, although their locations differ. Reductions in required points were observed, with 1D domains requiring half, 2D domains requiring one fourth, and 3D domains requiring one eighth of the initial interpolatory rule.

Example \ref{sec:expsino} showcases an exponential-sinusoidal function in a 3D domain. For this case, the CECM reduced the number of points from 133 of the original DECM rule, to 38. This example also showcased an scenario where the size of the matrices involved would render infeasible the computation of the SVD on a regular desktop computer. The SVD was computed using the sequential randomized SVD (SRSVD) algorithm presented also in this paper. The SRSVD allowed to efficiently compute the required orthogonal basis for matrices of up to 33 GB in size. 

Finally, the CECM was applied to an empirical interscale finite element (EIFE) example, for which the number of points was reduced from the original fine scale 42471 points to a 6-point rule, while incurring a negligible 
integration error of 0.005\%. 

Overall, the CECM algorithm presents an innovative approach to generating optimal integration rules. Its incorporation into frameworks that allow input of quadrature points' positions leads to substantial improvements in integral evaluation performance compared to standard interpolatory rules.

\section*{Acknowledgements}

This work is sponsored in part by  the Spanish Ministry of Economy and Competitiveness,
through the \emph{Severo Ochoa Programme for Centres of Excellence in R\&D} (CEX2018-000797-S)”. The authors also acknowledge the support of  the European High-Performance Computing Joint Undertaking (JU) under grant agreement No. 955558 (the JU receives, in turn,  support from the European Union’s Horizon 2020 research and innovation programme and Spain, Germany,  France, Italy, Poland, Switzerland, Norway).   J.A. Hern\'{a}ndez  expresseses gratitude by the  support of, on the one hand,   the   ``MCIN/AEI/10.13039/501100011033/y \emph{por FEDER una manera de hacer Europa}'' (PID2021-122518OB-I00), and, on the other hand, the European Union's Horizon 2020 research and innovation programme under Grant Agreement No. 952966 (project FIBREGY).  Lastly, both J.R. Bravo and S. Ares de Parge acknowledge the \emph{Departament de Recerca i Universitats de la Generalitat de Catalunya} for the financial
support through doctoral grants FI-SDUR 2020 and  FI SDUR-2021, respectively.

 \appendix 
 
 \section{Sequential Randomized SVD}
 \label{sec:SRSVD}
 
\subsection{Overview}

In this Appendix we explain and provide the implementational details (see Algorithm \ref{alg:006} ) of the proposed procedure  for computing the SVD of a partitioned matrix 
 $\A= [\A_1,\A_2 \ldots \A_p]$ ($\A_i \in \RRn{n}{m_i}$, with $m = \sum_{i=1}^p m_i$), alluded to in  Remark \ref{rem:001}.  The method is based on the same idea behind other \emph{randomized} algorithms \cite{halko2011finding}, namely, that the SVD of a matrix $\A$ ($\A = \U \S \V^T$)  can be alternatively computed   from the matrices of the SVD of $\L = \Q^T \A$,    $\Q \in \RRn{n}{r}$ being an arbitrary orthogonal basis matrix for the column space of $\A$ (here $r \le \min{(n,m)}$ denotes the rank of the matrix). More specifically,    given $\L = \Ubar \Sbar \Vbar^T$, then $\U = \Q \Ubar$, $\S = \Sbar$ and $\Vbar = \V$.   The proof of this property follows from expressing $\A$ as $\A = \Q (\Q^T \A)$ and replacing $\Q^T \A$ by its SVD: 
\begin{equation}
 \A = \Q \Q^T =  {(\Q \Ubar)}  \Sbar \Vbar^T. 
\end{equation}
Both $\Sbar$ and $\Vbar$ arise from an SVD, and therefore, are 
  diagonal with positive entries and orthogonal, respectively ---and $\Vbar$ is a basis matrix for the row space of $\L = \Q^T \A$, which is the same as the row space of $\A$.   Furthermore, $\U$ is   also an orthogonal matrix: 
\begin{equation}
 \U^T \U = (\Q \Ubar)^T (\Q \Ubar) = \Ubar^T (\Q^T \Q) \Ubar =  \ident. 
\end{equation}
Therefore, it follows from the uniqueness of the SVD (up to the signs of the left- and right- singular vectors)   that the
factorization  $\U \Sbar \Vbar^T$ is the   SVD of $\A$, as asserted. 

It can be readily shown that this property also holds for the case in which \emph{truncation} is introduced. Besides, since $\S = \Sbar$, the   truncation threshold for the SVD of $\L$ is the same as the input truncation threshold\footnote{This is not exactly true in the limiting case $\epsilon \rightarrow 0$, when the truncation criterion is established in terms of a machine-dependent precision parameter, which also depends on the size of the matrix being factorized (see Algorithm \ref{alg:006}, Line \ref{alg:006_finite}) } $\epsilon$, that is,  
\begin{equation}
 \normF{\L - \Ubar \Sbar \Vbar^T } \le \epsilon \normF{\L} =  \epsilon \normF{\A}.  
\end{equation}

This general strategy for computing the SVD of a matrix proves advantageous only when the following two conditions are met.  Firstly,    the rank of the matrix $r$ should be significantly  smaller\footnote{ This first condition may be relaxed by making $\A \approx \Q \Q^T \A$, that is, by determining an  {approximated} basis matrix for the column space of $\A$.   However, it should be noticed that, in this case,  the right-singular vector matrix $\V$ ceases to be  a basis matrix for the row space of $\A$, and therefore,  Eq. \ref{eq:follws1} would not hold.} than   the  number of columns and rows of the matrix  ---because otherwise  the SVD of $\L = \Q^T \A$ could become as costly as the original SVD. In the context of reduced-order models, this property is expected to hold ---and if it does not,   it means that the    parameterized boundary value problem we intend to solve might not be amenable to dimensionality reduction. The other condition  is that  the computation of the orthogonal basis matrix $\Q$ should be efficient, in the sense that it should be carried out by an algorithm in which     the asymptotic  count of floating point operations (flops)   is less than that required by the standard  SVD itself.

\begin{algorithm}[ht!]
\DontPrintSemicolon

\SetKwFunction{FMain}{SRSVD}
\SetKwFunction{RORTH}{SRORTH}
\SetKwProg{Fn}{Function}{}{}

\Fn{ $[\boldsymbol{U},  \boldsymbol{S}, \boldsymbol{V}]$   $\leftarrow$  \FMain{$[\A_1,\A_2,\ldots \A_p]$, $\epsilon$} }{ 
\KwData{$[\A_1,\A_2,\ldots \A_q]$:  Partitioned matrix (with $\boldsymbol{A}_i \in \mathbb{R}^{n \times m_i}$, and $m = \sum_{i=1}^p m_i$ ).  $\epsilon \in  [0,1] $  : relative error threshold  
} 

\KwResult{Truncated  Singular Value Decomposition (with relative truncation threshold $\epsilon$) of  $\A \in \RRn{n}{m}$, i.e.: $\A \approx \U  {\S} \V^T$, where  $\boldsymbol{U} \in \mathbb{R}^{n \times k}$ \ is the orthogonal matrix of left singular vectors,   $\boldsymbol{S} \in \mathbb{R}^{k \times k}$ is the diagonal matrix of positive singular values, and  $\boldsymbol{V} \in \mathbb{R}^{  m \times k}$ \ is the orthogonal matrix of right singular vectors. Here $k \le \min{(n,m)}$ denotes the number of retained singular vectors upon truncation, which is, by definition of SVD,  the lowest number of vectors such that $\normF{\boldsymbol{\A - \U  {\S} \V^T }} \le  \epsilon \normF{\A}$. 
} 
} 

\vspace{2mm} 

[$\boldsymbol{Q}$, $\L$] $\leftarrow$ \RORTH{$[\A_1,\A_2,\ldots \A_p]$} 
 \tcp{ Factorization $\A  = \Q \L$, where $\Q^T \Q = \ident$ \label{alg:006_ortho} ( see Algorithm \ref{alg:008})}

 $\epsilon_L \leftarrow   \epsilon \, \normF{\L}  $; $\; \mu_{mach} = \max{(n,m)} \, \texttt{eps}(\normF{\L}) $  \tcp{\label{alg:006_finite} Truncation thresholds (see definition of $\texttt{eps}$ in Algorithm \ref{alg:007}).   }

$[\tilde{\boldsymbol{U}}, \boldsymbol{S}, \boldsymbol{V}] \leftarrow$ \texttt{SVDT} ($\boldsymbol{L}$, $\epsilon_L$,$ \mu_{mach} $)  \tcp{Truncated SVD   (see Algorithm \ref{alg:007}). \label{alg:006_svd}.} 

$\boldsymbol{U} \leftarrow \boldsymbol{Q} \tilde{\boldsymbol{U}}$   \tcp{Matrix of left singular vectors \label{alg:006_u}}

\caption{ Sequential  Randomized Singular Value Decomposition (SRSVD) of a partitioned matrix $\A = [\A_1,\A_2,\ldots \A_p]$.  }
\label{alg:006}
\end{algorithm}

     \begin{algorithm}[!ht]
   \DontPrintSemicolon   
     \SetKwFunction{RORTH}{RORTH}
       \SetKwFunction{max}{max}
            \SetKwFunction{min}{min}
              \SetKwFunction{break}{break}
               \SetKwFunction{eps}{eps}
               \SetKwFunction{SVDT}{SVDT}
                    \SetKwFunction{svd}{svd}
    \SetKwProg{Fn}{Function}{:}{}
   \Fn{ $[\U,\S,\V] \leftarrow$ \SVDT{$\C,\epsilon,\mu_{mach}$} }{
    
    \KwData{ $\C \in \RRn{m}{n}$. Absolute tolerance threshold $\epsilon >0$. Machine-precision parameter $\mu_{mach} > 0$ (optional input, default value $\mu_{mach} = \max{n,m} \texttt{eps}(\normF{\C})$, where $\texttt{eps}(x)$ denotes the  machine-dependent floating-point relative precision associated to $x$.   }
    
    \KwResult{ Truncated SVD such that  $\normF{\C - \U  {\S} \V^T}  \le  \epsilon $ (if $\epsilon \le \mu_{mach}$, the truncation criterion changes to the one specified in  Line \ref{alg:007_rank}).   }

    $[\U,\S,\V] \leftarrow$ \svd{$\C$ }  \tcp{Thin SVD of $\C$  (In matlab's built-in function \texttt{svd($\C$,'econ')}}

 \uIf{$\epsilon \le \mu_{mach}$ \label{alg:007_mu}}{  $\;$  $k \leftarrow$ Smallest singular value  such that $\S(k,k) \ge \mu_{mach}$   \tcp{Numerical rank}  \label{alg:007_rank}}
     \Else{$k \leftarrow$ Find smallest $k$ such that $\sqrt{\sum_{i = k +1}^{n}  \S(i,i)^2}  \le \epsilon$  \tcp{Truncation level}  \label{alg:2:L9}}
     $\U \leftarrow \U(:,1:k)$, $\V \leftarrow \V(:,1:k)$, $\S \leftarrow  \S(1:k,1:k)$

   }
   \caption{Truncated Singular Value Decomposition (via standard, deterministic methods) \label{alg:007}} 
   \end{algorithm}

  \subsection{Sequential Randomized Orthogonalization }
  
  \subsubsection{Infinite-precision arithmetic}
  
  To meet this latter condition, we propose to determine $\Q$ by  the \emph{Sequential Randomized Orthogonalization} (SRORTH) invoked in Line \ref{alg:006_ortho} of Algorithm \ref{alg:006}, and with pseudo-code outlined in Algorithm \ref{alg:008}. The qualifier \emph{sequential}  refers to the   fact that the method only processes \emph{one block matrix at a time}, thus alleviating potential memory bottlenecks.   On the other hand, we call it  \emph{randomized}  because one of the factorizations employed by the algorithm is carried out by a modified version of the \emph{incremental randomized   SVD }proposed in Ref. \cite{martinsson2015randomized}.

  Let us describe first the overall structure of this orthogonalization procedure, without delving into  the randomized part of the algorithm, which will be treater later on, in Section \ref{sec:irsvd}.    In essence, the procedure   is a Gram–Schmidt orthogonalization which operates, rather than on single vectors,  on block matrices. Accordingly,   $\Q$ is constructed  as  the concatenation of  basis matrices (one for each block matrix $\A_i$):  
\begin{equation}
 \Q = \rowcuatro{\dQ_1}{\dQ_2}{\cdots}{\dQ_p}
\end{equation}
where $\dQ_i^T \dQ_j = \zero$ if $i \neq j$, and $\dQ_i^T \dQ_i = \ident$ ($i=1,2 \ldots p$). In turn, these orthogonal submatrices are computed by the recursion 
  \begin{align}
  \nonumber  \dQ_1  = \texttt{ORTH}(\A_1); \phantom{- \Q \Q^T \A_1} \hspace{2.4cm}   \Q^{(1)} = \dQ_1;   \hspace{1.9cm}  \P_1 & = {\Q^{(1)}}^T \A_1  \\
  \nonumber \dQ_2  = \texttt{ORTH}(\A_2 - \dQ_{1}  {\dQ_{1}}^T \A_2);  \hspace{1.30cm}  \Q^{(2)} = [\Q^{(1)},\dQ_2]; \hspace{0.85cm}  \P_2 & = \Q^{(2)^T} \A_2 \\
    \nonumber  \vdots   \hspace{2cm}  \vdots  \hspace{2cm} &   \\
   \dQ_i  = \texttt{ORTH}(\A_i - \Q^{(i-1)} \Q^{(i-1)^T} \A_i);  \hspace{1cm}  \Q^{(i)} = [\Q^{(i-1)},\dQ_i]; \hspace{0.65cm}  \P_i & = \Q^{(i)^T} \A_i \label{eq:recurses} \\  
   \nonumber   \vdots   \hspace{2cm}  \vdots  \hspace{2cm} &   \\
  \nonumber \dQ_p  = \texttt{ORTH}(\A_p - \Q^{(p-1)} \Q^{(p-1)^T} \A_p);  \hspace{1cm}  \Q  = [\Q,\dQ_p]; \hspace{1.55cm}  \P_p & = \Q^T  \A_p
  \end{align}
 (here $\texttt{ORTH}(\bullet)$ symbolizes the function  that determines an orthogonal basis matrix for the column space of its input). Notice that this procedure for determining $\Q$ need not store in the computer's main memory the entire matrix, but  just one block matrix at a time, as asserted earlier. The other matrix in the factorization, $\L = \Q^T \A$, can be also constructed incrementally as the algorithm progresses\footnote{Note that this statement would cease to be true in the case of approximated basis matrices for the range of $\A$ ($\A \approx \Q \Q^T \A$). }.  Indeed, by exploiting that $\dQ_j^T \A_i = \zero$ if $j>i$, we have that  
  \begin{equation}
  \label{eq:LL1}
 \L = \begin{bmatrix}
    \dQ_1^T\A_1 & \dQ_1^T \A_2 & \cdots & \dQ_1^T  \A_p \\  
    \dQ_2^T\A_1 & \dQ_2^T \A_2  &\cdots &\dQ_2^T \A_p \\   
   \vdots & \vdots & \ddots & \vdots  \\ 
    \dQ_p^T\A_1 & \dQ_p^T \A_1 & \cdots & \dQ_p^T \A_p \\   
    \end{bmatrix}  
    =  \begin{bmatrix}
    \dQ_1^T\A_1 & \dQ_1^T \A_2 & \cdots & \dQ_1^T  \A_p \\  
    \zero & \dQ_2^T \A_2  &\cdots &\dQ_2^T \A_p \\   
   \vdots & \vdots & \ddots & \vdots  \\ 
    \zero & \zero  & \cdots & \dQ_p^T \A_p \\   
    \end{bmatrix}.   
\end{equation}
Inspection of the nonzero entries in each column of the right-most matrix of the above equation  shows that these entries are expressible in terms of the matrices $\P_i$ ($i=1,2 \ldots p$) appearing in the recursion formulas \ref{eq:recurses} as follows: 
\begin{equation}
 \label{eq:LL2}
 \dQ_1^T\A_1 = \P_1; \hspace{0.5cm}  \coldos{\dQ_1^T }{\dQ_2^T} \A_2 = \P_2; \hspace{0.5cm}  \ldots \hspace{0.5cm}  \colcuatro{\dQ_1^T }{\dQ_2^T}{\vdots}{\dQ_i^T} \A_i = \P_i \hspace{0.5cm}   \ldots \hspace{0.5cm}  \colcuatro{\dQ_1^T }{\dQ_2^T}{\vdots}{\dQ_p^T} \A_p = \P_p,
\end{equation}
  
\subsubsection{Finite-precision arithmetic}

The preceding recursive scheme would in principle  work seamlessly in an ideal infinite-precision arithmetic scenario. Yet the devil is in the details, and when moving to the real case of finite-precision arithmetic, its performance  is seriously afflicted by sensitivity to rounding errors and loss of orthogonality over multiple steps ---as it occurs with  the classical Gram-Schmidt orthogonalization \cite{golub2012matrix}. The computational implementation described in Algorithm  \ref{alg:008} incorporates ingredients that mitigate these   deleterious effects.  Re-orthogonalization is carried out in Line  \ref{alg:008} by determining the component of $\dQ$ which is orthogonal to the current basis $\Q$, and then applying the SVD again, setting the corrected $\dQ$ equal to the matrix of left-singular vectors of such a  decomposition. The effect of rounding errors, on the other hand, is treated by computing the  orthogonal basis matrix $\dQ$ as the left-singular vectors of the truncated SVD of $\dA = \A_i - \Q \Q^T \A_i$  (see Line \ref{alg:008_rsvd}), but using as truncation tolerance  a machine-dependent precision parameter based on the norm of $\A_i$: 
\begin{equation}
 \mu_{mach} = \MAX{n,m_i} \, \texttt{eps}(\normF{\A_i})
\end{equation}
(see Line \ref{alg:008_mu}) rather than the default option,  which would be in terms of the norm of the input matrix  $\dA$ ( the definition of $\texttt{eps}(x)$ is given in the description of Algorithm \ref{alg:007}).

\begin{algorithm}[!ht]
\DontPrintSemicolon

\SetKwFunction{FMain}{SRORTH}
\SetKwFunction{Ftol}{ToleranceComputation}
\SetKwProg{Fn}{Function}{}{}

\Fn{ $[\boldsymbol{Q}$, $ \L]$  $\leftarrow$  \FMain{$[\A_1,\A_2, \ldots \A_p]$ } }{

\KwData{
$[\A_1,\A_2,\ldots \A_p]$:  Partitioned matrix (with $\boldsymbol{A}_i \in \mathbb{R}^{n \times m_i} )$, and $m = \sum_{i=1}^p m_i$)   
} 

\KwResult{ Factorization $\A = \Q \L$, where $\Q \in \RRn{n}{r}$ is an orthogonal basis matrix for the column space of  $\A$ ($r \le \MIN{n,m}$ is the numerical rank of $\A$), while  $\L = \Q^T \A$.      } 

} 

\vspace{2mm} 

$\boldsymbol{Q}  \leftarrow \emptyset$, \hspace{0.1cm} $\P_i \leftarrow \emptyset$, \hspace{0.1cm}  $i=1,2 \ldots p$  \tcp{Initializations}


\For{$i = 1$  $\texttt{to}$ $p$ }{
$\dQ \leftarrow \emptyset$  \tcp{Incremental basis matrix}

\eIf{
$ i = 1 $
}{
$\dA \leftarrow \boldsymbol{A}_i;  \; \; \;   r \leftarrow \texttt{ceil}(0.01 \, \MIN{n,m_i})$  \tcp{$r$: Estimation rank $\A_i$  \label{alg:008_ini}}
} %
{
$\dA \leftarrow \A_i -\Q \Q^T \A_i  $  \tcp{Component of  $\A_i$ orthogonal to the column space of $\Q$ }
}

$ \mu_{mach} \leftarrow \MAX{n,m_i} \, \texttt{eps}(\normF{\A_i})$ \tcp{Machine-dependent precision parameter   \label{alg:008_mu}}


\If{ $\normF{\dA} >$ $\mu_{mach}$ }{
    [$\dQ$,$\bullet$,$\bullet$ ] $\leftarrow$ \texttt{RSVDinc}($\dA,0,\mu_{mach},r$) \tcp{\label{alg:008_rsvd} Incremental randomized SVD, see Algorithm \ref{alg:009} }
    \lIf{$i > 1$}{  $[\dQ , \bullet, \bullet ] \leftarrow$  \texttt{SVDT}($(\dQ - \Q \Q^T \dQ),0$) \tcp{Re-orthogonalization  }  } 
    
    $\boldsymbol{Q} \leftarrow [\boldsymbol{Q},  {\dQ}]$ \tcp{Basis matrix augmented with the columns of $\dQ$}
    
}
%
%

    $\P_i \leftarrow \Q^T \A_i$

$r \leftarrow \texttt{ncol}(\dQ)$   \tcp{ \label{alg:008_est} Number of columns of $\dQ$ (estimation for the rank of the next submatrix in Line \ref{alg:008_rsvd})}

}

$\L \leftarrow     $ Use $\P_1,\P_2 \ldots \P_p$ to construct $\L$ according to expressions \refpar{eq:LL1} and \refpar{eq:LL2}

\caption{ Sequential Randomized Orthogonalization of a partitioned matrix $\A = [\A_1,\A_2, \ldots \A_p]$ (employed in Line \ref{alg:006_ortho} of Algorithm \ref{alg:006}). }
\label{alg:008}
\end{algorithm}

\subsection{Incremental randomized SVD}
\label{sec:irsvd}

Let us focus now on the randomized ingredient of the method, which is the \emph{Incremental Randomized SVD}  ($\texttt{RSVDinc}$) invoked in     Line \ref{alg:008_rsvd} of Algorithm \ref{alg:008} for determining the basis matrix for   $\dA$.   As commented previously, this randomized SVD is partially based on the adaptive randomized    algorithm proposed in Ref. \cite{martinsson2015randomized},  and it is the actual ingredient that   renders the proposed scheme faster than the standard ``deterministic'' SVD  when the rank of the input matrix is significantly smaller than   the minor dimension of the matrix.  The reason is that, as argued in Refs. \cite{martinsson2015randomized} (see also \cite{martinsson2015randomized,martinsson2020randomized} ),   the  asymptotic cost of this type of randomization algorithms  is $\ORDER{ n m  r}$, $r$ being the rank of the matrix ---as opposed to the standard SVD, whose cost is independent of the rank and scales quadratically with its minor dimension\footnote{For a detailed account of the asymptotic costs of classical SVD, the reader is referred to Ref. \cite{trefethen1997numerical}, Lecture 31.}.

    \begin{algorithm}[!ht]
   \DontPrintSemicolon   
   \SetKwFunction{FMain}{RSVDinc}

     \SetKwFunction{RORTH}{RORTHinc}
       \SetKwFunction{max}{max}
            \SetKwFunction{min}{min}
              \SetKwFunction{break}{break}
               \SetKwFunction{eps}{eps}
               \SetKwFunction{RSVDT}{RSVDT}
                    \SetKwFunction{svd}{svd}
    \SetKwProg{Fn}{Function}{:}{}
   \Fn{ $[\dQ,\S,\V] \leftarrow$ \FMain{$\dA,\epsilon$, {$\mu$,r}} }{
    
    \KwData{ $\dA \in \RRn{n}{m_i}$, tolerance threshold $\epsilon \ge 0$. OPTIONAL ARG.:  ${\mu}>0$: machine precision parameter; ${r}:$ estimation of $\rank{\dA}$ (default $r =\textrm{ceil}(0.05 \MAX{n,m_i})$)}
    
    \KwResult{ Truncated SVD such that  $\normF{\dA - \dQ  {\S} \V^T}  \le \epsilon$   }

    $[\H,\B] \leftarrow$ \RORTH{$\dA$,$\mu$,r}   \tcp{$\H$: Basis matrix for the column space of $\dA$, see Algorithm \ref{alg:010}}  \label{alg:009_orth}

     $\epsilon_B \leftarrow   \epsilon \, \normF{\B}  $; $\; \mu_{mach} = \MAX{n,m_i}  \, \texttt{eps}(\normF{\B})   $

    $[\dQ,\S,\V] \leftarrow$ \SVDT{$\B$,$\epsilon_B$,$\mu_{mach}$}      \tcp{Truncated SVD of $\B = \H^T\dA$, see Algorithm \ref{alg:007}}  \label{alg:009_svd}
    
    $\dQ \leftarrow \H \, \dQ$;  \hspace{0.5cm}

   }
   \caption{Incremental Randomized SVD  \label{alg:009}} 
   \end{algorithm}

The pseudo-code of this randomized SVD is described in Algorithm \ref{alg:009}.  The basic steps  are   identical to the ones employed in  the  \texttt{SRSVD} of Algorithm \ref{alg:006}, namely, 1)  determination  of an orthogonal basis matrix $\H$ for the column space (range) of the input matrix $\dA$ in Line \ref{alg:009_orth};  and 2) Truncated SVD of  $\B = \H^T \dA$ in Line \ref{alg:009_svd}.  The computation of the basis matrix  in the first step also shares common features with the one employed in the    \texttt{SRSVD} of Algorithm \ref{alg:006}. Indeed, as can be inferred from the  the pseudo-code  of Algorithm \ref{alg:010}  of  the function   devoted to this task    ($\texttt{RORTHinc}$  ), the desired orthogonal basis matrix $\dH$  is    built iteratively  in a Gram-Schmid orthogonalization fashion:
\begin{equation}
 \dH =  \rowcuatro{\dH_1}{\dH_2}{\cdots}{\dH_s}, \quad  \textrm{where} \quad \dH_i^T \dH_j = \zero \; (i \neq j), \hspace{0.5cm}  \dH_i^T \dH_i = \ident. 
\end{equation}

The actual difference with \texttt{SRORTH} in Algorithm  \ref{alg:008} is that, at  a given iteration $i$, the corresponding  orthogonal matrix $\dH_i$ is not determined from a column partition of the input matrix, but rather as an orthogonal basis matrix  for the range of the matrix defined by 
\begin{equation}
\label{eq:0909323}
 \Y_i = \dfrac{1}{\sqrt{n \Delta R_i}} \OmegaB_i \C_i , \hspace{1cm} i = 1,2 \ldots s
\end{equation}
(see Lines \ref{alg:010_rand} and \ref{alg:010_orth} in Algorithm \ref{alg:010}). Here,    $\C_i \in \RRn{n}{m}$ is the residual matrix at   iteration $u$, whereas $\OmegaB_i$   is a  $n \times \Delta R_i$ standard Gaussian test matrix (a \emph{random} matrix whose entries are independent standard normal variables).  The distinguishing feature of our algorithm with respect to the original scheme put forward in Ref. \cite{martinsson2015randomized} is that, in our case, the number of columns $\Delta R_i$ of this random matrix changes during the iterations. In the   first iteration ($i=1$), we set $\Delta R_1 = R$, $R$ being the number of columns of the incremental basis matrix of the previous block matrix (see Line \ref{alg:008_est} in Algorithm \ref{alg:008}).  As argued in Ref. \cite{halko2011finding},   if this initial guess is well above the rank of input matrix $\dA$,  it is highly probable that   the basis matrix for the range of $\Y_1$ in Eq. \ref{eq:0909323} is  the required orthogonal matrix $\dH$. Numerical experience shows  that when the submatrices $\A_{k-1}$ and $\A_{k}$  ($k=2,3 \ldots$ p) correspond to input parameters that are close in parameter space,  then this estimation is normally a reliable upper bound, and therefore, only one iteration  is required. 

If this first iteration is not sufficient to reduce the norm of the residual matrix $\C_i$  below the prescribed error threshold $\mu$ (see Line \ref{alg:010_while} in Algorithm \ref{alg:010} ), then it is necessary to calculate a guess for  the number of columns of the random matrix  in the next iteration.  Our proposal in this respect is to use the logarithmic estimation displayed in Line \ref{alg:010_estl} of Algorithm \ref{alg:010}.  This estimation is based on the observation that, in most physical problems amenable to dimensionality reduction, the singular values of the integrand matrix decay in an exponential manner.   Nevertheless, to avoid situation in which the estimated increments $\Delta R_i$ are either too large or too small, the minimum and maximum sizes of the increment, $\dRmin$ and $\dRmax$ respectively,  can be also specified as optional arguments. 

Lastly, it should be noted that    this randomized factorization is also  subject to the vagaries of finite precision arithmetics. To address this,  Algorithm \ref{alg:010} includes a     re-orthogonalization step in Line \ref{alg:010_re}.

  \begin{algorithm}[!ht]
   \DontPrintSemicolon   
     \SetKwFunction{RORTH}{RORTHinc}
    \SetKwFunction{SVDT}{SVDT}
       \SetKwFunction{max}{max}
            \SetKwFunction{min}{min}
              \SetKwFunction{break}{break}
    \SetKwProg{Fn}{Function}{:}{}
     \Fn{ $[\H,\B,c] \leftarrow$ \RORTH{$\C,\mu,R, {\dRmin}, {\dRmax}$} }{
    
    \KwData{ $\C \in \RRn{n}{m}$, tolerance $\mu >0$;  rank estimate $R < q$ (where $q =\MIN{n,m}$). Optional arguments $ {\dRmin}$, $ {\dRmax}$ bounds for  rank increment $\Delta R$  (default  $\CEIL{0.01 q }$ and  $\CEIL{0.25 q }$)   }
    
    \KwResult{ Orthogonal matrix $\H$ such that  $\normF{\C - \H \B}  \le \mu$, where $\B = \H^T \C$  }

 $\H \leftarrow \emptyset$;  $\B \leftarrow \emptyset$; $\Delta R \leftarrow R \;$; $c \leftarrow \mu +1 \;$; $i \leftarrow 1 \;$;     $\; q \leftarrow \MIN{n,m}$
  
    \While{$c\ge \mu$ \label{alg:010_while}}{

   $\Omega \leftarrow \randn(n,\Delta R)/\sqrt{n  \Delta R}$  \tcp{Draw a random $n \times \Delta R$   matrix }  \label{alg:010_rand}
   
   $[\dH,\cdot,\cdot ]\leftarrow $ \SVDT{$\C \Omega$,0} \tcp{ Orthogonal basis matrix for $\range{\C \, \Omega}$ }  \label{alg:010_orth}
   
   \lIf{$\dH = \emptyset$}{\break \tcp{Exiting the loop} } \vspace{-0.2cm}  
   
    \lIf{$i>1$}{$[\dH,\cdot,\cdot] \leftarrow $ \SVDT{$\dH - \H(\H^T \dH)$,0}   \tcp{Re-orthogonalization}  \label{alg:010_re}}  
    \vspace{-0.2cm}
    $\dB \leftarrow \dH^T\C$;  $\;\C \leftarrow \C - \dH \dB; $  $ \;c \leftarrow \normF{\C}$ \tcp{Residual update} \label{alg:1:L7}
     $\H \leftarrow \rowdos{\H}{\dH}$;  $\B \leftarrow \rowdos{\B^T}{\dB^T}^T$;   \tcp{Basis matrix update}

     { $\bar{R} \leftarrow R_i + \dfrac{R - R_i}{log{c} - log{c_i}}(log{\mu} -log{c_i})$;   $\; \bar{R} \leftarrow$ \min{q,$\bar{R}$}  \tcp{Logarithmic estimation rank}  \label{alg:010_estl}
     
       $\; \Delta R \leftarrow$ \min{$\dRmax$,$\ceil(\bar{R}-R)$}; $\; \Delta R \leftarrow$ \max{$\dRmin$,$\Delta R$} ;      $\;\bar{R} \leftarrow R + \Delta R$  
     
   }     
  $i \leftarrow i +1;$  $\;c_i \leftarrow c;$  $\; R_i \leftarrow R;  \;$  $R \leftarrow \bar{R}$ 
   }

   }
   
   \caption{Incremental Randomized orthogonalization, with steps of varying size \label{alg:010}} 
   \end{algorithm}

\subsection{Numerical study}
\label{sec:numst}

To compare the performance of the standard SVD and the proposed SRSVD, we use the convergence study  presented in Table \ref{tab:4} of  Section \ref{sec:expsino} for determining  an orthogonal basis functions for the parameterized, vector-valued function   of \refeq{eq:param333}. It can be gleaned from this table that the proposed SRSVD clearly outperforms the standard SVD, both in terms of computing time and memory requirements. For instance, for the 16 GB matrix, the SRSVD turns out to be almost 7 times faster than the standard SVD, and for the largest matrix of 33 GB, the standard SVD   simply cannot handle the operation because it exhausts the memory capabilities of the  employed 64 GB RAM computer. Furthermore, we can see that, in passing from the matrix of 16 GB to the matrix of 33 GB, the computing time of the SRSVD   increases by a factor slightly below 2, a fact that indicates that the cost scales approximately linearly with the number of columns ---as opposed  to the standard SVD, whose asymptotic costs scales quadratically with  its minor dimension\cite{trefethen1997numerical}. It is noteworthy also that, as the number of partitions increases, the number of iterations required by the incremental randomized orhtogonalization of Algorithm \ref{alg:010} tends to one. This  indicates that, as conjectured in Section \ref{sec:irsvd} of Appendix \ref{sec:SRSVD},   the rank of a given block   matrix is a reliable upper bound for the rank of the orhtogonal  complement of the  next block matrix in the sequence.  Incidentally, this may explain in part why the asymptotic costs of the SQRSV appears to scale linearly with its minor dimension, for the standard randomized SVD  itself exhibits  this desirable feature \cite{martinsson2015randomized}.  Last but not least, we show in Table \ref{tab:4} the  relative difference between the singular values computed by both approaches. The results reveal that the difference is negligible, a fact that supports the theoretical claim made at the outset of this Appendix, according to which\emph{ the proposed SRSVD is not an approximate method} for computing the truncated SVD of a matrix, but rather an \emph{alternative method to compute the exact factorization} ---one that  exploits   linear correlations existing between   blocks of  the matrix.

\bibliographystyle{elsarticle-harv}

\bibliography{Bibliography_no_url}

\end{document}